\documentclass[11pt,letterpaper,twoside,reqno,nosumlimits]{amsart}

\allowdisplaybreaks

\usepackage[usenames,dvipsnames]{xcolor}
\usepackage{fancyhdr}
\usepackage{amsmath,amsfonts,amsbsy,amsgen,amscd,amssymb,amsthm}
\usepackage{url}

\usepackage{mathtools}

\usepackage[font=small,margin=0.25in,labelfont={bf},labelsep={colon}]{caption}

\usepackage{subcaption}

\usepackage{tikz}
\usepackage{microtype}
\usepackage{enumitem}

\definecolor{dark-gray}{gray}{0.3}
\definecolor{dkgray}{rgb}{.4,.4,.4}
\definecolor{dkblue}{rgb}{0,0,.5}
\definecolor{medblue}{rgb}{0,0,.75}
\definecolor{rust}{rgb}{0.5,0.1,0.1}

\usepackage[colorlinks=true]{hyperref}

\hypersetup{urlcolor=rust}
\hypersetup{citecolor=dkblue}
\hypersetup{linkcolor=dkblue}

\usepackage{setspace}

\usepackage{graphicx}
\usepackage{booktabs,longtable,tabu} 
\usepackage{multirow} 

\usepackage{float}

\usepackage[full]{textcomp}

\usepackage[scaled=.98,sups,osf]{XCharter}
\usepackage[scaled=1.04,varqu,varl]{inconsolata}
\usepackage[type1]{cabin}
\usepackage[charter,vvarbb,scaled=1.07]{newtxmath}
\usepackage[cal=boondoxo]{mathalfa}
\usepackage[T1]{fontenc}

\usepackage{bm} 

\graphicspath{{figures/}}

\theoremstyle{definition}

\numberwithin{equation}{section} 
\numberwithin{figure}{section}
\numberwithin{table}{section}

\floatstyle{plaintop}
\newfloat{recipe}{thp}{lor}
\floatname{recipe}{Recipe}
\numberwithin{recipe}{section}

\providecommand{\mathbold}[1]{\bm{#1}}  
                                
\newcommand{\diff}[1]{\mathrm{d}{#1}}
\newcommand{\idiff}[1]{\, \diff{#1}}

\newcommand{\vct}[1]{\mathbold{#1}}

\newcommand{\triplenorm}[1]{{\left\vert\kern-0.25ex\left\vert\kern-0.25ex\left\vert #1
    \right\vert\kern-0.25ex\right\vert\kern-0.25ex\right\vert}}

\usepackage[numbers]{natbib}

\newcommand{\om}{\omega}
\newcommand{\vom}{\vct{\omega}}
\newcommand{\vu}{\vct{u}}

\begin{document}

\title[Potential singularity of $3$D Euler equations]{Potential singularity of the $3$D Euler equations in the interior domain}

\author[T. Y. Hou]{Thomas Y. Hou}
\address{Applied and Computational Mathematics, California Institute of Technology, Pasadena, CA 91125, USA}
\email{hou@cms.caltech.edu}

\date{\today}

\begin{abstract}
Whether the $3$D incompressible Euler equations can develop a finite time singularity from smooth initial data is one of the most challenging problems in nonlinear PDEs. In this paper, we present some new numerical evidence that the $3$D axisymmetric incompressible Euler equations with smooth initial data of finite energy develop a potential finite time singularity at the origin. This potential singularity is different from the blow-up scenario revealed by Luo-Hou in  \cite{luo2014potentially,luo2014toward}, which occurs on the boundary. 
Our initial condition has a simple form and shares several attractive features of a more sophisticated initial condition constructed by Hou-Huang in \cite{Hou-Huang-2021,Hou-Huang-2022}. 
One important difference between these two blow-up scenarios is that the solution for our initial data has a one-scale structure instead of a two-scale structure reported in \cite{Hou-Huang-2021,Hou-Huang-2022}. More importantly, the solution seems to develop nearly self-similar scaling properties that are compatible with those of the $3$D Navier--Stokes equations. We will present numerical evidence that the $3$D Euler equations seem to develop a potential finite time singularity. Moreover, the nearly self-similar profile seems to be very stable to the small perturbation of the initial data.

\end{abstract}

\maketitle

\section{Introduction}

The question regarding the global regularity of the $3$D incompressible Euler equations with smooth initial data of finite energy is one of the most important fundamental questions in nonlinear partial differential equations. 
The main difficulty associated with the global regularity of the $3$D Euler equations is due to the presence of vortex stretching \cite{majda2002vorticity}. In \cite{luo2014potentially,luo2014toward}, Luo-Hou presented some strong numerical evidence that the $3$D axisymmetric Euler equations develop a finite time singularity on the boundary. The presence of the boundary and the odd-even symmetry properties of the solution seem to play an essential role in generating a stable blow-up of the $3$D Euler equations.

In this paper, we present some new numerical evidence that the $3$D axisymmetric incompressible Euler equations with smooth initial data of finite energy seem to develop a potential finite time singularity at the origin. Our initial condition has a very simple analytic expression and is purely driven by large swirl initially. This potential singularity is different from the blow-up scenario revealed by Luo-Hou in  \cite{luo2014potentially,luo2014toward}, which occurs on the boundary. It is also different from the two-scale traveling wave scenario considered by Hou-Huang in \cite{Hou-Huang-2021,Hou-Huang-2022} although the two scenarios share some common features. One important difference between these two blow-up scenarios is that the solution for our initial data has a one-scale structure instead of a two-scale structure reported in \cite{Hou-Huang-2021,Hou-Huang-2022}. More importantly, the solution seems to develop nearly self-similar scaling properties that are compatible with those of the $3$D Navier--Stokes equations. This property is critical for the potentially singular behavior of the $3$D Navier--Stokes equations using our new initial data.



We consider the $3$D axisymmetric Euler and Navier--Stokes equations in a periodic cylindrical domain. We impose a no-flow boundary condition at $r=1$ for the Euler equations and a no-slip no-flow boundary condition at $r=1$ for the Navier--Stokes equations. We use a periodic boundary condition in the axial variable $z$ with period $1$ for both Euler and Navier--Stokes equations.
Let $u^\theta$, $\omega^\theta$, and $\psi^\theta$ be the angular components of the velocity, the vorticity, and the vector stream function, respectively. 
Following \cite{hou2008dynamic}, we make the following change of variables: 
\[
u_1= u^\theta/r, \quad \omega_1=\omega^\theta/r, \quad \psi_1 = \psi^\theta/r.
\]
Then the Navier--Stokes equations can be expressed by the following equivalent system
\begin{subequations}\label{eq:axisymmetric_NSE_0}
\begin{align}
u_{1,t}+u^ru_{1,r}+u^zu_{1,z} &=2u_1\psi_{1,z} + \nu \Delta u_1,\label{eq:as_NSE_0_a}\\
\om_{1,t}+u^r\om_{1,r}+u^z\om_{1,z} &=2u_1u_{1,z} + \nu \Delta \om_1,\label{eq:as_NSE_0_b}\\
 -\left(\partial_r^2+\frac{3}{r}\partial_r+\partial_z^2\right)\psi_1 &= \om_1,\label{eq:as_NSE_0_c}
\end{align}
\end{subequations}
where $u^r=-r\psi_{1,z}, \; u^z =2\psi_1 + r\psi_{1,r}$, and 
$\Delta = \partial_r^2+\frac{3}{r}\partial_r+\partial_z^2$ is the five-dimensional diffusion operator. 

\vspace{-0.1in}
\subsection{Major features of the potential singularity of the Euler equations}

Although the angular vorticity is set to zero initially, the large swirl and the oddness of $u_1$ as a function of $z$ induce a large odd angular vorticity dynamically. The oddness of angular vorticity induce two antisymmetric vortex dipoles, which generate a hyperbolic flow structure near the symmetry axis $r=0$. The antisymmetric vortex dipoles produce a strong shear layer for the axial velocity, which transports the solution toward $z=0$. At the same time, they also induce an antisymmetric local convective circulation that pushes the solution near $z=0$ toward the symmetry axis. Moreover, the oddness of $u_1$ in $z$ generates a large positive gradient $u_{1z}$ dynamically, which induces a rapid growth of $\omega_1$. The rapid growth of $\omega_1$ in turn feeds back to the rapid growth of $\psi_{1,z}$, forming a positive feedback loop. 

After a relatively short transition period, the solution of the $3$D Euler equations seems to develop nearly self-similar scaling properties. Denote by $(R(t),Z(t))$ the location in the $rz$ plane at which $u_1$ achieves its maximum. If we introduce $\xi = (r-R(t))/Z(t)$ and $\zeta=z/Z(t)$ as the dynamically rescaled variables, we observe nearly self-similar profiles of the rescaled solutions in terms of $(\xi, \zeta)$. Moreover, $Z(t)$ seems to scale like $O((T-t)^{1/2})$, which is consistent with potential blow-up scaling of the $3$D Navier--Stokes equations. In our computation, we observe that $R(t)$ and $Z(t)$ are roughly of the same order, but they have not settled down to a stable scaling relationship without viscous regularization.
On the other hand, our numerical results seem to indicate that the maximum vorticity grows like $O((T-t)^{-1})$ and $\int_0^t \|\vom (s)\|_{L^\infty}ds $ seems to grow without bound. The Beale-Kato-Majda blow-up criterion \cite{beale1984remarks} implies that the $3$D Euler equations would develop a finite time singularity. Furthermore, we observe that the maximum velocity grows like $O((T-t)^{-1/2})$, $\|u_1\|_{L^\infty}$ and $\|\psi_{1z}\|_{L^\infty}$ scale like $O((T-t)^{-1})$, and $\|\omega_1\|_{L^\infty}$ scales like $
O((T-t)^{-3/2})$. These scaling properties are consistent with the scaling property $Z(t) \sim (T-t)^{1/2}$.

One important feature is that $\psi_{1,z}(t,r,z)$ is relatively flat in a local region near the origin $\{(r,z) \in [0,0.9R(t)]\times [0,0.5Z(t)]\}$. This important property that was not observed in \cite{Hou-Huang-2021,Hou-Huang-2022}. We observe that $\psi_{1,z}(t,r,z)$ drops quickly beyond this local region and becomes negative near the tail region.  The large value of $\psi_{1,z}$ in the local region near the origin generates a large growth of $u_1$ through the vortex stretching term $2 \psi_{1,z} u_1$ in \eqref{eq:as_NSE_0_a}. On the other hand, the small or negative value of $\psi_{1,z}$ in the tail region leads a relatively slower growth rate of $u_1$. The difference in the growth rate in the local region and the tail region produces a one-scale traveling wave solution propagating toward the origin, overcoming the upward transport along the $z$-direction.

Our velocity field shares some important features observed in \cite{Hou-Huang-2021,Hou-Huang-2022}. In particular, we observe that the $2$D velocity field $(u^r(t),u^z(t))$ in the $rz$-plane forms a closed circle right above $(R(t),Z(t))$ and the corresponding streamline is trapped in the circle region in the $rz$-plane. This local circle structure has the desirable property of keeping the bulk parts of the $u_1,\om_1$ profiles near the most singular region without being transported away by the upward advection. The flow in this local circle region spins rapidly around the symmetry axis. As we get closer to the symmetry axis, the streamlines induced by the velocity field travel upward along the vertical direction and then move outward along the radial direction. The local blow-up solution resembles the structure of a tornado. On the other hand, we do not observe the formation of a no-spinning region in our blow-up scenario. The formation of a no-spinning region is an important signature of the two-scale traveling singularity reported in \cite{Hou-Huang-2021,Hou-Huang-2022}.

\vspace{-0.05in}
\subsection{Potentially singular behavior of the 3D Navier--Stokes equations}
Given that the solution of the $3$D Euler equations has scaling properties that are compatible with those of the Navier--Stokes equations, it is natural to investigate whether the $3$D Navier--Stokes equations would develop a potential finite time singularity using the same initial condition. We have performed some preliminary study by solving the Navier--Stokes equations with a relatively large viscosity $\nu=5\cdot 10^{-3}$.
Surprisingly, the viscous regularization enhances nonlinear alignment of vortex stretching. We observe a relatively long stable phase of strong nonlinear alignment of vortex stretching. The solution of the $3$D Navier--Stokes equations develops nearly self-similar scaling properties that are similar to what we have observed for the solution of the $3$D Euler equations. Moreover, the maximum vorticity has increased by a factor of $10^7$. To the best of our knowledge, such a large growth rate of maximum vorticity has not been reported in the literature for the $3$D Navier--Stokes equations. We refer to the companion paper \cite{Hou-nse-2021} on the potentially singular behavior of the Navier--Stokes equations for more discussion.

\vspace{-0.05in}
\subsection{Comparison with the two-scale traveling wave singularity} 

Our initial condition shares several attractive features of a more sophisticated initial condition constructed in 
\cite{Hou-Huang-2021,Hou-Huang-2022}.
However, there are also some important differences between our new initial condition and the one considered in 
\cite{Hou-Huang-2021,Hou-Huang-2022}. First, the solution of the $3$D Euler equations studied in \cite{Hou-Huang-2021,Hou-Huang-2022} has a three-scale structure. The smallest scale characterized by the thickness of the sharp front does not seem to settle down to a stable scaling relationship. That is why it is essential to apply degenerate viscosity coefficients of order $O(r^2)+ O(z^2)$ to select a two-scale solution structure. In comparison, the solution of the $3$D Euler equations using our initial data has essentially a one-scale structure with scaling properties compatible with those of the $3$D Navier--Stokes equations. This property is critical for us to observe potentially singular behavior of the $3$D Navier-Stokes equations. In comparison, the maximum vorticity of the $3$D Navier--Stokes solution with a constant viscosity $\nu=10^{-5}$ reported in \cite{Hou-Huang-2021} has grown by a factor less than $2$. 

A second important difference is that we roughly have $\|\psi_1 \|_{L^\infty} \sim (T-t)^{-1/2}$ for our new initial data while for the solution in 
\cite{Hou-Huang-2021,Hou-Huang-2022} we have $\|\psi_1 \|_{L^\infty} = O(1)$, which is consistent with scaling property $Z(t) \sim (T-t)$. 
Finally, the solution of the Navier--Stokes equations with degenerate viscosity coefficients in \cite{Hou-Huang-2021,Hou-Huang-2022} develops strong shearing instability in the tail region. We need to apply some low pass filtering to stabilize this shearing instability. In comparison, our solutions are very stable and do not suffer from the shearing instability in the tail region, thus 
there is no need to apply any low pass filtering.

\vspace{-0.075in}
\subsection{Numerical Methods}
We use a similar adaptive mesh strategy developed in \cite{Hou-Huang-2021} by constructing two adaptive mesh maps for $r$ and $z$ explicitly. For the Navier--Stokes equations, our solutions are much smoother than those considered in \cite{Hou-Huang-2021,Hou-Huang-2022} due to the relative large viscous regularization. For the Euler equations, we need to allocate a large percentage of the adaptive mesh to resolve the sharp front. The elliptic problem for the stream function becomes very ill-conditioned in the late stage. Moreover the interpolation from one adaptive mesh to another adaptive mesh introduces some high frequency errors. To alleviate this difficulty, we apply a second order numerical diffusion with $\nu=1/n_1^2$ for the $u_1$ and $\omega_1$ equations, here $n_1$ is the number of mesh points along $z$ direction (we use a square grid).
%

We use a second order finite difference method to discretize the spatial derivatives, and a second order explicit Runge--Kutta method to discretize in time. An adaptive time-step size is used according to the standard time-stepping stability constraint with the smallest time-step size of order $O(10^{-15})$. As in \cite{Hou-Huang-2021,Hou-Huang-2022}, we adopt the second order B-spline based Galerkin method developed in \cite{luo2014potentially,luo2014toward} to solve the Poisson equation for the stream function.  The overall method is second order accurate. We have performed careful resolution study and confirm that our method gives at least second order accuracy in the maximum norm.

\vspace{-0.05in}
\subsection{Review of previous works}

There have been a number of theoretical developments for the $3$D incompressible Euler equations. These include the well-known non-blowup criterion due to Beale--Kato--Majda \cite{beale1984remarks}, the geometric non-blowup criteria due to Constantin--Fefferman--Majda \cite{cfm1996} and its Lagrangian analog due to Deng-Hou-Yu \cite{dhy2005}. Recently Elgindi-Jeong \cite{elgindi2021} proved finite time singularity formation for incomprerssible $3$D Euler equations with bounded and piecewise smooth vorticities using Kida symmetry flow. There was a recent breakthrough due to Elgindi \cite{elgindi2019finite} (see also \cite{Elg19}) who proved that the $3$D axisymmetric Euler equations develop a finite time singularity for a class of $C^{1,\alpha}$ initial velocity with no swirl. There have been a number of interesting theoretical results inspired by the Hou--Luo blow-up scenario \cite{luo2014potentially,luo2014toward}, see e.g. \cite{kiselev2014small,choi2015,choi2014on,kryz2016,chen2019finite2,chen2020finite,chen2021finite3}
and the excellent survey article \cite{kiselev2018}.

There has been a number of previous attempts to search for potential Euler singularities numerically. These include the work of Grauer--Sideris \cite{gs1991} for the $3$D axisymmetric Euler equations, the work of E and Shu \cite{es1994} for the $2$D Boussinesq equations, the two anti-parallel vortex tube computation by Kerr in \cite{kerr1993} and a related work by Hou-Li in\cite{hl2006}, the work of  Boratav and Pelz in \cite{bp1994} for Kida's high-symmetry initial data (see also \cite{hl2008}), and a more recent work of Luo-Hou for $3$D axisymmetric Euler equations \cite{luo2014potentially,luo2014toward}. There is also an interesting proposal for potential Euler singularity by Brenner--Hormoz--Pumir in \cite{brenner2016euler}. We refer to a review article \cite{gibbon2008} for more discussions on potential Euler singularities.

In \cite{vasseur2020}, Vasseur and Vishik showed that blow-up solutions to $3$D Euler are hydrodynamically unstable (see also \cite{vasseur2021} for the axisymmetric Euler). In a recent preprint \cite{chenhou2022}, we use a Riccati eqution, the inviscid Burgers equation, the axisymmetric Euler equations with $C^{1,\alpha}$ velocity and boundary to demonstrate that the notion of stability proposed in \cite{vasseur2020} may not be suitable to study the stability of a blow-up solution. See Section \ref{blowup-stability} for a summary of the main findings in \cite{chenhou2022}.


The rest of the paper is organized as follows. In Section \ref{sec:setup}, we describe the setup of the problem. In Section \ref{sec:euler}, we report the potential finite time singularity of the $3$D Euler equations. In particular, we will study the rapid growth of the solution, the velocity field and the dipole structure. We also perform careful resolution study and scaling analysis. 
Some concluding remarks are made in Section \ref{sec:conclude}. The technical details regarding the construction of our adaptive mesh for the $3$D Euler equations will be deferred to the Appendix. 


\vspace{-0.1in}
\section{Description of the Problem}\label{sec:setup}
We consider the $3$D incompressible Navier--Stokes equations:
\begin{equation}\label{eq:NSE_vc}
\begin{split}
\vu_t + \vu\cdot \nabla\vu = -\nabla p + \nu \Delta \vu \;,\quad
\nabla \cdot \vu = 0,
\end{split}
\end{equation}
where $\vu = (u^x,u^y,u^z)^T$ is the $3$D velocity vector, $p$ is the scalar pressure, $\nabla = (\partial_x,\partial_y,\partial_z)^T$ is the gradient operator in $\mathbb{R}^3$, and $\nu$ is a constant diffusion coefficient.  When the diffusion is absent (i.e. $\nu = 0$), equations \eqref{eq:NSE_vc} reduce to the $3$D Euler equations. 
In this paper, we will study the potential singularity formulation for $3$D axisymmetric Euler and Navier--Stokes equations. In the axisymmetric scenario, it is convenient to rewrite equations~\eqref{eq:NSE_vc} in cylindrical coordinates. Consider the change of variables 
\[x = r\cos\theta,\quad y = r\sin\theta,\quad z = z.\]
We decompose the radially symmetric velocity field as follows 
\[\vct{u}(t,r,z) = u^r(t,r,z)\vct{e}_r + u^\theta(t,r,z)\vct{e}_\theta + u^z(t,r,z)\vct{e}_z,\]
\[\vct{e}_r = \frac{1}{r}(x,y,0)^T,\quad \vct{e}_\theta = \frac{1}{r}(-y,x,0)^T,\quad \vct{e}_z = (0,0,1)^T.\]
Define $\vom = \nabla\times \vu$ as the $3$D vorticity vector. The vorticity can be represented in cylindrical coordinates as follows:
\[\vom (t,r,z) = -(u^\theta)_z\vct{e}_r + \omega^\theta(t,r,z)\vct{e}_\theta + \frac{1}{r}(r u^\theta)_r\vct{e}_z.\]
%
Let $\psi^\theta$ be the angular stream function.
In \cite{hou2008dynamic}, Hou and Li introduced the variables 
\[u_1 = u^\theta/r,\quad \om_1=\om^\theta/r,\quad \psi_1 = \psi^\theta/r \]
and derived an equivalent system of equations for the axisymmetric Navier--Stokes equations as follows:
\begin{subequations}\label{eq:axisymmetric_NSE_1}
\begin{align}
u_{1,t}+u^ru_{1,r}+u^zu_{1,z} &=2u_1\psi_{1,z} + \nu\left(u_{1,rr} + \frac{3}{r}u_{1,r}\right) + \nu u_{1,zz},\label{eq:as_NSE_1_a}\\
\om_{1,t}+u^r\om_{1,r}+u^z\om_{1,z} &=2u_1u_{1,z} + \nu\left(\om_{1,rr} + \frac{3}{r}\om_{1,r}\right) + \nu \om_{1,zz}\label{eq:as_NSE_1_b}\\
 -\left(\partial_r^2+\frac{3}{r}\partial_r+\partial_z^2\right)\psi_1 &= \om_1,\label{eq:as_NSE_1_c}\\
 u^r=-r\psi_{1,z},\quad u^z&=2\psi_1 + r\psi_{1,r}. \label{eq:as_NSE_1_d}
\end{align}
\end{subequations}
This reformulation has the advantage of removing the $1/r$ singularity from the cylindrical coordinates. 

Our smooth initial condition has a very simple form and is given below:
\begin{equation}
\label{eq:initial-data}
u_1 (0,r,z) =\frac{12000(1-r^2)^{18}
\sin(2 \pi z)}{1+12.5(\sin(\pi z))^2}, \quad \om_1(0,r,z)=0.
\end{equation}
The nontrivial part of the initial data lies only in the angular velocity $u^\theta = r u_1$. The other two velocity components are set to zero initially. Thus the flow is completely driven by large swirl initially.
An important property of this initial condition is that $u_1$ is an odd function of $z$ and decays rapidly as $r$ approaches the boundary $r=1$. The specific form of the denominator in $u_1$ is also important. It breaks the symmetry along the $z$-direction and generates an appropriate bias toward $z=0$. The maximum (or the minimum) of $u_1$ is located at $r=0$ and is closer to $z=0$ than to $z=1/2$ (or $z=-1/2$). This specific form of the initial condition is designed in such a way that the solution has comparable scales along the $r$ and $z$ directions, leading to a one-scale traveling solution moving toward the origin.  

Such initial condition does not favor nonlinear alignment of vortex stretching initially. In fact, the maximum of $u_1$ first decreases in time due to the large negative value of $\psi_{1z}$ in the very early stage. Then the solution propagates away from the symmetry axis $r=0$ and generates some favorable solution structure. The location of maximum $u_1$ soon turns around and approaches the origin. From that time on, the solution generates a positive feedback loop to maintain strong nonlinear alignment of vortex stretching. 

\vspace{-0.1in}
\subsection{Settings of the solution}\label{sec:settings} We will solve the transformed equations~\eqref{eq:axisymmetric_NSE_1}-\eqref{eq:initial-data} in the cylinder
\[\mathcal{D} = \{(r,z): 0\leq r\leq 1\}.\]
We will impose a periodic boundary condition in $z$ with period $1$ and the odd symmetry of $u_1$, $\omega_1$ and $\psi_1$ as a function of $z$. According to \cite{liuwang2006}, $u^\theta,\om^\theta,\psi^\theta$ must be an odd function of $r$, which implies that $u_1,\om_1,\psi_1$ must be an even function of $r$. Thus, we impose the following pole conditions:
\begin{equation}\label{eq:even_r}
u_{1,r}(t,0,z) = \om_{1,r}(t,0,z) = \psi_{1,r}(t,0,z) = 0.
\end{equation}
For both the Euler and Navier--Stokes equations, the velocity satisfies a no-flow boundary condition on the solid boundary $r=1$:
\begin{equation}\label{eq:no-flow}
\psi_1(t,1,z) = 0\quad \text{for all $z$}.
\end{equation}
For the Navier--Stokes equations, we further impose a no-slip boundary condition at $r=1$:
\begin{equation}\label{eq:no-slip}
u^\theta(t,1,z) = u^z(t,1,z) = 0,\quad \text{for all $z$}.
\end{equation}
Using \eqref{eq:as_NSE_1_d} and \eqref{eq:no-flow}, we obtain $\psi_{1,r}(t,1,z) = 0$. Therefore, the no-slip boundary is reduced to the following boundary condition for $u_1$ and $\omega_1$ at $r=1$:
\begin{equation}\label{eq:no-slip1}
u_1(t,1,z) = 0,\quad \om_1(t,1,z) = -\psi_{1,rr}(t,1,z),\quad \text{for all $z$}.
\end{equation}
The second no-slip boundary condition will be implmeneted as a numerical vorticity boundary condition for $\omega_1$ by using a ghost grid point and enforcing $\psi_{1,r}(t,1,z) = 0$.
By the periodicity and the odd symmetry of the solution, we only need to solve equations~\eqref{eq:axisymmetric_NSE_1} in the half-period domain 
\[\mathcal{D}_1 = \{(r,z): 0\leq r\leq 1, 0\leq z\leq 1/2\}.\] 
Moreover, we have
\[u^r = -r \psi_{1,z} = 0 \quad \text{on $r=0,1$}\quad \text{and}\quad u^z = 2\psi_1+r\psi_{1,r} = 0\quad \text{on $z=0,1/2$},\]
due to the periodicity and the odd symmetry of the solution. Thus 
the boundaries of $\mathcal{D}_1$ behave like ``impermeable walls''. 
To numerically compute the potential singularity formation of the equations \eqref{eq:axisymmetric_NSE_1}-\eqref{eq:initial-data}, we adopt the numerical methods developed in my recent joint work with Dr. De Huang \cite{Hou-Huang-2021} except that there is no need to apply any low pass filtering in our computation. The detailed descriptions of our numerical methods can be found in Appendix A in \cite{Hou-Huang-2021}. 


A key step in our numerical method is the construction of the adaptive mesh.
The construction of the adaptive mesh will be different for the $3$D Euler and Navier-Stokes equations since the Euler solution is more singular than the Navier--Stokes solution. We need to allocate more grid points to resolve the sharp front for the $3$D Euler equations in the late stage. We will provide more details how to construct the adaptive mesh for the $3$D Euler and Navier--Stokes equations in the Appendix.

\section{Potential finite time singularity of the 3D Euler equations}
\label{sec:euler}
In this section, we will investigate potential finite time singularity of the $3$D Euler equations at the origin. Due to the fast decay of the initial data near the boundary, the boundary at $r=1$ essentially has no effect on the singularity formation. This is very different from the Hou-Luo blow-up scenario in which the boundary plays an essential role in generating a stable finite time singularity.

\subsection{Numerical Results: First Sign of Singularity}\label{sec:first_sign}
We have numerically solved the $3$D axisymmetric Euler equations \eqref{eq:axisymmetric_NSE_1}-\eqref{eq:initial-data} with $\nu=0$ on the half-period cylinder $\mathcal{D}_1=\{(r,z):0\leq r\leq 1, 0\leq z\leq 1/2\}$ using meshes of size $(n_1,n_2) = (256p,256p)$ for $p = 2, 3, \dots, 6$. In this subsection, we first present the major features of the potential Euler singularity revealed by our computation. In Section \ref{sec:performance_study_euler}, we carry out a careful resolution study of the numerical solutions. Then we provide some qualitative study of the scaling properties in Section \ref{sec:scaling_study_euler}.

\subsubsection{Profile evolution}\label{sec:profile_evolution_euler}
In this subsection, we investigate how the profiles of the solution evolve in time. We will use the numerical results computed on the adaptive mesh of size $(n_1,n_2) = (1536,1536)$. We have computed the numerical solution up to time $t=0.00227648$ when it is still well resolved. 

\begin{figure}[!ht]
\centering
    \includegraphics[width=0.32\textwidth]{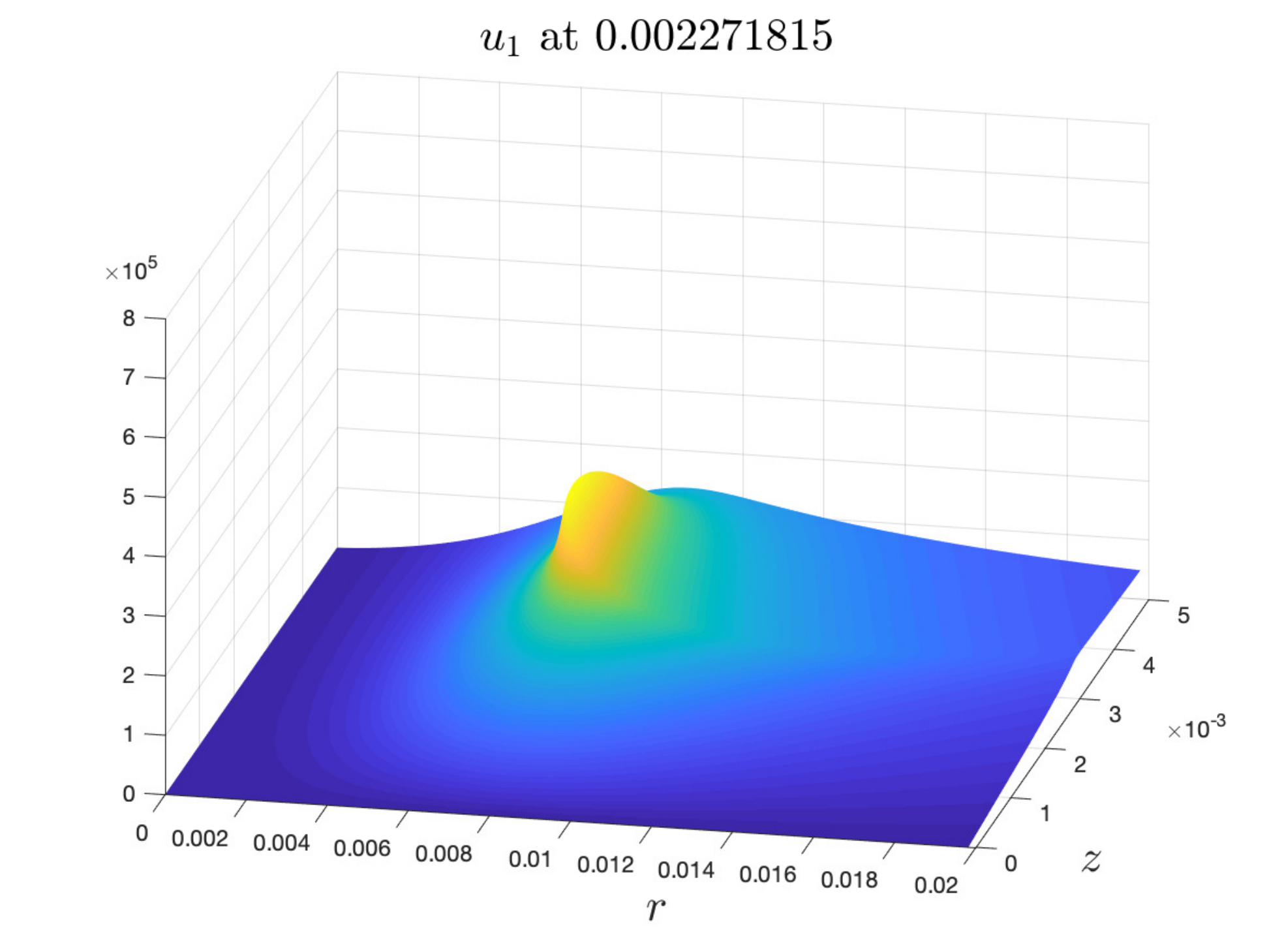}
    \includegraphics[width=0.32\textwidth]{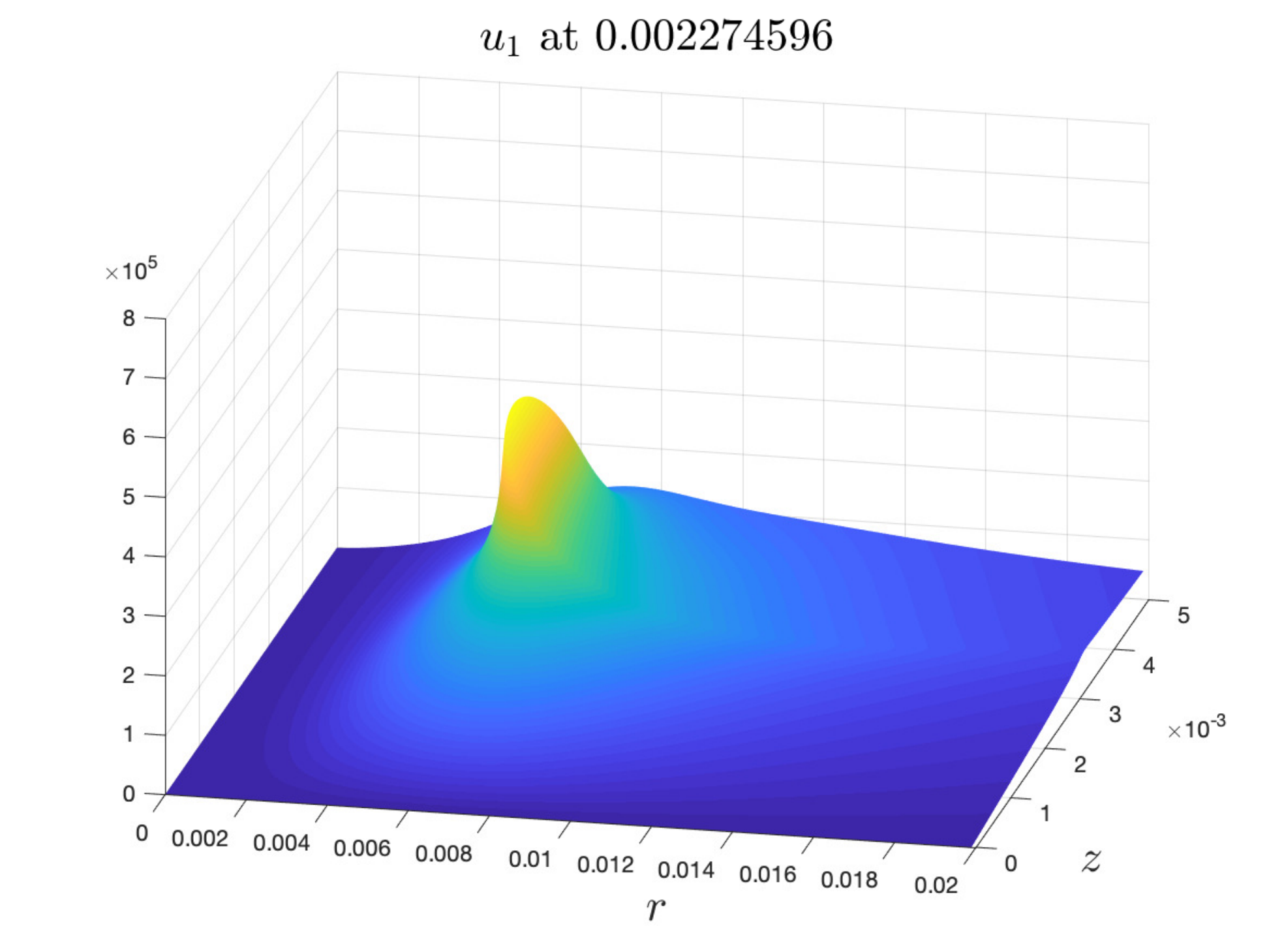}
    \includegraphics[width=0.32\textwidth]{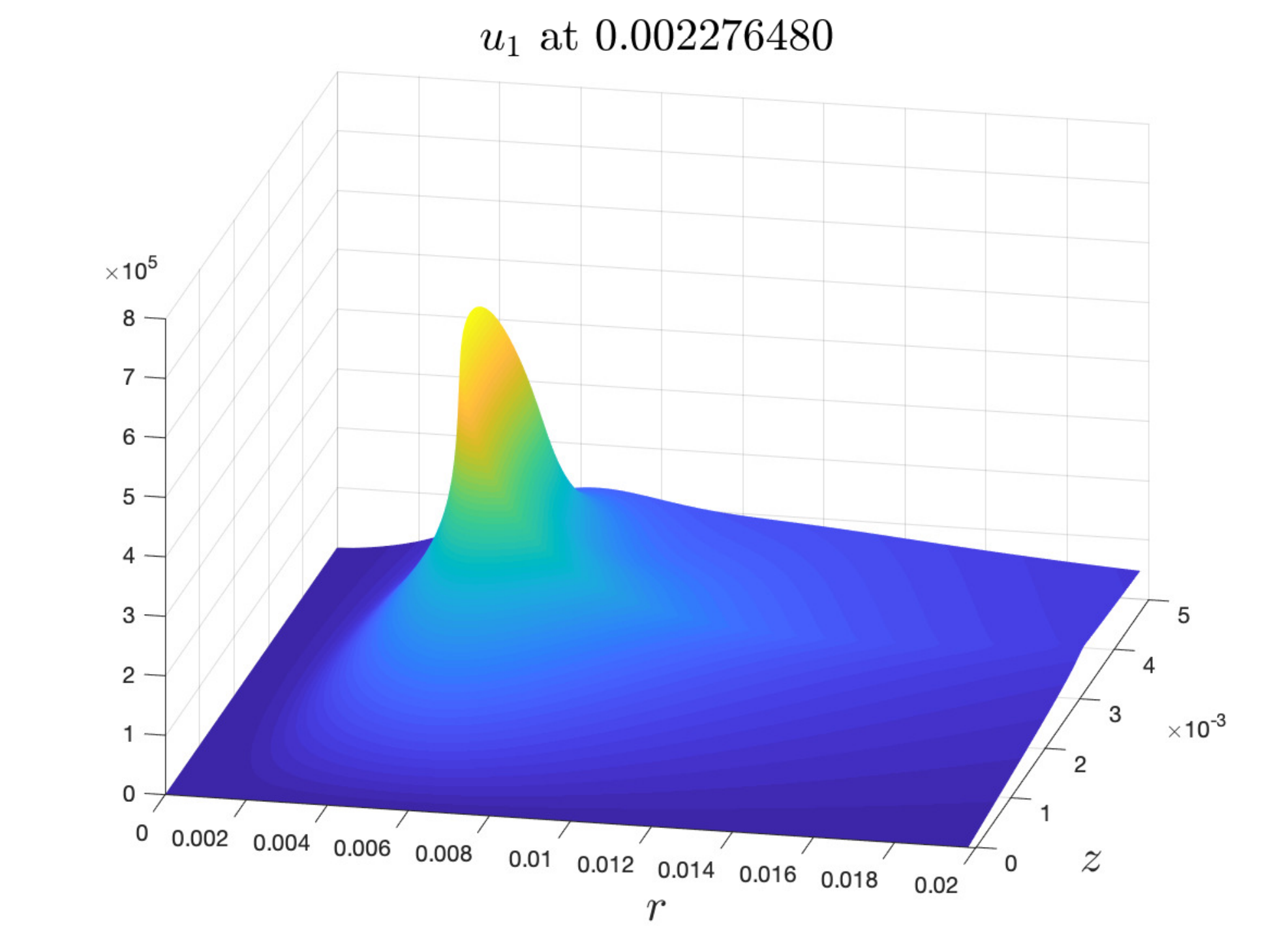}
   \includegraphics[width=0.32\textwidth]{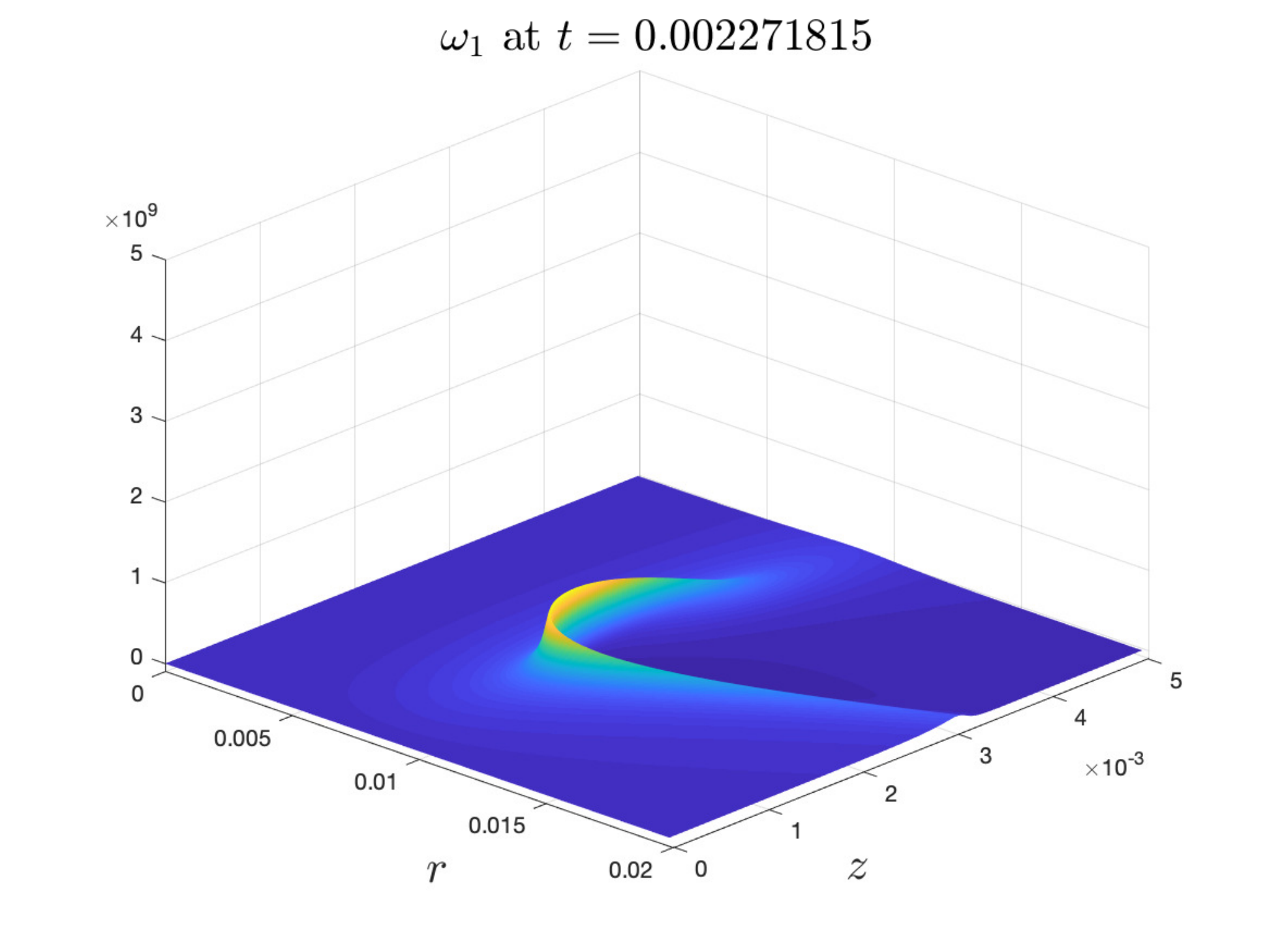}
    \includegraphics[width=0.32\textwidth]{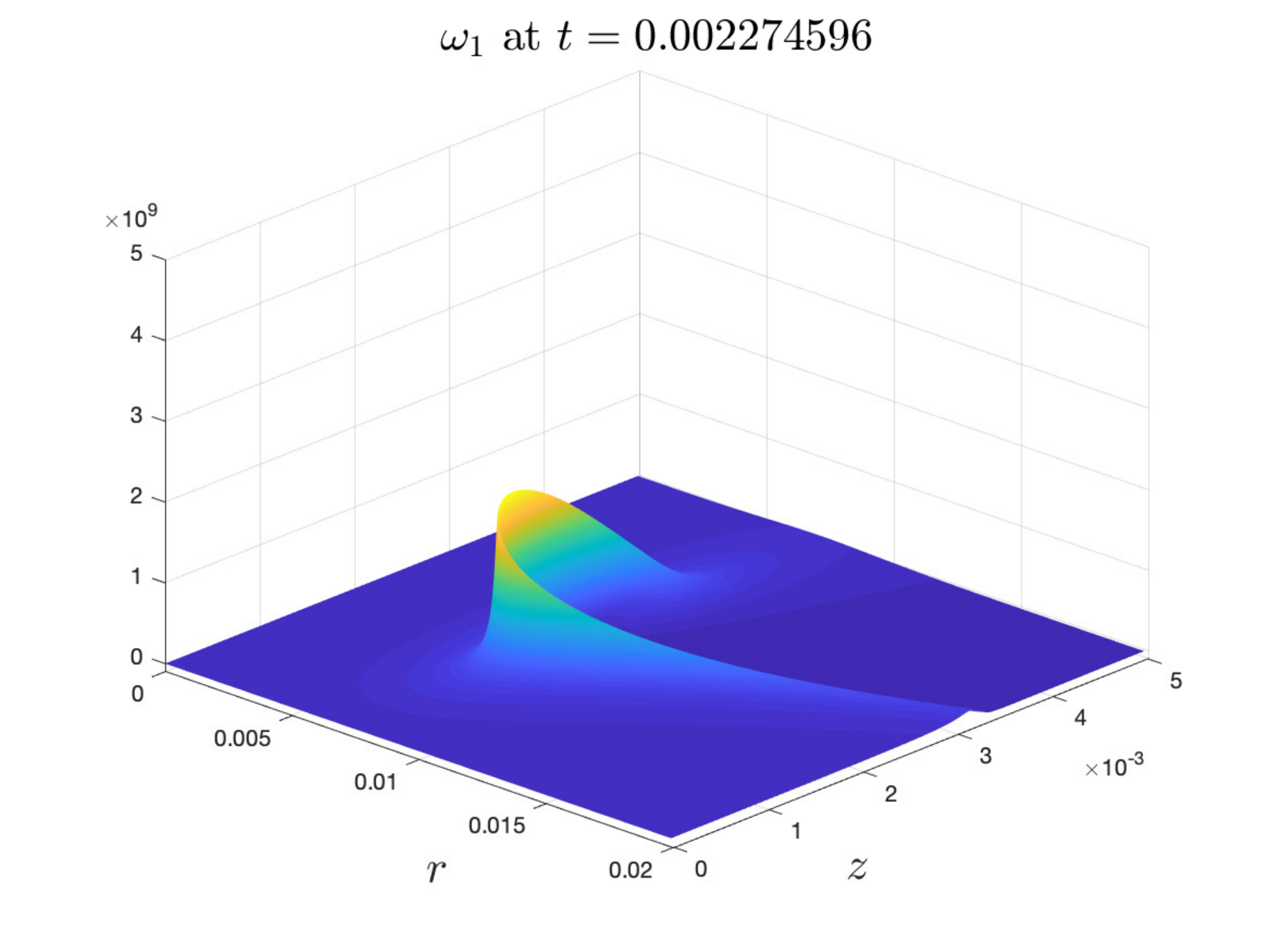}
    \includegraphics[width=0.32\textwidth]{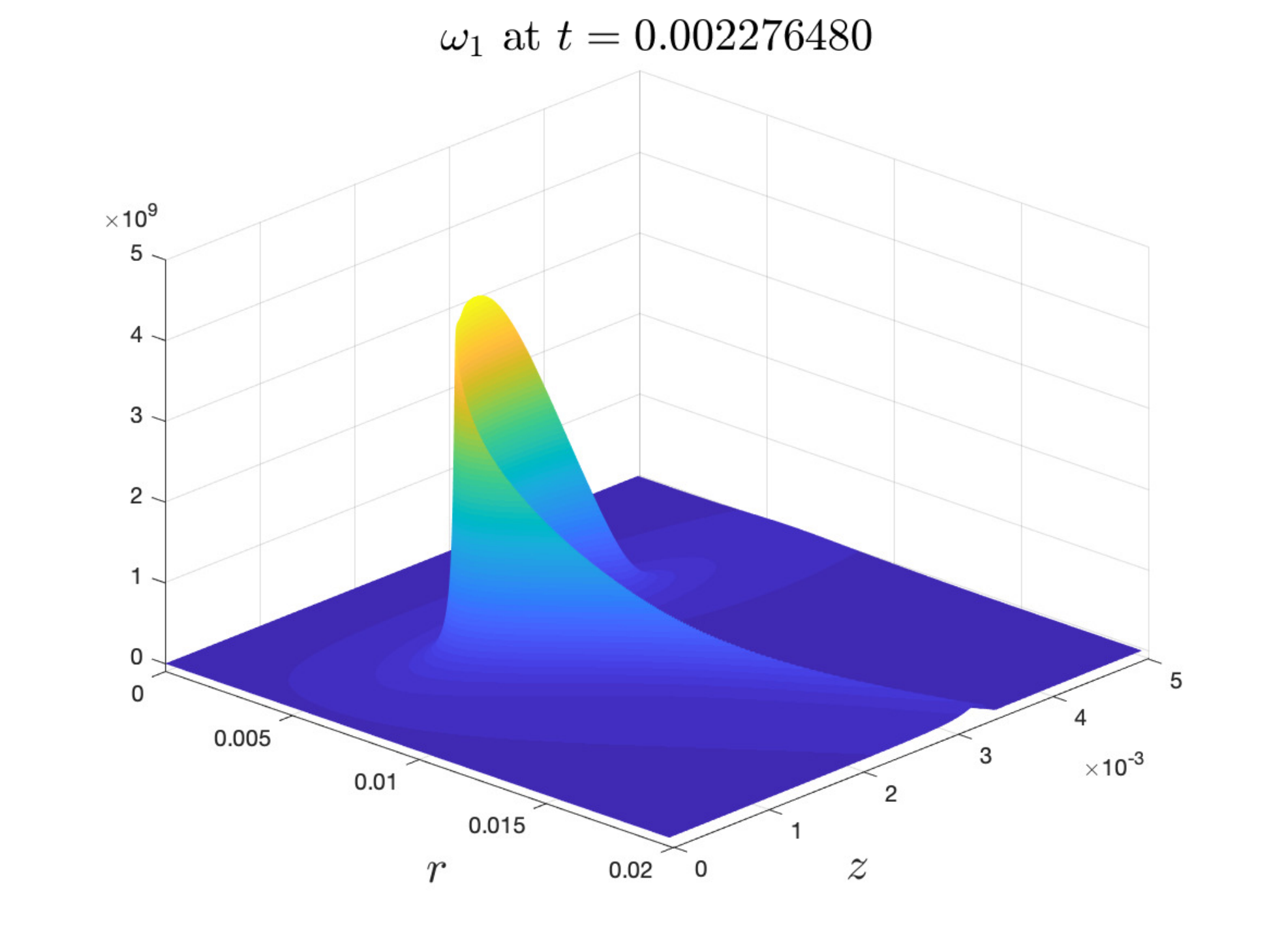}
    \caption[Profile evolution]{The evolution of the profiles of $u_1$ (row $1$) and $\om_1$ (row $2$) for the $3$D Euler equations at three different times.}  
    \label{fig:profile_evolution}
\end{figure}

Figure \ref{fig:profile_evolution} illustrates the evolution of $u_1,\om_1$ by showing the solution profiles at $3$ different times $t = 0.002271815, \;0.002274596, \;0.002276480$. We can see that the magnitudes of $u_1,\om_1$ grow in time. The profiles travel toward the origin and the singular support shrinks in space. The solution $u_1$ develops sharp gradients in both directions. 
Moreover, we observe that $\om_1$ is essentially supported along a thin curved region. Both $u_1$ and $\om_1$ form a tail part that decays rapidly into the far field. 
 
\begin{figure}[!ht]
\centering
    \includegraphics[width=0.38\textwidth]{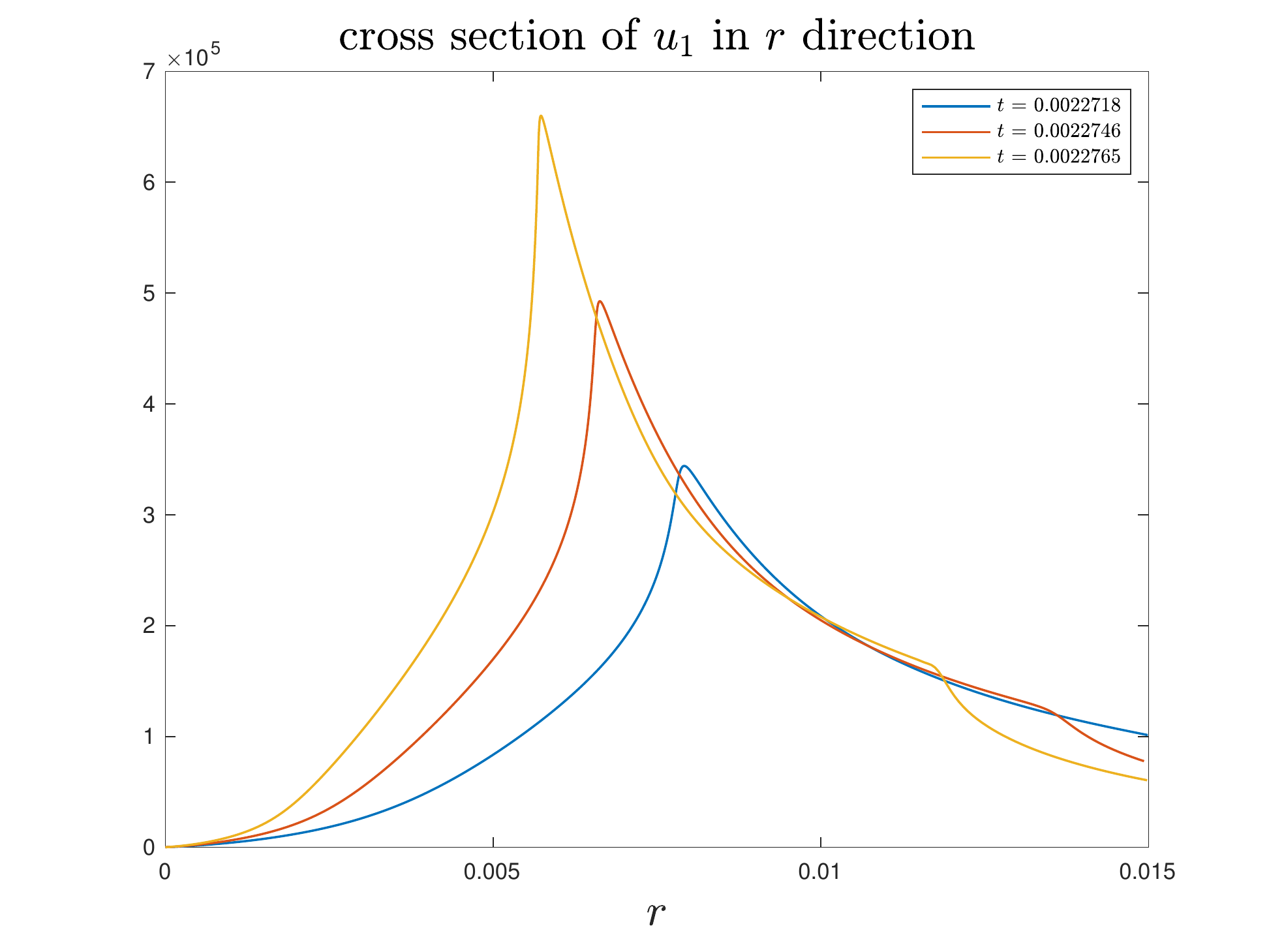}
    \includegraphics[width=0.38\textwidth]{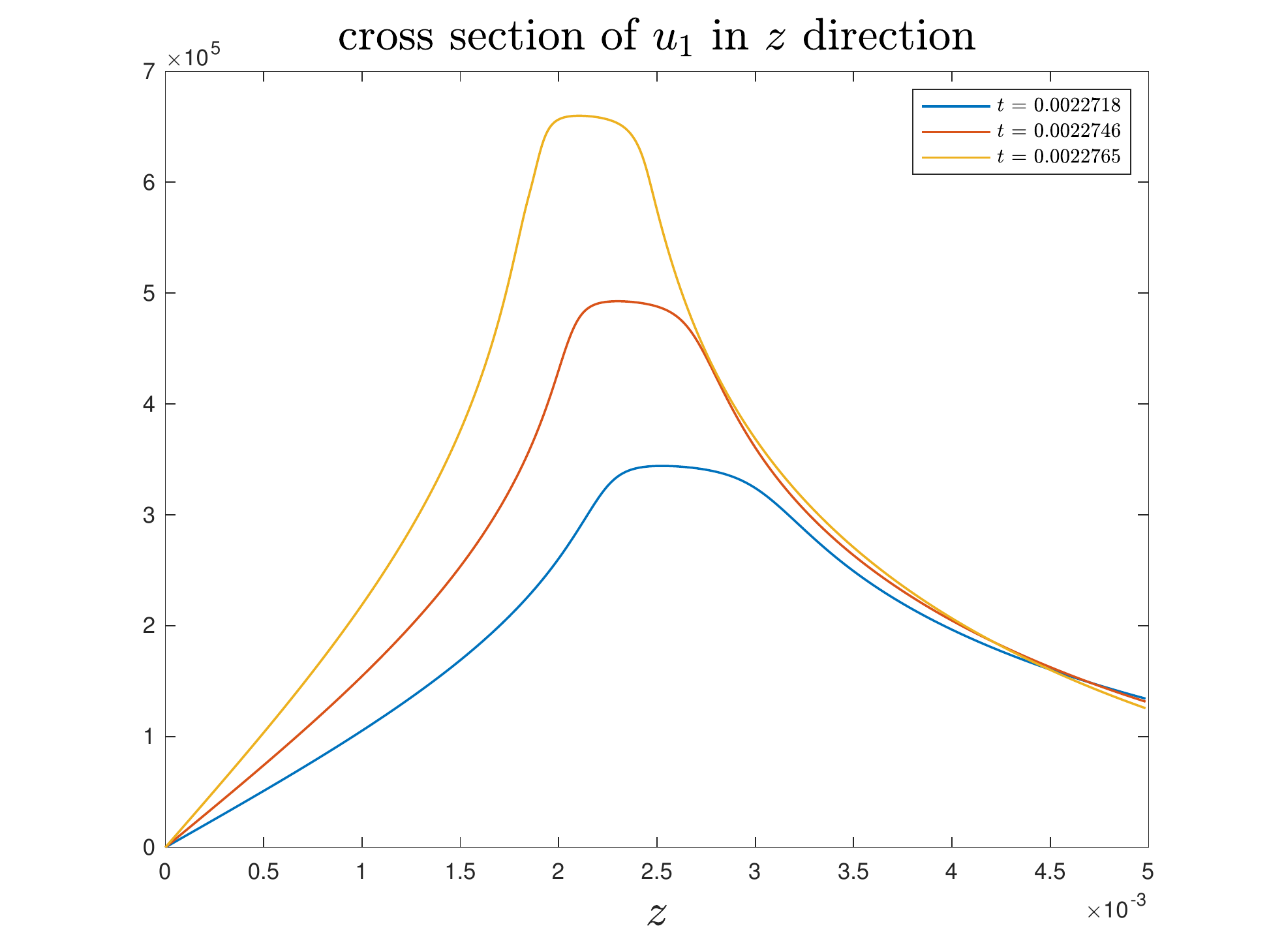}
    \caption[Cross section]{Cross sections of $u_1$ in both directions at $t=0.002271815$, $t=0.002274596$ and $t=0.002276480$, respectively.}  
    \label{fig:cross_section}
\end{figure}

Let $(R(t),Z(t))$ denote the maximum location of $u_1(t,r,z)$. We will always use this notation throughout the paper. In Figure \ref{fig:cross_section}, we plot the cross sections of $u_1$ going through the point $(R(t),Z(t))$ in both directions, i.e. $u_1(t,r,Z(t))$ versus $r$ and $u_1(t,R(t),z)$ versus $z$, respectively. We can see more clearly that $u_1$ develops sharp gradients in both directions and $u_1$ develops a sharp front along the $r$-direction. Unlike the two-scale travel wave solution reported in \cite{Hou-Huang-2021,Hou-Huang-2022}, we do not observe the formation of a no-spinning region between the sharp front and $r=0$ and $u_1$ does not form a compact support along the $z$-direction that is shrinking toward $z=0$. 

\subsection{Trajectory and alignment}\label{sec:two_scale} 
\begin{figure}[!ht]
\centering
    \begin{subfigure}[b]{0.38\textwidth}
        \centering
        \includegraphics[width=1\textwidth]{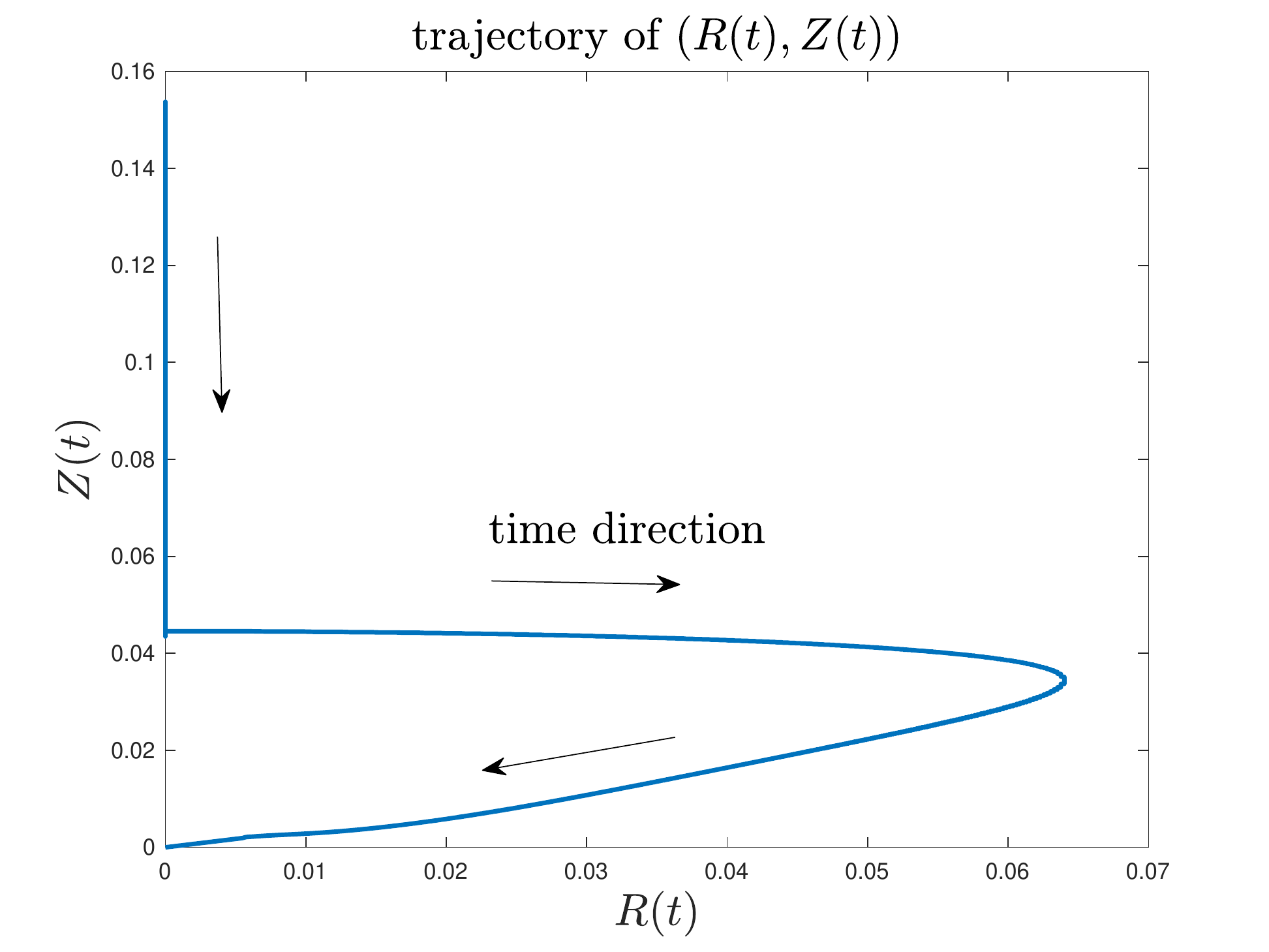}
        \caption{Trajectory $(R(t),Z(t))$}
    \end{subfigure}
    \begin{subfigure}[b]{0.38\textwidth}
        \centering
        \includegraphics[width=1\textwidth]{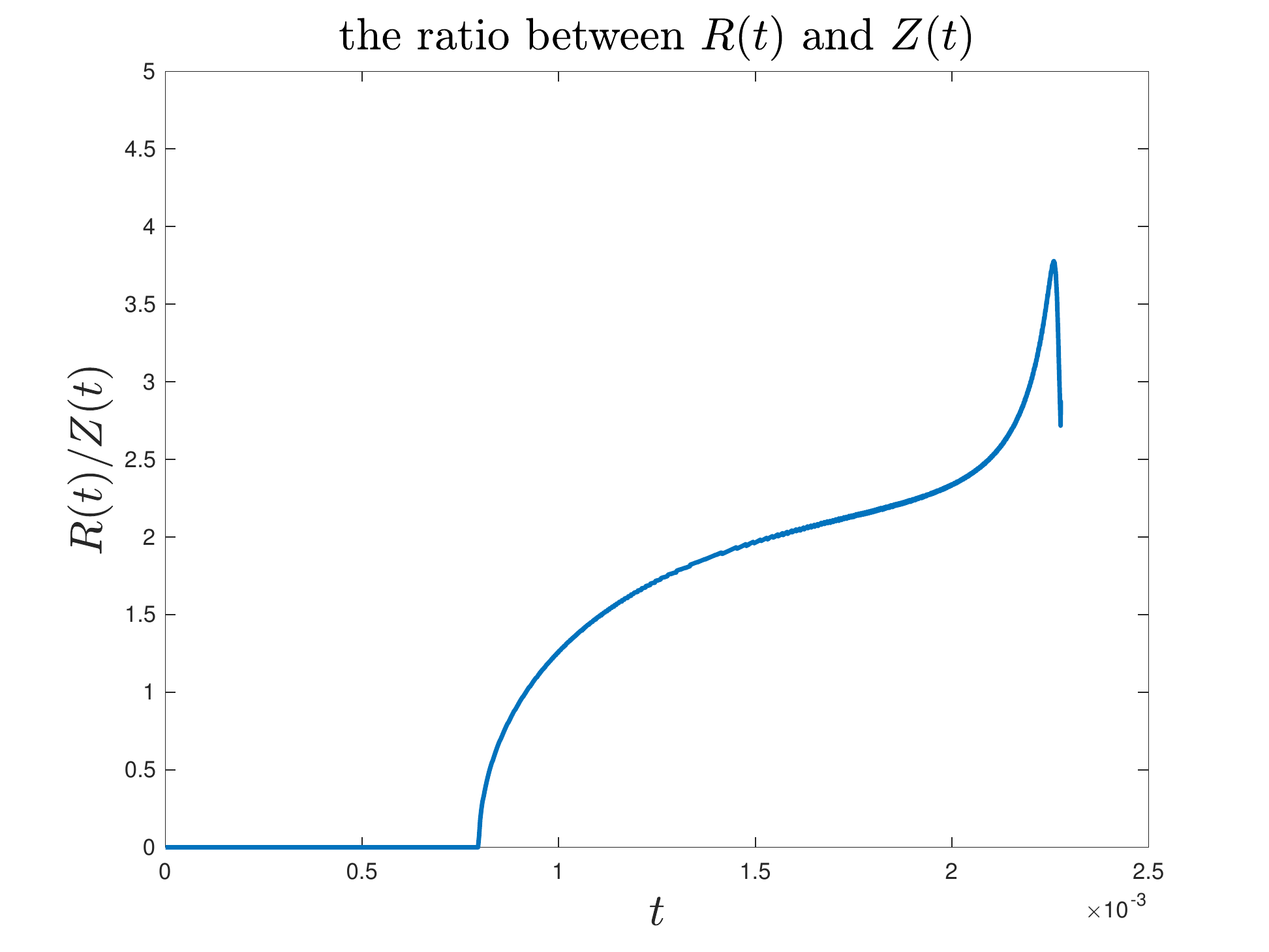}
        \caption{The ratio $R(t)/Z(t)$}
    \end{subfigure}
    \vspace{0.1in}
  \begin{subfigure}[b]{0.38\textwidth}
        \centering
        \includegraphics[width=1\textwidth]{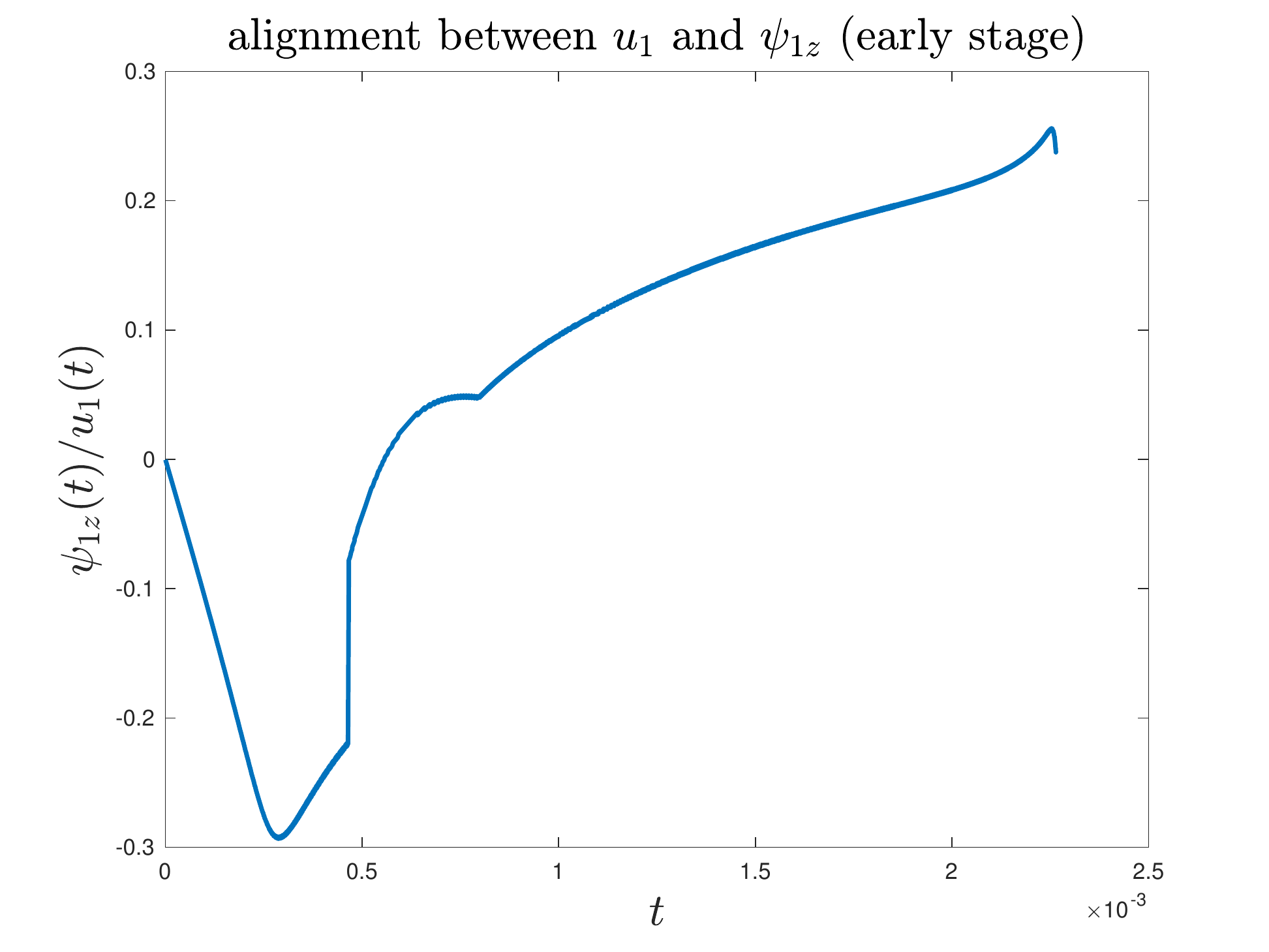}
        \caption{alignment $\psi_{1z}/u_1$ early stage}
    \end{subfigure}
    \begin{subfigure}[b]{0.38\textwidth}
        \centering
        \includegraphics[width=1\textwidth]{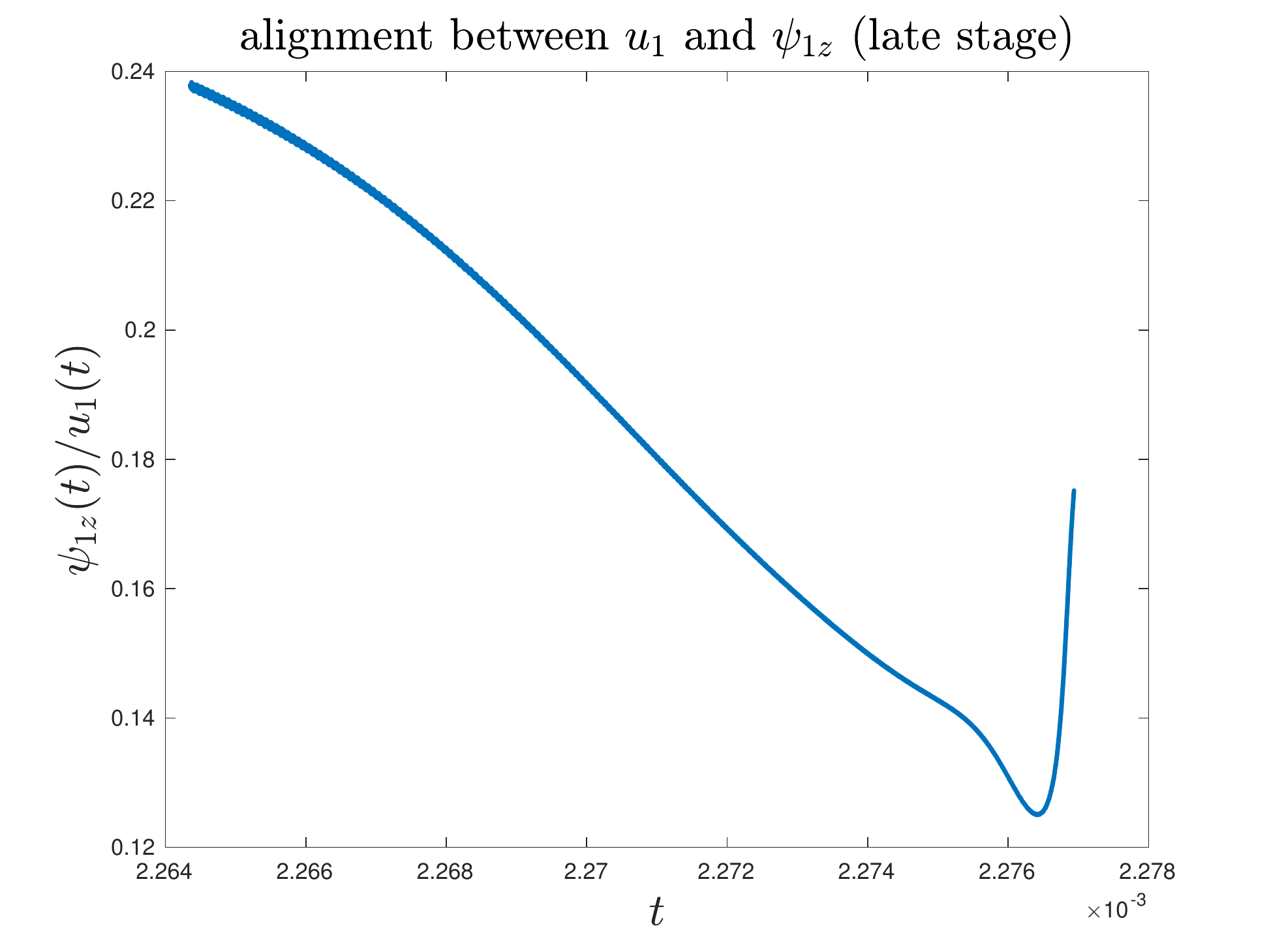}   
          \caption{alignment $\psi_{1z}/u_1$ late stage}
    \end{subfigure}
    \caption[Trajectory]{(a) the trajectory of $(R(t),Z(t))$, (b) the ratio $R(t)/Z(t)$ as a function of time ($0 \leq t \leq 0.002276938$.  (c) the alignment between $\psi_{1z}$ and $u_1$ at $(R(t),Z(t))$ in the early stage with $0 \leq t \leq 0.002264353$. (d) the late stage with $ 0.002264353 \leq t \leq 0.002276938$.}  
     \label{fig:trajectory}
        \vspace{-0.05in}
\end{figure}

Figure \ref{fig:trajectory} (first column) shows the trajectory of the maximum location $(R(t),Z(t))$ of $u_1(t,r,z)$. In the early time, the maximum of $u_1$ lies on the symmetry axis and travels downward along the symmetry axis.  There is a negative alignment between $\psi_{1z}$ and $u_1$. After a short time, $(R(t),Z(t))$ moves away from the symmetry axis almost horizontally and the alignment continues to be negative for a while. Then $(R(t),Z(t))$ begins to turn around and propagates toward the origin $(r,z)=(0,0)$. From this time on, we have a positive alignment between $\psi_{1z}$ and $u_1$ and the alignment becomes stronger for some period of time, see Figure \ref{fig:trajectory} (c). As shown in the Figure \ref{fig:trajectory} (b), the ratio $R(t)/Z(t)$ has a modest growth during this time but is still of $O(1)$, indicating a one-scale solution structure. Then we see a phase transition. The ratio $R(t)/Z(t)$ begins to drop and the alignment begins to decrease. During this second phase, the cross section of $u_1$ develops a sharp front and $\psi_{1z}(t,r,Z(t))$ develops a sharp drop along the $r$-direction just in front of the sharp front. This explains the sharp drop in the alignment. However, the alignment eventually becomes stabilized and begins to grow again in the final stage of our computation. We will revisit this point in Section \ref{sec:mechanism}.

\subsubsection{Rapid growth}\label{sec:rapid_growth_euler} 
\begin{figure}[!ht]
\centering
    \includegraphics[width=0.32\textwidth]{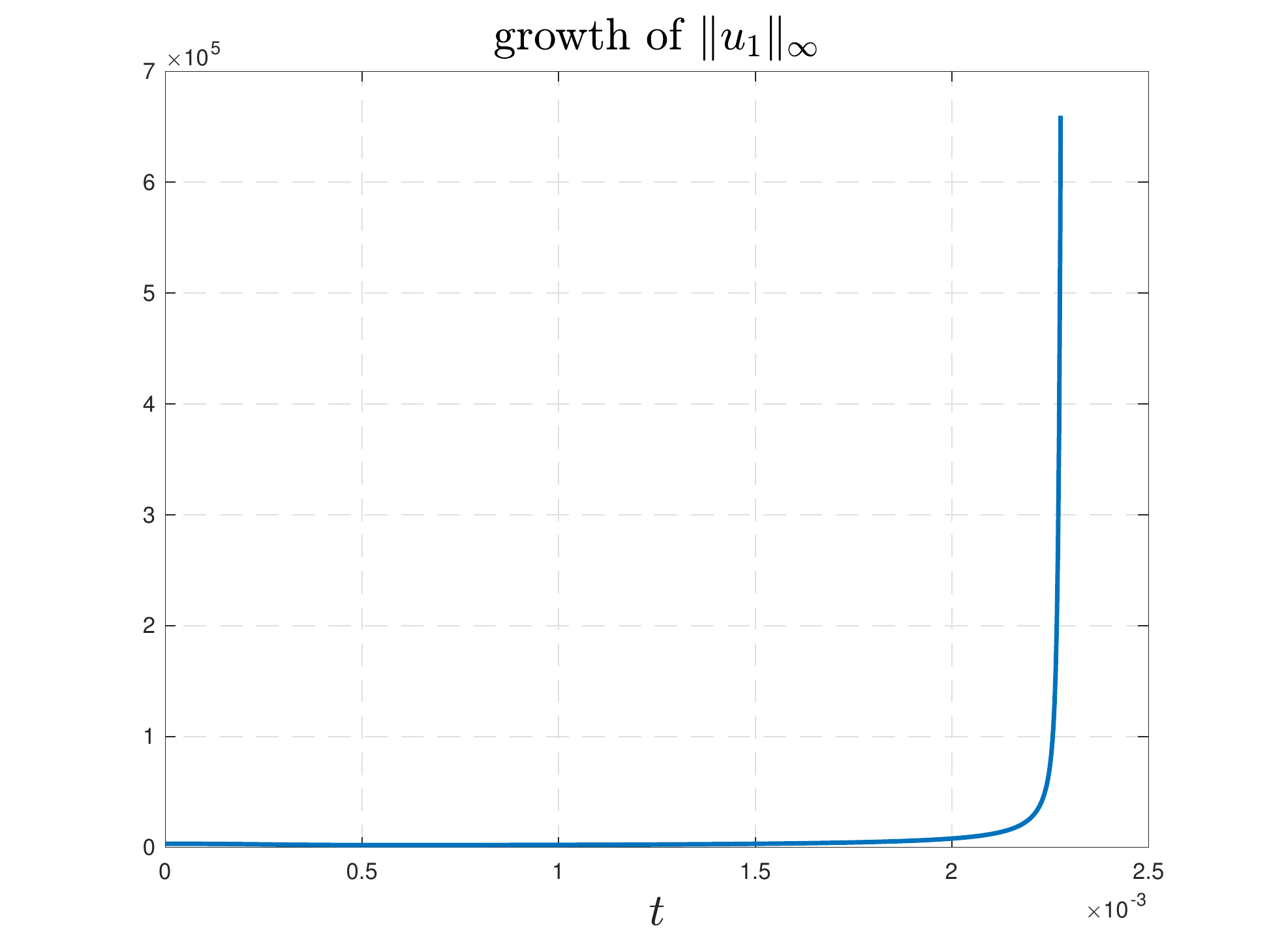}
    \includegraphics[width=0.32\textwidth]{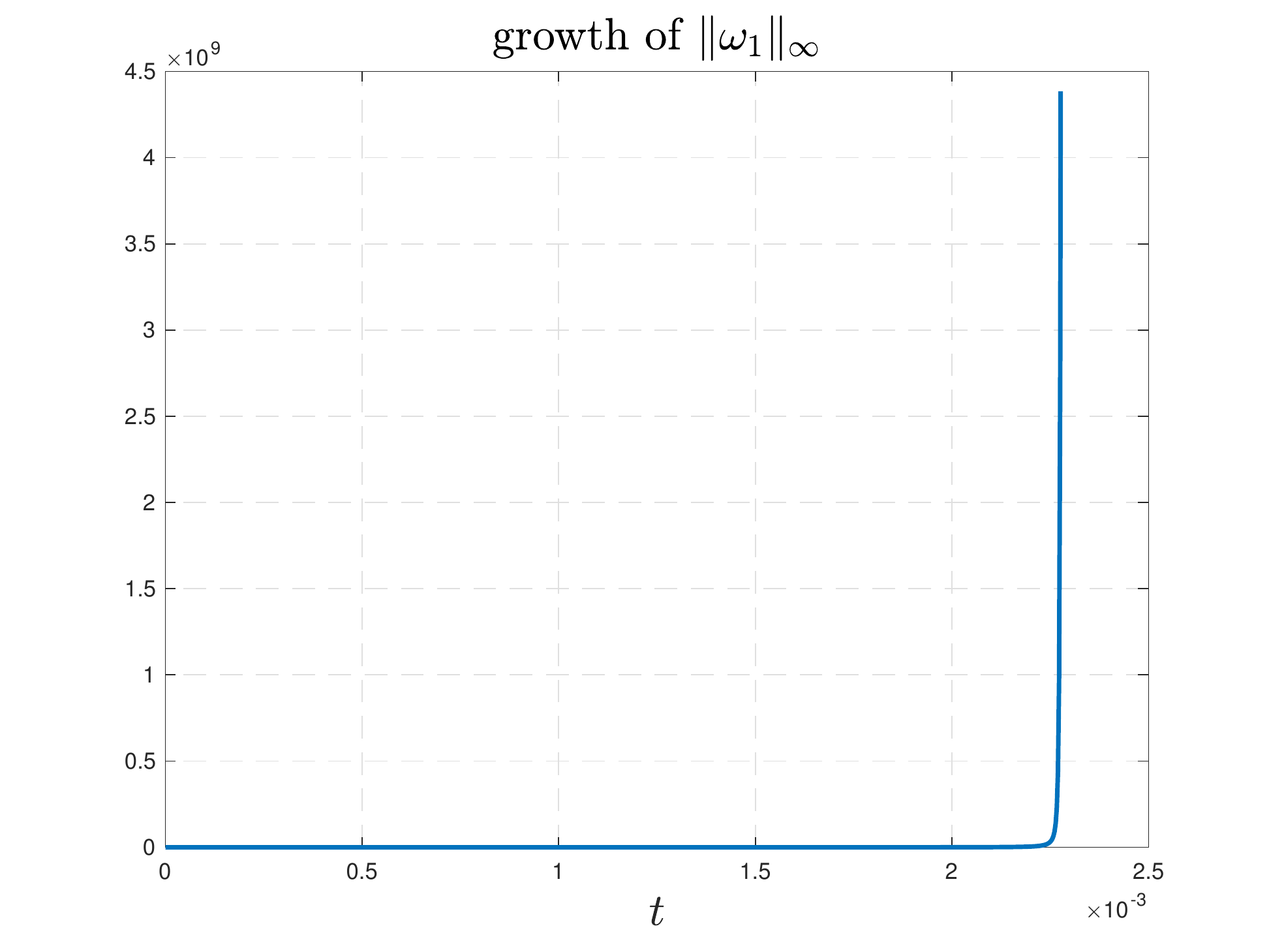} 
    \includegraphics[width=0.32\textwidth]{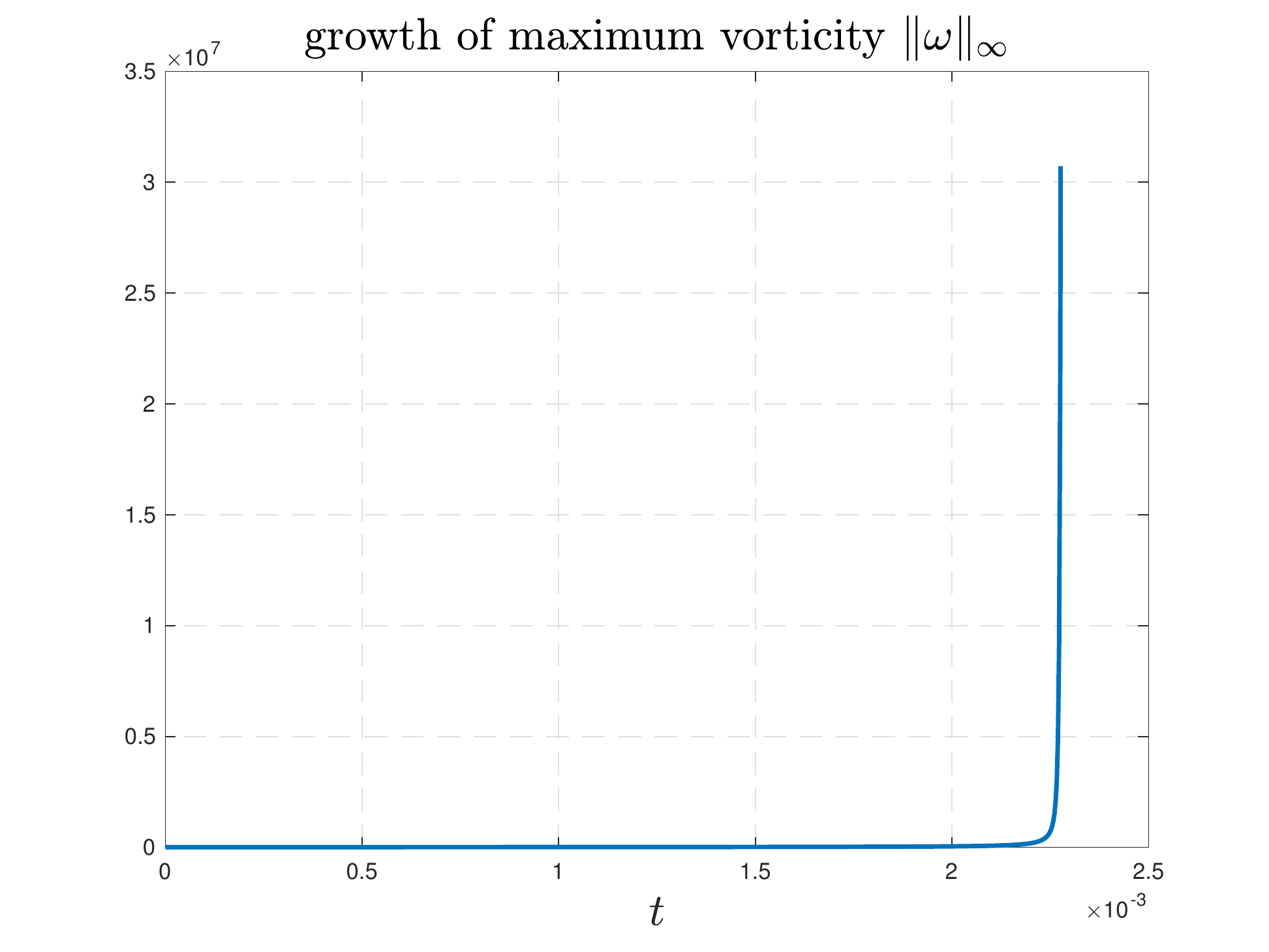}
    \includegraphics[width=0.32\textwidth]{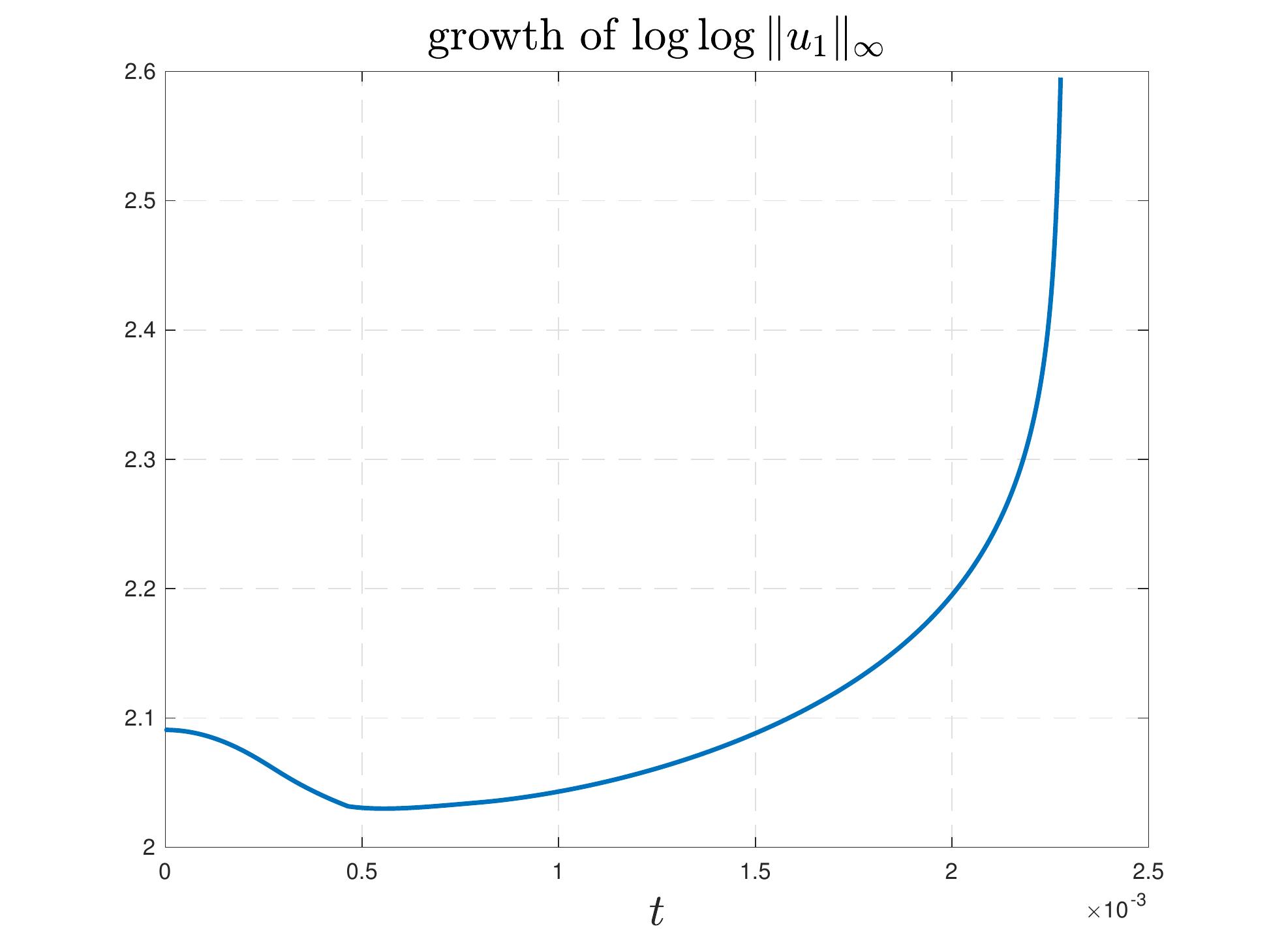}
    \includegraphics[width=0.32\textwidth]{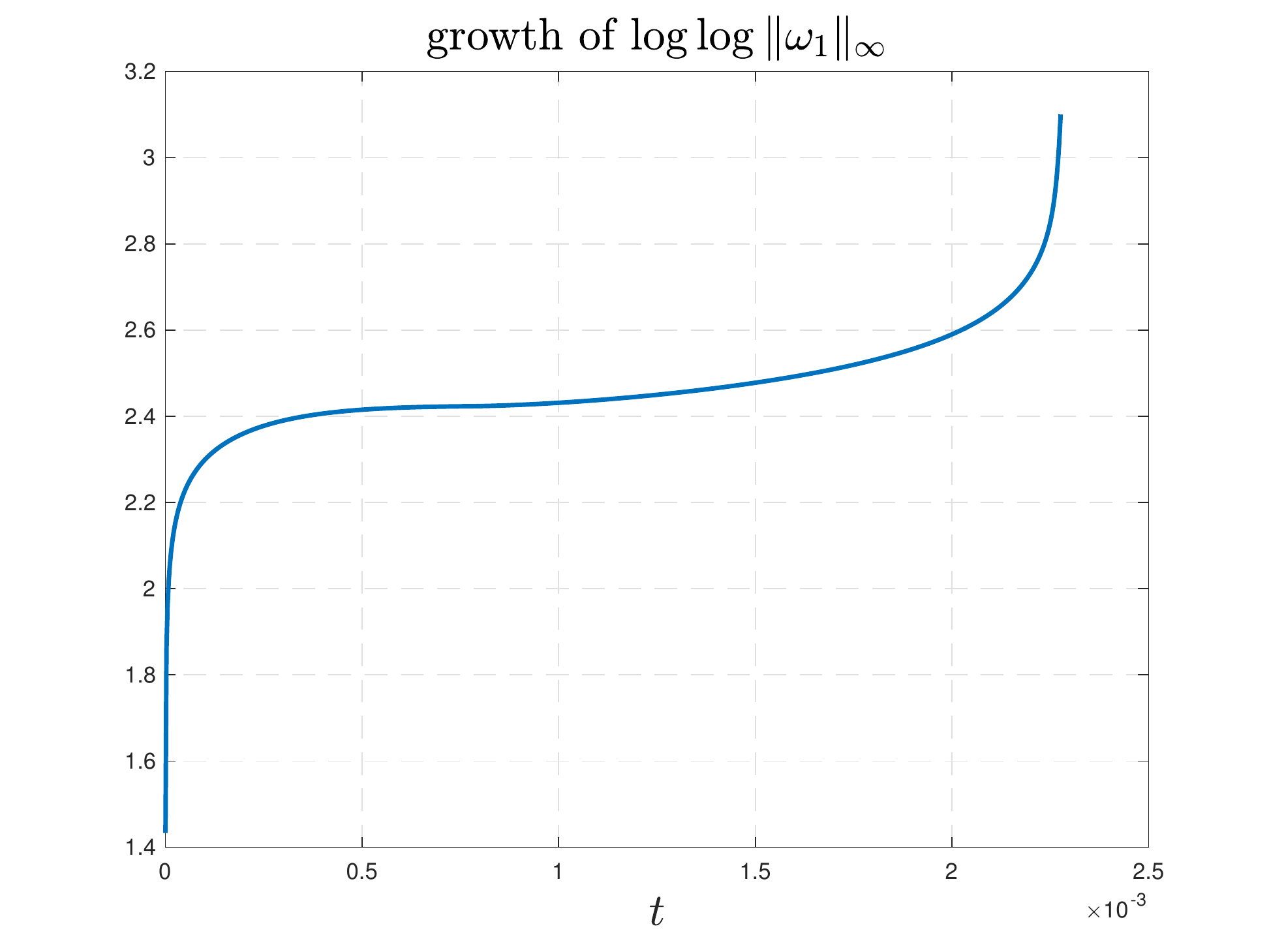}
    \includegraphics[width=0.32\textwidth]{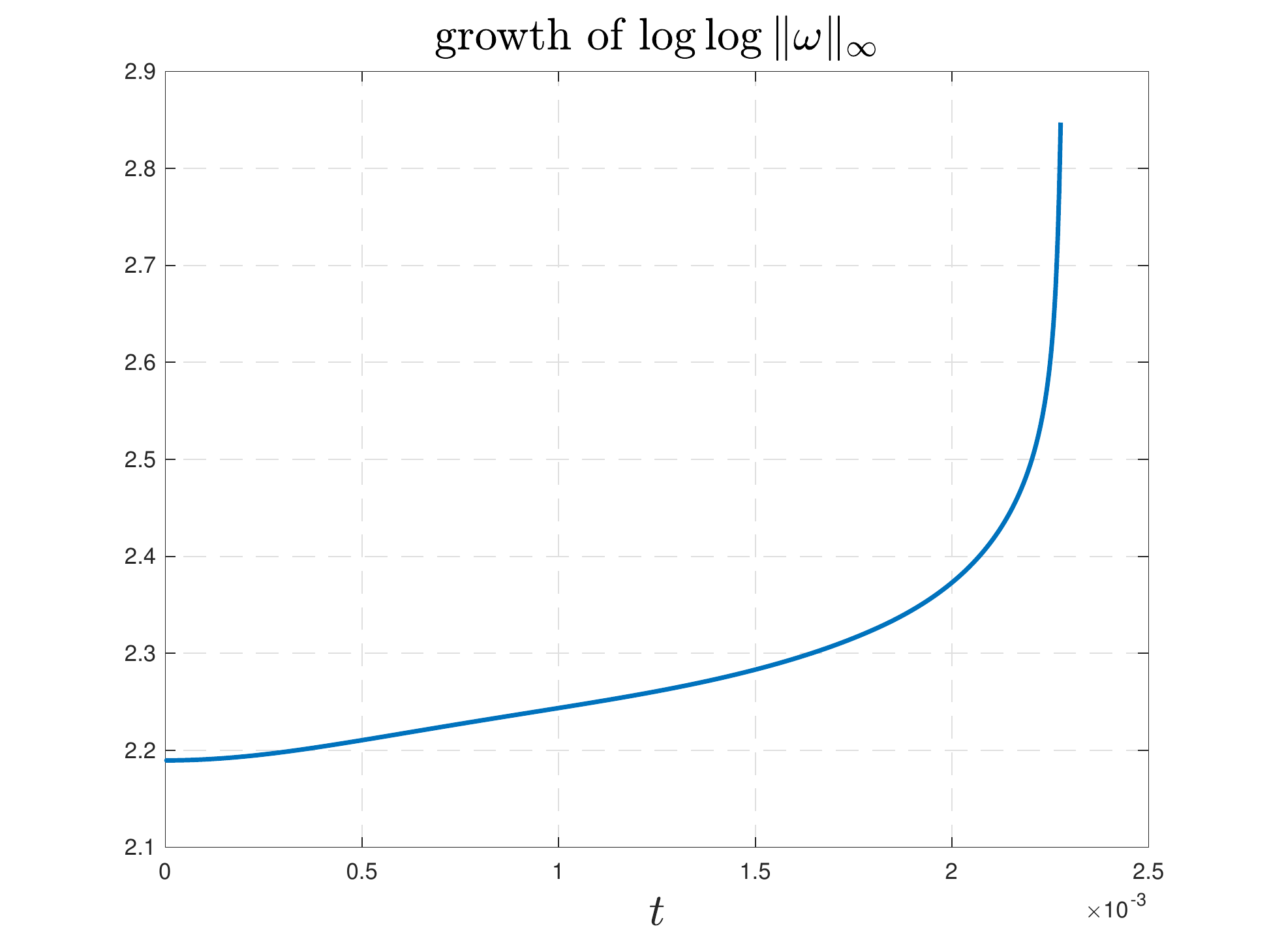}
    \caption[Rapid growth]{First row: the growth of $\|u_1\|_{L^\infty}$, $\|\om_1\|_{L^\infty}$ and $\|\vom\|_{L^\infty}$ as functions of time. Second row: $\log\log\|u_1\|_{L^\infty}$, $\log\log\|\om_1\|_{L^\infty}$ and $\log\log\|\vom\|_{L^\infty}$.} 
    \label{fig:rapid_growth}
\end{figure}

In this subsection, we study the rapid growth of the solution. First, we define the vorticity vector:
\[\vom = (\om^\theta,\om^r,\om^z)^T = (r\om_1\,,\,-r u_{1,z}\,,\,2u_1+ru_{1,r})^T\]
  and its magnitude
\[|\vom| = \sqrt{(\om^\theta)^2+(\om^r)^2 + (\om^z)^2}.\]
In Figure \ref{fig:rapid_growth}, we plot $\|u_1\|_{L^\infty},\|\om_1\|_{L^\infty}$ and $\|\vom\|_{L^\infty}$ as functions of time. 
We can see that these variables grow rapidly in time. Moreover, from the second row in Figure \ref{fig:rapid_growth}. we conclude that the solution grows much faster than a double-exponential rate. 

\begin{figure}[!ht]
\centering
    \includegraphics[width=0.38\textwidth]{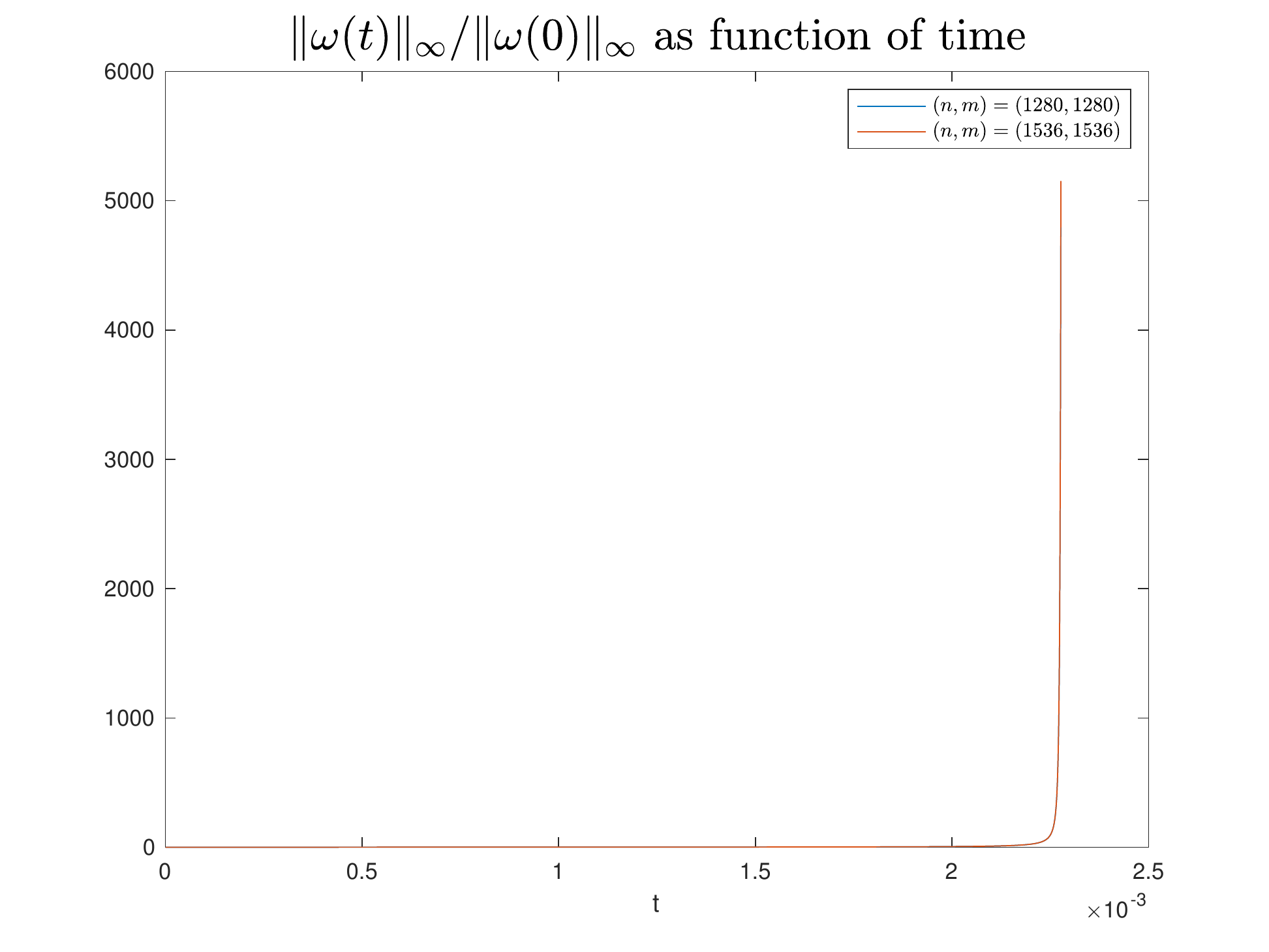}
    \includegraphics[width=0.38\textwidth]{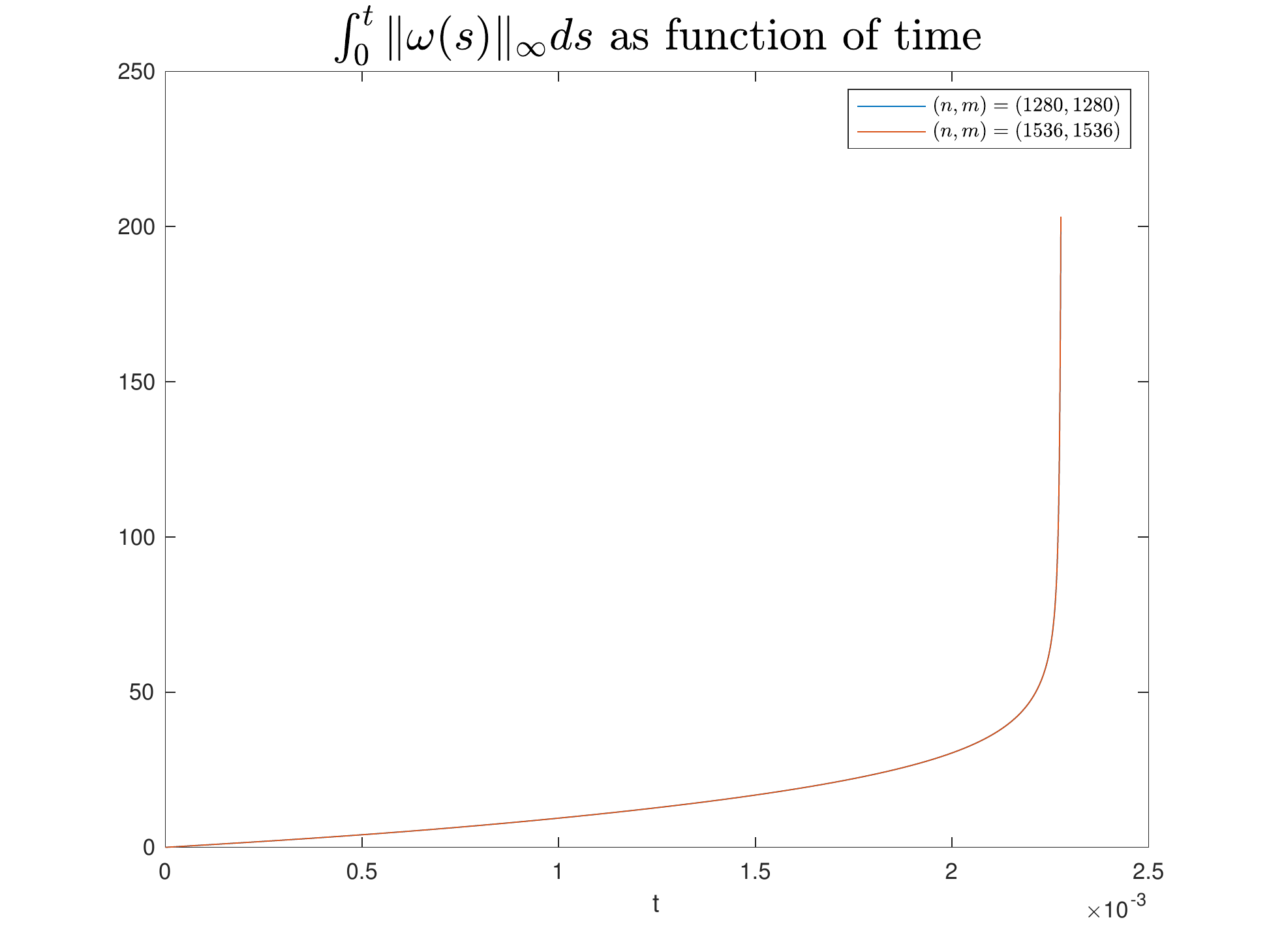} 
    \caption[Rapid growth]{ Left plot: the amplification of maximum vorticity relative to its initial maximum vorticity, $\|\vom (t)\|_{L^\infty}/\|\vom (0)\|_{L^\infty}$ as a function of time. Right plot: the time integral of maximum vorticity,  $\int_0^t \|\vom (s)\|_{L^\infty}ds$ as a function of time. Two resolutions with $(n_1,n_2)=(1280,1280)$ and $(n_1,n_2)=(1536,1536)$ are used. They are almost indistinguishable. The final time instant is $t=0.002276938$.} 
    \label{fig:rapid_growth2}
\end{figure}

The rapid growth of the maximum vorticity $\|\vom\|_{L^\infty}$ is an important indicator of a finite time singularity. The well-known Beale--Kato--Majda criterion \cite{beale1984remarks} states that the solution to the $3$D Euler equations ceases to exist in some regularity class $H^s$ (for $s\geq 3$) at some finite time $T_*$ if and only if
\begin{equation}\label{eq:BKM}
\int_0^{T}\|\vom(t)\|_{L^\infty}\idiff t = +\infty. 
\end{equation}
We will demonstrate in Section~\ref{sec:scaling_study_euler} that the growth of $\|\vom\|_{L^\infty}$ has a very good fitting of the form
\[\|\vom(t)\|_{L^\infty} \approx (T-t)^{-1}\;.\]
This is only a qualitative fitting.
If this qualitative scaling relationship holds, it implies that the solution would develop a potential finite time singularity. 
Moreover, we also compute the relative growth of maximum vorticity $\|\vom (t)\|_{L^\infty}/\|\vom (0)\|_{L^\infty}$ and $\int_0^{t}\|\vom(s)\|_{L^\infty}\idiff s$ in Figure \ref{fig:rapid_growth2}. The final time of this computation is at $t=0.002276938$. We use two resolutions to compute this quantity ($1280\times1280$ vs $1536\times1536$). The two solutions are almost indistinguishable. We observe that the maximum vorticity has grown more than $5000$ compared with the maximum of the initial vorticity vector. The rapid growth of $\int_0^{t}\|\vom(s)\|_{L^\infty}\idiff s$ provides further evidence that the $3$D Euler equations seem to develop a finite time singularity.

\subsubsection{Velocity field} In this subsection, we investigate the velocity field. We first study the $3$D velocity field $\vu = u^r\vct{e}_r + u^\theta \vct{e}_\theta + u^z \vct{e}_z$ (denoted by $(u^r,u^\theta,u^z)$) by studying the induced streamlines. A streamline $\{\Phi(s;X_0)\}_{s\geq0}\subset \mathbb{R}^3$ is determined by the background velocity $\vu$ and the initial point $X_0 = (x_0,y_0,z_0)^T$ through the initial value problem
\[\frac{\partial}{\partial s}\Phi(s;X_0) = \vu(\Phi(s;X_0)),\quad s\geq 0;\quad  \Phi(0;X_0) = X_0.\]

We will generate different streamlines with different initial points $X_0 = (r_0\cos(2\pi \theta),r_0\sin(2\pi\theta),z_0)^T$. Due to the  axisymmetry of the velocity field $\vu$,  it is sufficient to prescribe $(r_0,z_0)$ as a starting point. We will plot the streamlines with different angular parameter $\theta$ to illustrate the rotational symmetry of the streamlines.  

\subsubsection{A tornado singularity} 
\begin{figure}[!ht]
\centering
\vspace{-0.7in}
    \begin{subfigure}[b]{0.40\textwidth}
        \centering
        \includegraphics[width=1\textwidth]{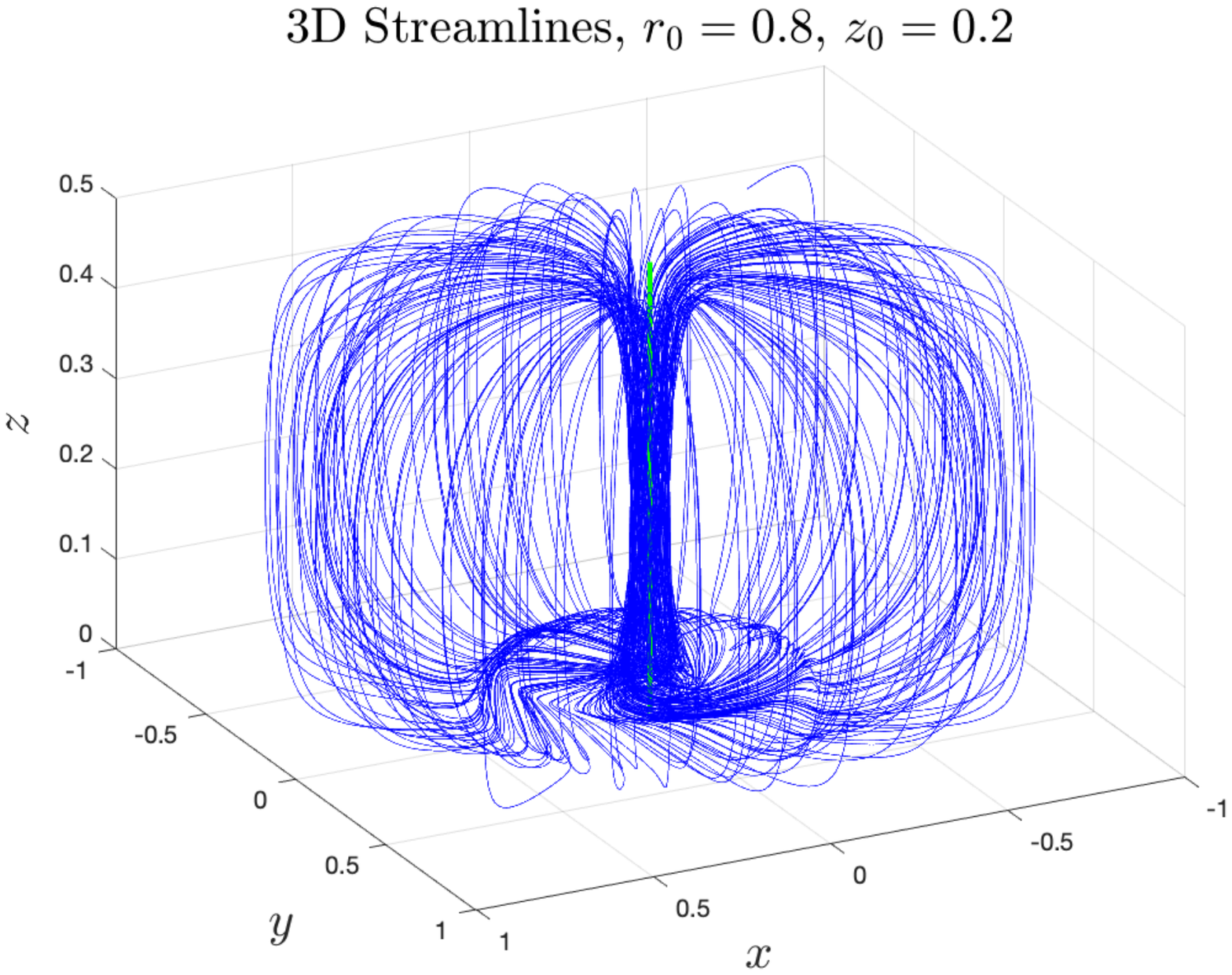}
        \vspace{-0.9in}
        \caption{$r_0 = 0.8$, $z_0 = 0.2$}
    \end{subfigure}
    \begin{subfigure}[b]{0.40\textwidth}
        \centering
        \includegraphics[width=1\textwidth]{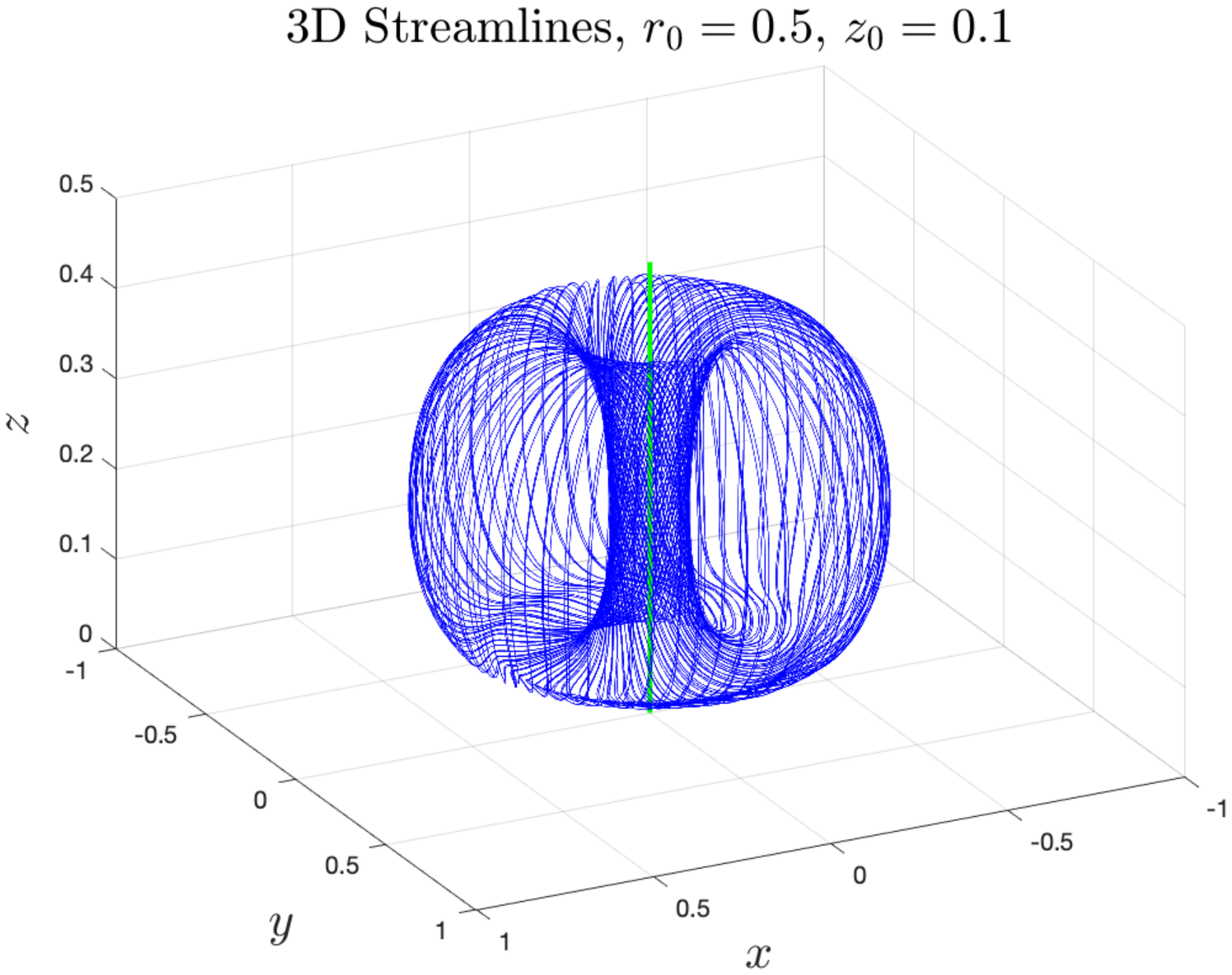}
        \vspace{-0.9in}
        \caption{$r_0 = 0.5$, $z_0 = 0.1$}
    \end{subfigure}
    \begin{subfigure}[c]{0.40\textwidth}
        \centering
        \vspace{-0.7in}
        \includegraphics[width=1\textwidth]{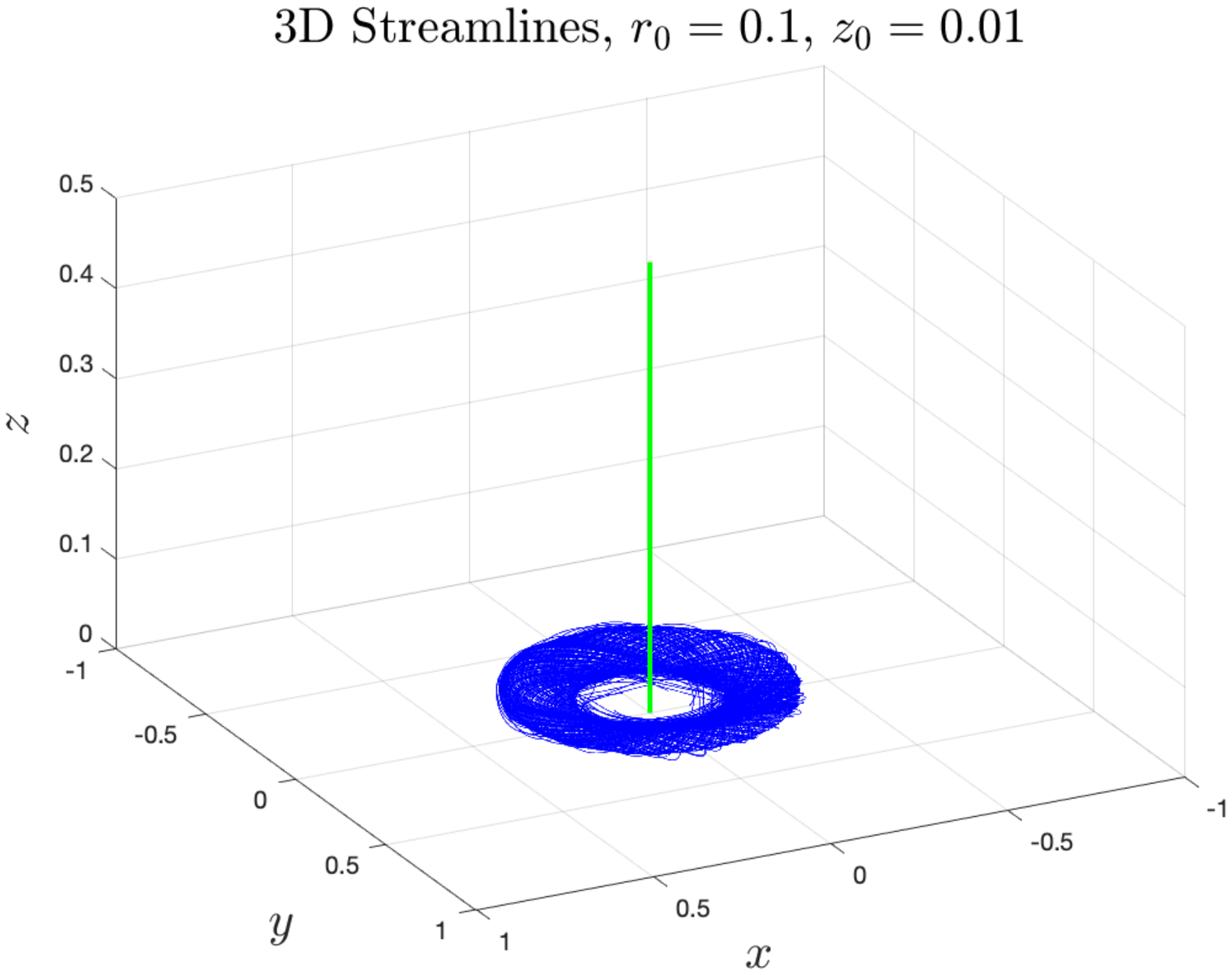}
        \vspace{-0.9in}
        \caption{$r_0 = 0.1$, $z_0 = 0.01$}
    \end{subfigure}
    \begin{subfigure}[d]{0.40\textwidth}
        \centering
         \vspace{-0.7in}
        \includegraphics[width=1\textwidth]{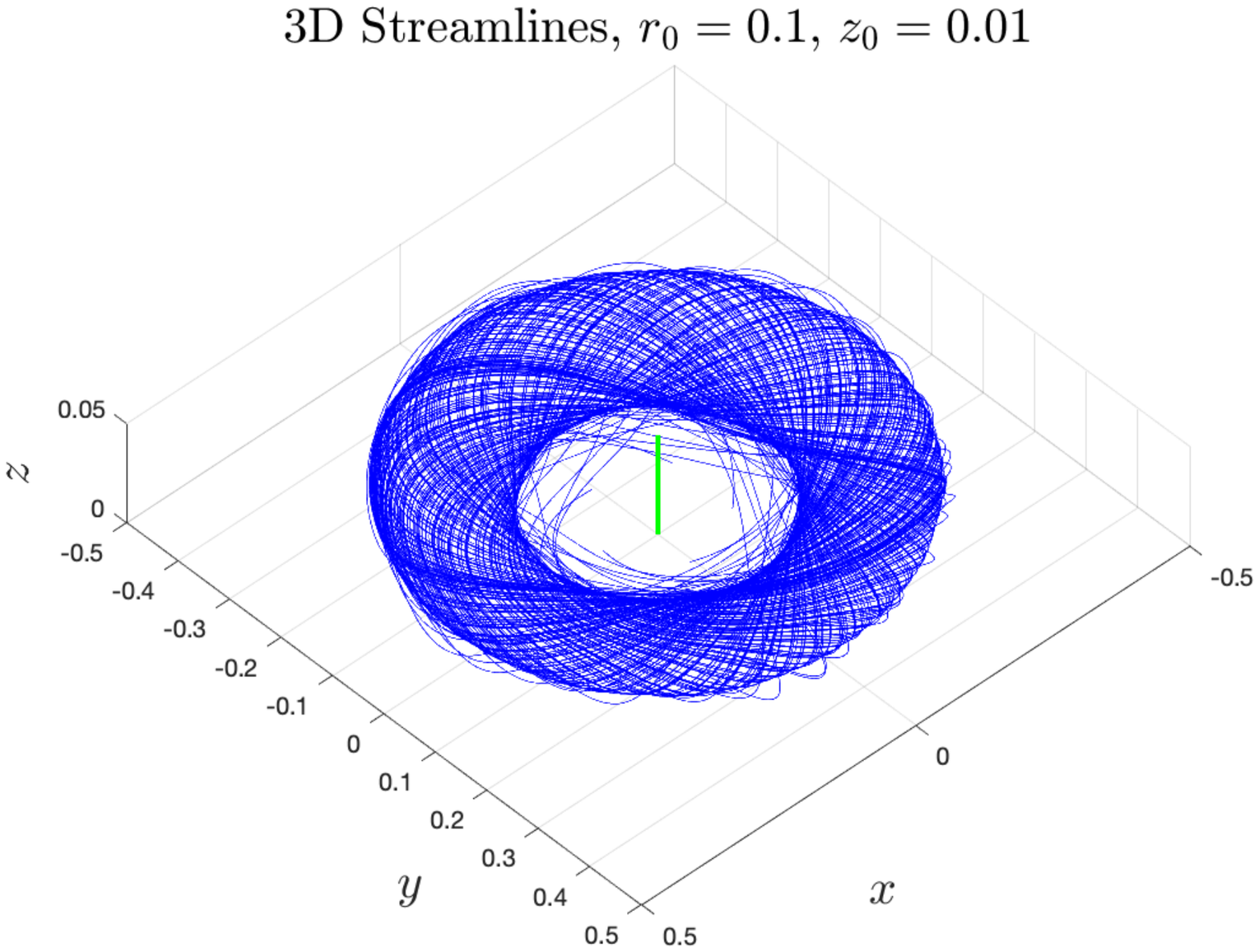}
        \vspace{-0.9in}
        \caption{same as (c), a zoom view}
    \end{subfigure}
    \caption[Local streamline]{The streamlines of $(u^r(t),u^\theta(t),u^z(t))$ at time $t=0.002276480$ with initial points given by (a) $(r_0,z_0) = (0.8,0.2)$, (b) $(r_0,z_0) = (0.5,0.1)$, (c) $(r_0,z_0) = (0.1,0.01)$ ($3$D view), (d) $(r_0,z_0)$ is the same as (c), a zoom view. The green pole is the symmetry axis $r=0$.}  
     \label{fig:streamline_3D_global}
        \vspace{-0.05in}
\end{figure}

In Figure~\ref{fig:streamline_3D_global}, we plot  the streamlines induced by the velocity field $\vu(t)$ at $t = 0.00227648$ in a macroscopic scale (the whole cylinder domain $\mathcal{D}_1\times [0,2\pi]$) for different initial points with (a) $(r_0,z_0) = (0.8,0.2)$, (b) $(r_0,z_0) = (0.5,0.1)$, and (c)-(d) $(r_0,z_0) = (0.1,0.01)$ . We can see that the velocity field generates a tornado like structure spinning around the symmetry axis (the green pole). In  Figure~\ref{fig:streamline_3D_global}(a), we observe that the streamlines first travel toward the symmetry axis, then move upward toward $z=1/2$, and at last turn outward away from the symmetry axis. If we move $z_0$ a bit lower toward $z=0$, we observe in Figure~\ref{fig:streamline_3D_global}(b) that the streamlines first approach the symmetry axis, move upward along the symmetry axis, and then turn outward as they approach $z=1/2$, and move downward. At the same time, they slowly circle around the symmmetry axis. On the other hand, if the initial point is very close to $z=0$ as in Figure~\ref{fig:streamline_3D_global}(c)-(d), the streamlines will just spin around the symmetry axis and stay near $z=0$ instead of moving upward.

\begin{figure}[!ht]
\centering
\vspace{-0.7in}
    \begin{subfigure}[b]{0.40\textwidth}
        \centering
        \includegraphics[width=1\textwidth]{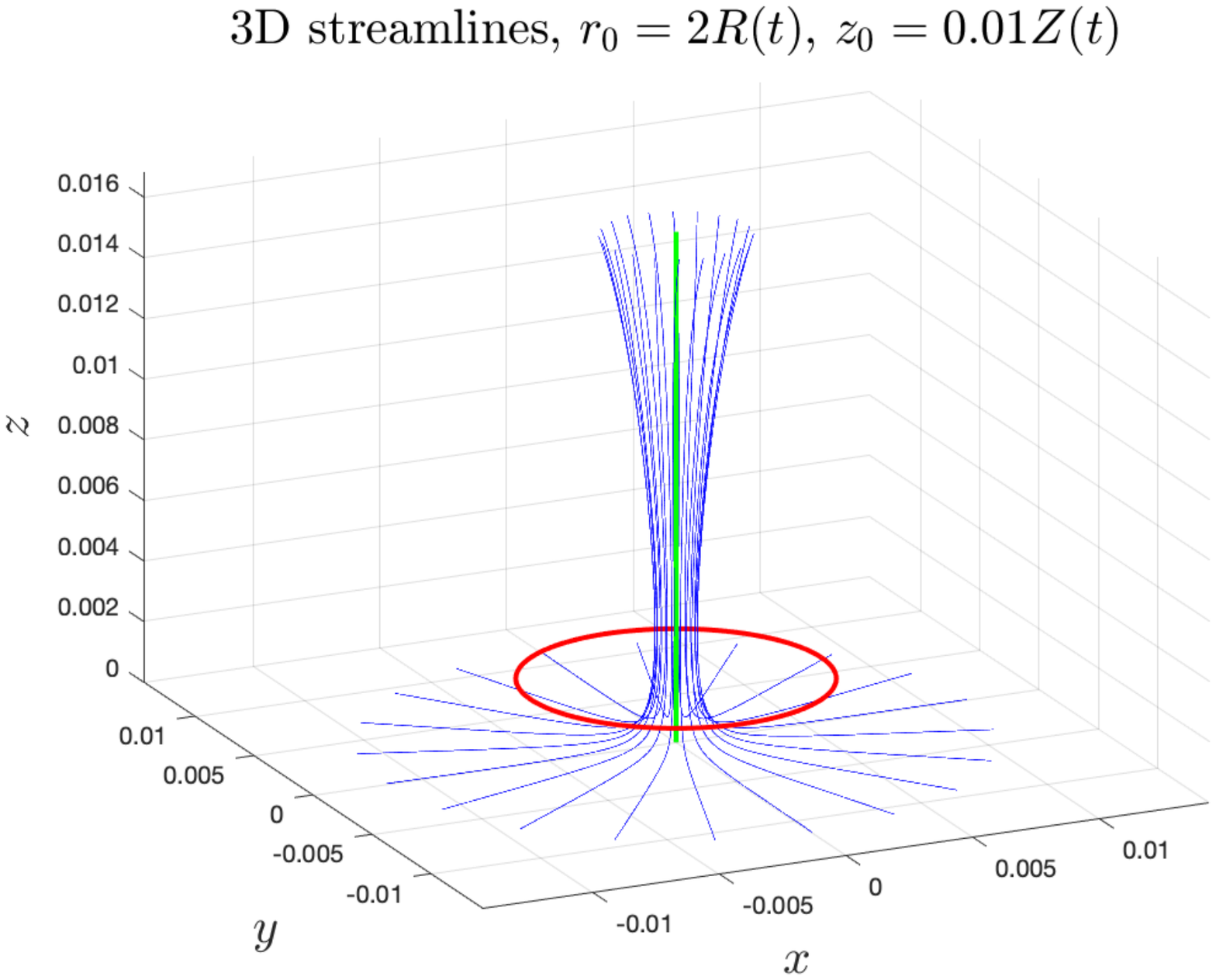}
        \vspace{-0.9in}
        \caption{$r_0 = 2R(t)$, $z_0 = 0.01Z(t)$}
    \end{subfigure}
    \begin{subfigure}[b]{0.40\textwidth}
        \centering
        \includegraphics[width=1\textwidth]{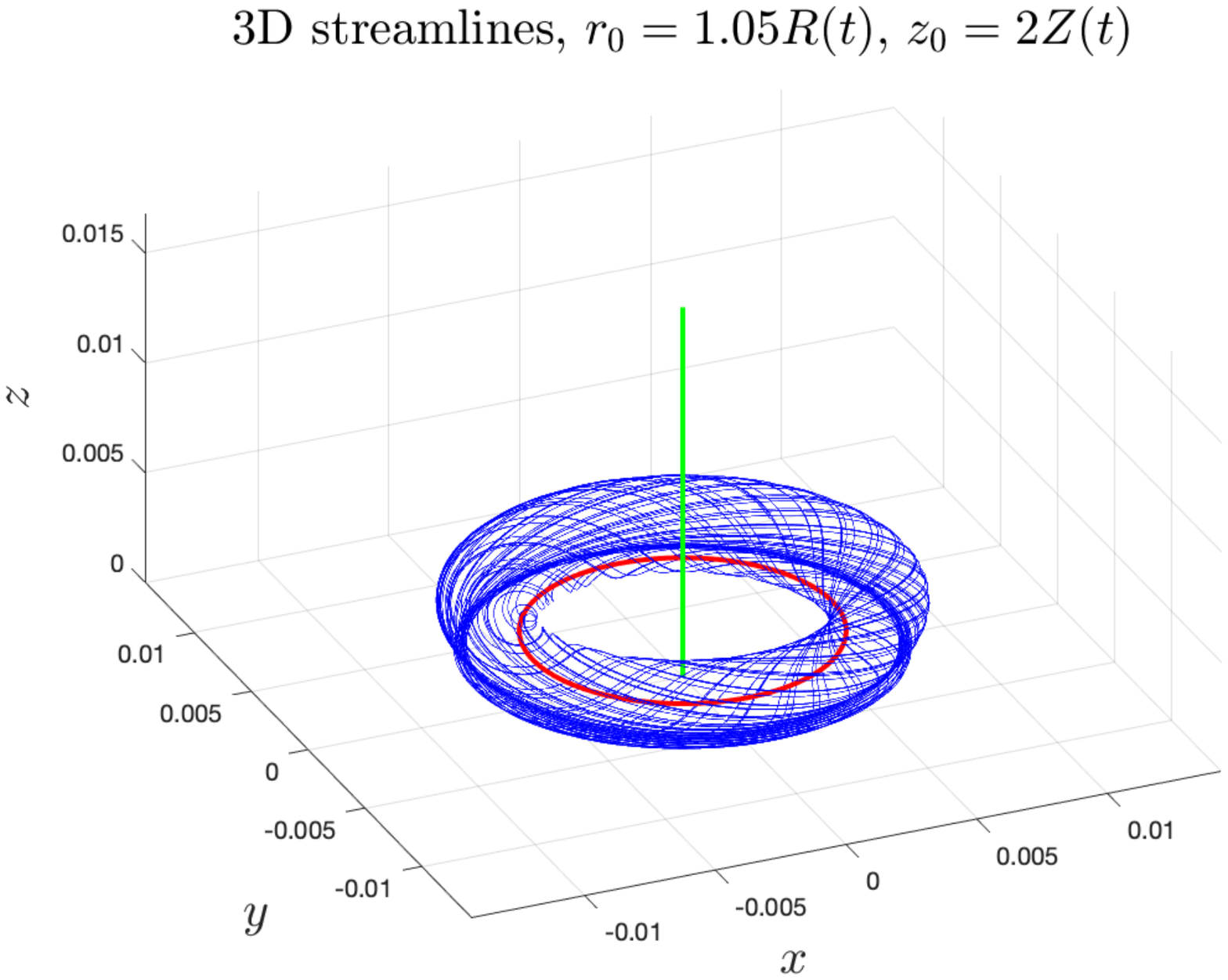}
        \vspace{-0.9in}
        \caption{$r_0 = 1.05R(t)$, $z_0 = 2Z(t)$}
    \end{subfigure}
    \begin{subfigure}[c]{0.40\textwidth}
        \centering
        \vspace{-0.7in}
        \includegraphics[width=1\textwidth]{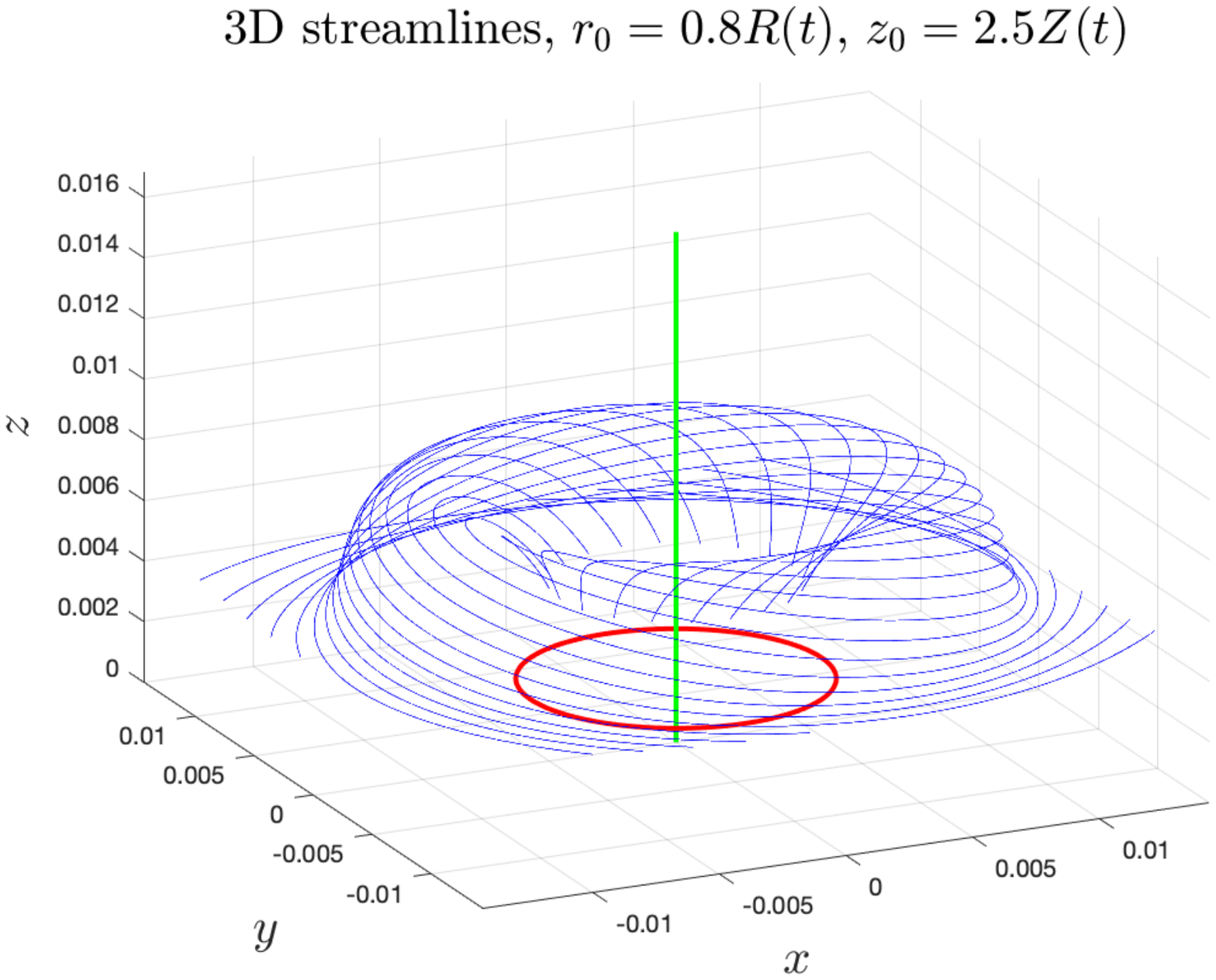}
        \vspace{-0.9in}
        \caption{$r_0 = 0.8R(t)$, $z_0 = 2.5Z(t)$, $3$D view}
    \end{subfigure}
    \begin{subfigure}[d]{0.40\textwidth}
        \centering
         \vspace{-0.7in}
        \includegraphics[width=1\textwidth]{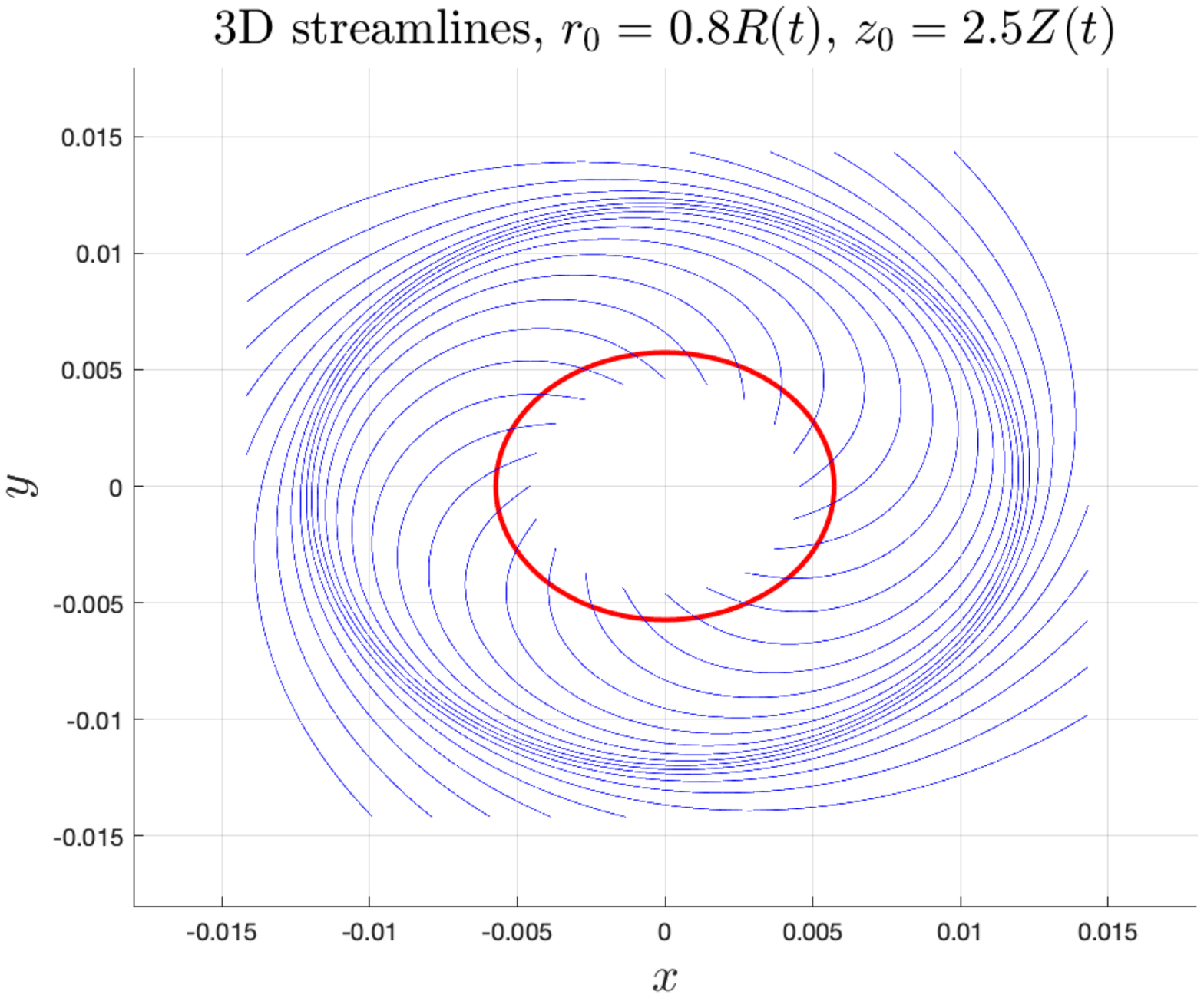}
        \vspace{-0.9in}
        \caption{same as (c), a top view}
    \end{subfigure}
    \caption[Local streamline]{The streamlines of $(u^r(t),u^\theta(t),u^z(t))$ at time $t=0.002276480$ with initial points given by (a) $(r_0,z_0) = (2R(t),0.01Z(t))$, (b) $(r_0,z_0) = (1.05R(t),2Z(t))$, (c) $(r_0,z_0) = (0.8R(t),2.5Z(t))$ ($3$D view), (d) $(r_0,z_0)$ is the same as (c), a top view. $(R(t),Z(t))$ is the maximum location of $u_1(t)$, indicated by the red ring. The green pole is the symmetry axis $r=0$.}  
     \label{fig:streamline_3D_zoomin}
        \vspace{-0.05in}
\end{figure}

Next, we demonstrate the $3$D velocity in a local region near $(R(t),Z(t))$. In Figure~\ref{fig:streamline_3D_zoomin}, we plot the streamlines at time $t=0.00227648$ for different initial points near $(R(t),Z(t))$. The red ring represents the location of $(R(t),Z(t))$, and the green pole is the symmetry axis $r=0$. The $3$ different starting points $(r_0,z_0)$ are given as follows.
\begin{itemize}
\item[(a)] $(r_0,z_0) = (2R(t),0.01Z(t))$. The streamlines start near $z=0$ and below the red ring $(R(t),Z(t))$. They first travel toward the symmetry axis and then travel upward away from $z=0$. There is almost no spinning since $u^\theta = ru_1$ is small in the region near $r=0$.
\item[(b)] $(r_0,z_0) = (1.05R(t),2Z(t))$. The streamlines start some distance above the ring $(R(t),Z(t))$. They get trapped in a local region, oscillating and spinning around the symmetry axis periodically. The spinning is strong since $u^\theta = ru_1$ is quite large near $(r_0,z_0)$.
\item[(c)-(d)] $(r_0,z_0) = (0.8R(t),2.5Z(t))$.  The streamlines start even higher than $z=2Z(t)$, but is inside the ring $(R(t),Z(t))$. They first spin upward and outward, then travel downward away from the blow-up region. We have computed several other starting points and obtained the same qualitative result as long as $z_0$ is larger than $2 Z(t)$. We can take $r_0$  smaller or larger than $R(t)$.
\end{itemize}

\subsubsection{The $2$D flow} 
\begin{figure}[!ht]
\centering
    \begin{subfigure}[b]{0.38\textwidth}
        \centering
        \includegraphics[width=1\textwidth]{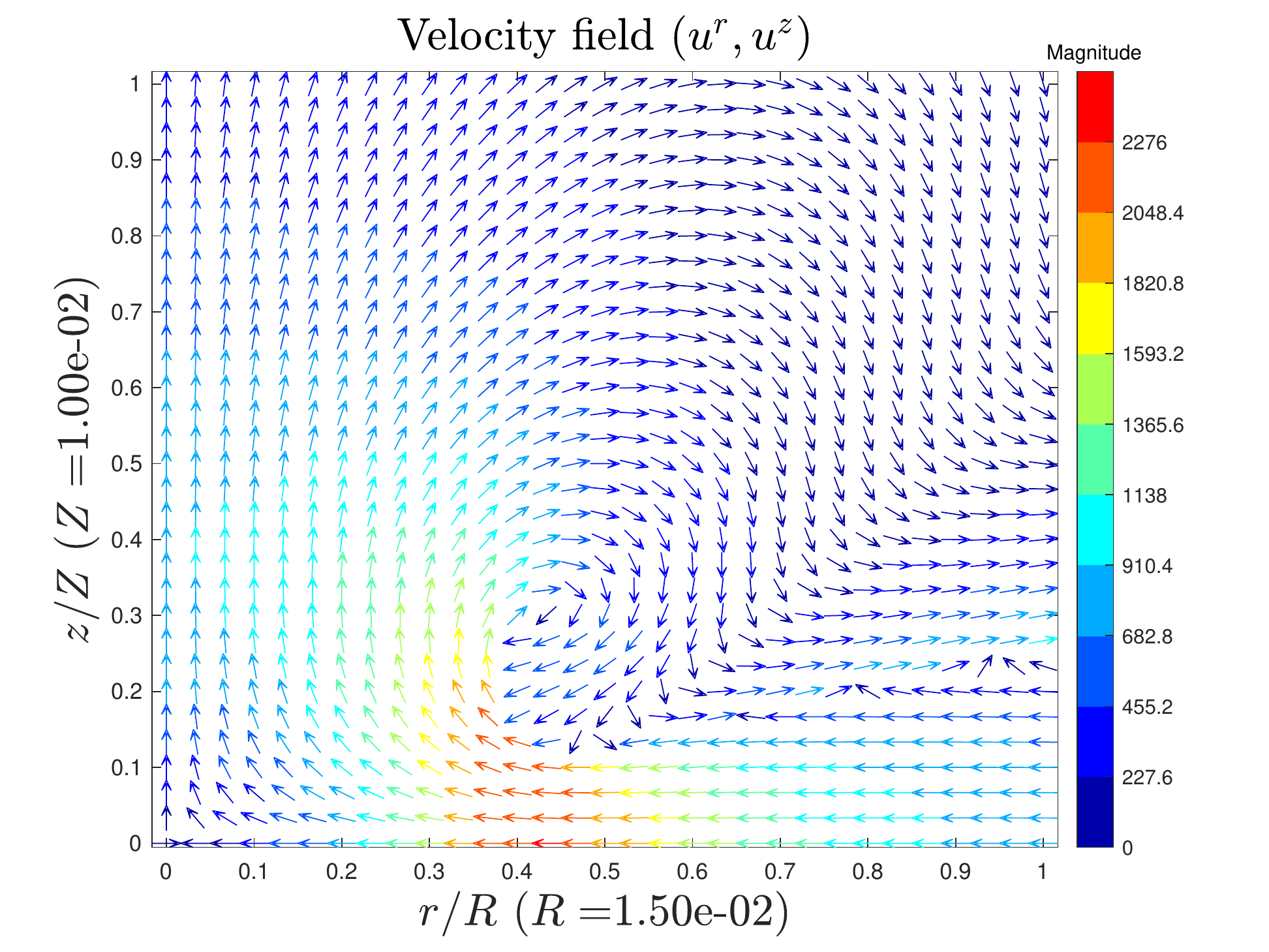}
        \caption{the velocity field $(u^r,u^z)$}
    \end{subfigure}
    \begin{subfigure}[b]{0.38\textwidth}
        \centering
        \includegraphics[width=1\textwidth]{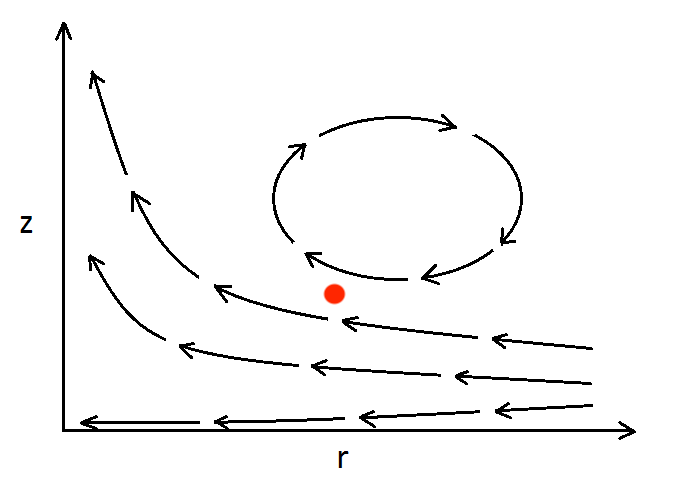} 
        \caption{a schematic}
    \end{subfigure}
    \caption[Global streamline]{(a) The velocity field $(u^r(t),u^z(t))$ near the maximum location $(R(t),Z(t))$ of $u_1(t)$ at $t=0.00227648$. The color corresponds to the amplitude of the velocity field. The size of the domain has been rescaled. (b) A schematic of the vector field near the point $(R(t),Z(t))$.}  
     \label{fig:velocity_field}
        \vspace{-0.05in}
\end{figure}

To understand the phenomena in the blow-up region as shown in Figure~\ref{fig:streamline_3D_zoomin}, we look at the $2$D velocity field $(u^r,u^z)$ in the computational domain $\mathcal{D}_1$. In Figure~\ref{fig:velocity_field}(a), we show the vector field of $(u^r(t),u^z(t))$ in a local microscopic domain $[0,R_b]\times [0,Z_b]$, where $R_b = 0.015$ and $Z_b = 0.01$. The color corresponds to the amplitude of the velocity field. We have rescaled the domain in the figure for better visualization. Figure~\ref{fig:velocity_field}(b) is a schematic for the vector field in Figure~\ref{fig:velocity_field}(a). 

We observe a hyperbolic flow structure. The streamlines below $(R(t),Z(t))$ first travel toward $r=0$ and then move upward away from $z=0$, going around the sharp front near $(R(t),Z(t))$ as if it were an obstacle. As the flow gets close to $r=0$, the strong axial velocity $u^z$ transports $u_1$ from near $z=0$ upward along the $z$ direction. Due to the odd symmetry of $u_1$, the angular velocity $u^\theta = ru^1$ is almost $0$ in the region near $z=0$. As a consequence, the value of $u_1$ becomes very small near the symmetry axis $r=0$. This explains why the streamlines almost do not spin around the symmetry axis in this region, as illustrated in Figure~\ref{fig:streamline_3D_zoomin}(a). 

Moreover, we observe that the velocity field $(u^r(t),u^z(t))$ forms a closed circle right above $(R(t),Z(t))$ as illustrated in Figure~\ref{fig:velocity_field}(b). The corresponding streamline is trapped in the circle region in the $rz$-plane. Due to the large positive value of $u_1$ in this region, the induced streamlines oscillate and spin fast inside a $3$D torus surrounding the symmetry axis, see also Figure~\ref{fig:streamline_3D_zoomin}(b). This local circle structure of the $2$D velocity field plays an essential role in stabilizing the blow-up process since it retains a large portion of $u_1,\om_1$ within this circle structure instead of being pushed upward. Moreover, $u_1$ and $\omega_1$ within this circle region develop some favorable alignment that leads to strong nonlinear alignment between $u_1$ and $\psi_{1z}$ and a sustainable blow-up mechanism. We will revisit this point in the next subsection when we study the mechanism for the potential blow-up of the $3$D Euler equations.

\begin{figure}[!ht]
\centering
    \includegraphics[width=0.38\textwidth]{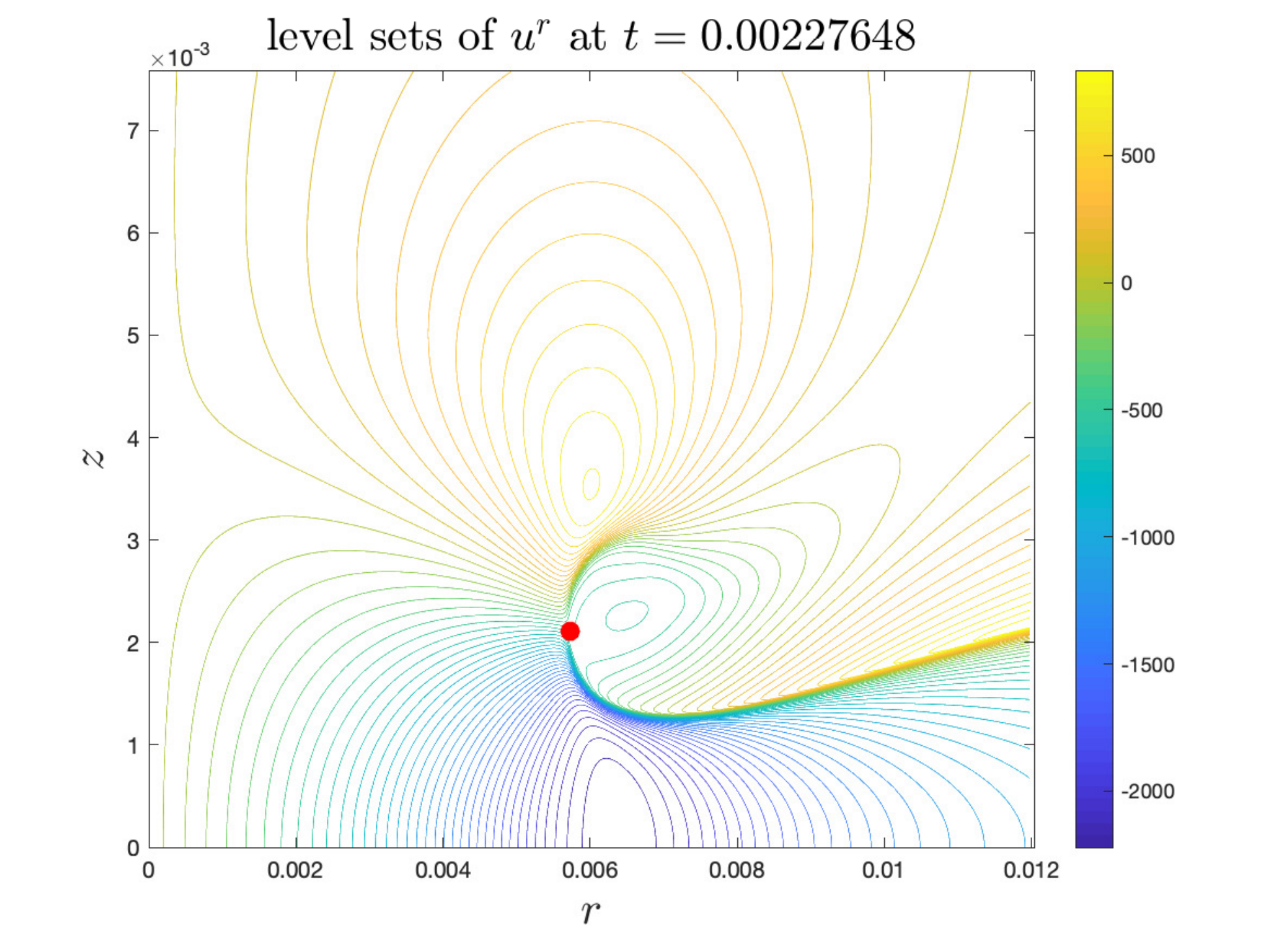}
    \includegraphics[width=0.38\textwidth]{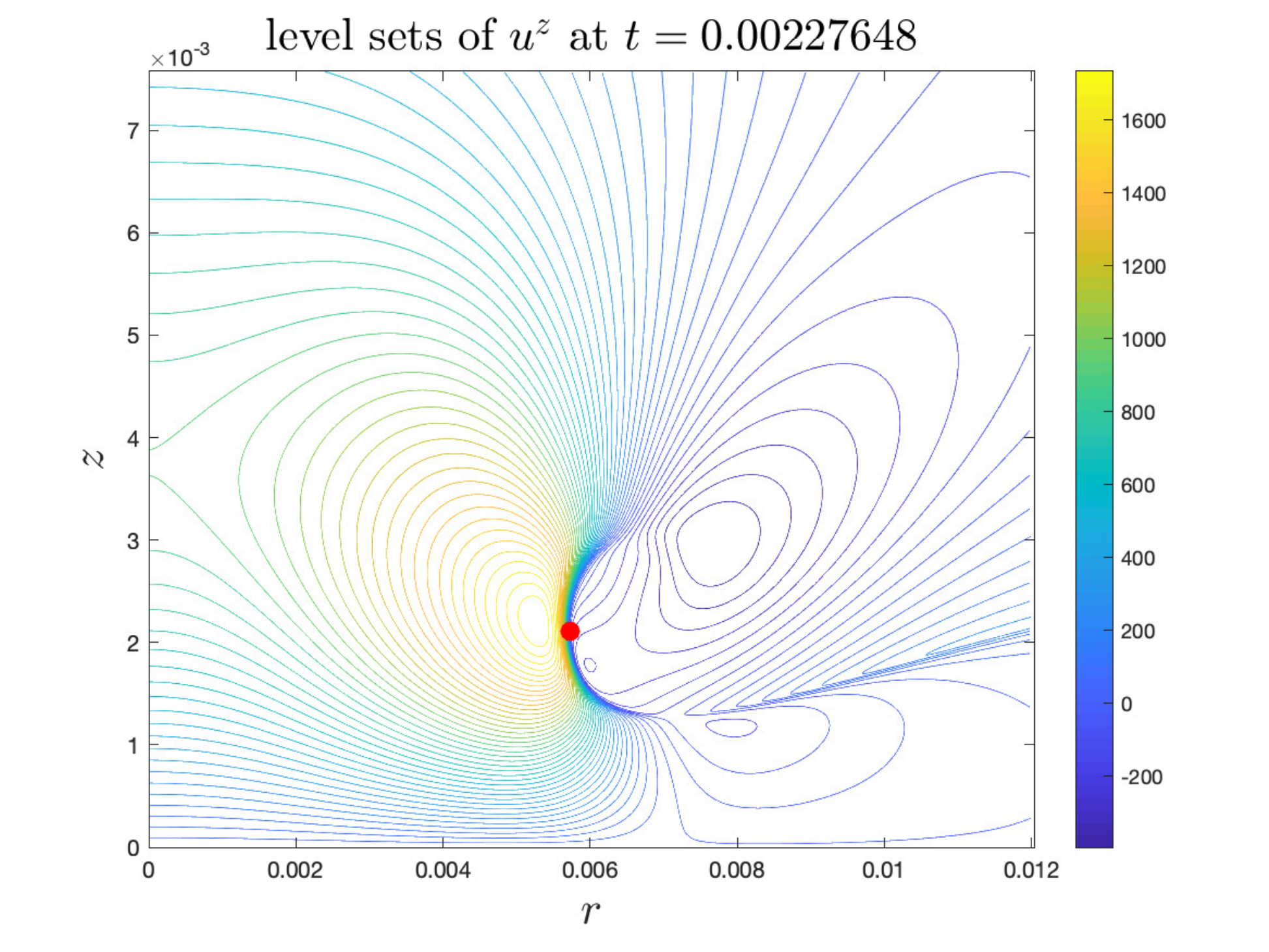} 
    \caption[Velocity level sets]{The level sets of $u^r$ (left) and $u^z$ (right) at $t = 0.00227648$. The red point is the maximum location $(R(t),Z(t))$ of $u_1(t)$.}  
     \label{fig:velocity_levelset}
        \vspace{-0.05in}
\end{figure}

The velocity field $(u^r(t),u^z(t))$ also explains the sharp structures of $u_1,\om_1$ in their local profiles. Figure~\ref{fig:velocity_levelset} shows the level sets of $u^r,u^z$ at $t = 0.00227648$. One can see that the radial velocity $u^r$ has a strong shearing layer below $(R(t),Z(t))$ (the red point). This shearing layer contributes to the sharp gradient of $u_1$ in the $z$-direction. Similarly, the axial velocity $u^z$ also has a strong shearing layer close to the point $(R(t),Z(t))$. This shearing layer contributes to the sharp front of $u_1$ in the $r$-direction. 

 Despite the sharp front developed near $(R(t),Z(t))$ in the solution as shown in the contours of $u^r$ and $u^z$ in Figure \ref{fig:velocity_levelset}, we tend to believe that the solution will not develop a shock like singularity before it reaches the symmetry axis $r=0$. The sharp front that we observe is due to the strong shearing layers developed along the $r$ and $z$ directions, which is very different from a shock formation in the compressible Euler equations. The potential singularity of the $3$D Euler equations most likely will occur at the origin, which is a stagnation point of the flow where the effect of advection is minimized. 

\subsection{Understanding the blow-up mechanism}\label{sec:mechanism} In this subsection, 
we further investigate the blow-up mechanism by examining several important factors that lead to a sustainable blow-up solution. 

\subsubsection{Vortex dipoles and hyperbolic flow} First of all, the dipole structure induced by the antisymmetric angular vorticity plays an important role in generating a favorable flow structure. To better visualize the dipole structure, we plot the velocity field in the whole period $\{(r,z); 0\leq r\leq 1, -1/2\leq z\leq 1/2\}$. Since $u_1$ and $\psi_1$ are even function of $r$, the $2$D velocity field $(u^r, u^\theta)$ is also an even function of $r$. We can extend it to the negative $r$ plane as an even function of $r$. The odd symmetry (in $z$) of $\omega_1$ contributes to a dipole structure of the angular vorticity $\om^\theta$. The dipole structure induces a hyperbolic flow in the $rz$-plane and generates a pair of antisymmetric (with respect to $z$) local circulations. This pair of antisymmetric convective circulations has the desirable property of pushing the solution near $z=0$ toward $r=0$.

\begin{figure}[!ht]
\centering
    \includegraphics[width=0.38\textwidth]{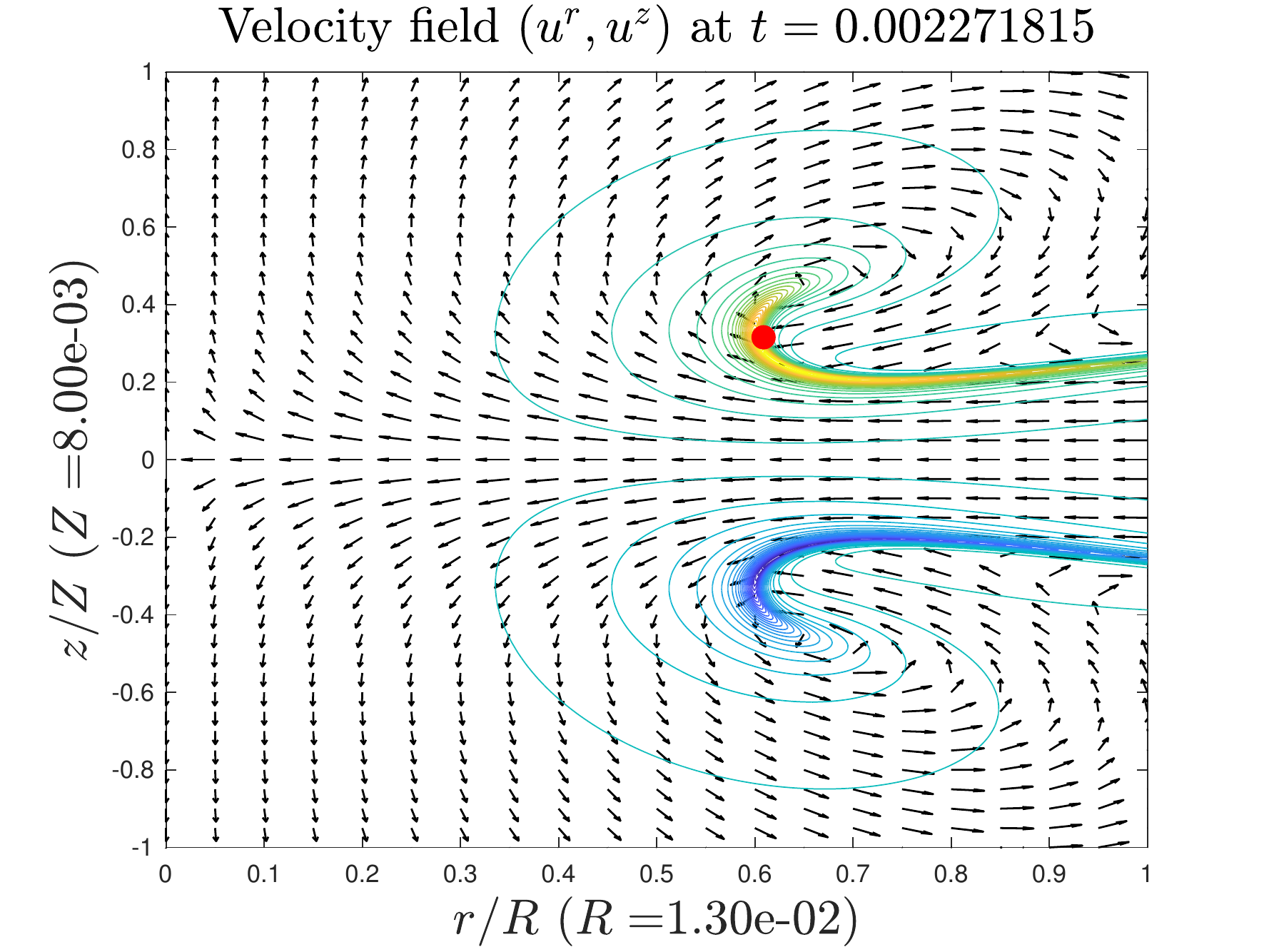}
    \includegraphics[width=0.38\textwidth]{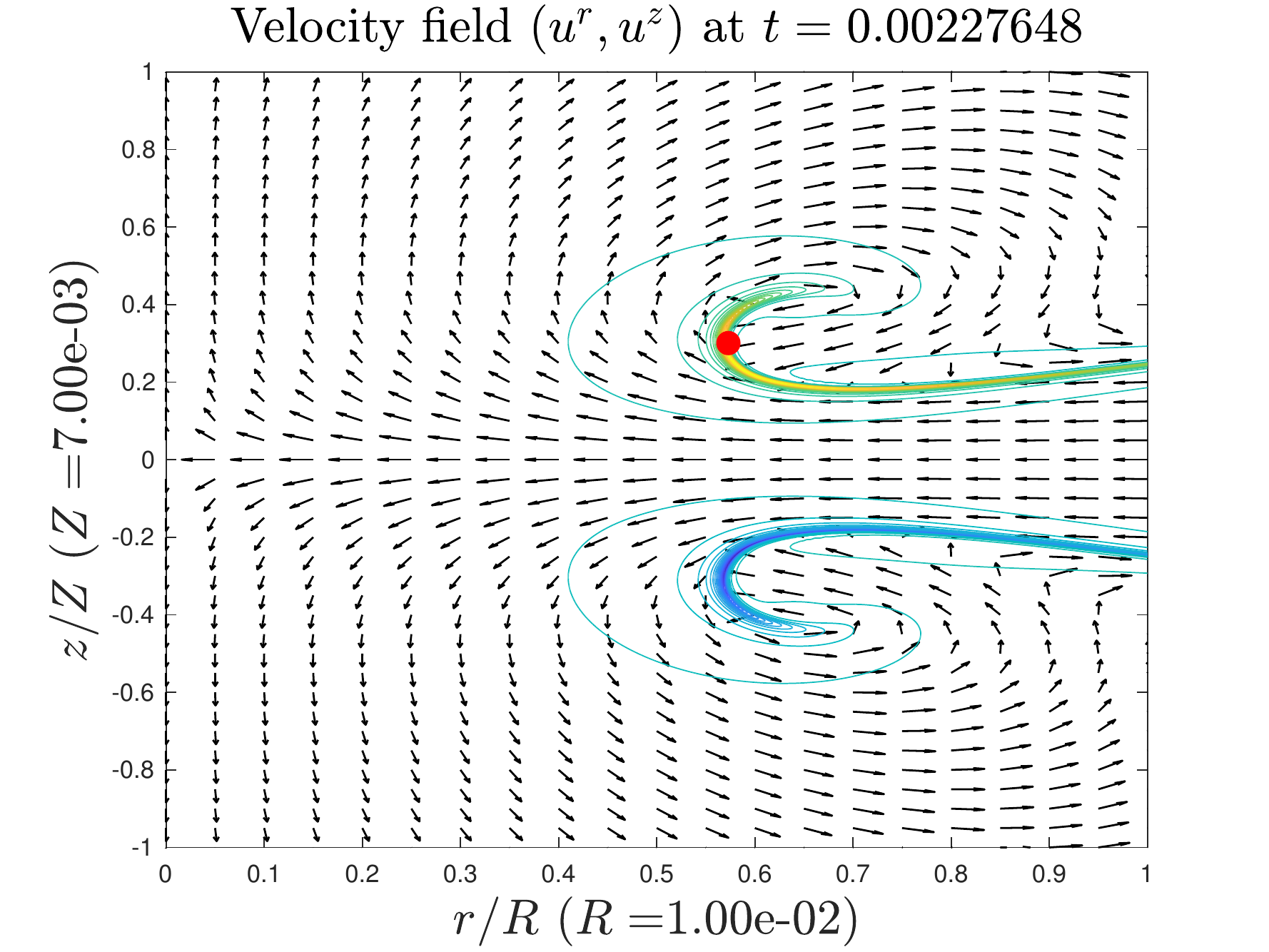} 
    \caption[Dipole]{The dipole structure of $\omega_1$  and the induced local velocity field at two different times, $t=0.002271815$ (left plot) and $t=0.00227648$ (right plot). The red point is the maximum location $(R(t),Z(t))$ of $u_1(t)$.}  
     \label{fig:dipole}
        \vspace{-0.05in}
\end{figure}

In Figure~\ref{fig:dipole}, we plot the dipole structure of $\om_1$ in a local symmetric region with the hyperbolic velocity field induced by the dipole structure in the background. The antisymmetric vortex dipoles generate a negative radial velocity near $z=0$, which pushes the solution toward $r=0$. This is one of the driving mechanisms for a potential singularity on the symmetry axis. 

\subsubsection{The odd symmetry and sharp gradient of $u_1$} 

\begin{figure}[!ht]
\centering
	\begin{subfigure}[b]{0.35\textwidth}
    \includegraphics[width=1\textwidth]{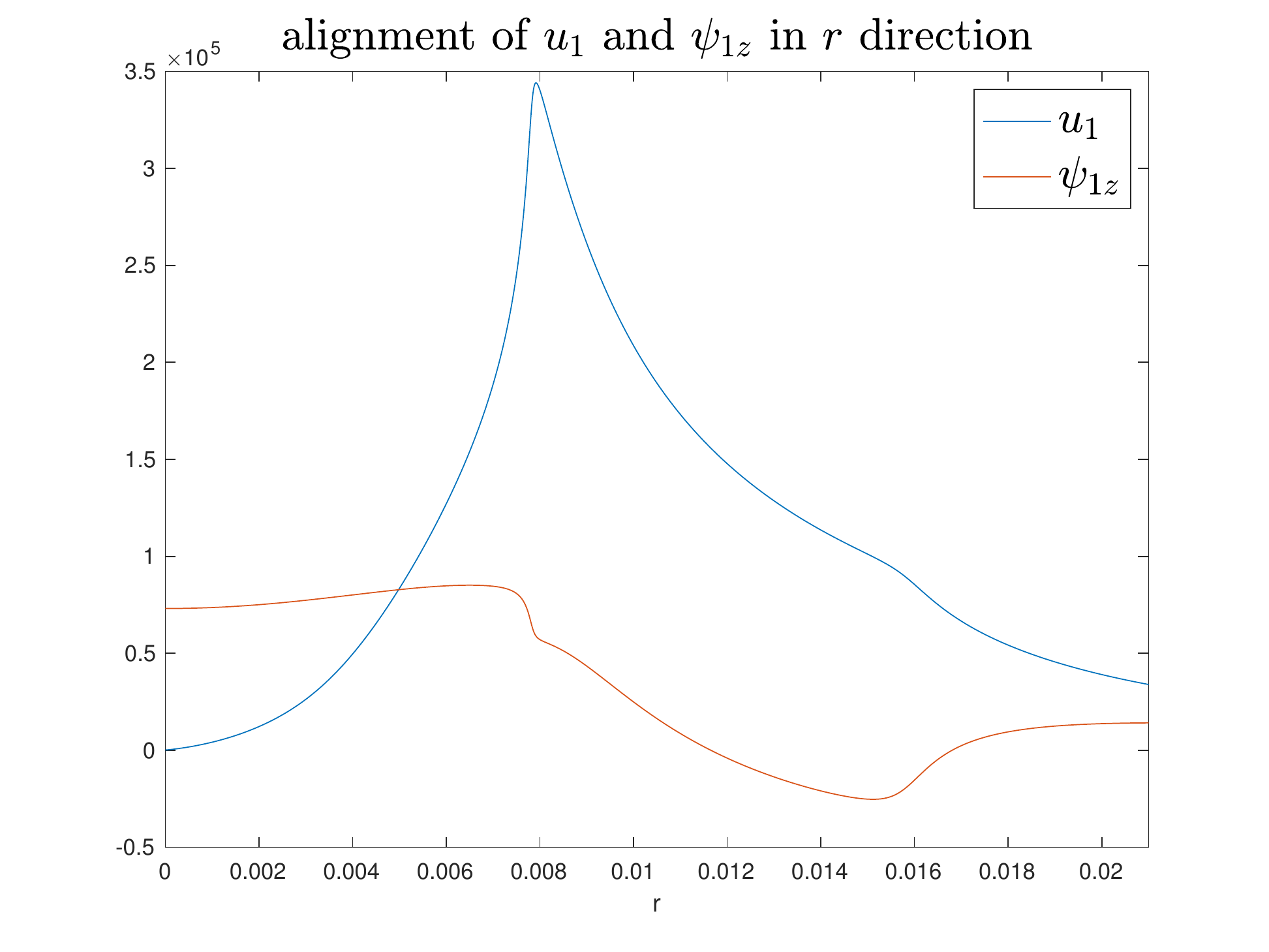}
    \caption{$r$ cross sections of $u_1,\psi_{1,z}$}
    \end{subfigure}
    \begin{subfigure}[b]{0.35\textwidth}
    \includegraphics[width=1\textwidth]{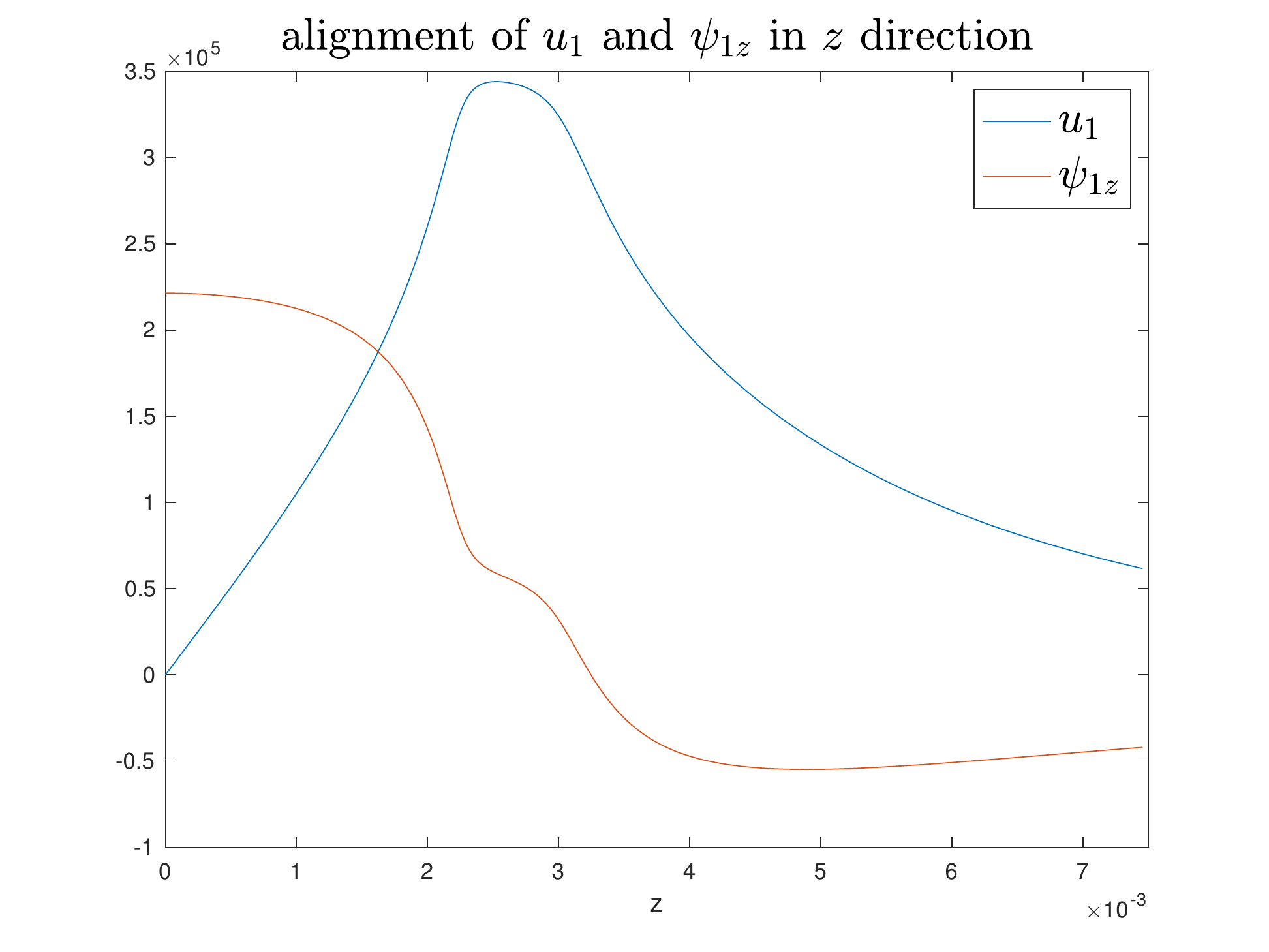}
    \caption{$z$ cross sections of $u_1,\psi_{1,z}$}
    \end{subfigure}
  	\begin{subfigure}[b]{0.35\textwidth}
    \includegraphics[width=1\textwidth]{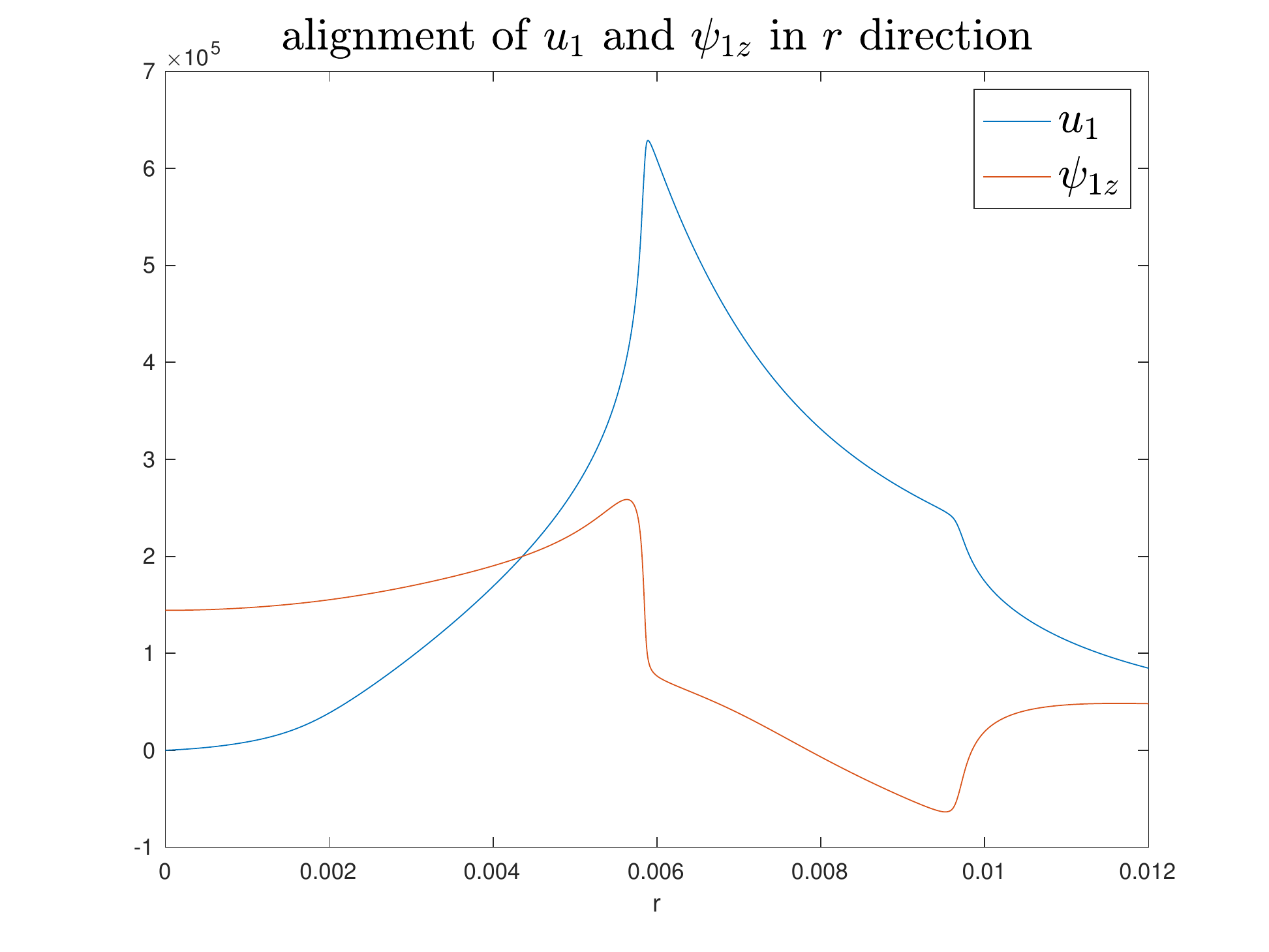}
    \caption{$r$ cross sections of $u_1,\psi_{1,z}$}
    \end{subfigure}
     \begin{subfigure}[b]{0.35\textwidth}
    \includegraphics[width=1\textwidth]{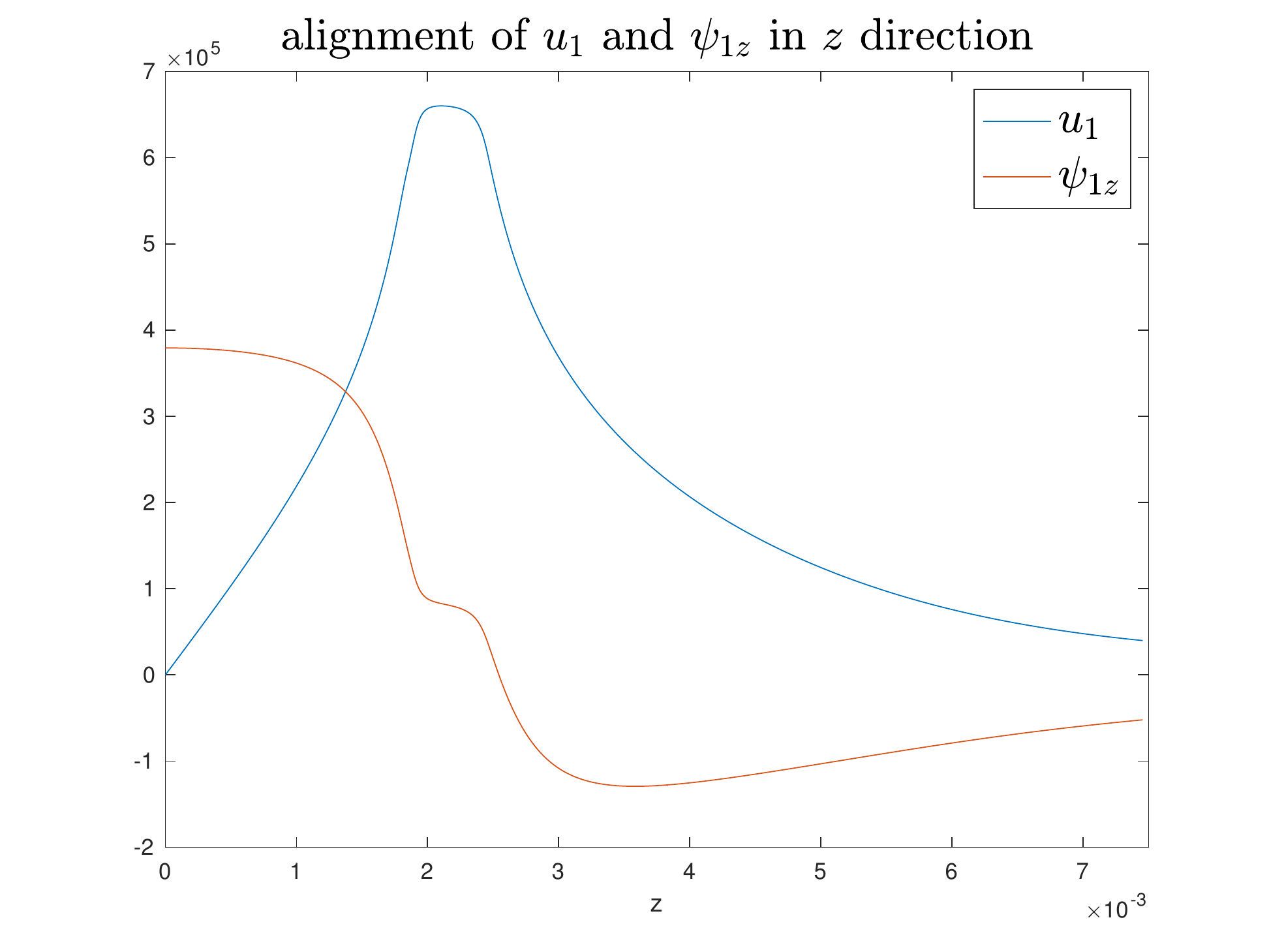}
    \caption{$z$ cross sections of $u_1,\psi_{1,z}$}
    \end{subfigure}
    \caption[Alignment]{The alignment between $u_1$ and $\psi_{1,z}$. (a) and (b): cross sections of $u_1$ and $\psi_{1,z}$ through the point $(R(t),Z(t))$ at $t=0.002271815$. (c) and (d): cross sections of $u_1$ and $\psi_{1,z}$ through the point $(R(t),Z(t))$ at $t=0.00227648$.}  
     \label{fig:alignment}
        \vspace{-0.05in}
\end{figure}

One important feature of our initial condition is that $\psi_{1,z}(r,z,t)$ is large, positive and relatively flat in a local rectangular domain $\{ (r,z) \; | \; 0 \leq r \leq 0.9 R(t), \; 0 \leq z \leq 0.5Z(t) \}$. This is an important property that was not observed in the blow-up scenario considered in \cite{Hou-Huang-2021,Hou-Huang-2022}. Moreover, $\psi_{1z}$ decays quickly outside this local domain and becomes negative in the tail region. Through the vortex stretching term $2 \psi_{1z} u_1$ in the $u_1$-equation, $u_1$ grows rapidly in this local domain while the growth outside his local domain (especially in the tail region) is relatively slow. The difference in the growth rate in the local region and the tail region produces a one-scale traveling wave solution propagating toward the origin.  In particular, this implies that $(R(t),Z(t))$ will propagate toward the origin. Since the close circle structure of the flow is located right above $(R(t),Z(t))$ and there is still a relatively strong nonlinear alignment within the circle region, see Figure \ref{fig:alignment}, the induced traveling wave also pushes the close circle region toward the origin.

The fact that $\psi_{1z}$ achieves its maximum at $z=0$ can be explained by evaluating the Poisson equation \eqref{eq:as_NSE_1_c} at $z=0$ using the odd symmetry of $\psi_1$, i.e. $\psi_{1,zz} = -\om_1$ at $z=0$. This reduced equation is still approximately valid near $z=0$. Due to the oddness of $\omega_1$, $\omega_1$ vanishes at $z=0$. Thus $z=0$ is a critical point of $\psi_{1z}$. Moreover, we have $\omega_1 > 0$ for $z>0$ and $\omega_1 <0$ for  $z<0$. This further implies that $\psi_{1z}$ achieves a local maximum at $z=0$. 

The odd symmetry of $u_1$ as a function of $z$ and the the fact that $\psi_{1z}(t,R(t),z)$ is large, positive and monotonically decreasing near $z=0$ induces a large positive gradient of $u_1^2$ in the $z$-direction between $z = Z(t)$ and $z=0$. The vortex stretching term $2(u_1^2)_z$ then induces a rapid growth for $\om_1$. The growth of $\om_1$ generates a rapid growth of $\psi_{1,z}$. This in turn contributes to the rapid growth of $u_1$ through the vortex stretching term $2\psi_{1,z}u_1$ in the $u_1$-equation. 
As $Z(t)$ approaches to $z=0$, $u_1^2$ develops an even sharper gradient in the $z$-direction, leading to an even faster growth of $\omega_1$. The whole coupling mechanism described above forms a positive feedback loop, 
\begin{equation}\label{eq:mechanism}
(u_1^2)_z \uparrow \quad \Longrightarrow \quad\om_1 \uparrow \quad \Longrightarrow \quad \psi_{1,z} \uparrow \quad \Longrightarrow \quad u_1 \uparrow \quad \Longrightarrow \quad (u_1^2)_z\uparrow,
\end{equation}
leading to a sustainable blow-up solution at the origin. 

To ensure that this positive feedback loop is dynamically stable, it is important that the maximum location of $\om_1$ should align with the location where $u_{1,z}$ is positive and large, which is slightly below $z=Z(t)$. Although our initial condition of $\om_1$ is set to zero, the dynamically generated $\om_1$ automatically has this favorable alignment property between $u_1$ and $\om_1$. In \cite{Hou-Huang-2021,Hou-Huang-2022}, we need to design the initial data in such a way that this special alignment property is achieved initially and continues to hold at later time. It is quite amazing that our simple initial data have this desirable property developed dynamically although such alignment property between $u_1$ and $\om_1$ does not hold in the very early stage.

Figure \ref{fig:alignment} demonstrates the alignment between $\psi_{1,z}$ and $u_1$. Figure \ref{fig:alignment}(b) shows the cross section of $u_1,\psi_{1,z}$ in the $z$-direction through $(R(t),Z(t))$ at $t = 0.002271815$. We can see that $\psi_{1,z}(t,R(t),z)$ is monotonically decreasing for $z \in [0, 2Z(t)]$ and relatively flat for $z \in [0, 0.5 Z(t)]$. Moreover, $\psi_{1,z}$ is large, positive, and comparable to $u_1$ in magnitude, which leads to the rapid growth of $u_1$ and pushes $Z(t)$ moving toward $z=0$. 

\begin{figure}[!ht]
\centering
	\begin{subfigure}[b]{0.35\textwidth}
    \includegraphics[width=1\textwidth]{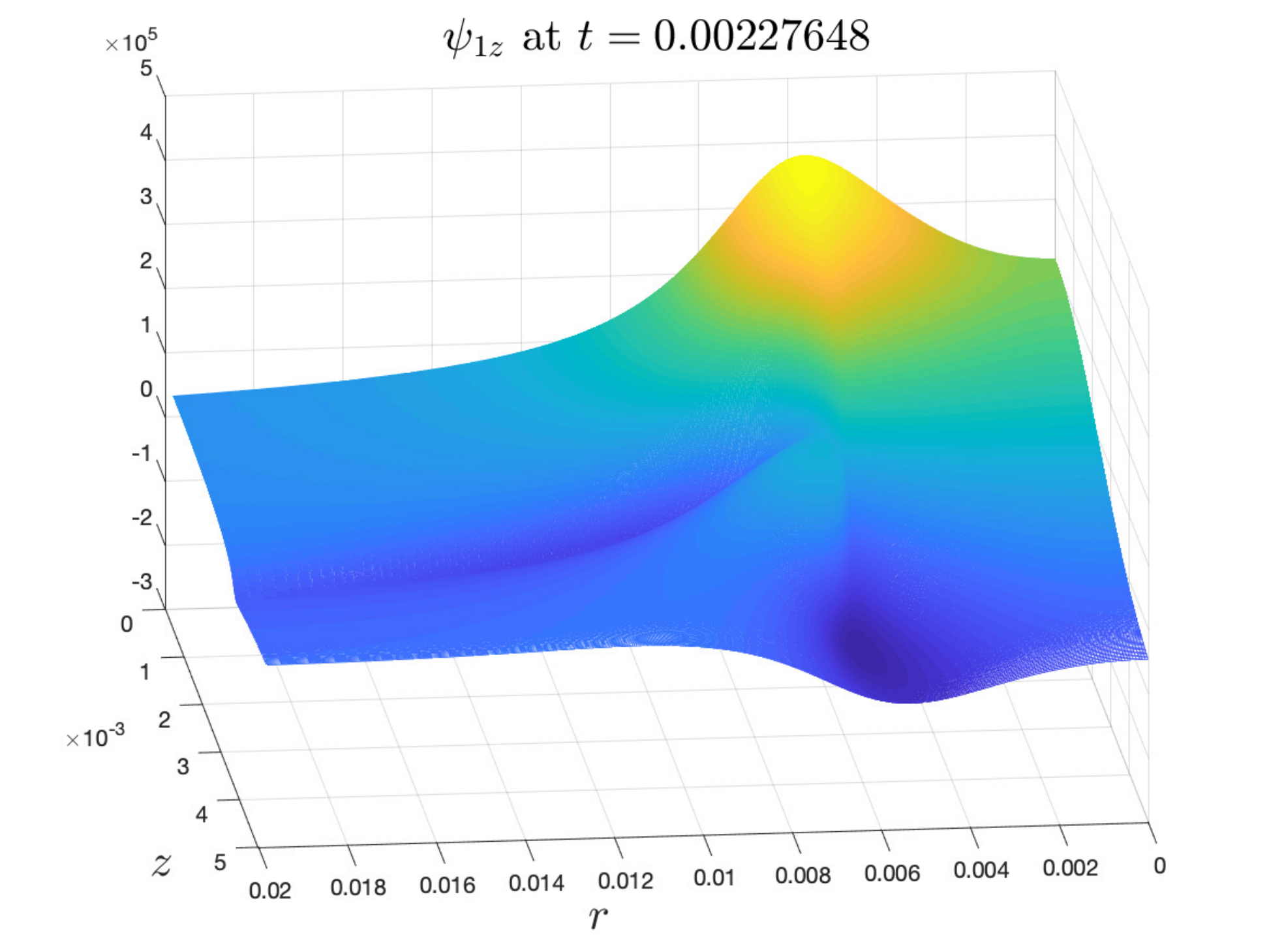}
    \caption{$3$D plot of $\psi_{1,z}$, back view}
    \end{subfigure}
    \begin{subfigure}[b]{0.35\textwidth}
    \includegraphics[width=1\textwidth]{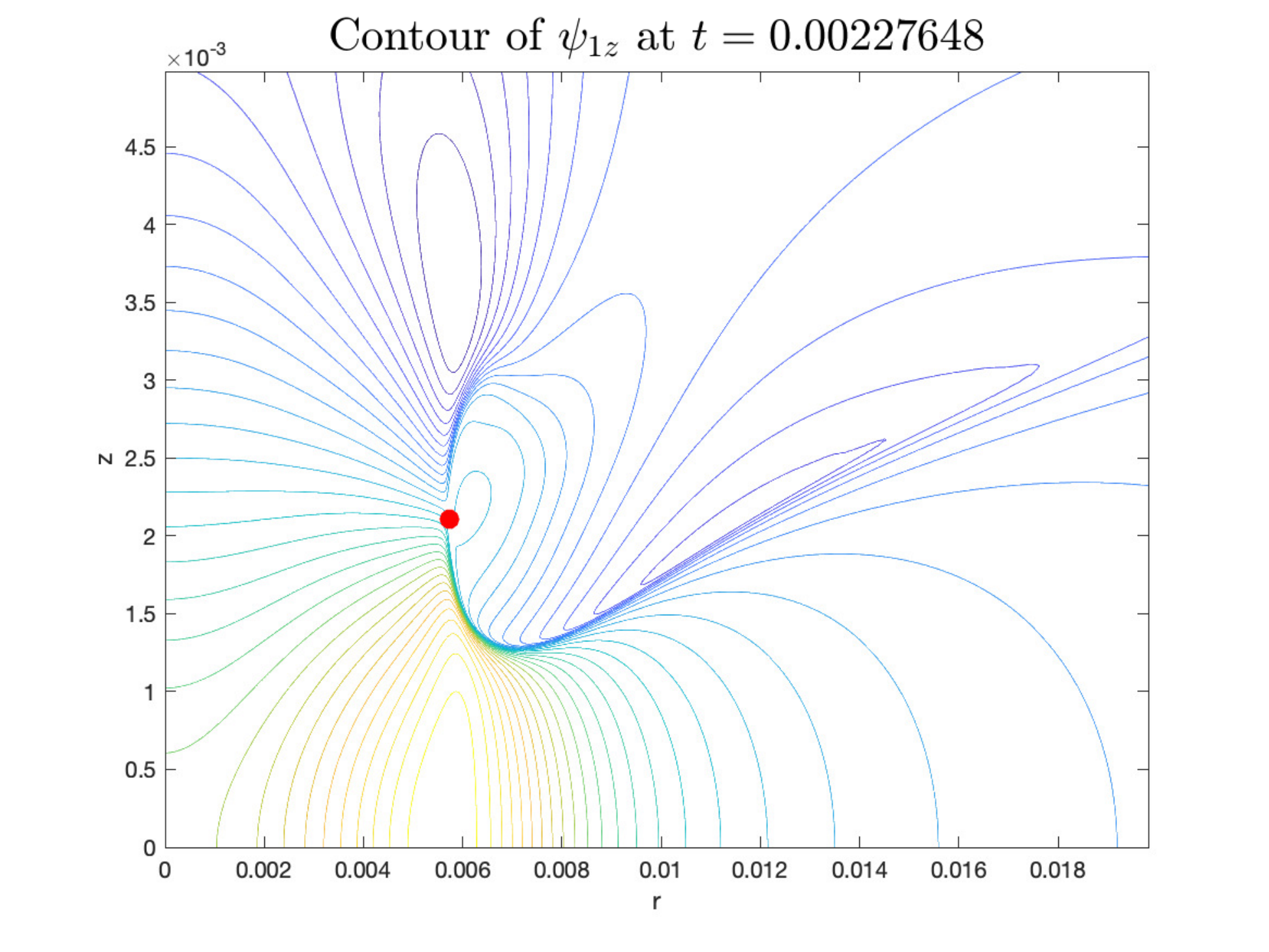}
    \caption{$z$ Contour of $\psi_{1,z}$}
    \end{subfigure}
    \caption[Alignment]{(a) The $3$D profile  of $\psi_{1,z}$ at $t=0.00227648$, back view. (b) The contour of $\psi_{1,z}$, the red dot corresponds to $(R(t),Z(t))$. }  
     \label{fig:psiz_profile}
        \vspace{-0.05in}
\end{figure}

Figure \ref{fig:alignment}(c) shows the cross section of $u_1,\psi_{1,z}$ in the $r$-direction through $(R(t),Z(t))$ at $t = 0.00227648$. We observe that $\psi_{1,z}(r,Z(t),t)$ is relatively flat for $r \in [0, 0.9 R(t)]$. This property is critical for $u_1$ remains large between the sharp front and $r=0$, thus avoiding the formation of a no-spinning region and a two-scale solution structure. We also observe that $\psi_{1z}(t,r,Z(t))$ experiences a sharp drop as a function of $r$ near the sharp front. In Figure \ref{fig:psiz_profile}(a)-(b), we plot the local $3$D profile of $\psi_{1z}$ and its contours. We see more clearly that there is a sharp drop and a plateau region for $\psi_{1z}$ behind $(R(t),Z(t))$. This explains why the alignment between $\psi_{1z}$ and $u_1$ experiences a sharp decline in the late stage, see Figure \ref{fig:trajectory}(d). 

We would like to emphasize that the blow-up scenario revealed in this paper is genuinely three dimensional. The most singular behavior and the strong nonlinear interaction occur away from the symmetry axis $r=0$ and the symmetry plane $z=0$. The angular velocity $u^\theta$ develops a sharp front and the angular vorticity $\omega^\theta$ develops a Delta function like structure and changes sign across the sharp front. The strong vortex stretching generates a traveling wave that propagates toward the origin. The geometric structure of the solution in the most singular region is quite complicated. Thus, we cannot model this blow-up phenomena by using the $1$D Hou-Li model along the symmetry axis $r=0$ \cite{hou2008dynamic}. The magical cancellation between the advection and vortex stretching in the Hou-Li model does not occur in our blow-up scenario. The positive  feedback loop that develops near $(R(t),Z(t))$ away from $r=0$ and $z=0$ makes such a blow-up scenario sustainable in the interior domain.

\subsection{Numerical Results: Resolution Study}\label{sec:performance_study_euler}
In this subsection, we perform resolution study and investigate the convergence property of our numerical methods. In particular, we will study 
\begin{itemize}
\item[(i)] the effectiveness of the adaptive mesh (Section \ref{sec:mesh_effectiveness_euler}), and 
\item[(ii)] the convergence of the solutions when $h_\rho,h_\eta\rightarrow 0$ (Section \ref{sec:resolution_study_euler}).
\end{itemize}

\subsubsection{Effectiveness of the adaptive mesh}\label{sec:mesh_effectiveness_euler} To resolve the potential finite time singularity of the $3$D Euler equations, it is essential to design the adaptive mesh to resolve the solution in the most singular region as well as in the far field. The detail of how to construct our adaptive mesh for the Euler equations will be given in Appendix~\ref{apdx:adaptive_mesh}. In this subsection, we study the effectiveness of our adaptive moving mesh. 

\begin{figure}[!ht]
\centering
	\begin{subfigure}[b]{0.40\textwidth}
    \includegraphics[width=1\textwidth]{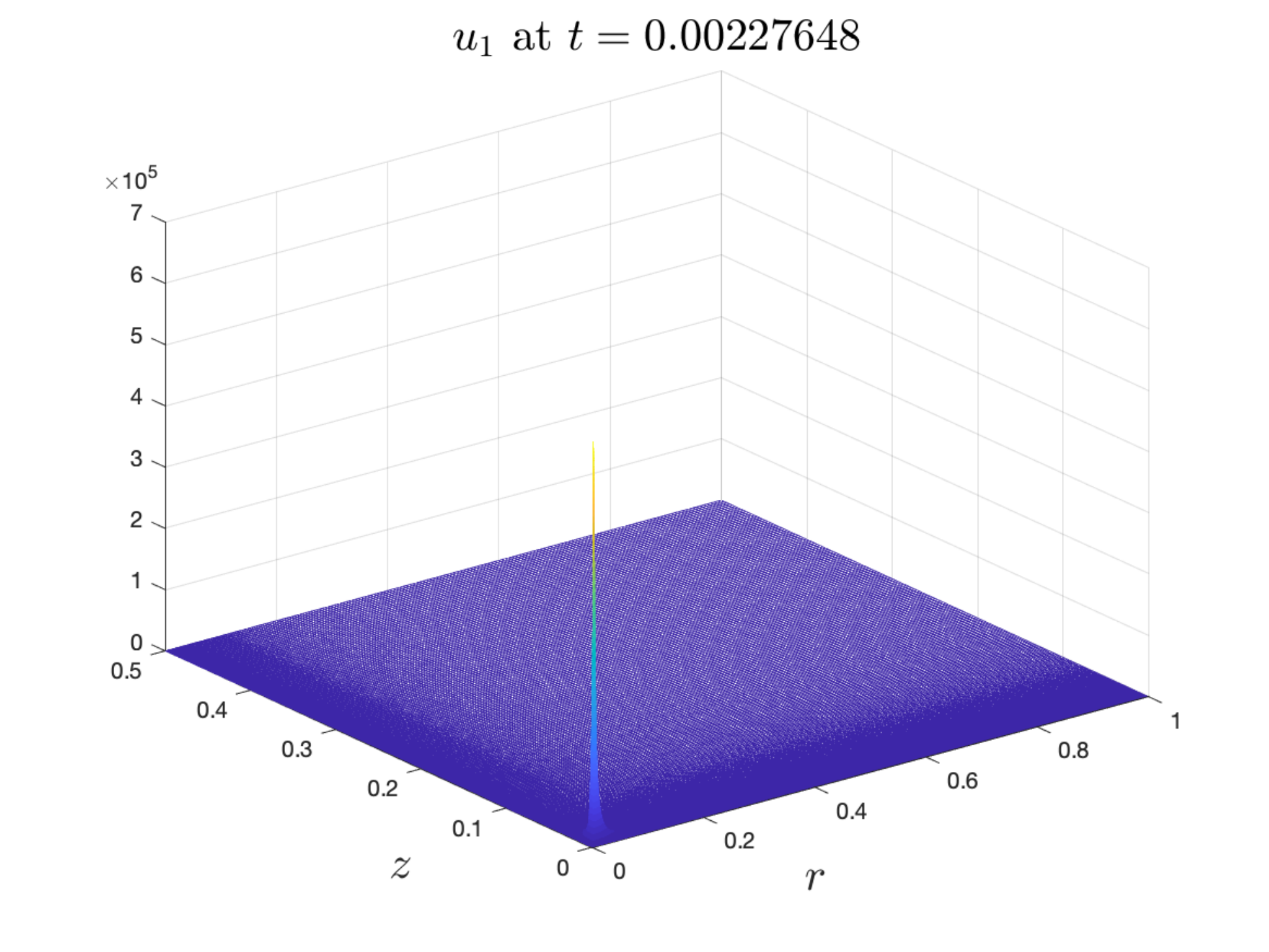}
    \caption{$u_1$ in the $rz$-plane}
    \end{subfigure}
  	\begin{subfigure}[b]{0.40\textwidth}
    \includegraphics[width=1\textwidth]{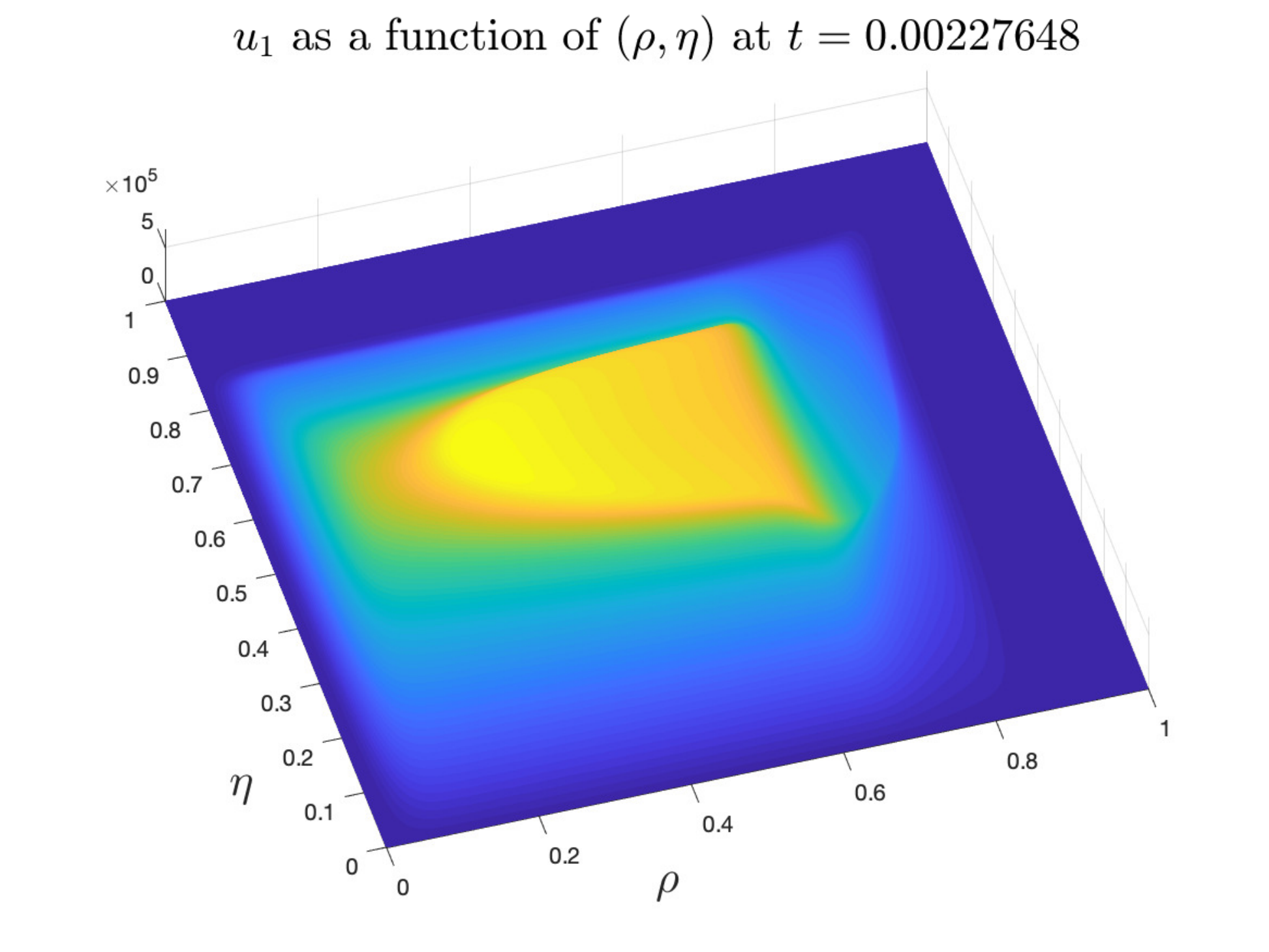}
    \caption{$u_1$ in the $\rho\eta$-plane}
    \end{subfigure}
    \caption[Mesh effectiveness]{The adaptive mesh resolves the solution in the $\rho\eta$-plane. The left subplot (a) shows the focusing singularity structure of $u_1$ at $t=0.00227648$ in the $rz$-plane on the whole computational domain $\mathcal{D}_1$. The right subplot (b) plots the profile of $u_1$ in the $\rho\eta$-plane.}  
     \label{fig:mesh_effective}
        \vspace{-0.05in}
\end{figure}

To see how well the adaptive mesh resolves the solution, we first visualize how it transforms the solution from the $rz$-plane to the $\rho\eta$-plane. Figure~\ref{fig:mesh_effective}(a) shows the function $u_1$ at $t=0.00227648$ in the original $rz$-plane. This plot suggests that the solution develops a focusing singularity at the origin. For comparison, Figure~\ref{fig:mesh_effective}(b) plots the profile of $u_1$ at the same time in the $\rho\eta$-plane. We can see that our adaptive mesh resolves the potential singular solution in the $(\rho,\eta)$ coordinates.

\begin{figure}[!ht]
\centering
 \vspace{-0.7in}
	\begin{subfigure}[b]{0.40\textwidth}
    \includegraphics[width=1\textwidth]{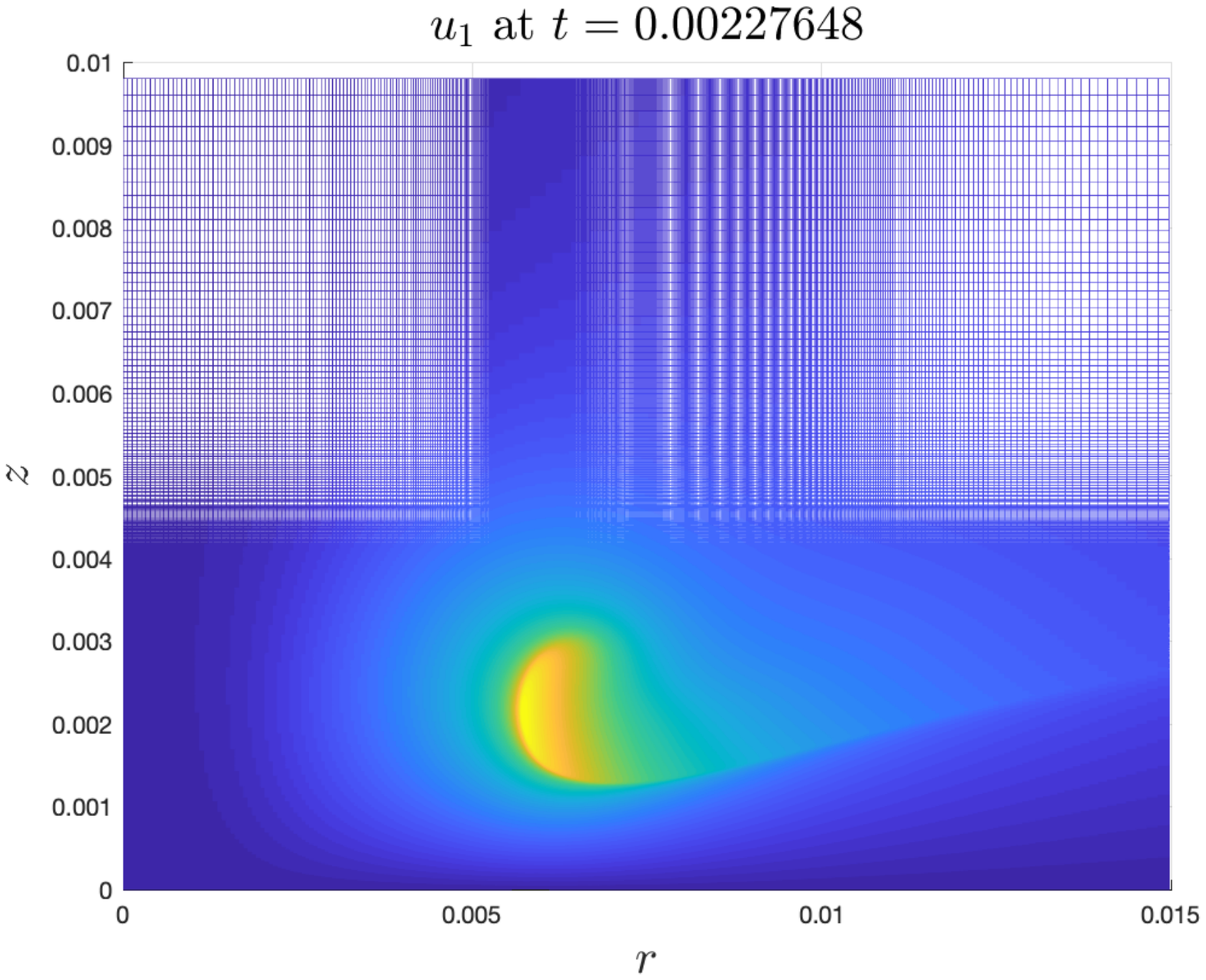}
    \end{subfigure}
  	\begin{subfigure}[b]{0.40\textwidth}
    \includegraphics[width=1\textwidth]{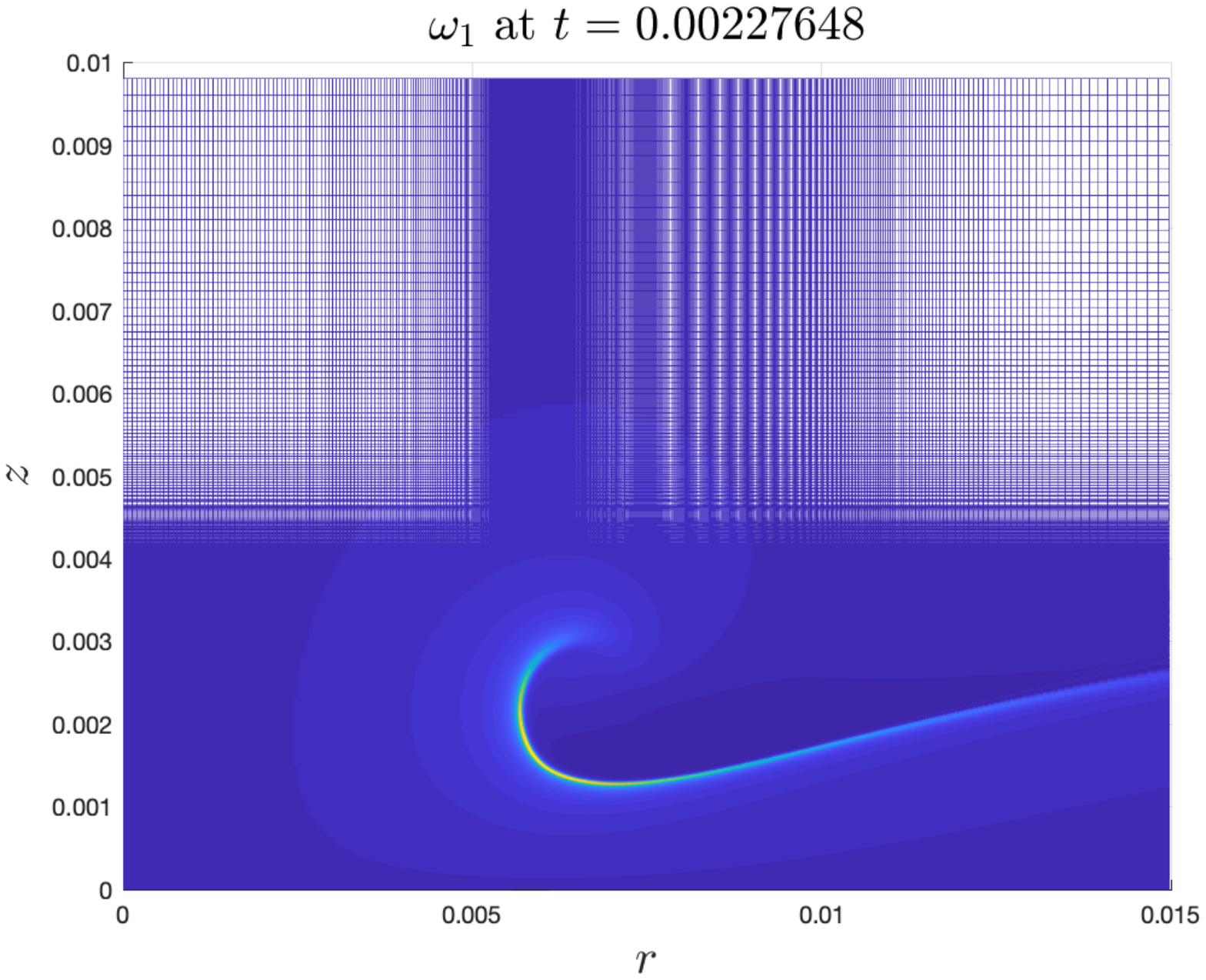}
    \end{subfigure} 
     \vspace{-0.7in}
    \caption[Mesh phases]{The adaptive mesh has different densities in different regions. Left: adaptive mesh for $u_1$. Right: adaptive mesh for $\omega_1$. }  
     \label{fig:mesh_phases}
        \vspace{-0.05in}
\end{figure}

In Figure~\ref{fig:mesh_phases}, we show the top views of the profiles of $u_1,\om_1$ in a local domain at $t=0.00227648$. This figure demonstrates how the mesh points are distributed in different phases of the adaptive mesh. We can see that a large number of the adaptive meshes concentrate in phase $1$ in both directions where solution is most singular.

Inspired by the recent work in \cite{Hou-Huang-2021}, we define the mesh effectiveness functions $ME_\rho(v),ME_\eta$ with respect to some solution variable $v$ as follows:
\[ME_\rho(v) = \frac{h_\rho v_\rho}{\|v\|_{L^\infty}} = \frac{h_\rho r_\rho v_r}{\|v\|_{L^\infty}},\quad ME_\eta(v) = \frac{h_\eta v_\eta}{\|v\|_{L^\infty}} = \frac{h_\eta r_\eta v_z}{\|v\|_{L^\infty} }.\]
We further define the corresponding mesh effectiveness measures (MEMs) as follows:
\[ME_{\rho,\infty}(v) = \|ME_\rho(v)\|_{L^\infty},\quad ME_{\eta,\infty}(v) = \|ME_\eta(v)\|_{L^\infty}.\]
The MEMs quantify the largest relative growth of a function $v$ in one single mesh cell.  The smallness of the MEMs measures the effectiveness of an adaptive mesh. 

\begin{figure}[!ht]
\centering
    \includegraphics[width=0.38\textwidth]{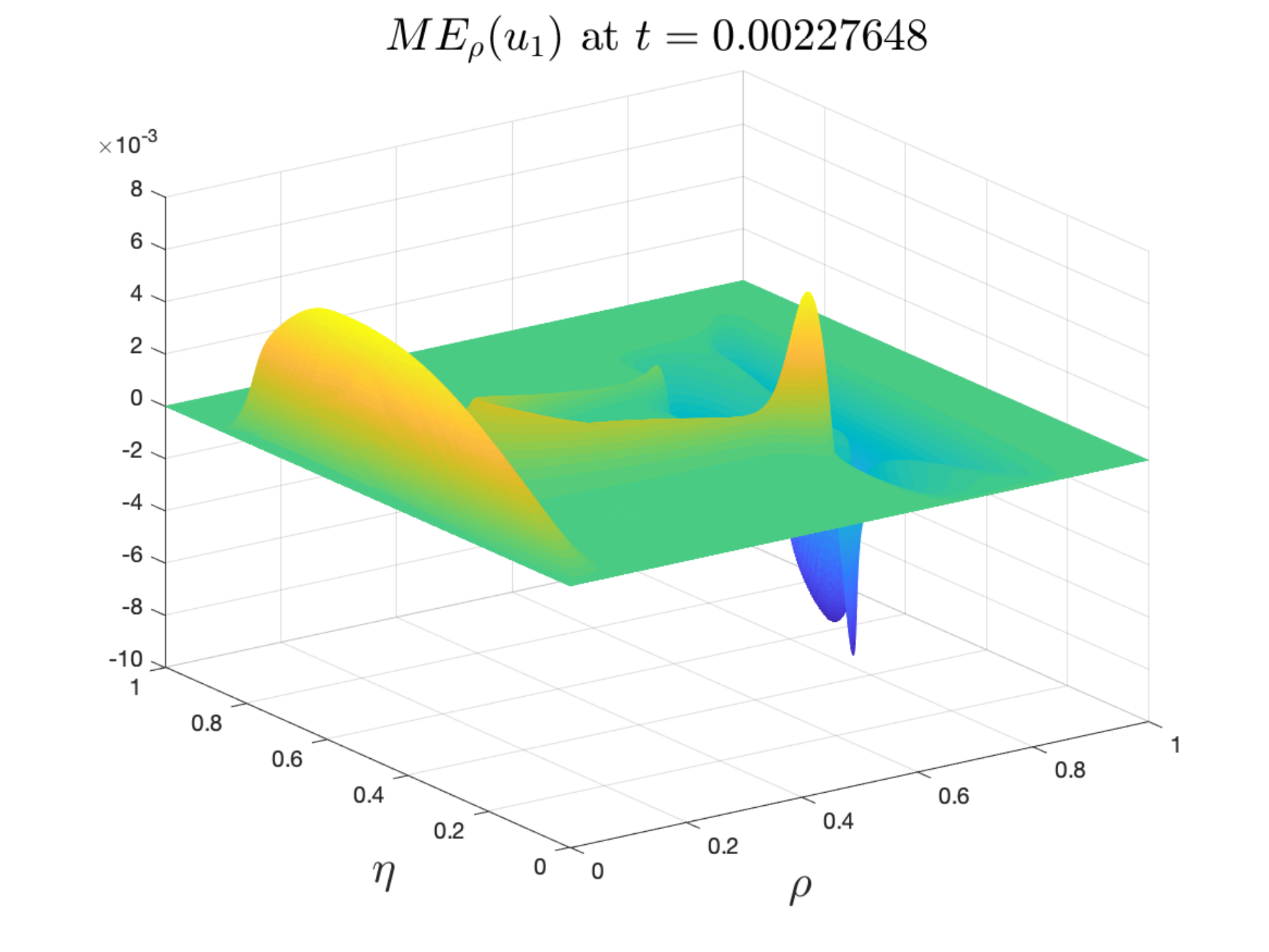}
    \includegraphics[width=0.38\textwidth]{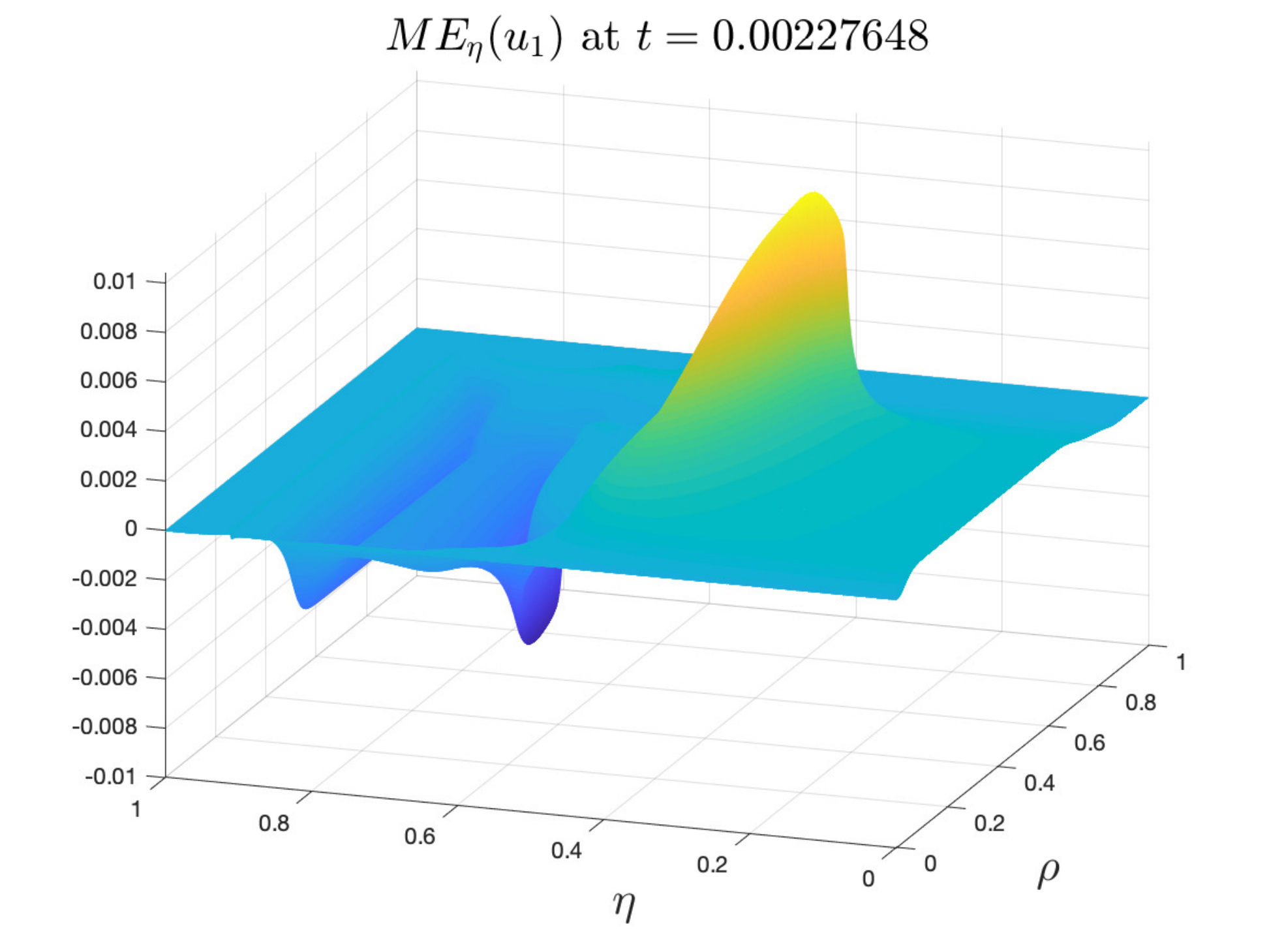} 
    \includegraphics[width=0.38\textwidth]{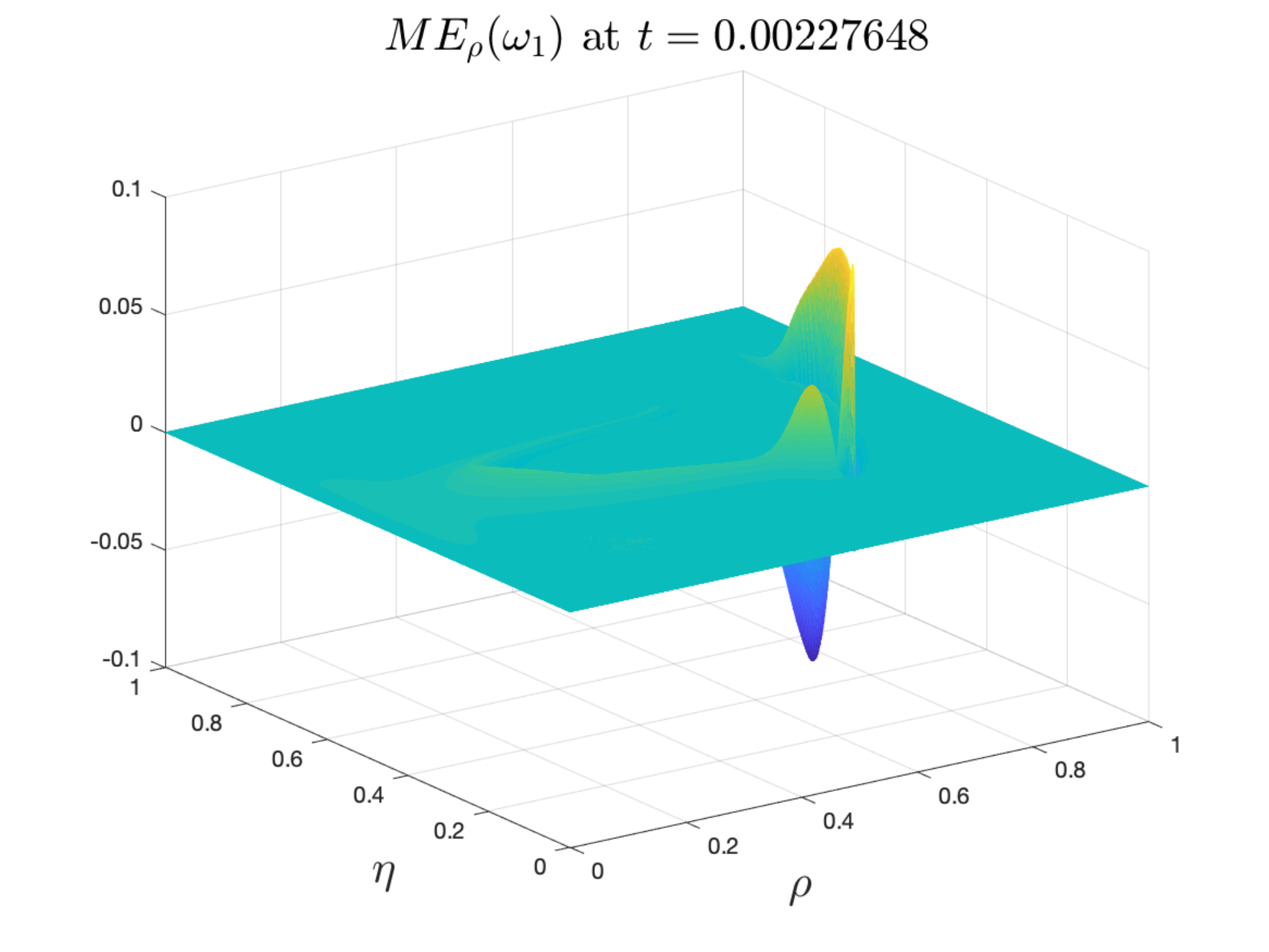}
    \includegraphics[width=0.38\textwidth]{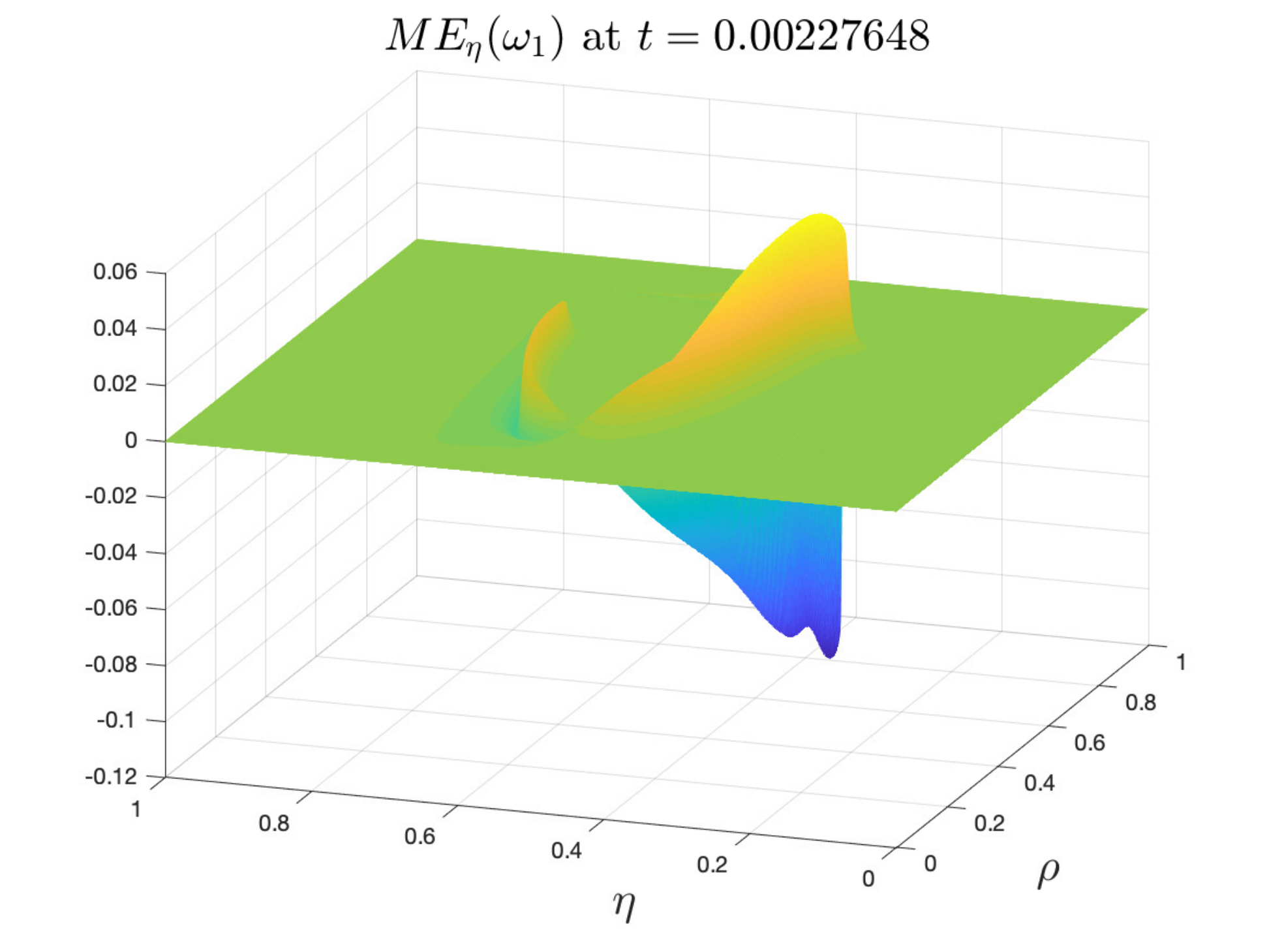}
    \caption[Mesh effectiveness functions]{First row: the mesh effectiveness functions of $u_1$ at $t=0.00227648$ with mesh dimension $(n_1,n_2) = (1536,1536)$. Second row: the mesh effectiveness functions of $\om_1$ in the same setting.}  
     \label{fig:MEF}
        \vspace{-0.05in}
\end{figure}

In Figure~\ref{fig:MEF}, we plot the mesh effectiveness functions of $u_1,\om_1$ at time $0.00227648$ on the mesh of size $(n_1,n_2) = (1536,1536)$. We can see that the mesh effectiveness functions of $u_1,\om_1$ are relatively small and uniformly bounded by a small number. 

Table~\ref{tab:MEM_mesh} reports the MEMs of $u_1,\om_1$ at $t=0.00227648$ on meshes of different sizes. We can see that the MEMs decrease as the grid sizes $h_\rho,h_\eta$ decrease since the MEMs are proportional to $h_\rho,h_\eta$. Table~\ref{tab:MEM_time} reports the MEMs of $u_1,\om_1$ at different times with the same mesh size $(n_1,n_2) = (1536,1536)$. We can see that the MEMs remain relatively small throughout this time interval, although there is an increasing trend in time. The above study implies that our adaptive mesh strategy is effective in resolving the potentially singular solution both in the most singular region and in the far field.

We remark that our mesh strategy can well resolve the solution up to $t=0.00227648$. With our highest resolution $(1536,1536)$, we can push the computation up to $t=0.002276938$ and still maintain a reasonable accuracy. To continue the computation beyond $t=0.002276938$, one may need to use a finer resolution or a more sophisticated adaptive mesh strategy to resolve the sharp front. 

\begin{table}[!ht]
\centering
\footnotesize
\renewcommand{\arraystretch}{1.5}
    \begin{tabular}{|c|c|c|c|c|}
    \hline
    \multirow{2}{*}{Mesh size} & \multicolumn{4}{c|}{MEMs on mesh at $t=0.00227648$} \\ \cline{2-5} 
    						   & $ME_{\rho,\infty}(u_1)$ & $ME_{\eta,\infty}(u_1)$ & $ME_{\rho,\infty}(\om_1)$ & $ME_{\eta,\infty}(\om_1)$  \\ \hline 
    $768\times 768$ & $0.016$ & $0.020$ & $0.129$ & $0.163$ \\ \hline 
    $1024\times 1024$ & $0.013$ & $0.016$ & $0.113$ & $0.148$  \\ \hline 
    $1280\times 1280$ & $0.010$ & $0.012$ & $0.096$ & $0.126$ \\ \hline 
    $1536\times 1536$ & $0.009$ & $0.010$ & $0.082$ & $0.108$  \\ \hline 
    \end{tabular}
    \caption{\small MEMs of $u_1,\om_1$ at $t=0.00227648$ on the meshes of different sizes.}
    \label{tab:MEM_mesh}
    \vspace{-0.3in}
\end{table}

\begin{table}[!ht]
\centering
\footnotesize
\renewcommand{\arraystretch}{1.5}
    \begin{tabular}{|c|c|c|c|c|}
    \hline
    \multirow{2}{*}{Time} & \multicolumn{4}{c|}{MEMs on mesh $(n_1,n_2) = (1536,1536)$} \\ \cline{2-5} 
    						   & $ME_{\rho,\infty}(u_1)$ & $ME_{\eta,\infty}(u_1)$ & $ME_{\rho,\infty}(\om_1)$ & $ME_{\eta,\infty}(\om_1)$  \\ \hline 
    $0.00226435$ & $0.006$ & $0.004$ & $0.014$ & $0.020$  \\ \hline 
    $0.002271815$ & $0.006$ & $0.006$ & $0.022$ & $0.031$ \\ \hline 
    $0.002274696$ & $0.007$ & $0.006$ & $0.042$ & $0.060$  \\ \hline 
    $0.002276480$ & $0.009$ & $0.010$ & $0.082$ & $0.108$ \\ \hline 
    $0.002276938$ & $0.01$ & $0.012$ & $0.093$ & $0.126$ \\ \hline 
    \end{tabular}
    \caption{\small MEMs of $u_1,\om_1$ at different times on the mesh of size $(n_1,n_2) = (1536,1536)$.}
    \label{tab:MEM_time}
    \vspace{-0.3in}
\end{table}

\subsubsection{Resolution study}\label{sec:resolution_study_euler}
In this subsection, we perform resolution study on the numerical solutions of the initial-boundary value problem \eqref{eq:axisymmetric_NSE_1} at various time instants $t$. We will estimate the relative error of a solution variable $f_p$ computed on the $256p\times 256p$ mesh by comparing it to a reference variable $\hat{f}$ that is computed on the finest mesh of size $1536\times 1536$ at the same time instant. If $f_p$ is a number, the absolute relative error is computed as 
$e_p = |f_p-\hat{f}|/|\hat{f}|$.
If $f_p$ is a spatial function, the reference variable $\hat{f}$ is first interpolated to the mesh on which $f$ is computed. Then the sup-norm relative error is computed as 
\begin{align*}
e_p &= \frac{\|f_p-\hat{f}\|_{\infty}}{\|\hat{f}\|_{\infty}}\quad \text{if $f$ is a scalar function,}\\
\text{and}\quad e_p &= \frac{\left\|\big|(f_p^\theta-\hat{f}_p^\theta,f_p^r-\hat{f}_p^r,f_p^z-\hat{f}_p^z)\big|\right\|_{\infty}}{\left\|\big|(\hat{f}_p^\theta,\hat{f}_p^r,\hat{f}_p^z)\big|\right\|_{\infty}}\quad \text{if $f$ is a vector function.}
\end{align*}
The numerical order of the error is computed as 
\[\beta_p = \log_{\frac{p}{p-1}}\left(\frac{e_{p-1}}{e_p}\right).\]


We first study the sup-norm error of the solution, which is the most important measure of accuracy for our numerical method. Tables~\ref{tab:sup-norm_error_1-1}--\ref{tab:sup-norm_error_1-4} report the sup-norm relative errors and numerical orders of different solution variables at times $t = 0.002264353$ and $t = 0.002271815$, respectively. The results confirm that our method is at least $2$nd-order accurate. 

\begin{table}[!ht]
\centering
\footnotesize
\renewcommand{\arraystretch}{1.5}
    \begin{tabular}{|c|c|c|c|c|c|c|}
    \hline
    \multirow{2}{*}{Mesh size} & \multicolumn{6}{c|}{Sup-norm relative error at $t=0.002264353$ for $3$D Euler equations} \\ \cline{2-7} 
    						   & $u_1$ & Order & $\omega_1$ & Order & $\psi_1$ & Order \\ \hline 
    $512\times512$ & $1.080\times10^{-1}$ & -- & $3.6465\times10^{-1}$ & -- & $3.0167\times10^{-2}$ & -- \\ \hline 
    $768\times768$ & $4.2066\times10^{-2}$ & $2.326$ & $1.5520\times10^{-1}$ & $2.107$ & $1.1469\times10^{-2}$ & $2.385$ \\ \hline 
    $1024\times1024$ & $1.7730\times10^{-2}$ & $3.003$ & $6.6567\times10^{-2}$ & $2.942$ & $4.8033\times10^{-3}$ & $3.025$ \\ \hline 
    $1280\times1280$ & $6.2685\times10^{-3}$ & $4.659$ & $2.3668\times10^{-2}$ & $4.634$ & $1.6947\times10^{-3}$ & $4.669$ \\ \hline 
    \end{tabular}
    \caption{\small Sup-norm relative errors and numerical orders of $u_1,\om_1,\psi_1$ at $t = 0.002264353$ for the $3$D Euler equations.}
    \label{tab:sup-norm_error_1-1}
    \vspace{-0.2in}
\end{table}

\begin{table}[!ht]
\centering
\footnotesize
\renewcommand{\arraystretch}{1.5}
    \begin{tabular}{|c|c|c|c|c|c|c|}
    \hline
    \multirow{2}{*}{Mesh size} & \multicolumn{6}{c|}{Sup-norm relative error at $t=0.002264353$ for $3$D Euler equations} \\ \cline{2-7} 
    						   & $u^r$ & Order & $u^z$ & Order & $\vom=(\om^\theta,\om^r,\om^z)$ & Order \\ \hline 
    $512\times512$ & $1.0761\times10^{-1}$ & -- & $2.4350\times10^{-1}$& -- & $3.6878\times10^{-1}$ & -- \\ \hline 
    $768\times768$ & $4.2298\times10^{-2}$ & $2.303$ &
    $9.6378\times10^{-2}$ & $2.286$ & $1.5573\times10^{-1}$ & $2.126$ \\ \hline 
    $1024\times1024$ & $1.7850\times10^{-2}$ & $2.999$ & $4.0709\times10^{-2}$ & $2.996$ &  $6.6667\times10^{-2}$ & $2.949$ \\ \hline 
    $1280\times1280$ & $6.3127\times10^{-3}$ & $4.658$ &
    $1.4403\times10^{-2}$ & $4.656$ & $2.3694\times10^{-2}$ & $4.636$ \\ 
    \hline
    \end{tabular}
    \caption{\small Sup-norm relative errors and numerical orders of $u^r,u^z,\vom$ at $t = 0.002264353$ for $3$D Euler equations.}
    \label{tab:sup-norm_error_1-2}
    \vspace{-0.2in}
\end{table}

\begin{table}[!ht]
\centering
\footnotesize
\renewcommand{\arraystretch}{1.5}
    \begin{tabular}{|c|c|c|c|c|c|c|}
    \hline
    \multirow{2}{*}{Mesh size} & \multicolumn{6}{c|}{Sup-norm relative error at $t=0.002271815$ for $3$D Euler equations} \\ \cline{2-7} 
    						   & $u_1$ & Order & $\omega_1$ & Order & $\psi_1$ & Order \\ \hline 
    $512\times512$ & $2.4869\times10^{-1}$ & -- & $7.8344\times10^{-1}$ & -- & $5.3814\times10^{-2}$ & -- \\ \hline 
    $768\times768$ & $1.1086\times10^{-1}$ & $1.993$ & $5.5565\times10^{-1}$ & $0.847$ & $2.0597\times10^{-2}$ & $2.369$ \\ \hline 
    $1024\times1024$ & $4.8370\times10^{-2}$ & $2.883$ & $2.7929\times10^{-1}$ & $2.391$ & $8.6479\times10^{-3}$ & $3.017$ \\ \hline 
    $1280\times1280$ & $1.7214\times10^{-2}$ & $4.630$ & $1.0298\times10^{-1}$ & $4.471$ & $3.0525\times10^{-3}$ & $4.667$ \\ \hline 
    \end{tabular}
    \caption{\small Sup-norm relative errors and numerical orders of $u_1,\om_1,\psi_1$ at $t = t=0.002271815$ for the $3$D Euler equations.}
    \label{tab:sup-norm_error_1-3}
    \vspace{-0.2in}
\end{table}

\begin{table}[!ht]
\centering
\footnotesize
\renewcommand{\arraystretch}{1.5}
    \begin{tabular}{|c|c|c|c|c|c|c|}
    \hline
    \multirow{2}{*}{Mesh size} & \multicolumn{6}{c|}{Sup-norm relative error at $t=0.002271815$ for $3$D Euler equations} \\ \cline{2-7} 
    						   & $u^r$ & Order & $u^z$ & Order & $\vom=(\om^\theta,\om^r,\om^z)$ & Order \\ \hline 
    $512\times512$ & $3.1741\times10^{-1}$ & -- & $6.3493\times10^{-1}$& -- & $8.0538\times10^{-1}$ & -- \\ \hline 
    $768\times768$ & $1.4943\times10^{-1}$ & $1.858$ &
    $3.1370\times10^{-1}$ & $1.739$ & $5.5844\times10^{-1}$ & $0.903$ \\ \hline 
    $1024\times1024$ & $6.5946\times10^{-2}$ & $2.843$ & $1.3984\times10^{-1}$ & $2.809$ &  $2.7905\times10^{-1}$ & $2.411$ \\ \hline 
    $1280\times1280$ & $2.3730\times10^{-2}$ & $4.580$ &
    $4.9944\times10^{-2}$ & $4.614$ & $1.0285\times10^{-1}$ & $4.473$ \\ 
    \hline
    \end{tabular}
    \caption{\small Sup-norm relative errors and numerical orders of $u^r,u^z,\vom$ at $t = t=0.002271815$ for $3$D Euler equations.}
    \label{tab:sup-norm_error_1-4}
    \vspace{-0.2in}
\end{table}

\begin{figure}[!ht]
\centering
    \includegraphics[width=0.35\textwidth]{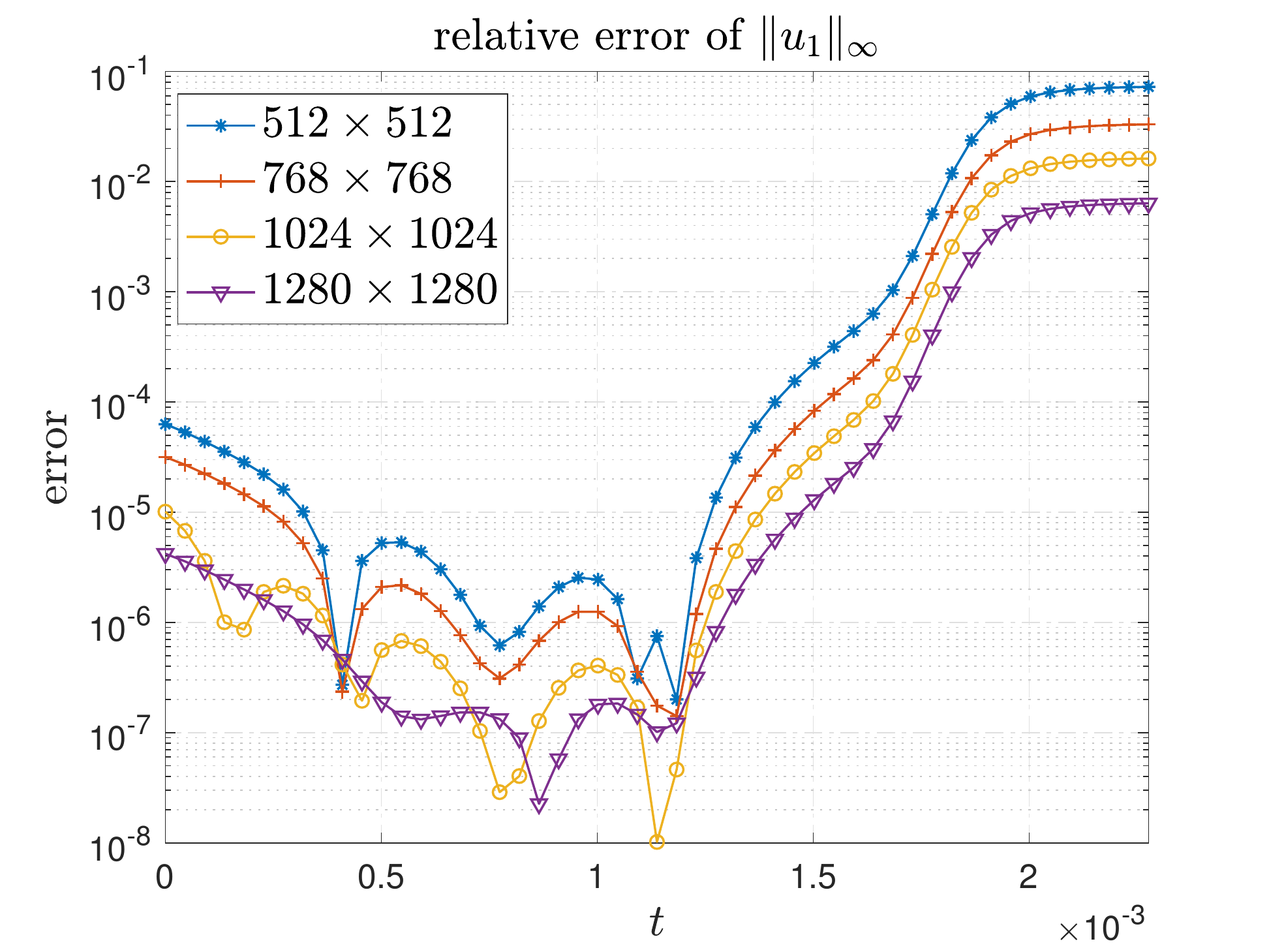}
\includegraphics[width=0.35\textwidth]{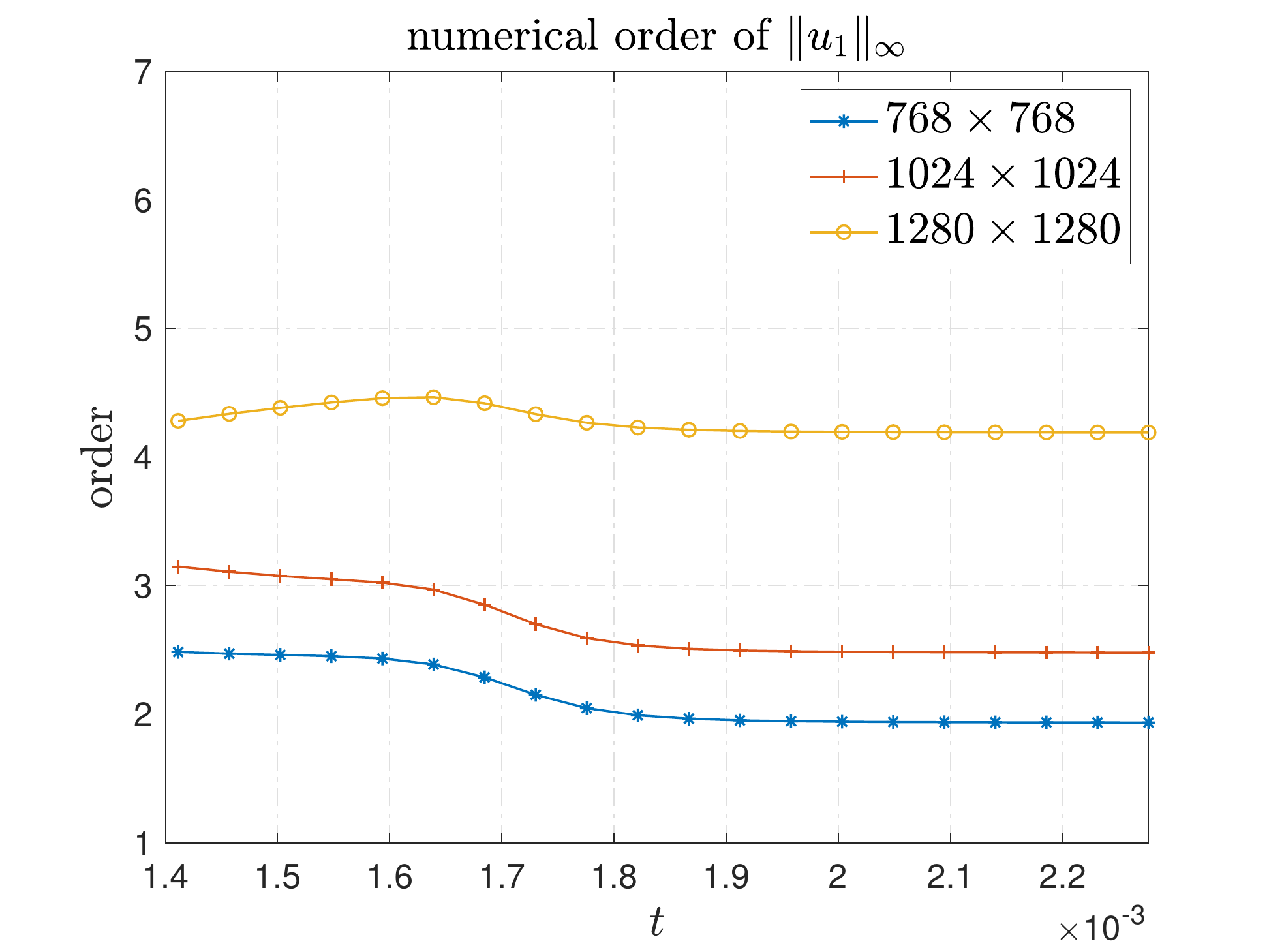}
    \includegraphics[width=0.35\textwidth]{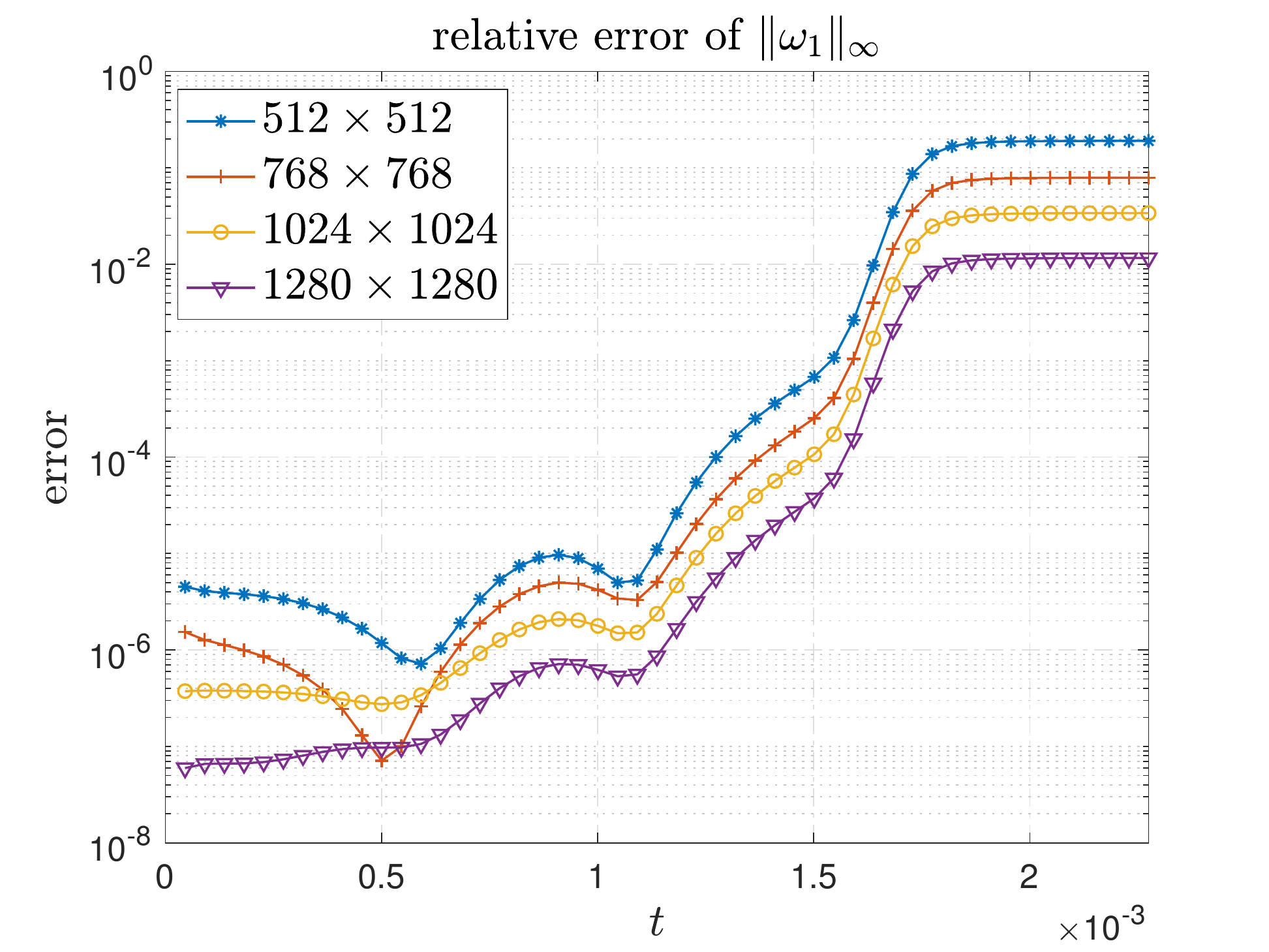}
    \includegraphics[width=0.35\textwidth]{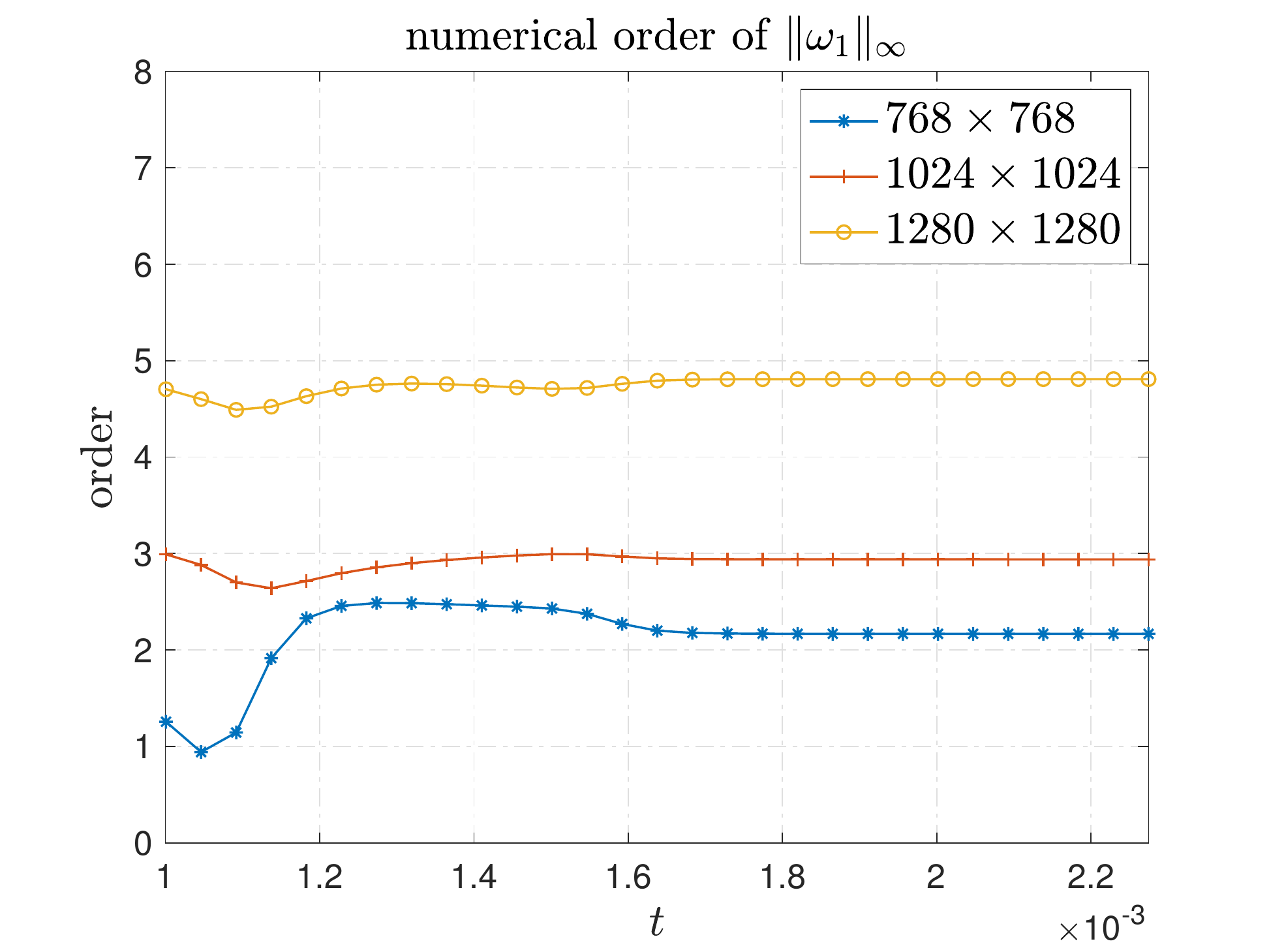}
    \caption[Relative error I]{First row: relative error and numerical order of $\|u_1(t)\|_{L^\infty}$. Second row: relative error and numerical order of $\|\om_1(t)\|_{L^\infty}$. The last time instant shown in the figure is $t= 0.00227648$ for the first row and $t= 0.002274596$ for the second row.}  
    \label{fig:relative_error_1}
\end{figure}

\begin{figure}[!ht]
\centering  
    \includegraphics[width=0.35\textwidth]{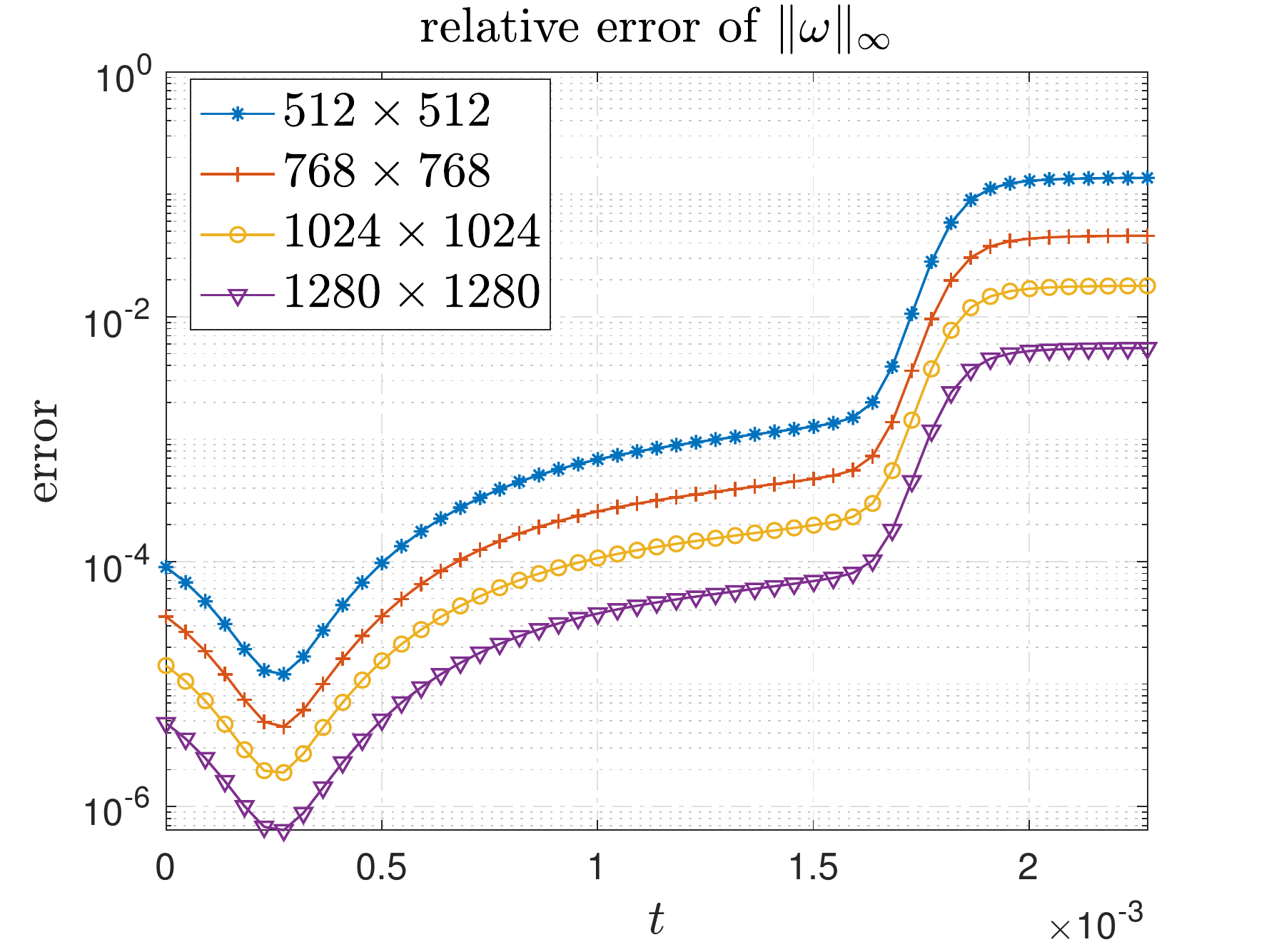}
    \includegraphics[width=0.35\textwidth]{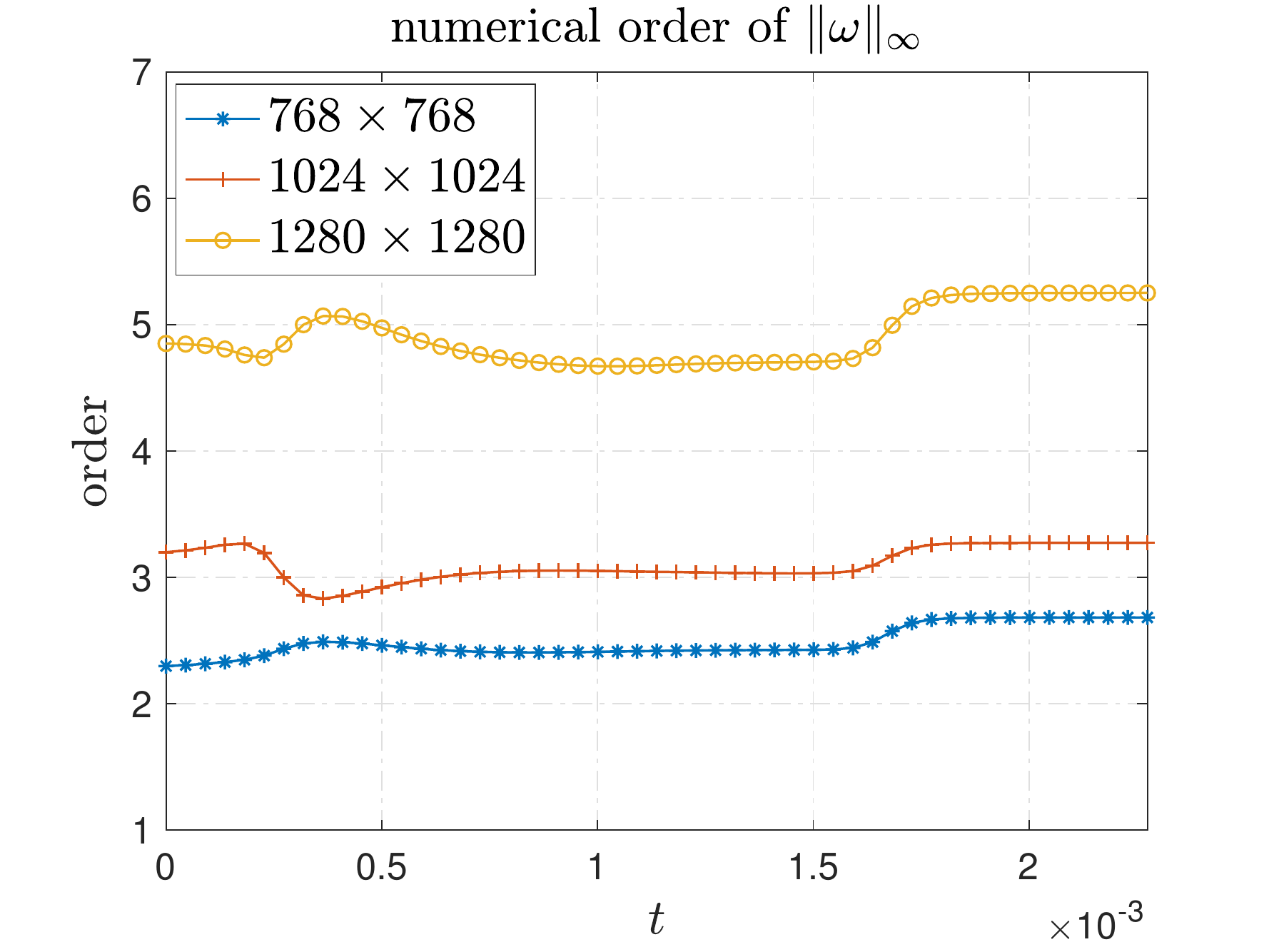} 
      \includegraphics[width=0.35\textwidth]{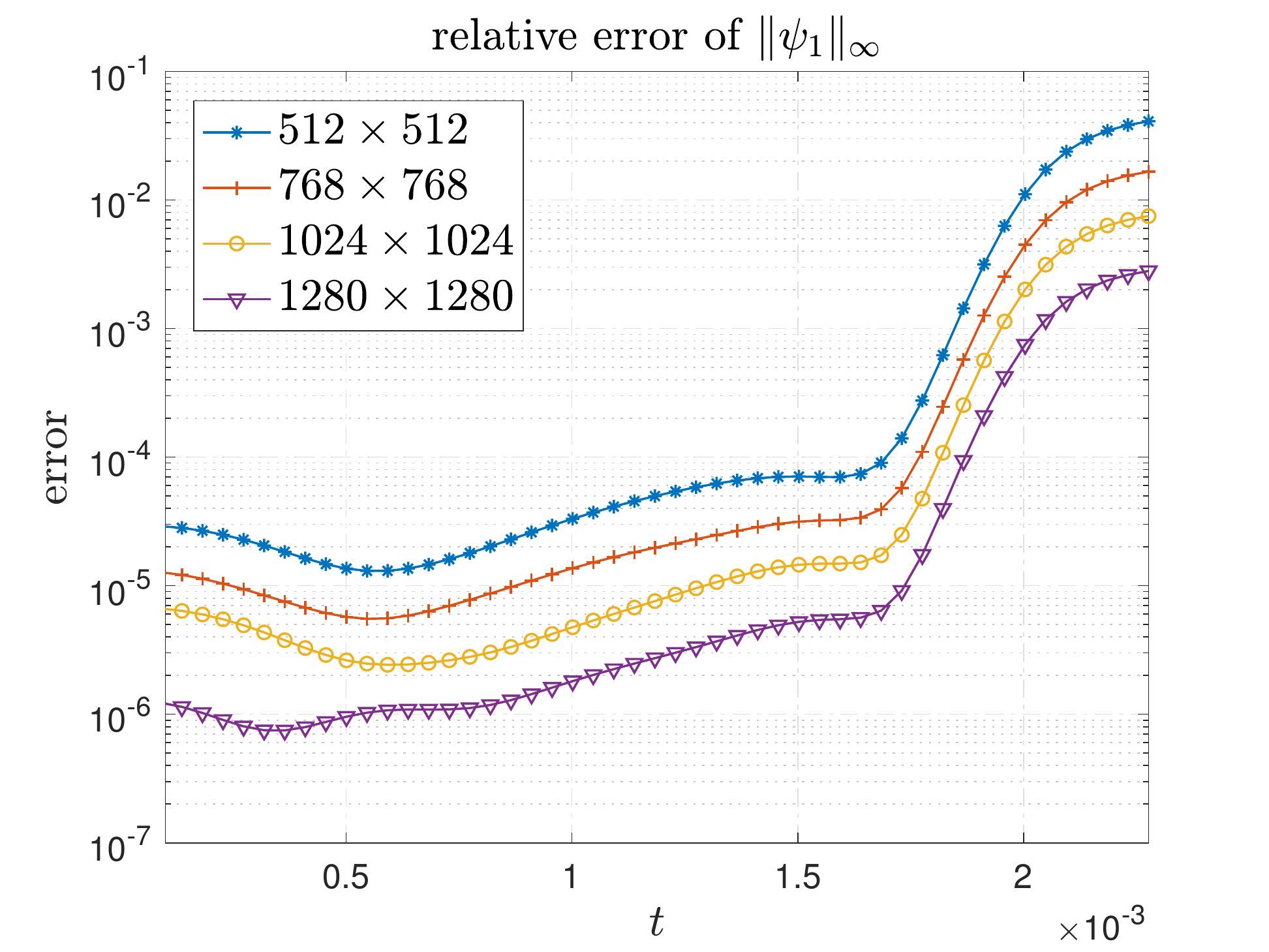}
    \includegraphics[width=0.35\textwidth]{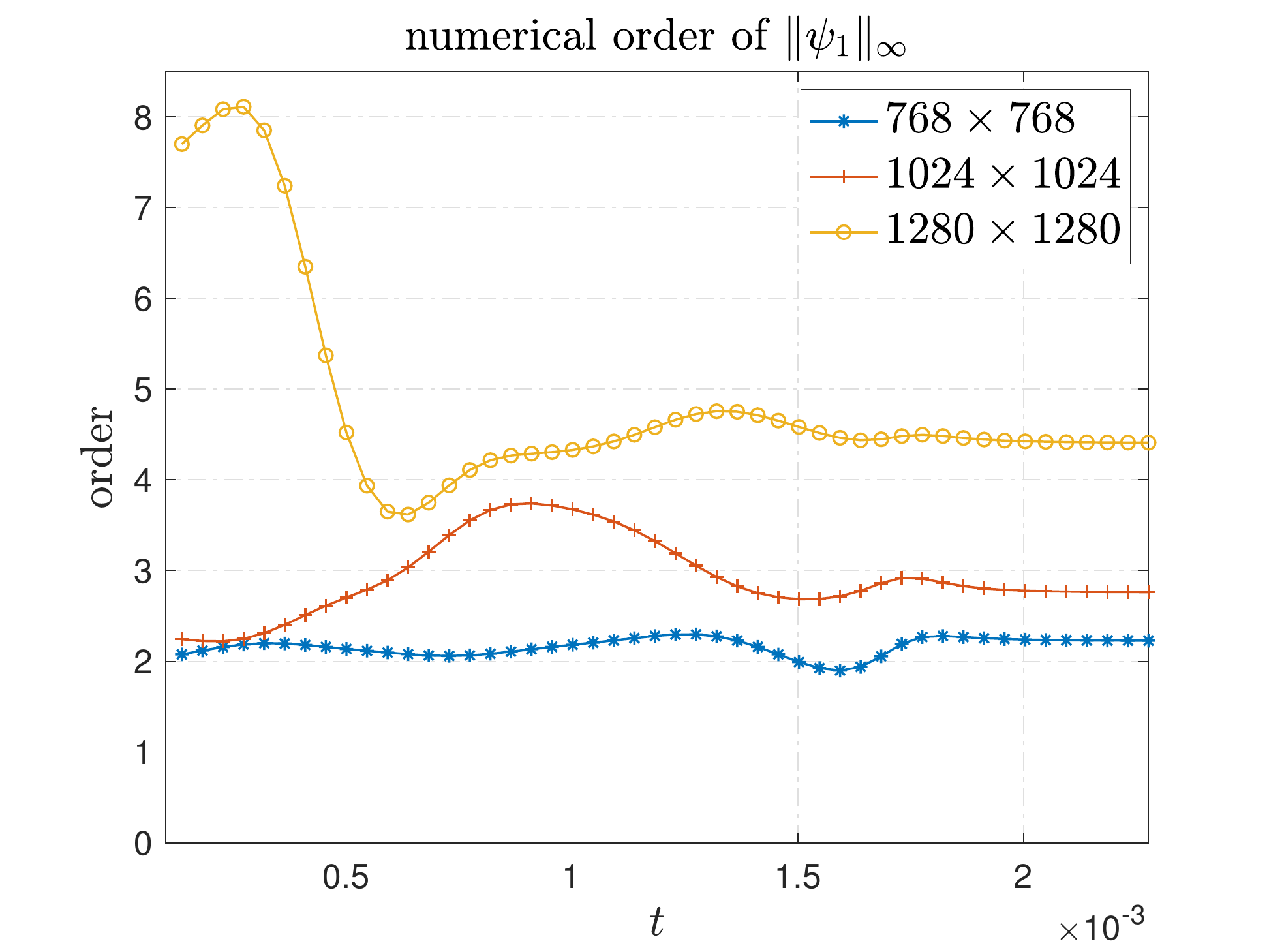} 
    \includegraphics[width=0.35\textwidth]{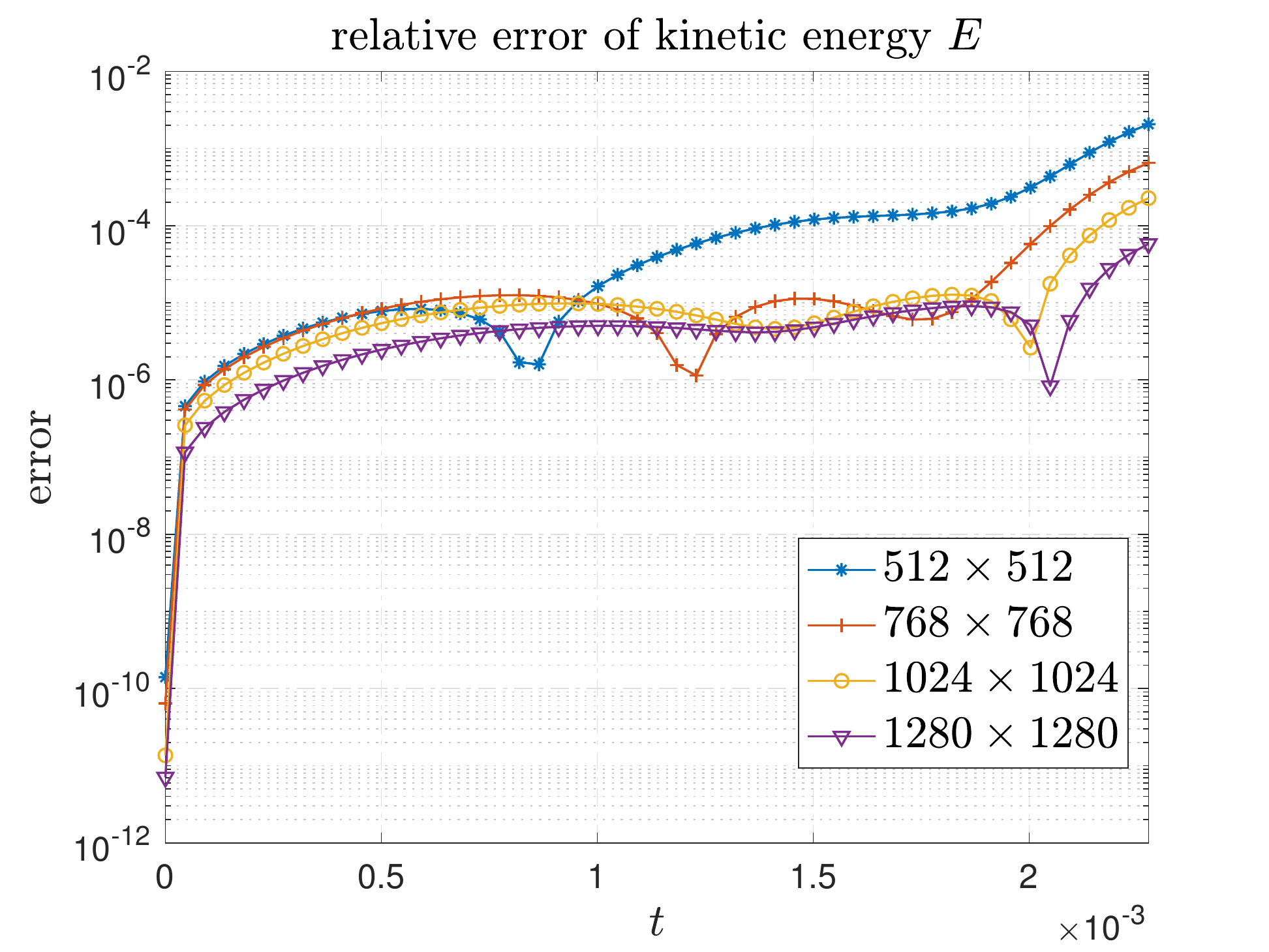}
    \includegraphics[width=0.35\textwidth]{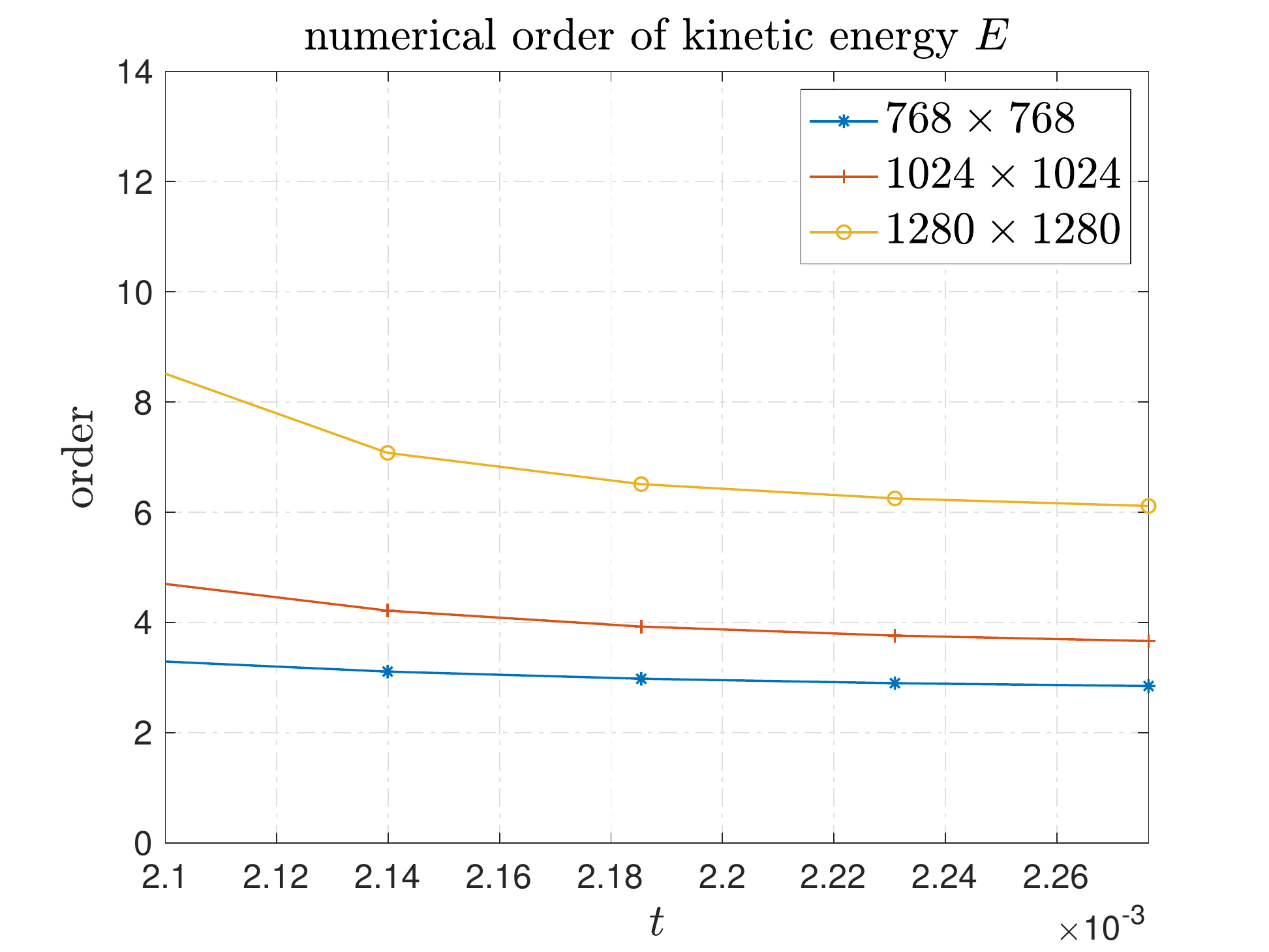}
    \caption[Relative error II]{First row: relative error and numerical order of $\|\vom(t)\|_{L^\infty}$. Second row: relative error and numerical order of $\|\psi_1(t)\|_{L^\infty}$. Third row: relative error and numerical order of $E(t)$. The last time instant shown in the first row is $t= 0.002274596$ and the last time instant in the second and third rows is $t=0.00227648$.}  
    \label{fig:relative_error_2}
       \vspace{-0.05in}
\end{figure}

Next, we study the convergence of some variables as functions of time. In particular, we report the convergence of the quantities $\|u_1\|_{L^\infty}$, $\|\om_1\|_{L^\infty}$, $\|\vom\|_{L^\infty}$, and the kinetic energy $E$. Here the kinetic energy $E$ is given by 
\[E := \frac{1}{2}\int_{\mathcal{D}_1}|\vu|^2 \idiff x = \frac{1}{2}\int_0^1\int_0^{1/2}\left(|u^r|^2+ |u^\theta|^2 + |u^z|^2\right)r\idiff r \idiff z. \]
For smooth solutions, the kinetic energy is a non-increasing function of time: $E(t_1)\leq E(t_2)$ for $t_2\geq t_1\geq0$. Figures~\ref{fig:relative_error_1} and \ref{fig:relative_error_2} plot the relative errors and numerical orders of these quantities as functions of time. The results further confirm that our method is at least $2$nd-order in $h_\rho,h_\eta$. 

We remark that in the early stage of the computation, the solution is quite smooth. As a result, the discretization error is small, the total error may be dominated by the interpolation error when we change from one adaptive mesh to another. The change of mesh tends to happen more frequently for the computation with a finer mesh. Thus it is possible that the total error on a coarser mesh may be smaller than that on a finer mesh in the early stage of the computation, as we can see in the first row of Figure~\ref{fig:relative_error_1}. We can only observe the expected order of accuracy when the discretization error dominates the interpolation error. On the other hand, we also observe an increasing trend in the relative errors of $\|u_1\|_{L^\infty}$, $\|\om_1\|_{L^\infty}$, and $\|\vom\|_{L^\infty}$. This is due to the rapidly decreasing thickness of the sharp front as $t$ approaches the potential singularity time. 

We also observe that the order of convergence seems to be better than the second order accuracy. One possible explanation for this super-convergence phenomenon is that the coefficient in the leading error term may vanish at the position at which $u_1$ or $\omega_1$ is most singular due to some symmetry property of the solution. As we refine the mesh, we seem to obtain a $4th$ order convergence. 

%
%
%

\subsection{Potential Blow-up Scaling Analysis}\label{sec:scaling_study_euler}
In this subsection, we will study the scaling properties of the Euler solution. 
We will provide some qualitative numerical evidences that the $3$D Euler equations seem to develop nearly self-similar scaling properties.

\subsubsection{Fitting of the growth rate}\label{sec:growth_fitting_euler}
We will perform linear fitting to the numerical solutions obtained in our computation. Since the solution has not settled down to a stable phase, we will not try to look for an optimal fitting exponent as we did in \cite{Hou-Huang-2021}. Instead, we will perform some qualitative fitting and use the $R$-Square values of the linear fitting as an indication on the goodness of the linear fitting.

\begin{figure}[!ht]
\centering
	\begin{subfigure}[b]{0.40\textwidth}
    \includegraphics[width=1\textwidth]{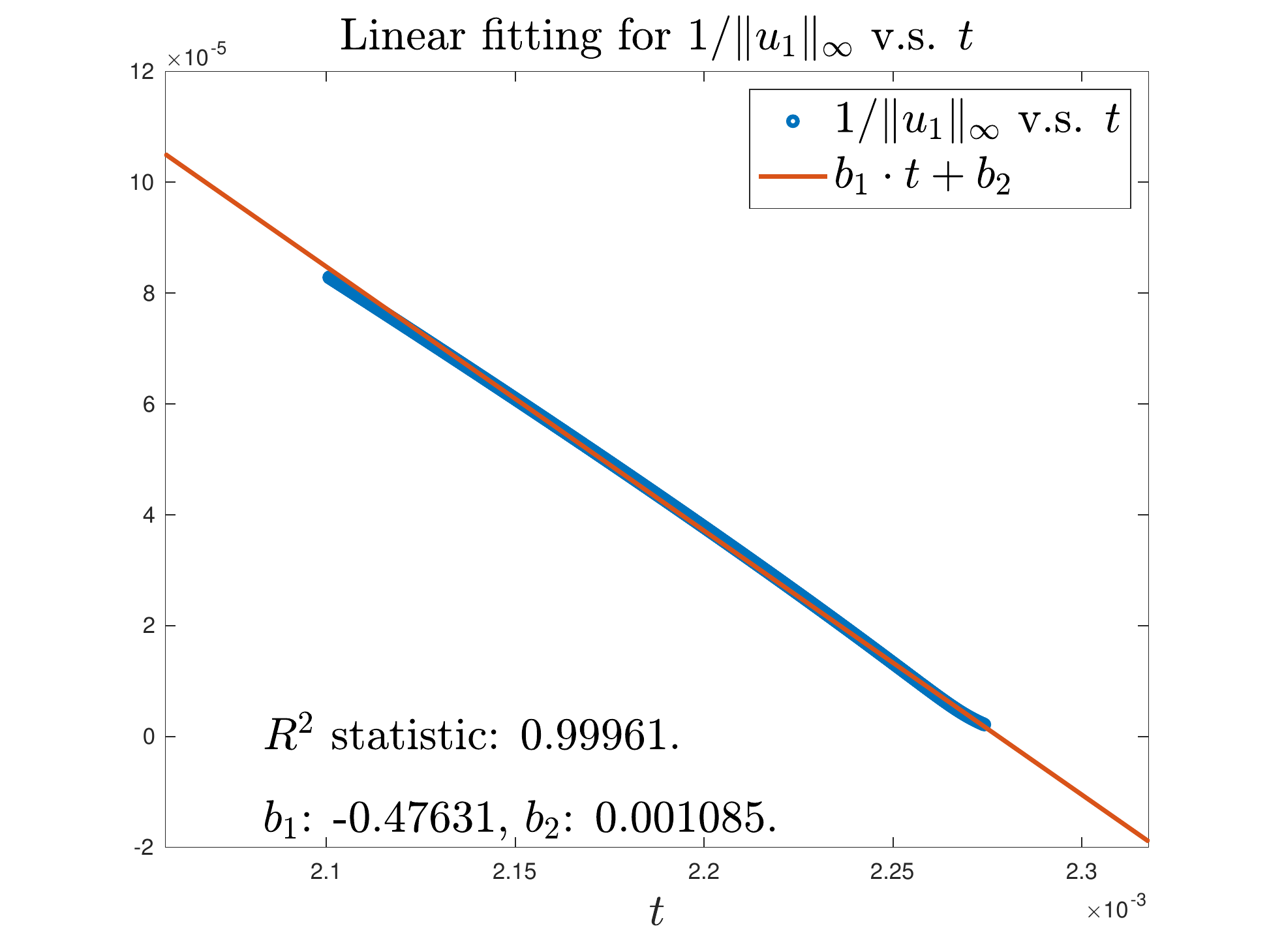}
    \caption{linear regression of $\|u_1(t)\|_{L^\infty}^{-1}$}
    \end{subfigure}
  	\begin{subfigure}[b]{0.40\textwidth} 
    \includegraphics[width=1\textwidth]{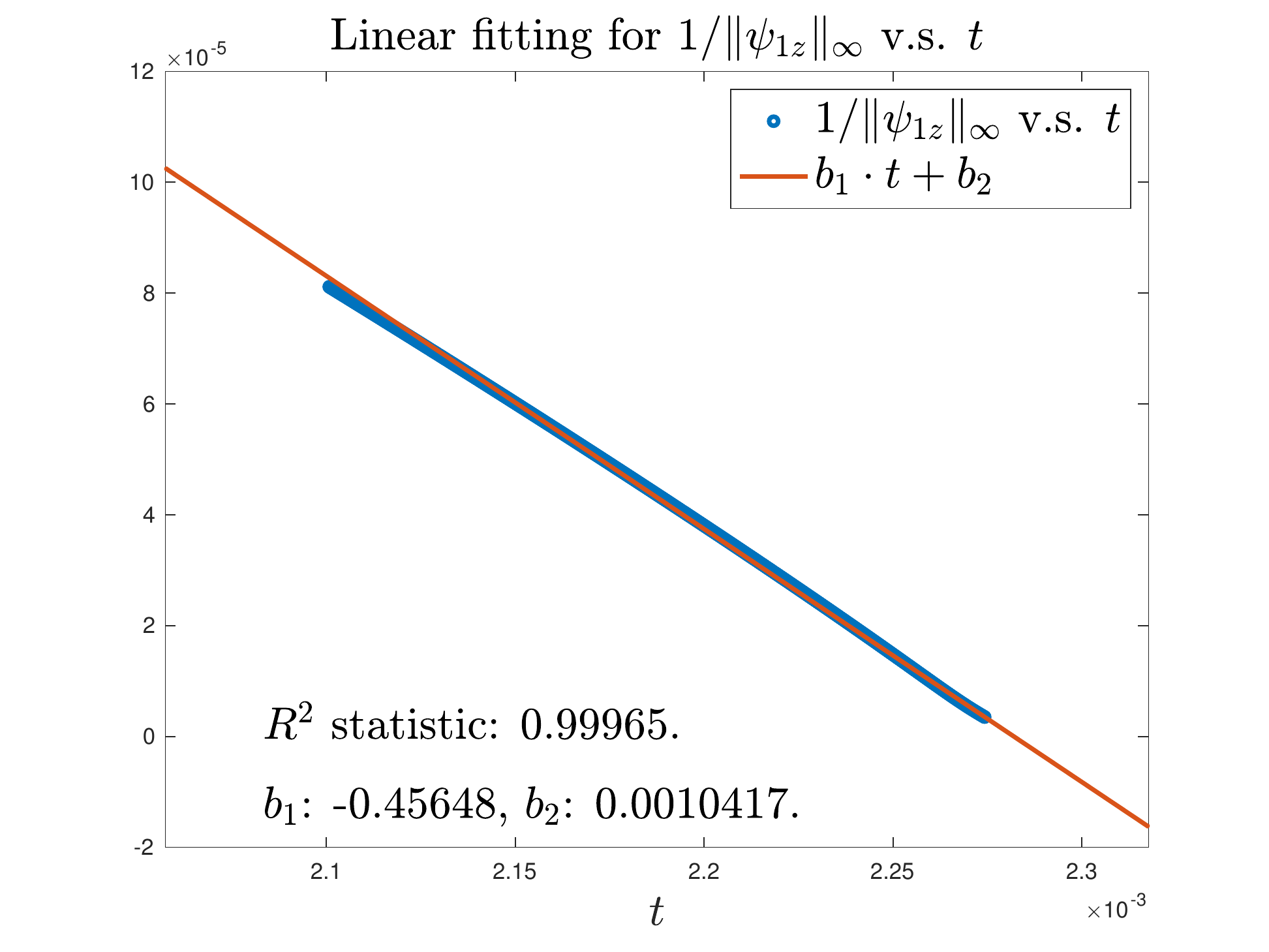}
    \caption{linear regression of $\|\psi_{1z}\|_{L^\infty}^{-1}$}
    \end{subfigure} 
    \caption[Linear regression $u_1$]{The linear regression of (a) $\|u_1\|_{L^\infty}^{-1}$ vs $t$, (b) $\|\psi_{1z}\|_{L^\infty}^{-1}$ vs $t$. 
    The blue points are the data points obtained from our computation, and the red lines are the linear models. }   
    \label{fig:linear_regression_u1}
       \vspace{-0.05in}
\end{figure}

We use the data obtained from $1536\times1536$ resolution. Figure\eqref{fig:linear_regression_u1} shows the fitting results for the quantity $\|u_1(t)\|_{L^\infty}$ and $\|\psi_{1,z}\|_{L^\infty} $ on the time interval $[t_1,t_2] = [0.0021007568,0.0022742813]$. 
We can see that the linear fitting 
$\|u_1(t)\|_{L^\infty}^{-1} \sim (T-t)$ and
$ \|\psi_{1,z}\|_{L^\infty}^{-1}  \sim (T-t)$ have excellent linear fitness with $R$-Square values (we use $R^2$ values for short in the relevant figures) very close to $1$. As we mentioned earlier, we observed a strong positive alignment between $\psi_{1,z}$ and $u_1$ around the maximum location $(R(t),Z(t))$ of $u_1$ in the early stage with $\psi_{1,z}(R(t),Z(t),t)\sim u_1(R(t),Z(t),t) $. Thus, the equation of $\|u_1(t)\|_{L^\infty}$ can be approximated by 
\[\frac{\diff \,}{\diff t} \|u_1(t)\|_{L^\infty} = 2\psi_{1,z}(R(t),Z(t),t)\cdot u_1(R(t),Z(t),t) \sim c_0 \|u_1(t)\|_{L^\infty}^2.\]
This would imply that 
$\|u_1(t)\|_{L^\infty} \sim  (T-t)^{-1} $
for some finite time $T$, which is consistent with our linear fitting results.

\begin{figure}[!ht]
\centering
    \begin{subfigure}[b]{0.38\textwidth}
    \includegraphics[width=1\textwidth]{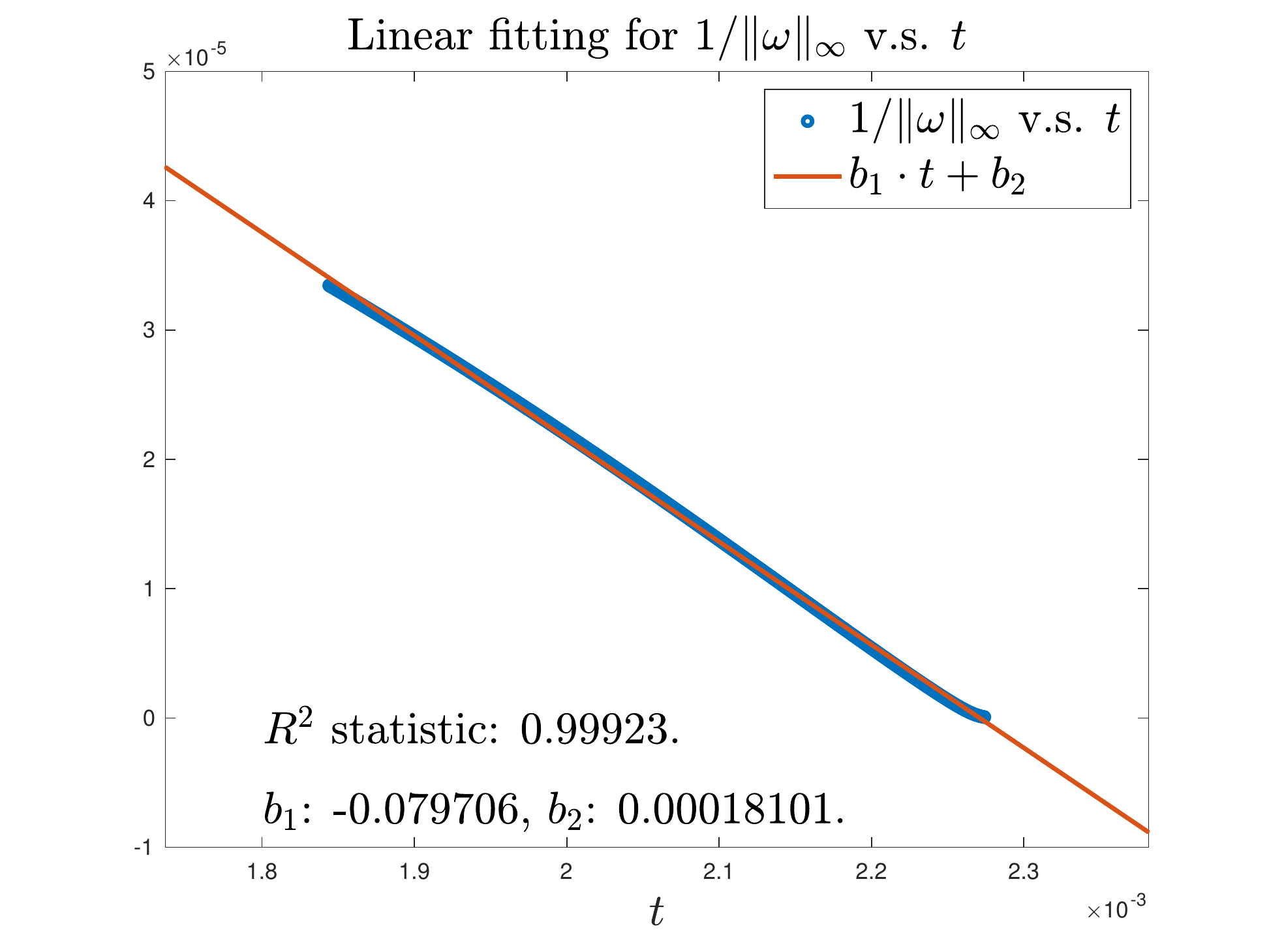}
    \caption{linear regression of $\|\omega (t)\|_{L^\infty}^{-1}$}
    \end{subfigure}
  	\begin{subfigure}[b]{0.38\textwidth} 
    \includegraphics[width=1\textwidth]{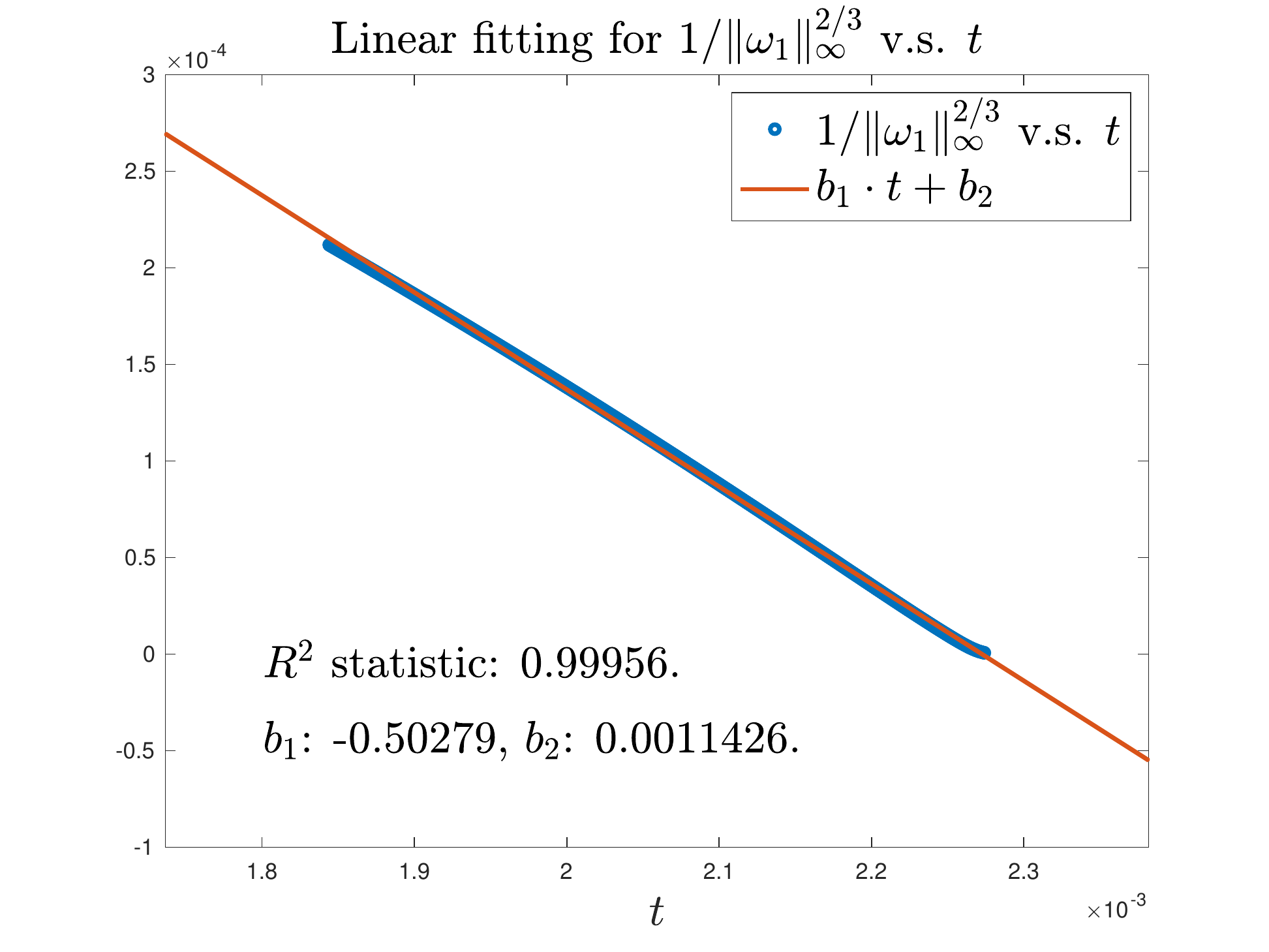}
    \caption{linear regression of $\|\omega_1(t)\|_{L^\infty}^{-2/3}$}
    \end{subfigure} 
    \caption[Linear regression $\omega$]{The linear regression of (a) $\|\omega(t)\|_{L^\infty}^{-1}$ vs $t$, (b) $\|\omega_1(t)\|_{L^\infty}^{-2/3}$ vs $t$.
    The blue points are the data points obtained from our computation, and the red lines are the linear models. }   
    \label{fig:linear_regression_w1}
       \vspace{-0.05in}
\end{figure}

\begin{figure}[!ht]
\centering
    \begin{subfigure}[b]{0.38\textwidth}
    \includegraphics[width=1\textwidth]{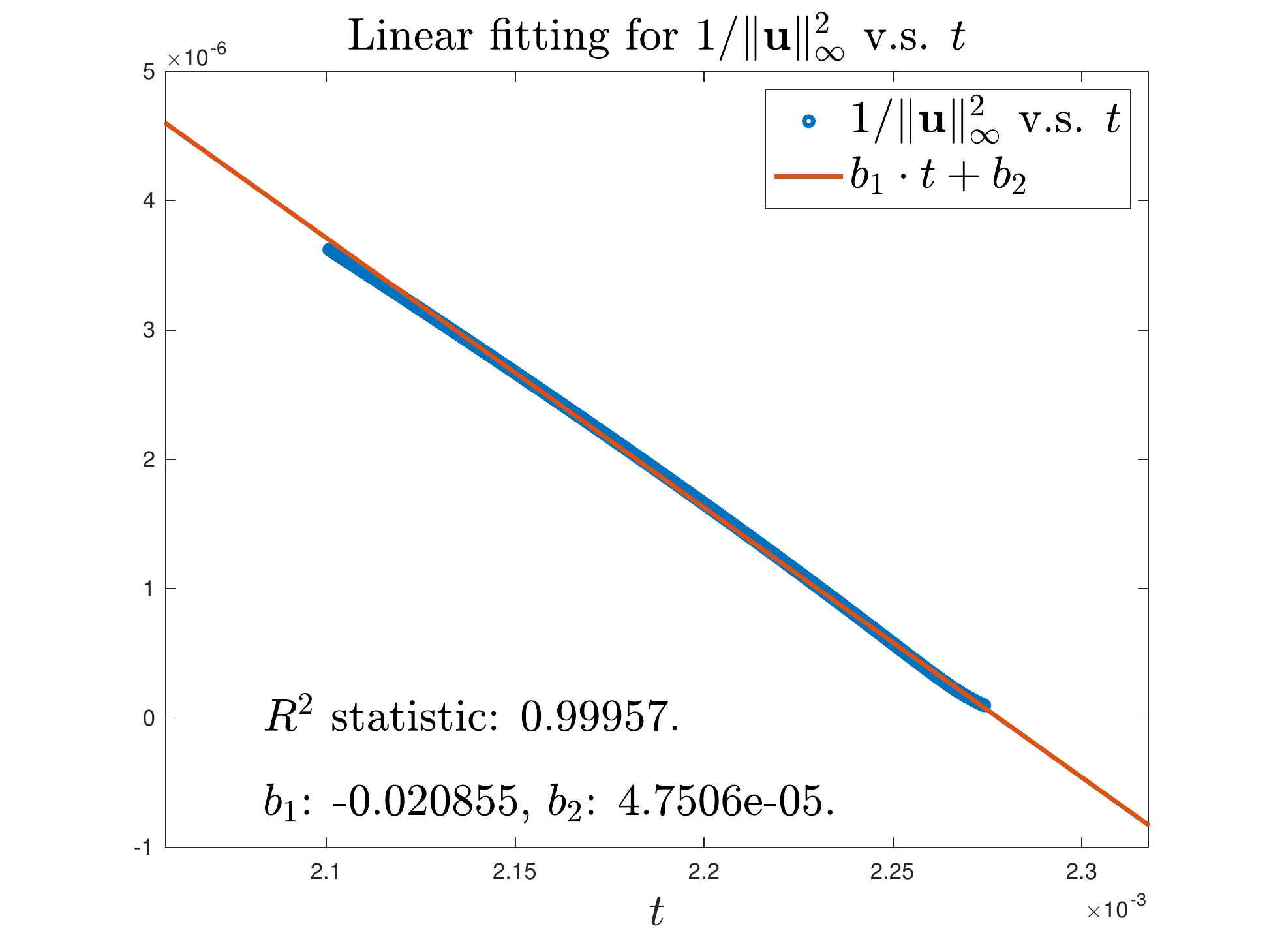}
    \caption{linear regression of $\|{\bf u} (t)\|_{L^\infty}^{-2}$}
    \end{subfigure}
  	\begin{subfigure}[b]{0.38\textwidth} 
    \includegraphics[width=1\textwidth]{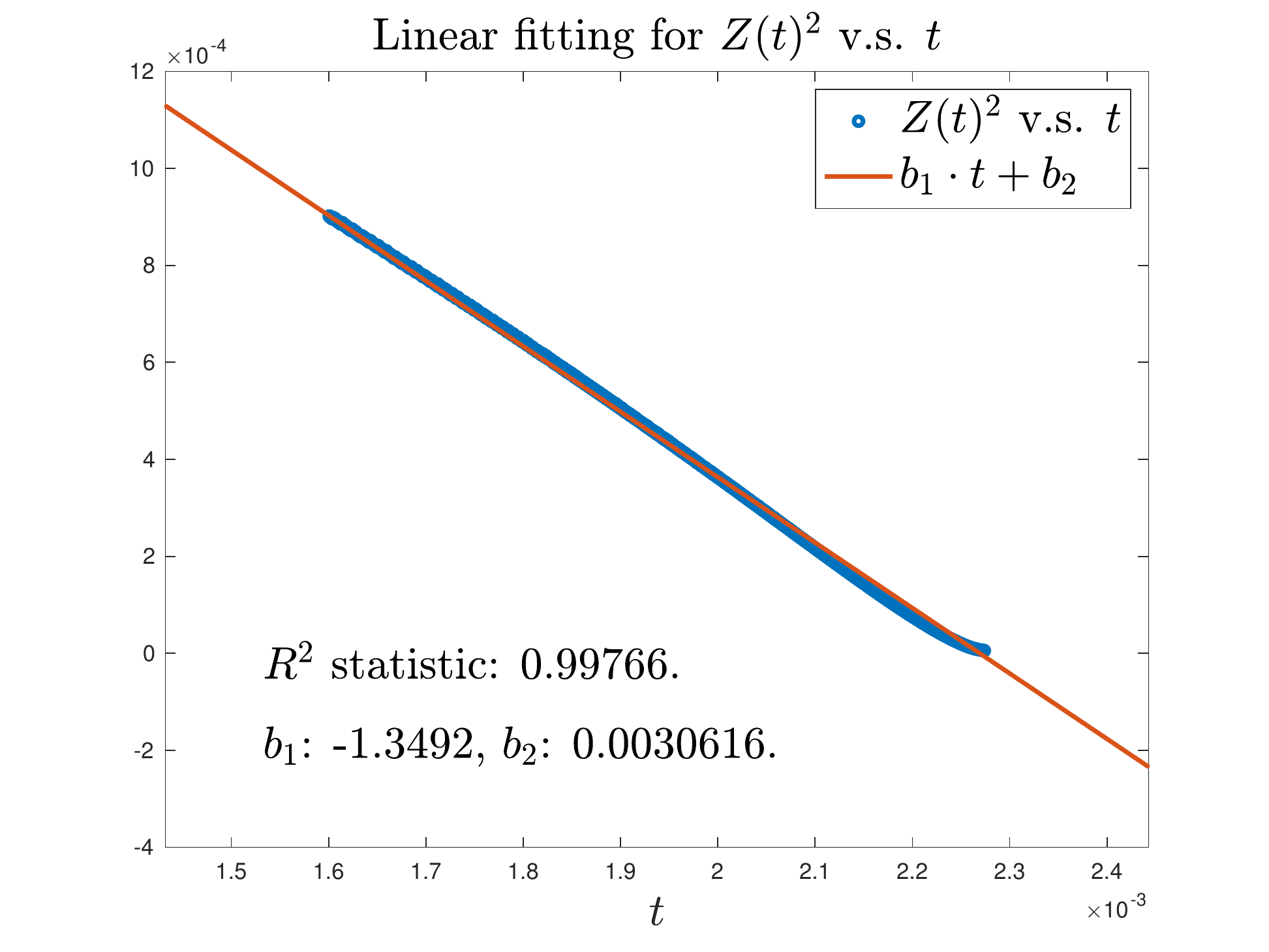}
    \caption{linear regression of $Z(t)^2$}
    \end{subfigure} 
    \caption[Linear regression ${\bf u}$]{The linear regression of (a) $\|{\bf u}(t)\|_{L^\infty}^{-2}$ vs $t$, (b) $Z(t)^2$ vs $t$.
    The blue points are the data points obtained from our computation, and the red lines are the linear models. }   
    \label{fig:linear_regression_vel}
       \vspace{-0.05in}
\end{figure}

Next, we study the growth of the maximum vorticity $\|\vom (t)\|_{L^\infty}$. Figure \ref{fig:linear_regression_w1} (a) shows the linear fitting of $\|\vom(t)\|_{L^\infty}^{-1}$ as a function of time on the time interval $[t_1,t_2] = [0.0018441297, 0.0022742812]$. 
Again, we see that $ \|\vom(t)\|_{L^\infty}^{-1} \sim (T-t)$ has good linear fitness with $R$-Square values very close to $1$, which provides evidences of the finite-time blow-up of $\|\vom (t)\|_{L^\infty}$ in the form of an inverse power law
\[\|\vom(t)\|_{L^\infty} \sim (T-t)^{-1}.\] 
This scaling property is consistent with the rapid growth of $\int_0^t \| \vom (s)\|_{L^\infty} ds $ as we demonstrated in Figure \ref{fig:rapid_growth2}(b).
Similarly, Figure \ref{fig:linear_regression_w1} (b) shows the linear fitting of $\|\omega_1(t)\|_{L^\infty}^{-2/3}$ on the same time interval. We see that $ \|\omega_1(t)\|_{L^\infty} = O(1/(T-t)^{3/2})$ has good fitness with $R$-Square values close to $1$.

To further illustrate the existence of a potential finite-time blow-up, we perform the linear fitting for the maximum velocity and $Z(t)$. For the fitting of the maximum velocity, the fitting time interval is the same as that for $u_1$ and $\psi_{1z}$, i.e. $[t_1,t_2] = [0.0021007568,0.0022742813]$.
For $Z(t)$, we use a larger time interval $[t_1,t_2] = [0.0016006384,0.0022742813]$. 
In Figure \ref{fig:linear_regression_vel} (a), we observe that $ \|{\bf u}(t)\|_{L^\infty} = O(1/(T-t)^{1/2})$ has good linear fitness with $R$-Square values very close to $1$. This implies that $\|{\bf u} (t)\|_{L^\infty}$ has an inverse power law
\[\|{\bf u}(t)\|_{L^\infty} \sim (T-t)^{-1/2}.\] 

The scaling properties of the maximum vorticity and maximum velocity seem to suggest that the spatial length scale of the potential blow-up, which is characterized by $Z(t)$, should be $Z(t) \sim (T-t)^{1/2}$. Indeed, Figure \ref{fig:linear_regression_vel}(b), we observe that $Z(t) = O((T-t)^{1/2})$ has reasonably good linear fitness although the $R$-Square values are not as close to $1$ as in the other quantities. This is due to the fact that the solution experiences a phase transition in the late stage as the thickness of the sharp front becomes smaller and smaller. We can also see that the deviation from the linear fitting mainly occurs near the end. The fitting in the early time is much better.

\begin{figure}[!ht]
\centering
  	\begin{subfigure}[b]{0.38\textwidth} 
    \includegraphics[width=1\textwidth]{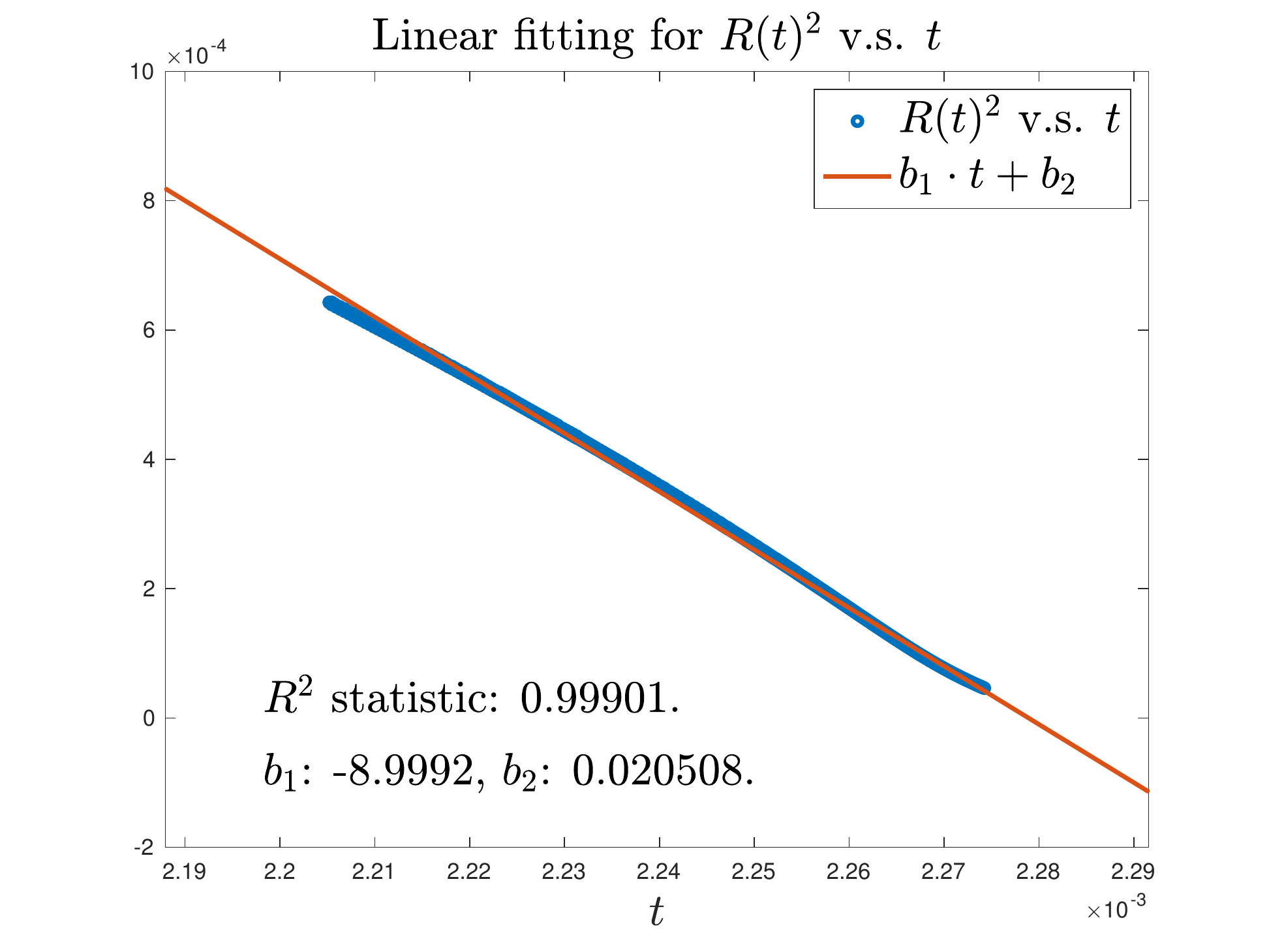}
    \caption{linear regression of $R(t)^2$}
    \end{subfigure} 
    \begin{subfigure}[b]{0.38\textwidth} 
    \includegraphics[width=1\textwidth]{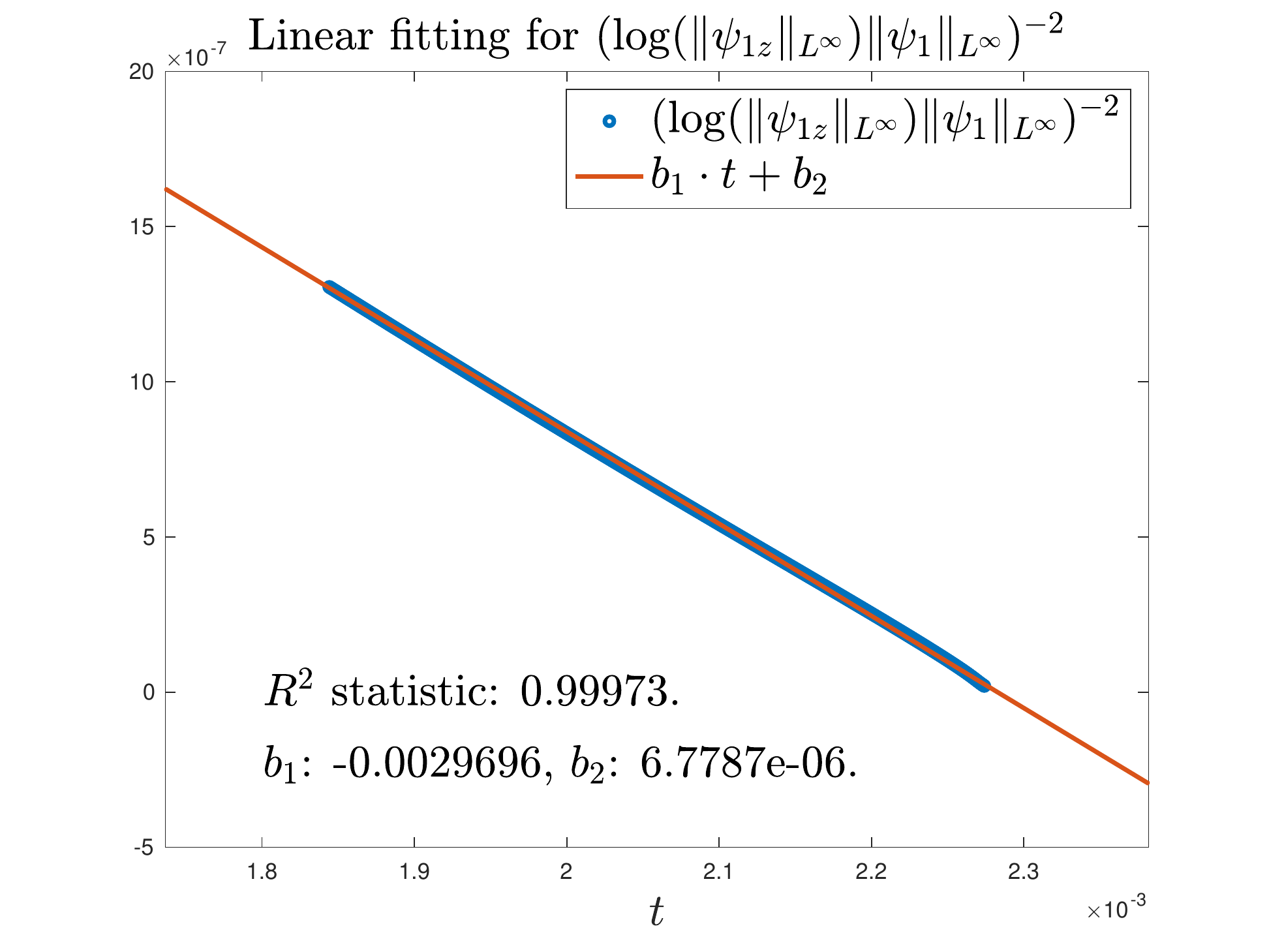}
    \caption{linear fitting of $(\log(\|\psi_{1z}\|_{L^\infty})\|\psi_1\|_{L^\infty})^{-2}$}
    \end{subfigure} 
    \caption[Linear regression ${\bf u}$]{The linear regression of (a) $R(t)^2$ vs $t$. (b) Linear regression of $(\log(\|\psi_{1z}\|_{L^\infty})\|\psi_1\|_{L^\infty})^{-2}$.
    The blue points are the data points obtained from our computation, and the red lines are the linear models. }   
    \label{fig:linear_regression_R}
       \vspace{-0.05in}
\end{figure}

We remark that we have not been able to obtain good linear fitting of $R(t)^2$ in the same time interval in which we obtained good linear fitting for $Z(t)^2$. This is possibly due to the fact that the solution of the $3$D Euler equations have not settled down to a stable phase.  In Figure \ref{fig:linear_regression_R} (a), we show the linear fitting of $R(t)^2$ in a smaller time interval, $[0.0022052019, 0.0022742813]$.
During this time interval, we do observe good linear fitting of $R(t)^2$ with $R$-Square values close to $1$. In Figure \ref{fig:linear_regression_R} (b), we plot the linear fitting for $(\log(\|\psi_{1z}\|_{L^\infty})\|\psi_1\|_{L^\infty})^{-2}$ in the same time interval as we did for the fitting of $\omega$ and $\omega_1$, i.e. $[t_1,t_2] = [0.0018441297, 0.0022742812]$. We can see that 
\[
\|\psi_1\|_{L^\infty} \sim \frac{1}{(T-t)^{1/2} |\log(T-t)|},
\]
gives good linear fitting with $R$-Square values close to $1$. The fact that $R(t)^2$ does not enjoy a good linear fitting on the same time interval as $Z(t)^2$ did and the logarithmic correction in the fitting of $\|\psi_1\|_{L^\infty}$ seem to suggest that the potential blow-up of the Euler equations is only nearly self-similar.

\vspace{-0.05in}
\subsubsection{Numerical evidence of locally self-similar profiles}\label{sec:evidence_of_self-similar_euler}
In addition to performing linear fitting of some variables in the previous subsection, we will further investigate the nearly self-similar solution by comparing the properly normalized profiles of the solution at different time instants. We will focus our study on the solution profile in a small neighborhood centered around $(R(t),Z(t))$.

As in \cite{Hou-Huang-2021}, we propose the following self-similar ansatz in the axisymmetric setting: 
\begin{subequations}\label{eq:self-similar_ansatz}
\begin{align}
u_1(r,z,t) &\sim (T-t)^{-c_{u}}\overline{U}\left(t, \frac{r-R(t)}{(T-t)^{c_l}}\,,\,\frac{z}{(T-t)^{c_l}}\right), \label{eq:self-similar_u1}\\
\om_1(r,z,t) &\sim (T-t)^{-c_{\om}}\overline{\Omega}\left(t,\frac{r-R(t)}{(T-t)^{c_l}}\,,\,\frac{z}{(T-t)^{c_l}}\right), \label{eq:self-similar_w1}\\
\psi_1(r,z,t) &\sim (T-t)^{-c_{\psi}}\overline{\Psi}\left(t,\frac{r-R(t)}{(T-t)^{c_l}}\,,\,\frac{z}{(T-t)^{c_l}}\right), \label{eq:self-similar_psi1}\\
R(t) &\sim (T-t)^{c_s}R_0. \label{eq:self-similar_R}
\end{align}
\end{subequations}
Here $\overline{U},\overline{\Omega},\overline{\Psi}$ denote the self-similar profiles of $u_1,\om_1,\psi_1$ respectively. By this ansatz, the solution develops locally nearly self-similar blow-up centered at the point $(R(t),0)$ that travels toward the origin with scaling $(T-t)^{c_l}$. 

By slightly modifying the asymptotic scaling analysis developed in \cite{Hou-Huang-2021} for variable diffusion coefficients to the Euler equations, one can obtain a set of governing equations for the evolution of the rescaled profiles 
$\overline{U},\overline{\Omega},\overline{\Psi}$ by further introducing a rescaled time $\tau$. From our numerical study of the scaling properties of the solution, we have roughly $c_l=1/2$. This in turns gives $c_u = 1$, $c_\omega = 3/2$, and $c_\psi = 1/2$. From these scaling properties, we obtain 
$\| \omega\|_{L^\infty} \sim \frac{1}{T-t}$, which is consistent with what we observed numerically. It is important to note that the scaling exponent $c_l=1/2$ is consistent with the scaling properties of the Navier--Stokes equations with a constant viscosity coefficient. 

We remark that such scaling analysis would not be able to account for the potential logarithmic correction in $\psi_1$. It is likely that the dynamic rescaling equations that govern the evolution of the rescaled profiles may not converge to a steady state as the rescaled time tends to infinity. Instead it may oscillate around an approximate self-similar profile. If this is the case, we will not have an asymptotically self-similar blow-up of the $3$D Euler equations. Instead, the solution of the dynamic rescaled Euler equations is close to an approximate self-similar profile with scaling properties qualitatively similar to those of an asymptotic self-similar blow-up. 

\begin{figure}[!ht]
\centering
  	\begin{subfigure}[b]{0.32\textwidth} 
    \includegraphics[width=1\textwidth]{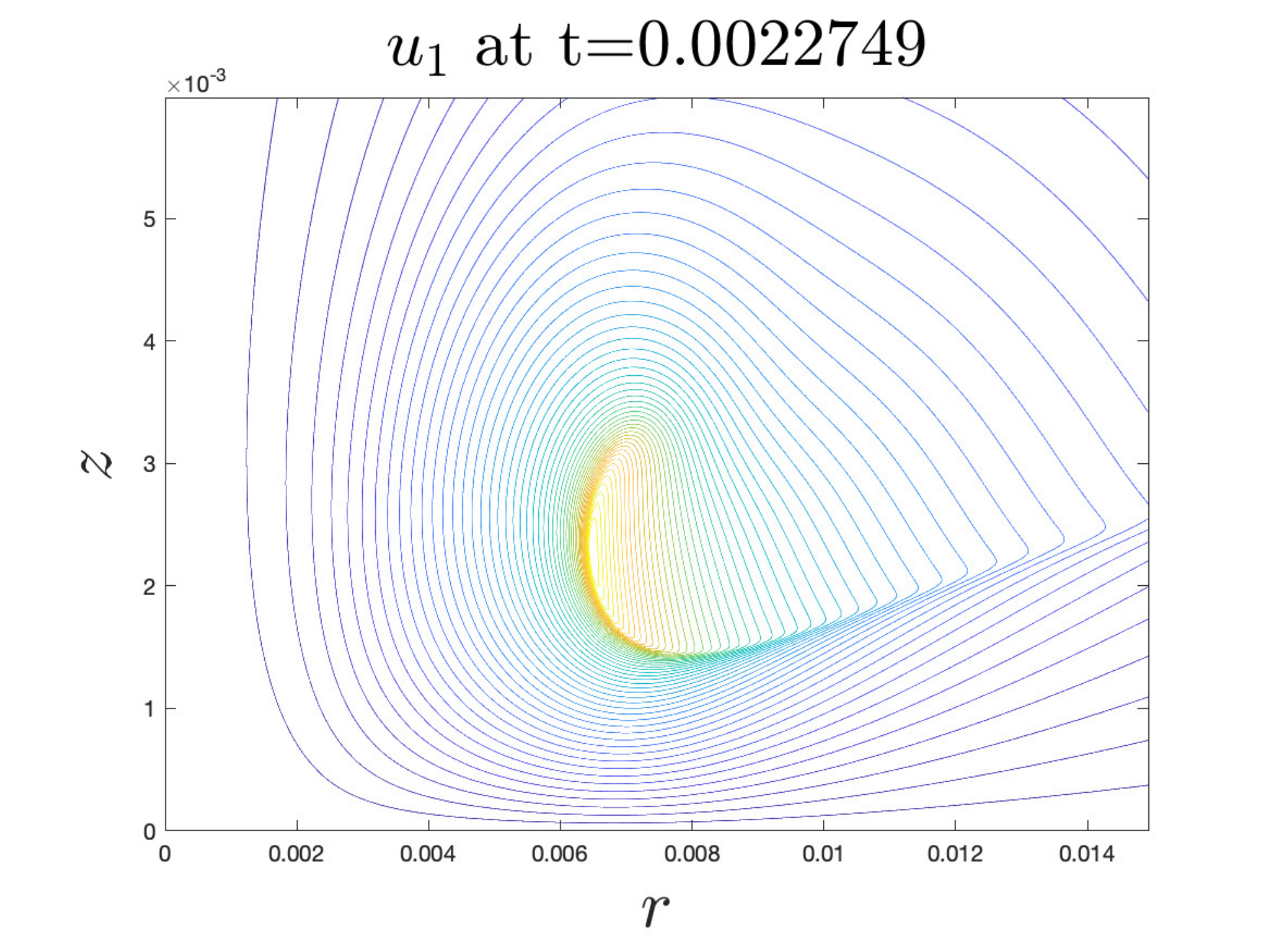}
    \end{subfigure} 
    \begin{subfigure}[b]{0.32\textwidth} 
    \includegraphics[width=1\textwidth]{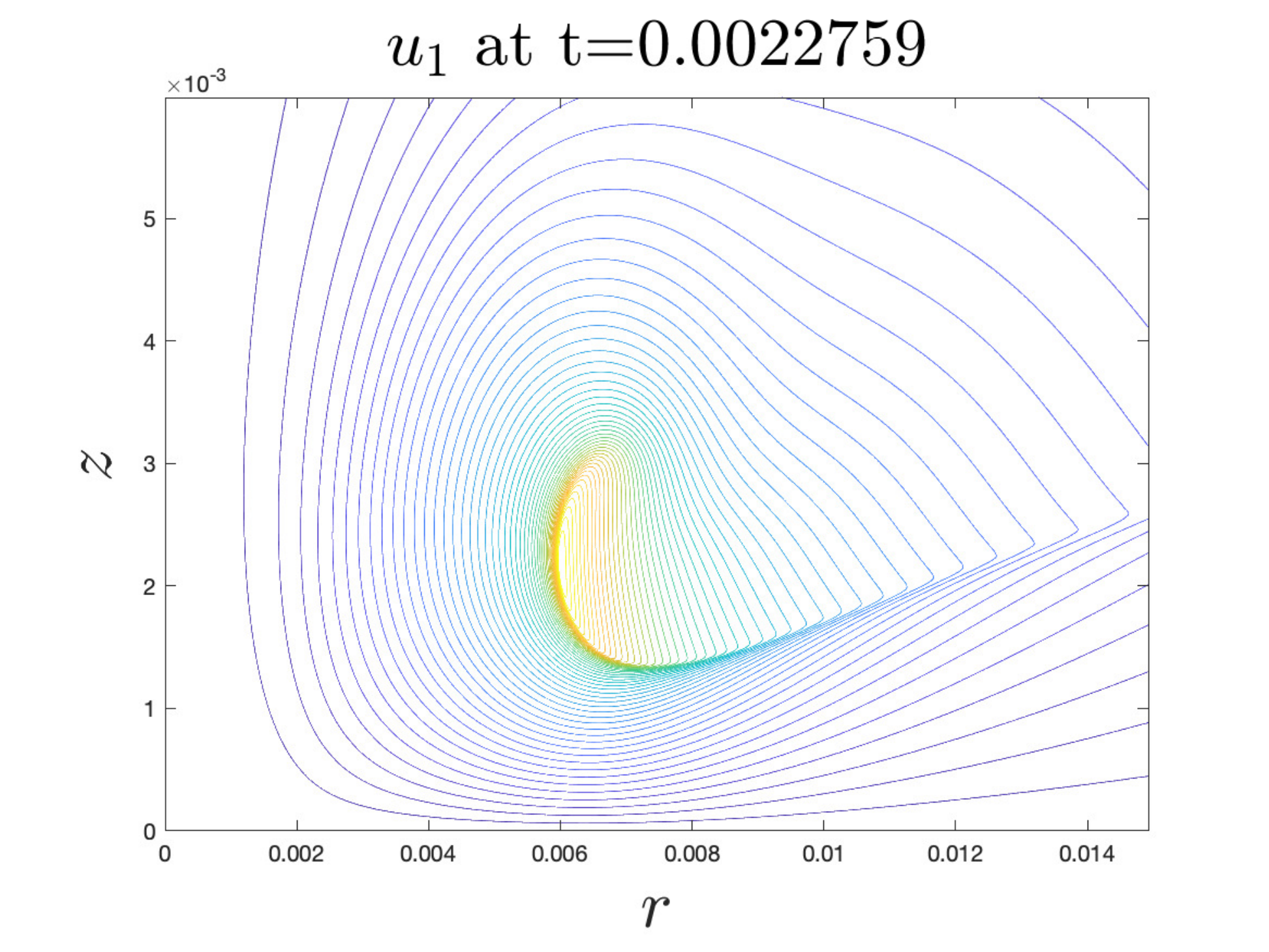}
    \end{subfigure} 
    \begin{subfigure}[b]{0.32\textwidth} 
    \includegraphics[width=1\textwidth]{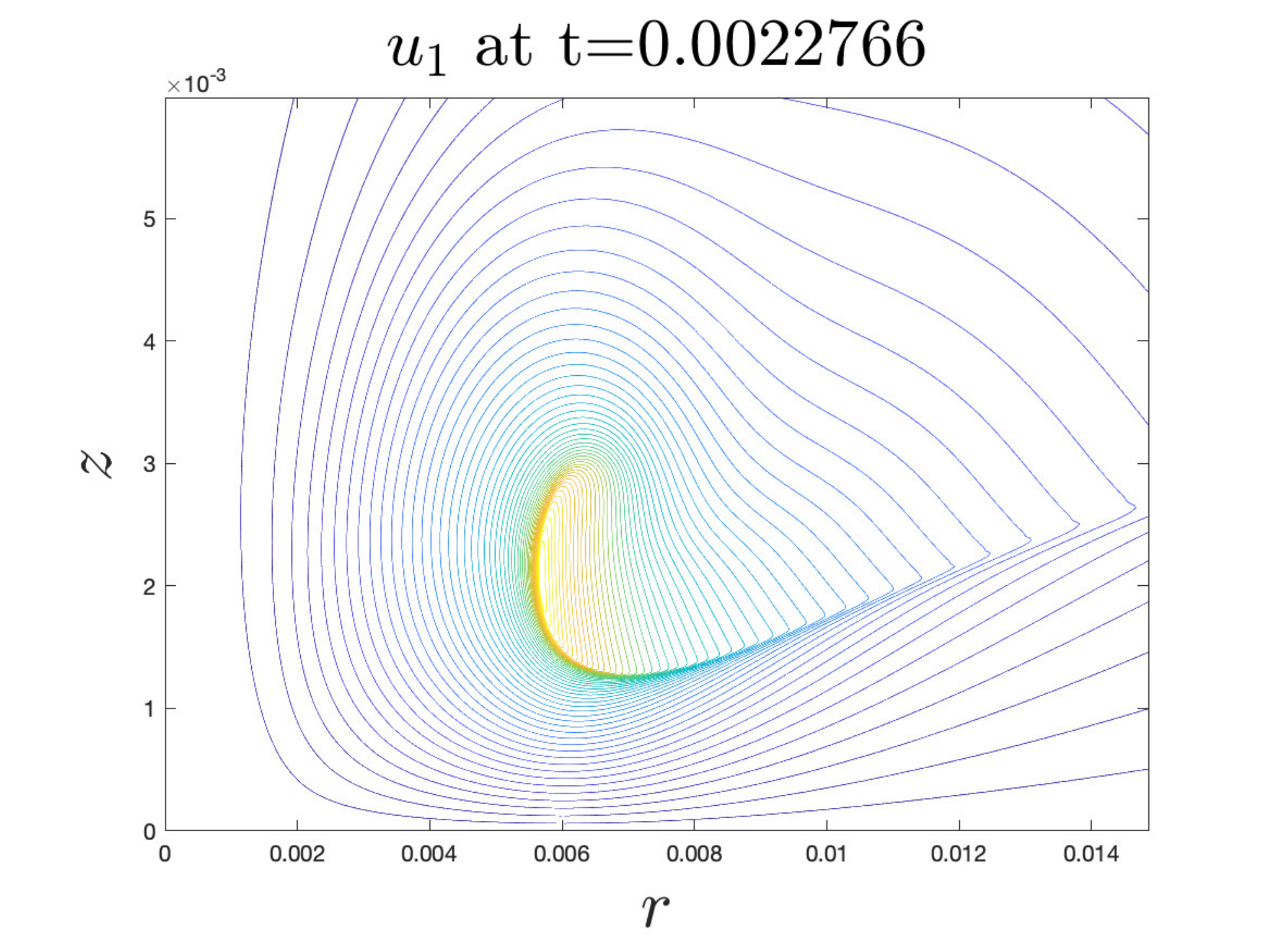}
    \end{subfigure} 
    \begin{subfigure}[b]{0.32\textwidth} 
    \includegraphics[width=1\textwidth]{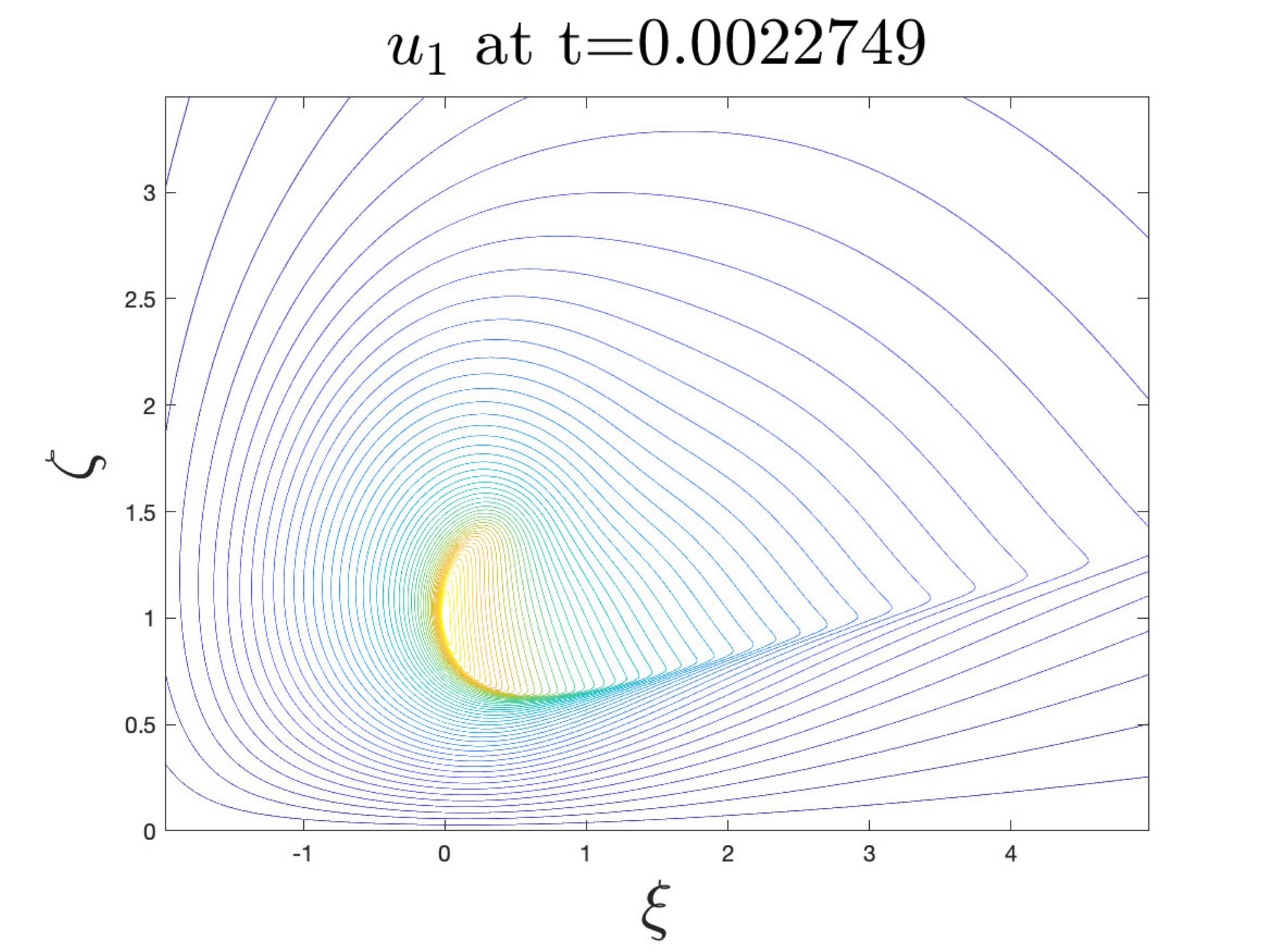}
    \end{subfigure} 
    \begin{subfigure}[b]{0.32\textwidth} 
    \includegraphics[width=1\textwidth]{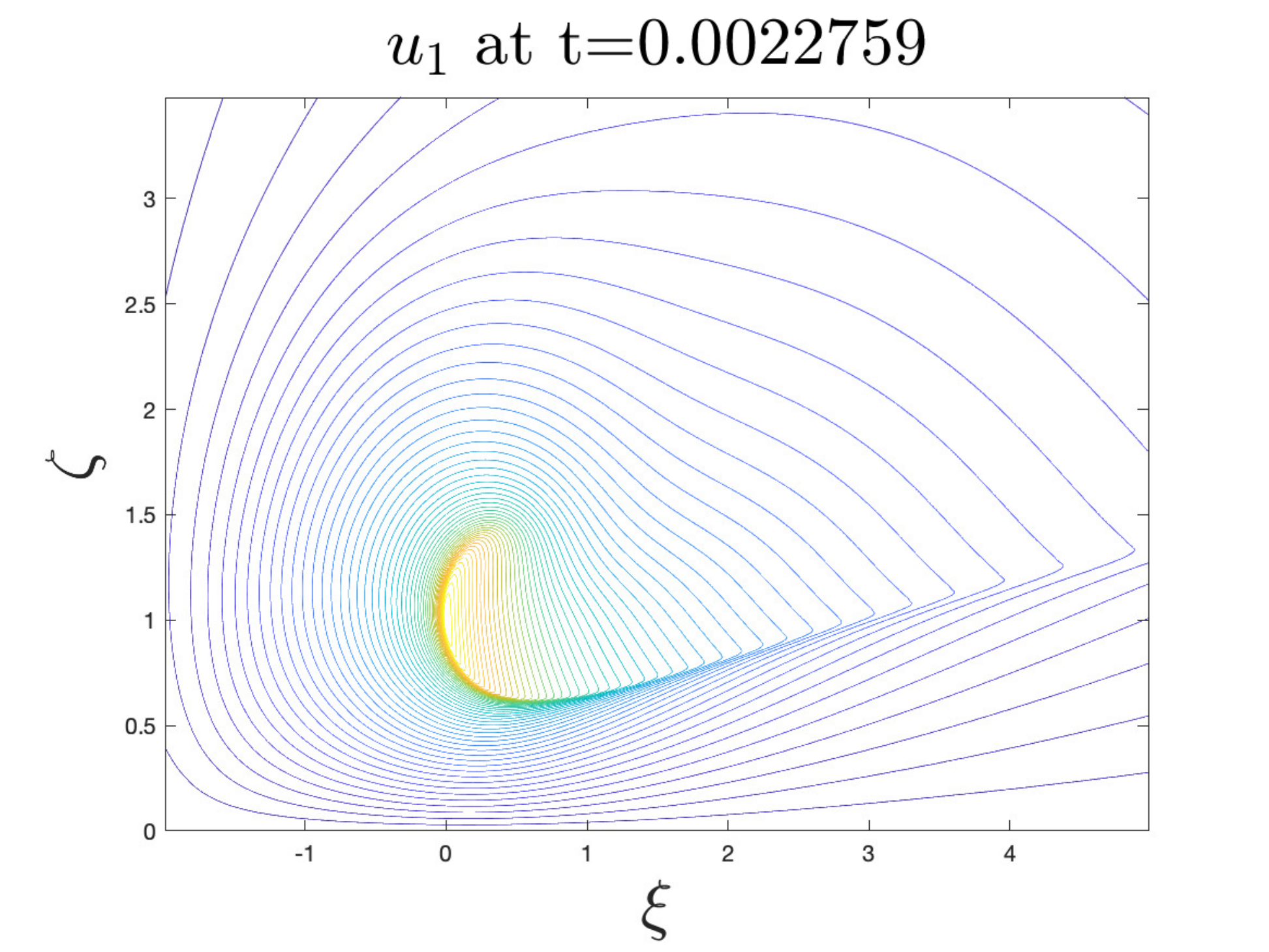}
    \end{subfigure} 
    \begin{subfigure}[b]{0.32\textwidth} 
    \includegraphics[width=1\textwidth]{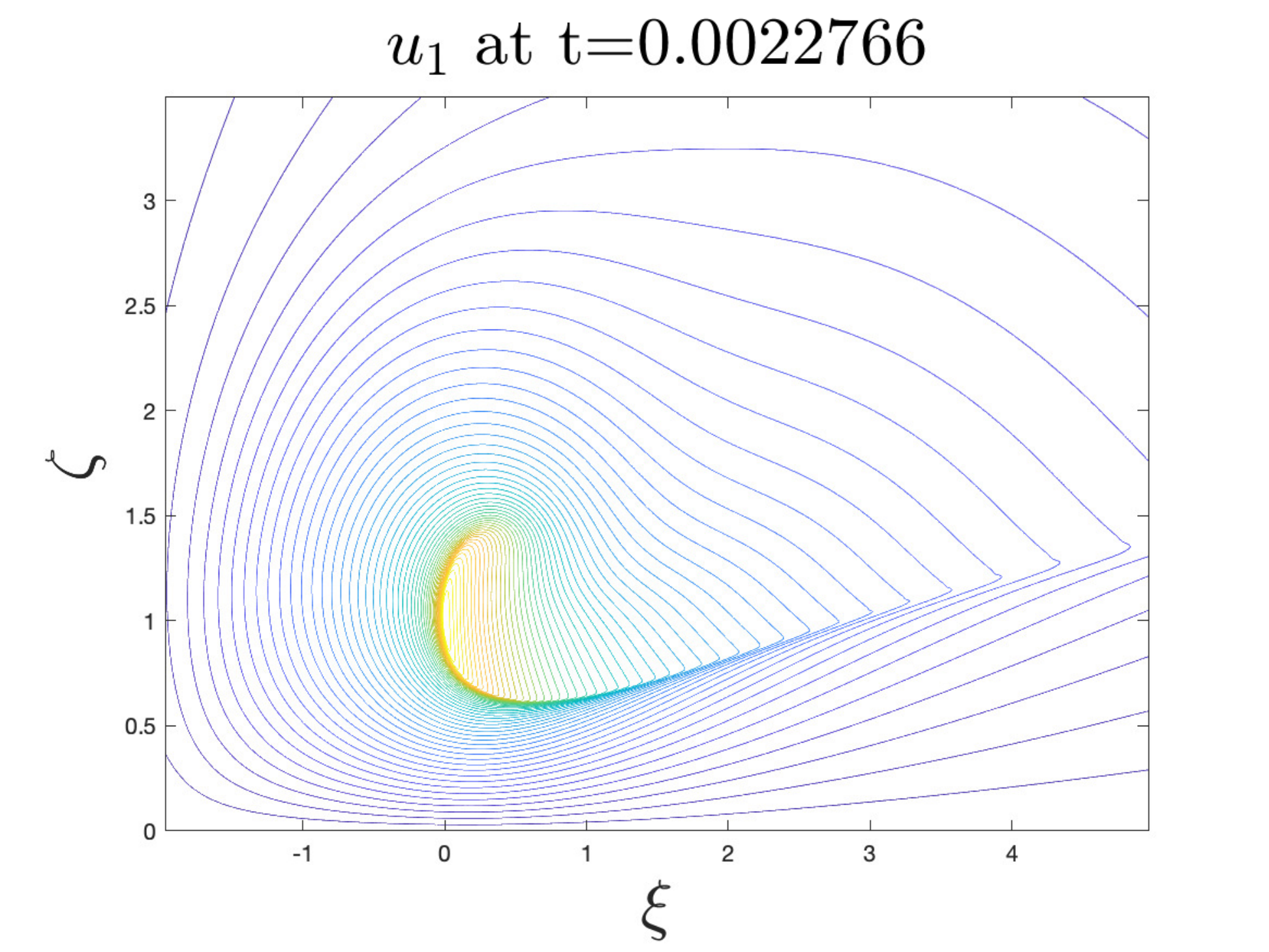}
    \end{subfigure} 
    \caption{Comparison of the level sets of $u_1$ at $t=0.0022749$,  $t=0.0022759$, and $t= 0.0022766$, respectively. First row: original level sets of $u_1$ in the domain $(r,z)$ in different times. Second row: rescaled level sets of $u_1$ as a function of $(\xi,\zeta)$ in the domain $(\xi,\zeta)$.}   
    \label{fig:levelset_compare_u1}
\end{figure}


In Figure \ref{fig:levelset_compare_u1}, we compare the level sets of $u_1$ at different time instants. In the first row of Figure \ref{fig:levelset_compare_u1}, we plot the level sets of $u_1$ in a local domain $(r,z)\in[0,0.015]\times[0,0.006]$. We plot the profiles in a short time interval from $t = 0.0022749$ to $t = 0.0022766$. The main part of the profile shrinks in space and travels toward the origin. However, if we plot the level sets of the spatially rescaled function 
\begin{equation}\label{eq:stretch_u1}
\tilde{u}_1(\xi,\zeta,t) = u_1(R(t)+Z(t)\xi ,Z(t)\zeta,t)
\end{equation}
as in the second row of Figure \ref{fig:levelset_compare_u1}, we can see that the rescaled profile of $\tilde{u}_1$ (in the $\xi\zeta$-plane) is almost indistinguishable in a small neighborhood near $(R(t),Z(t))$. Here 
\[\xi = \frac{r - R(t)}{Z(t)} ,\quad \zeta = \frac{z}{Z(t)} \]
are the dynamically rescaled variables motivated by the self-similar ansatz \eqref{eq:self-similar_ansatz}. This observation suggests that there exists an approximate self-similar profile $U(\xi,\zeta)$. 

\begin{figure}[!ht]
\centering
  	\begin{subfigure}[b]{0.32\textwidth} 
    \includegraphics[width=1\textwidth]{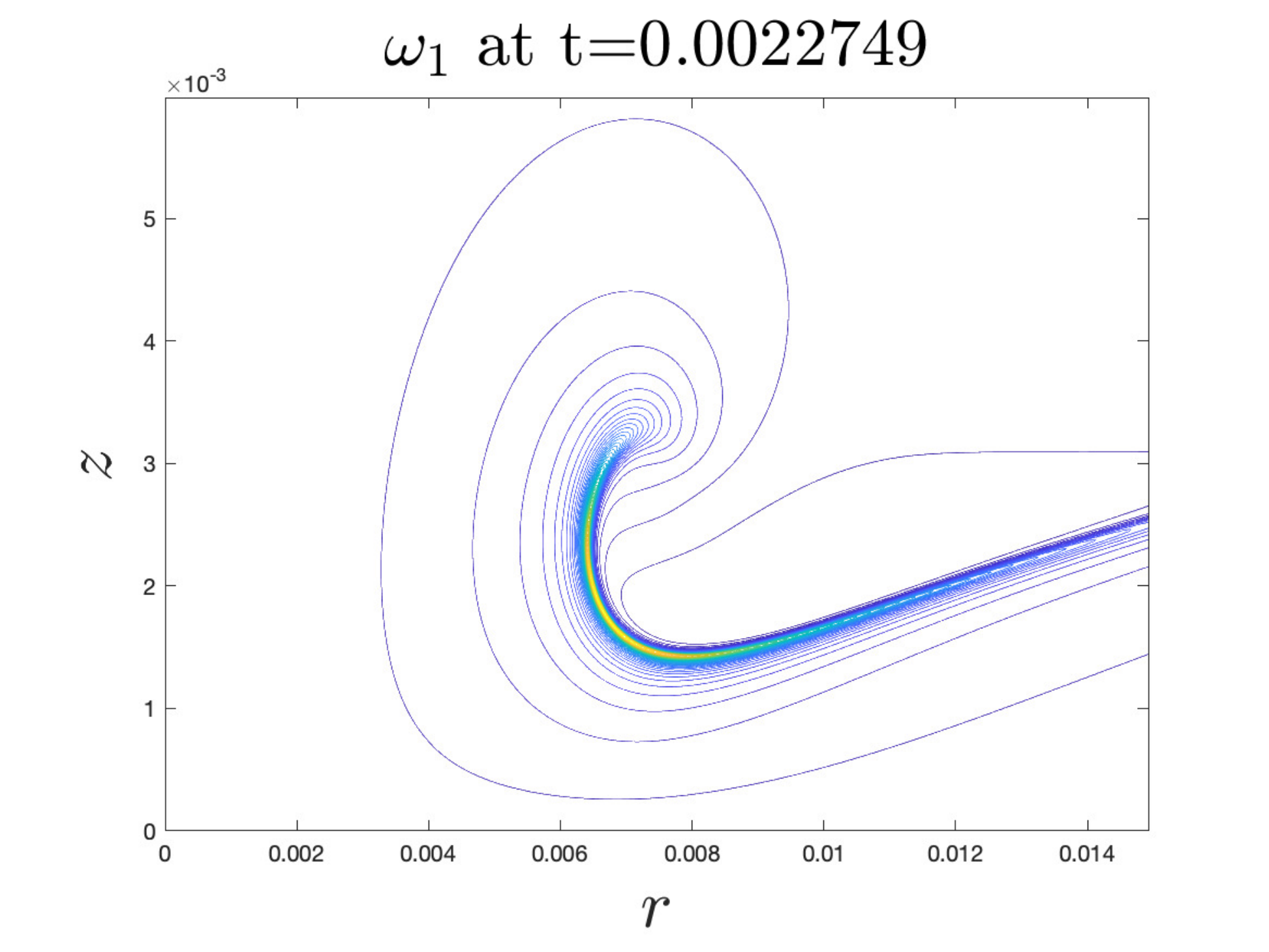}
    \end{subfigure} 
    \begin{subfigure}[b]{0.32\textwidth} 
    \includegraphics[width=1\textwidth]{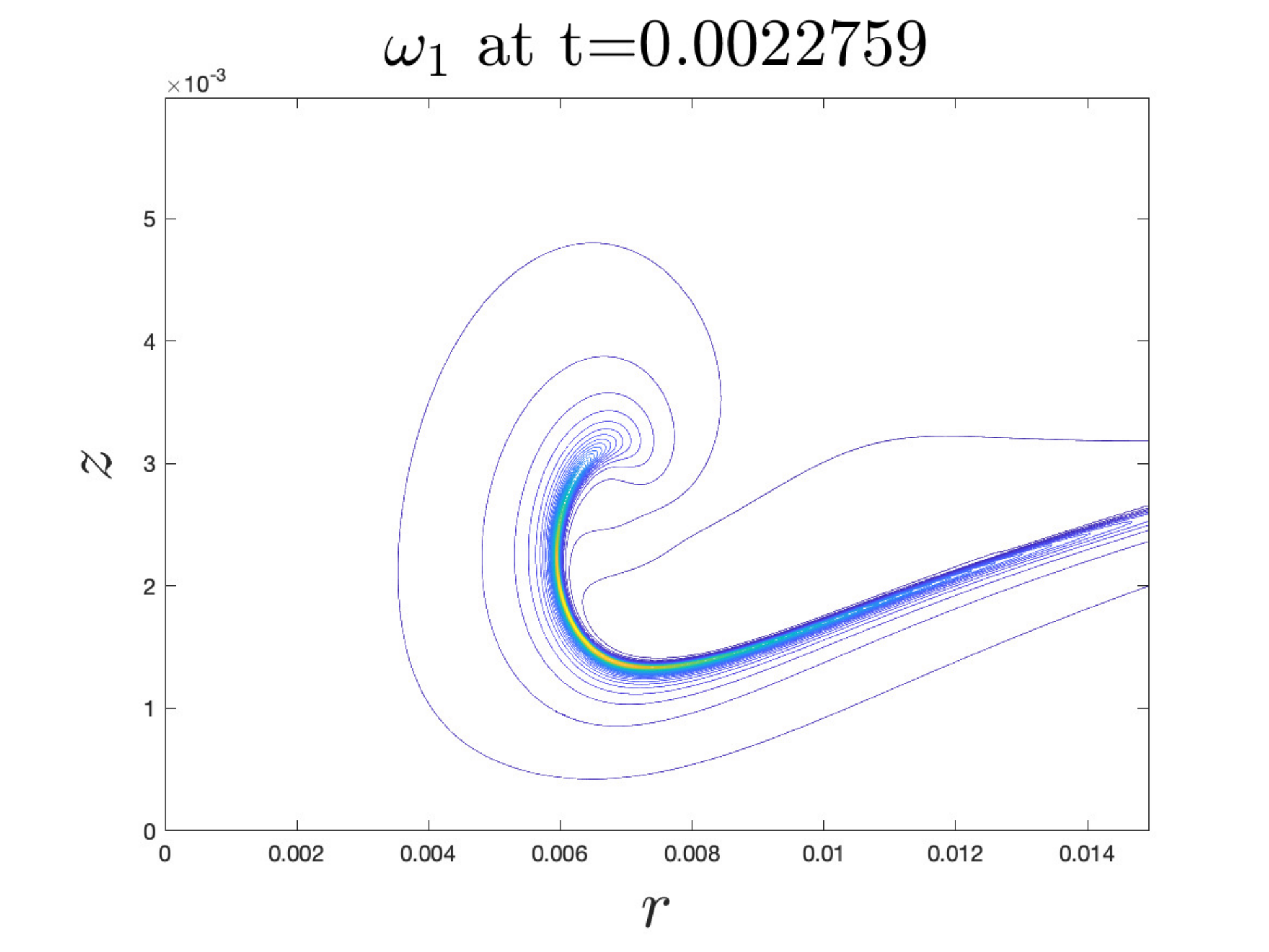}
    \end{subfigure} 
    \begin{subfigure}[b]{0.32\textwidth} 
    \includegraphics[width=1\textwidth]{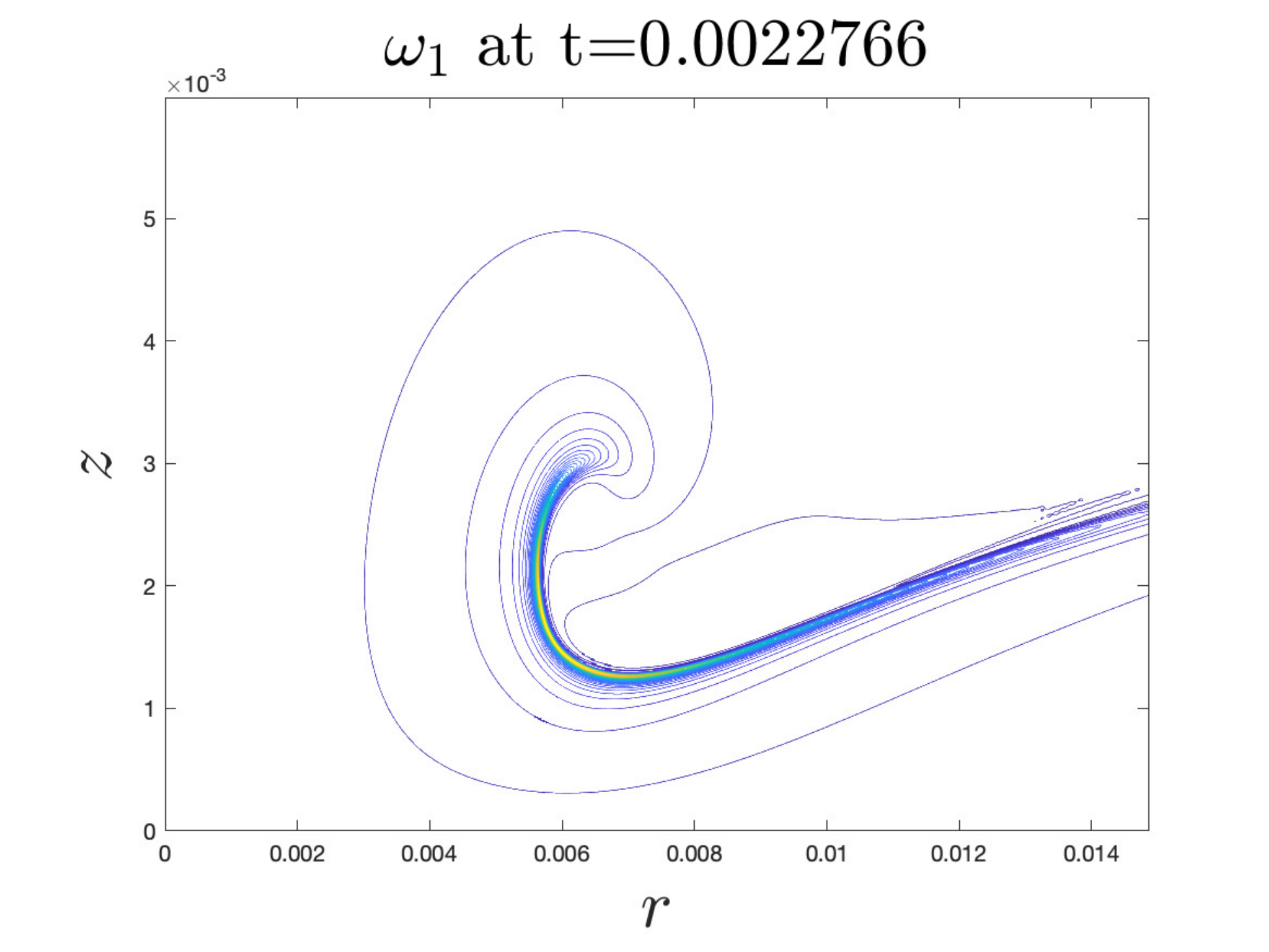}
    \end{subfigure} 
    \begin{subfigure}[b]{0.32\textwidth} 
    \includegraphics[width=1\textwidth]{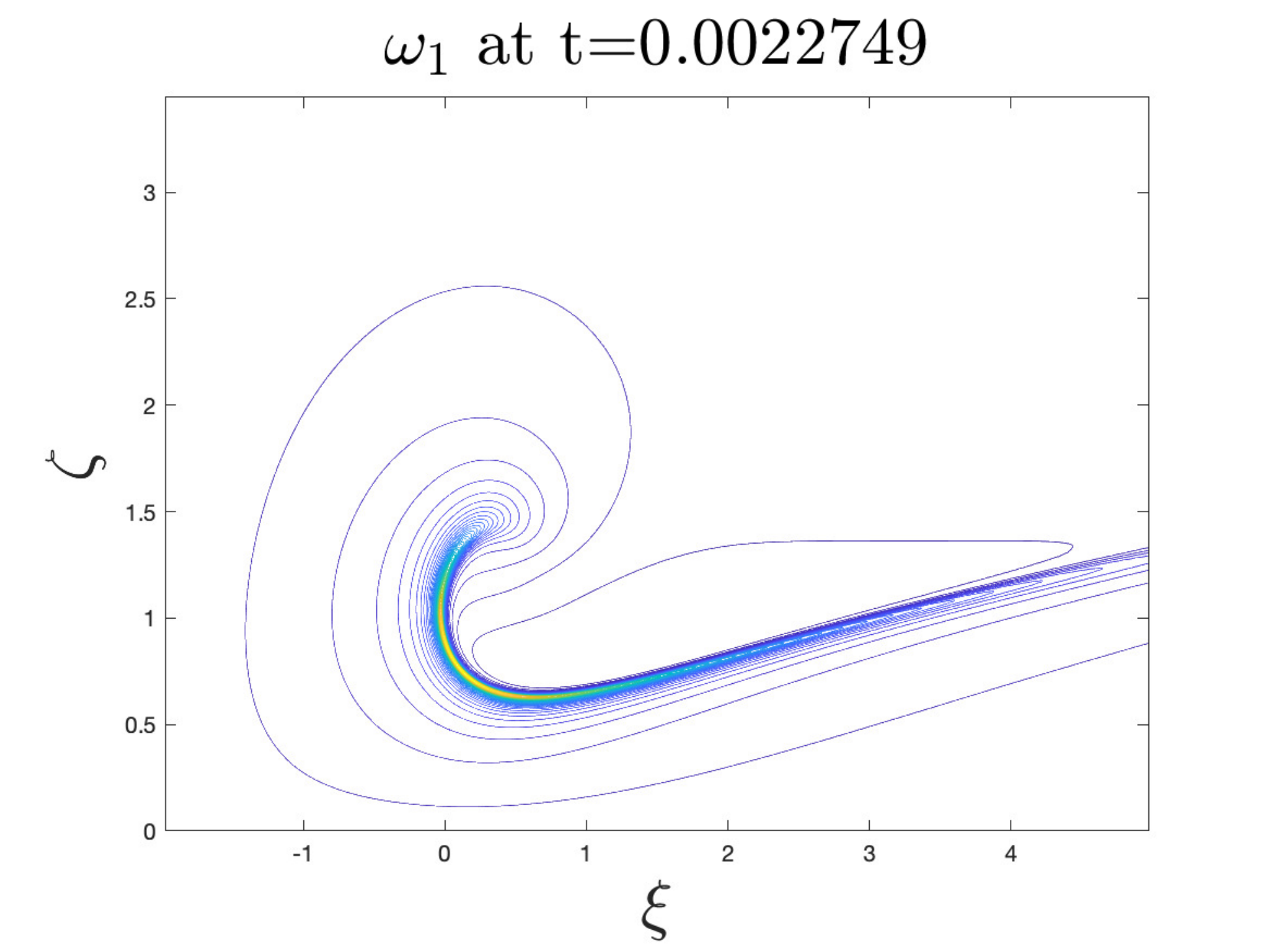}
    \end{subfigure} 
    \begin{subfigure}[b]{0.32\textwidth} 
    \includegraphics[width=1\textwidth]{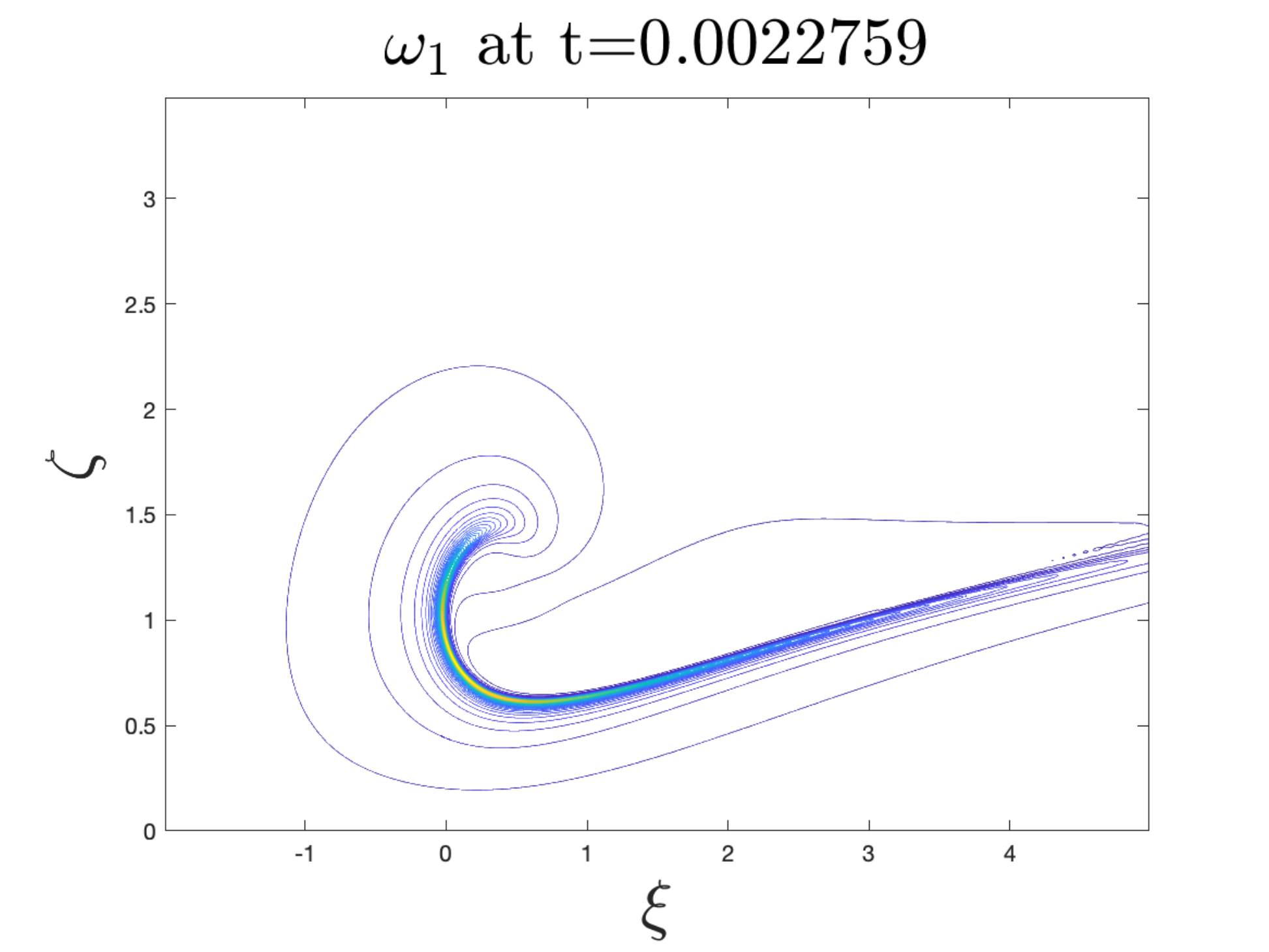}
    \end{subfigure} 
    \begin{subfigure}[b]{0.32\textwidth} 
    \includegraphics[width=1\textwidth]{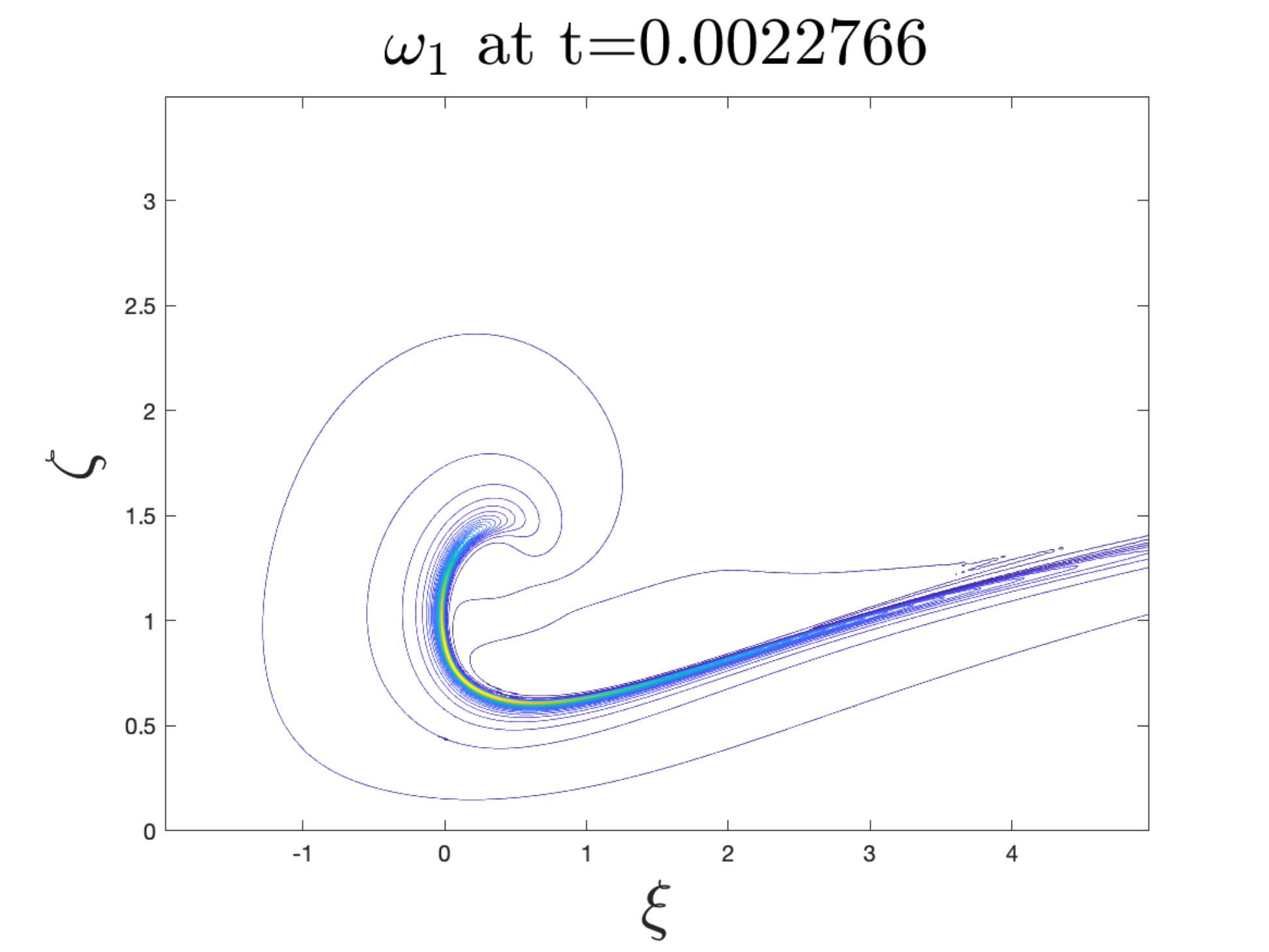}
    \end{subfigure} 
    \caption{Comparison of the level sets of $u_1$ at $t=0.0022749$,  $t=0.0022759$, and $t= 0.0022766$, respectively. First row: original level sets of $\omega_1$ in the domain $(r,z)$ in different times. Second row: rescaled level sets of $\omega_1$ as a function of $(\xi,\zeta)$ in the domain $(\xi,\zeta)$.}   
    \label{fig:levelset_compare_w1}
\end{figure}

%

In Figure \ref{fig:levelset_compare_w1}, we compare the level sets of $\om_1$ and the level sets of the spatially rescaled function
\begin{equation}\label{eq:stretch_w1}
\tilde{\om}_1(\xi,\zeta,t) = \om_1(R(t)+Z(t)\xi,Z(t)\zeta,t)
\end{equation}
in a similar manner. Again, we can see that although the profile of $\om_1$ has changed noticeably in the original physical space, the rescaled profile $\tilde{\om}_1$ seems to converge. This again suggests that there exists an approximate self-similar profile $\Omega(\xi,\zeta)$.

\begin{figure}[!ht]
\centering
	\begin{subfigure}[b]{0.38\textwidth}
    \includegraphics[width=1\textwidth]{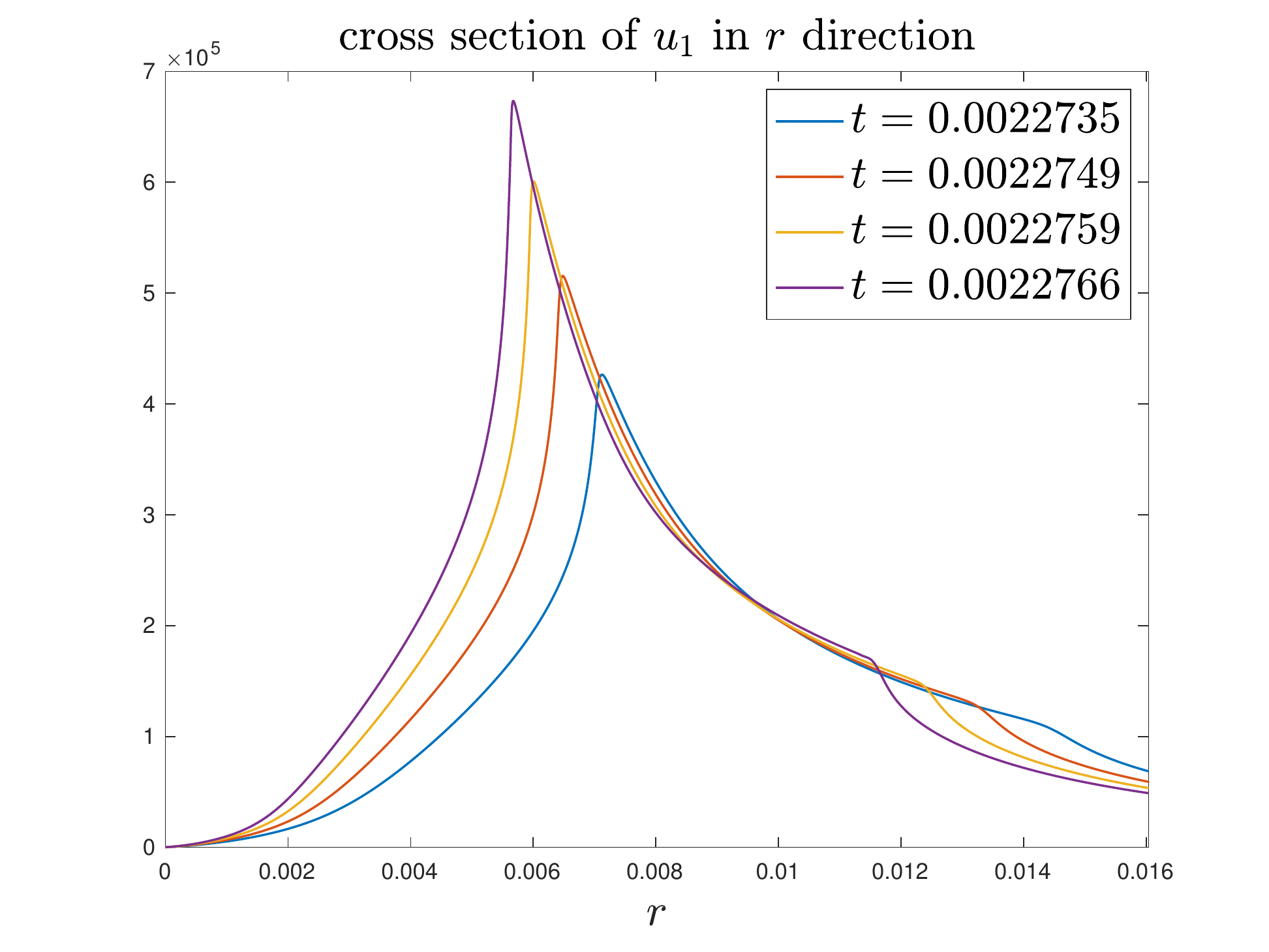}
    \caption{cross sections of $u_1$ in $r$}
    \end{subfigure}
  	\begin{subfigure}[b]{0.38\textwidth} 
    \includegraphics[width=1\textwidth]{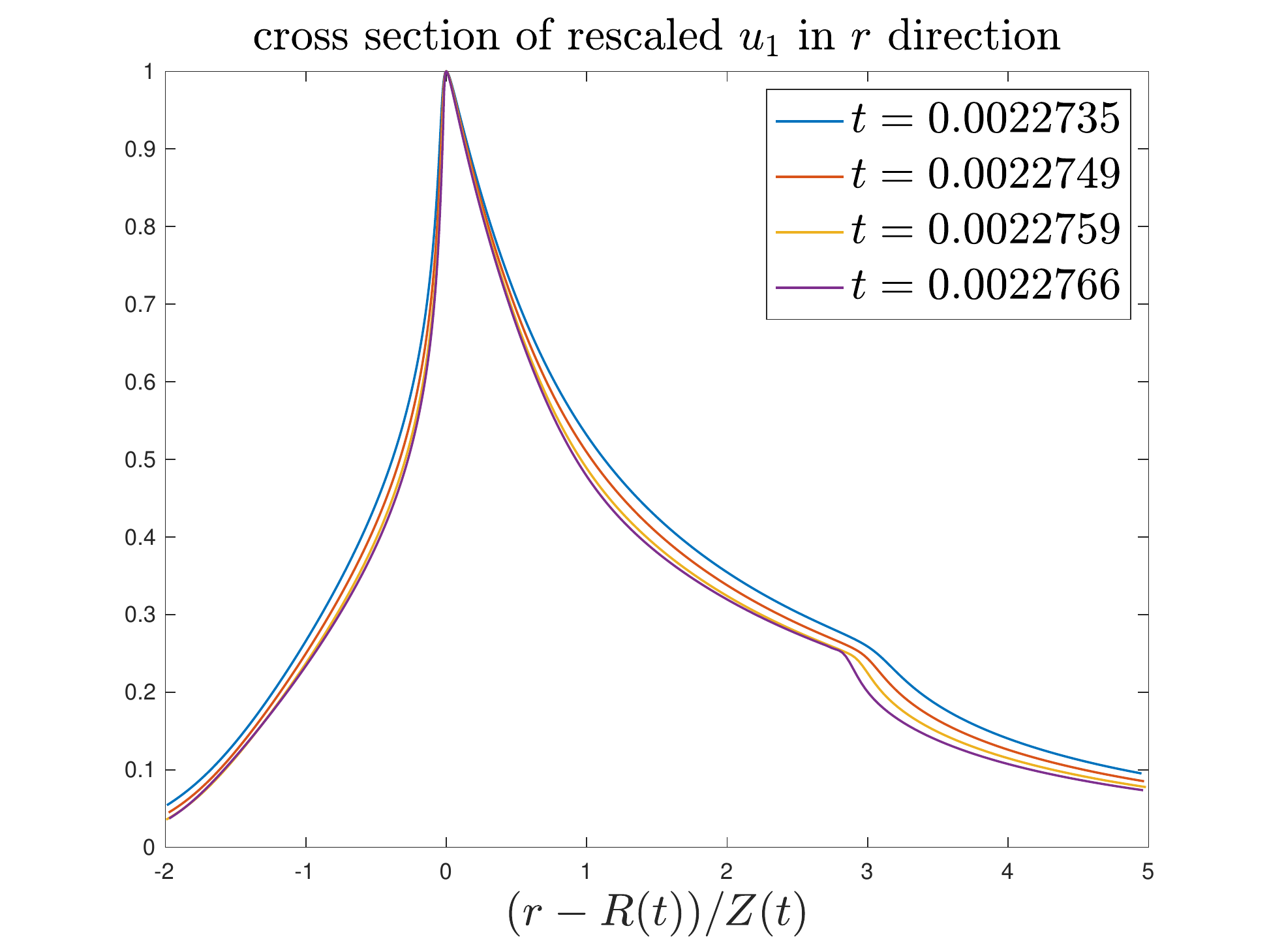}
    \caption{rescaled cross sections of $u_1$ in $r$}
    \end{subfigure} 
    \begin{subfigure}[b]{0.38\textwidth} 
    \includegraphics[width=1\textwidth]{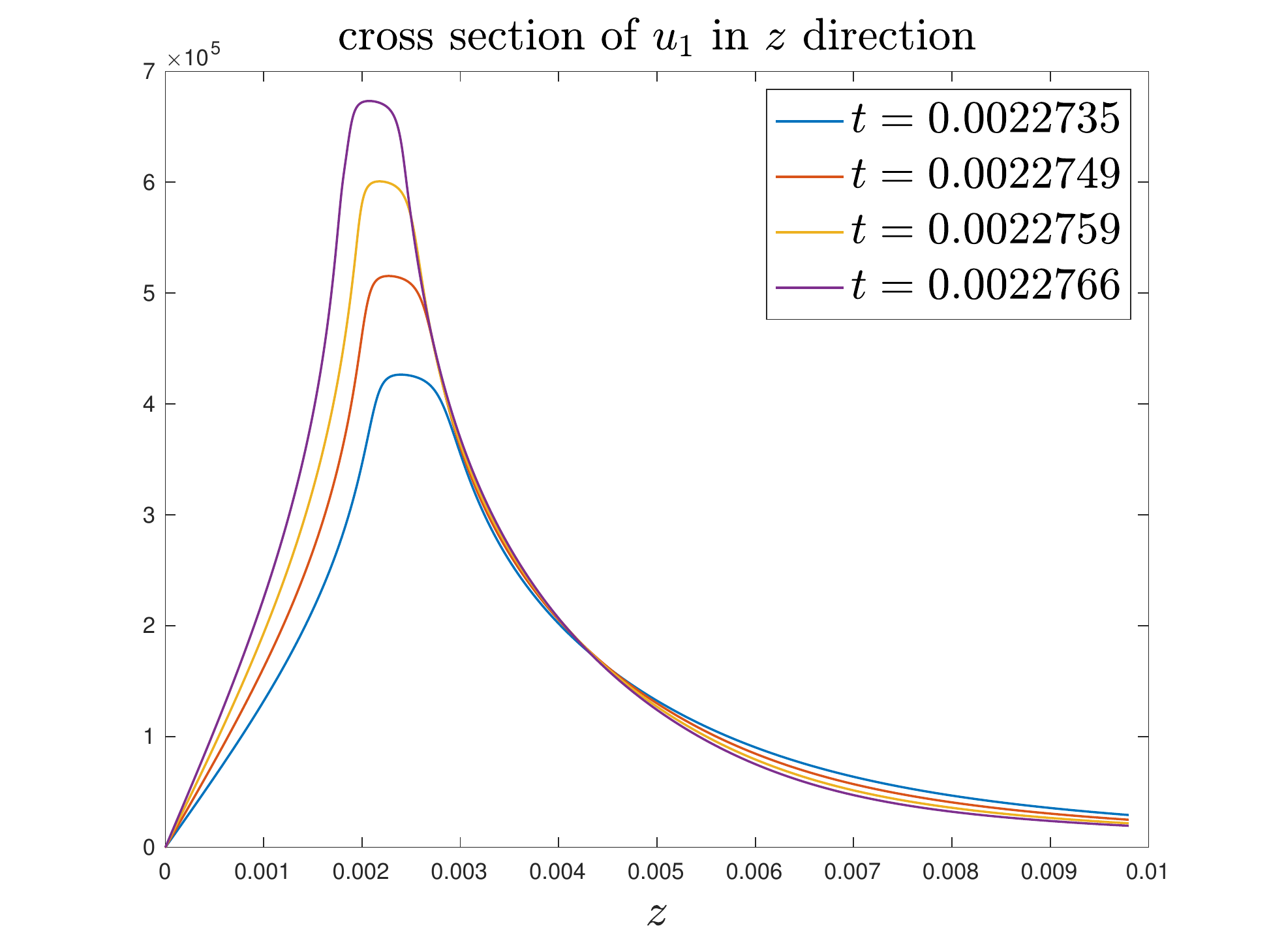}
    \caption{cross sections of $u_1$ in $z$}
    \end{subfigure}
  	\begin{subfigure}[b]{0.38\textwidth} 
    \includegraphics[width=1\textwidth]{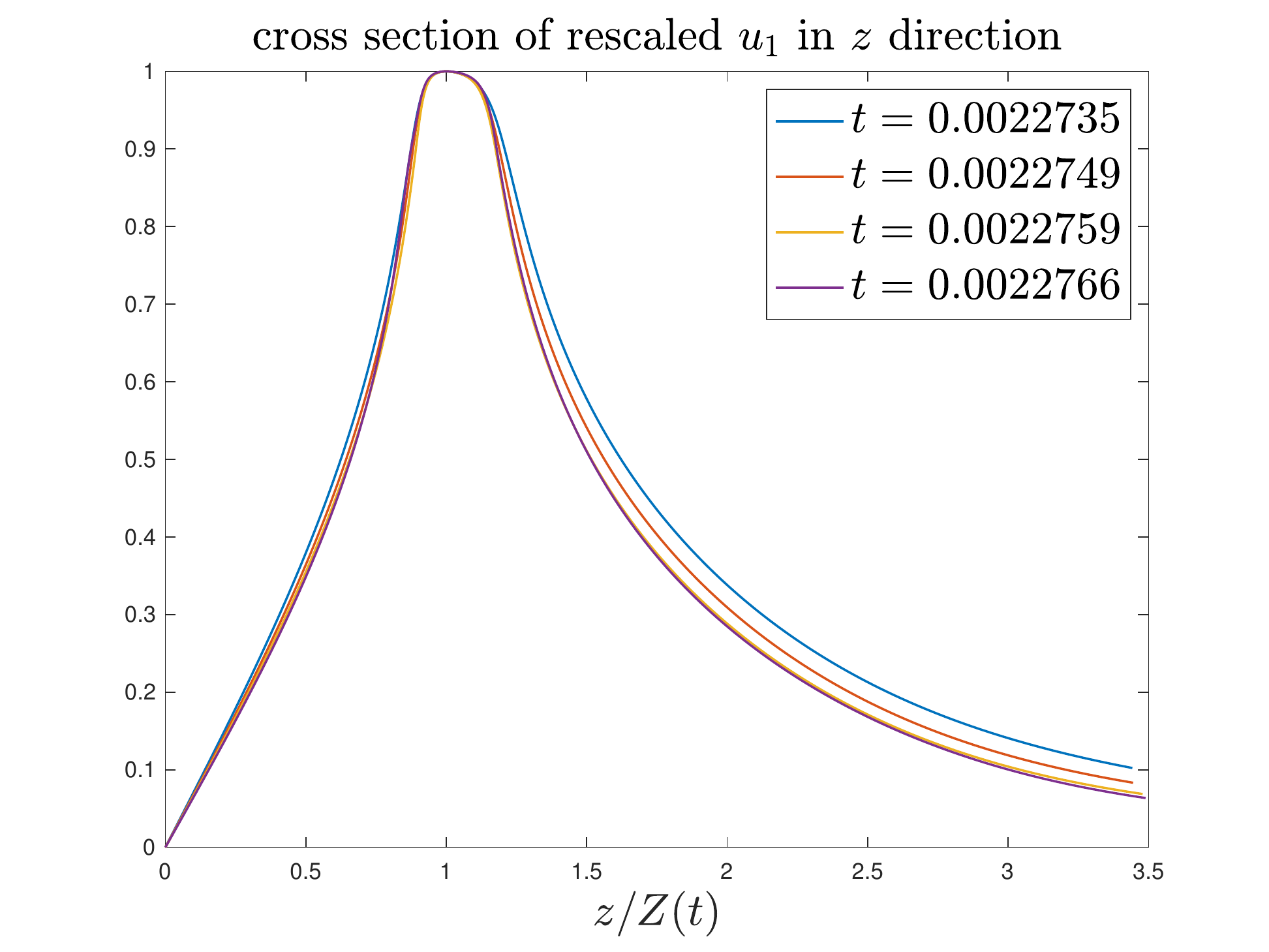}
    \caption{rescaled cross sections of $u_1$ in $z$}
    \end{subfigure} 
    \caption[Cross section compare]{Cross sections and rescaled cross sections of $u_1$ through the point $R(t),Z(t)$ in both directions at different time instants instants. (a) Cross sections in the $r$ direction. (b) Rescaled cross sections in the $r$ directions. (c) Cross sections in the $z$ direction. (d) Rescaled cross sections in the $z$ directions.}   
    \label{fig:cross_section_compare}
       \vspace{-0.05in}
\end{figure}

Finally, we compare the cross sections of the solution at different time instants to study the nearly self-similar blow-up. In Figure \ref{fig:cross_section_compare}(a) and (c), we present the evolution of the cross sections of $u_1$ through the point $R(t),Z(t)$ in both directions. The length scale of the profile shrinks in both directions, and the sharp front travels toward $r=0$. For comparison, in Figure \ref{fig:cross_section_compare}(b) and (d), we plot the corresponding cross sections of the normalized function $u_1/\|u_1\|_{L^\infty}$ in the dynamically rescaled variables $(\xi, \zeta)$. We can see that the normalized profiles for $u_1$ seem to converge to a limiting profile as time increases. These results further support the existence of an approximate self-similar profile of the solution near the reference location $(R(t),Z(t))$.

\begin{figure}[!ht]
\centering
	\begin{subfigure}[b]{0.38\textwidth}
    \includegraphics[width=1\textwidth]{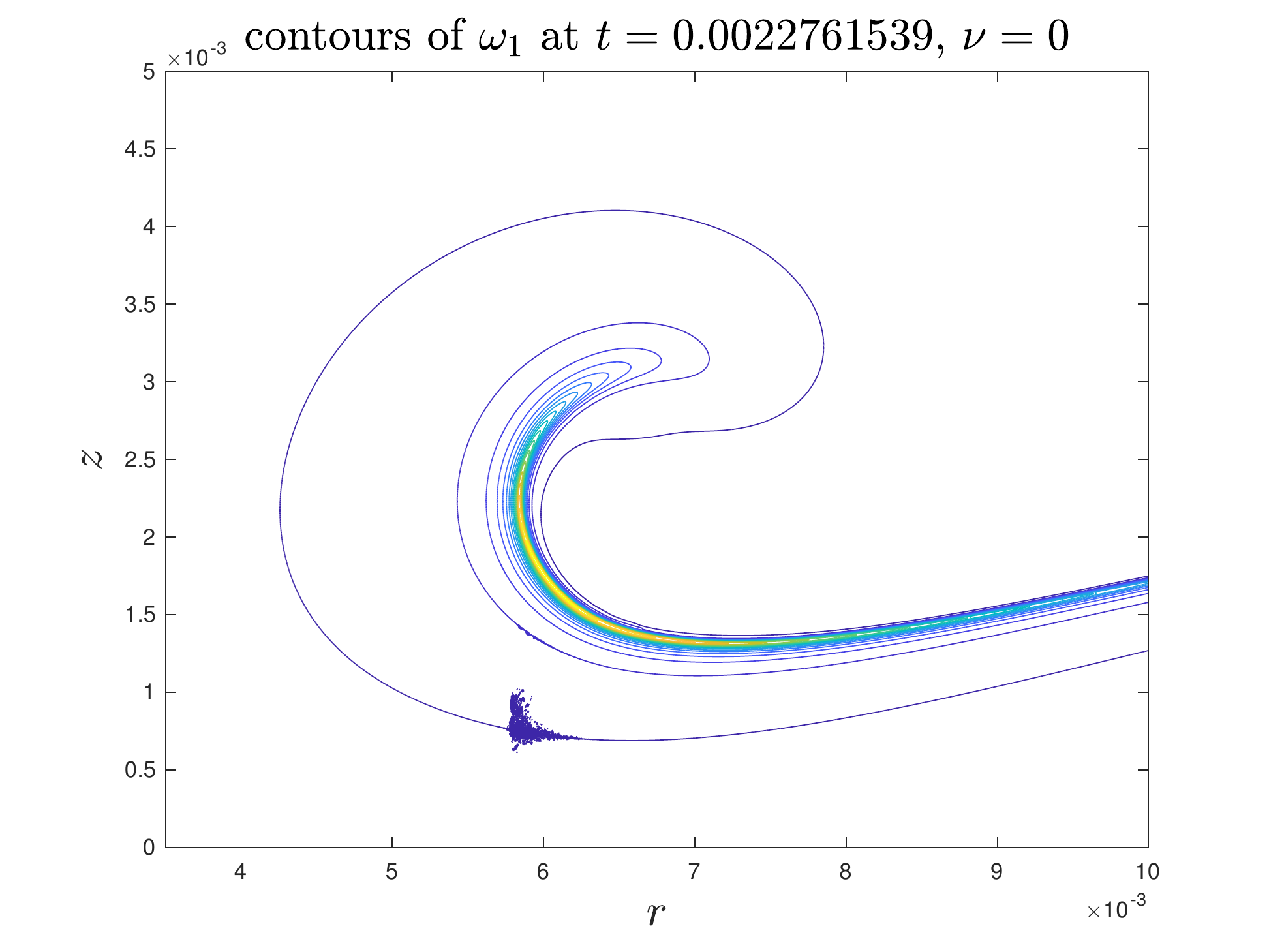}
    \caption{contours of $\omega_1$ with $\nu=0$}
    \end{subfigure}
  	\begin{subfigure}[b]{0.38\textwidth} 
    \includegraphics[width=1\textwidth]{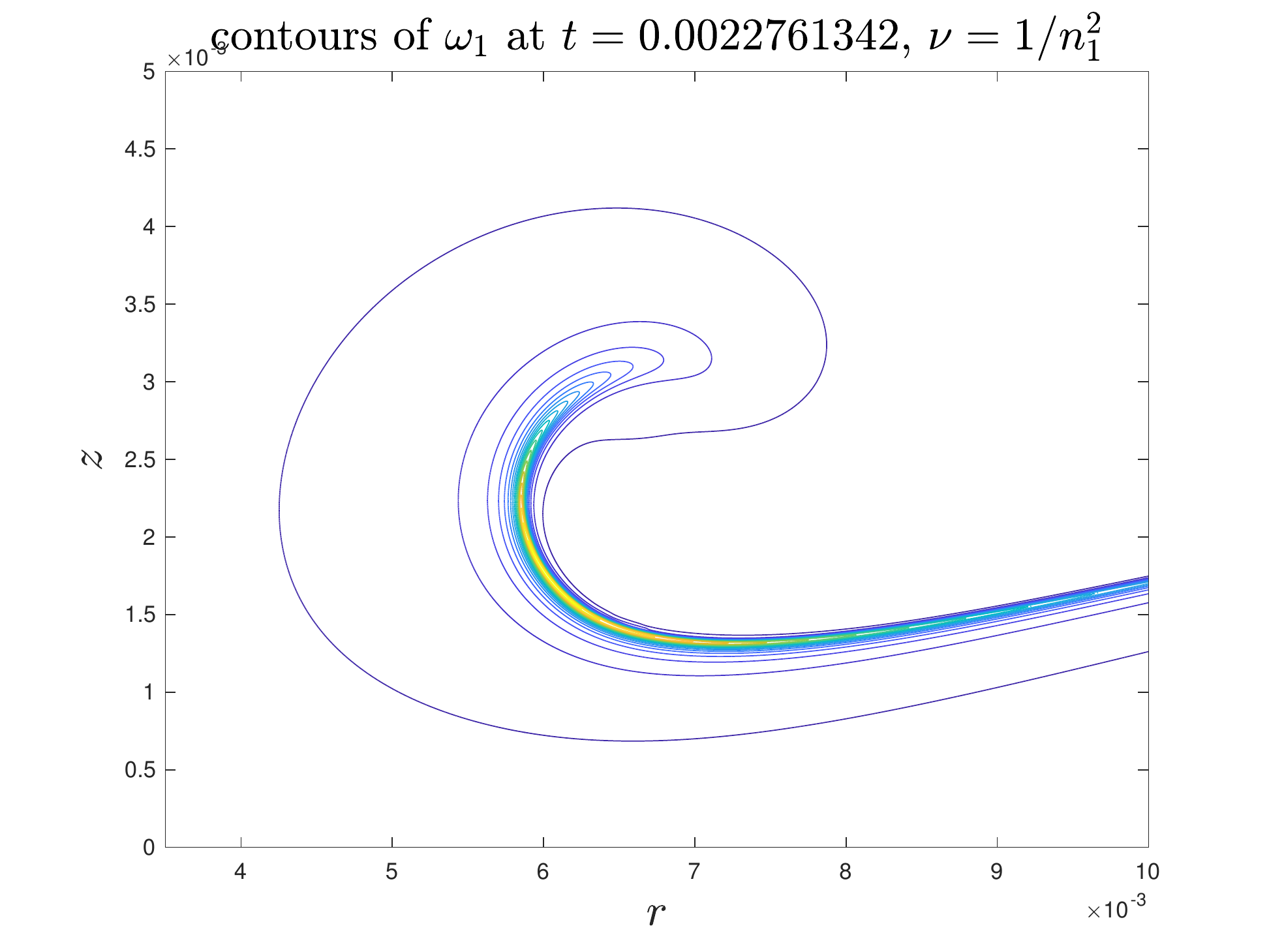}
    \caption{contours of $\omega_1$ with $\nu=1/n_1^2$}
    \end{subfigure} 
    \begin{subfigure}[b]{0.38\textwidth} 
    \includegraphics[width=1\textwidth]{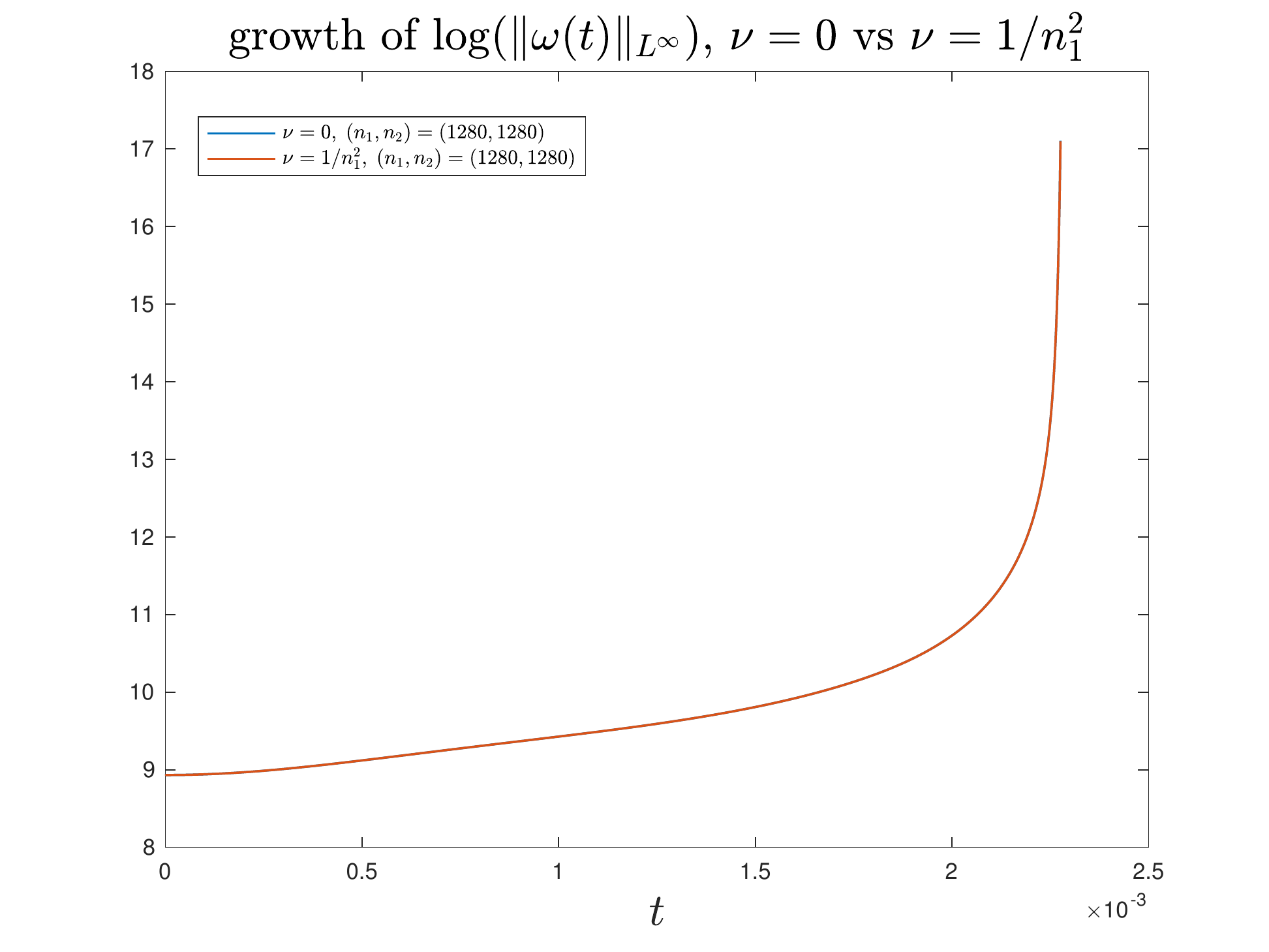}
    \caption{Comparison of $\log(\|\omega\|_\infty)$}
    \end{subfigure}
  	\begin{subfigure}[b]{0.38\textwidth} 
    \includegraphics[width=1\textwidth]{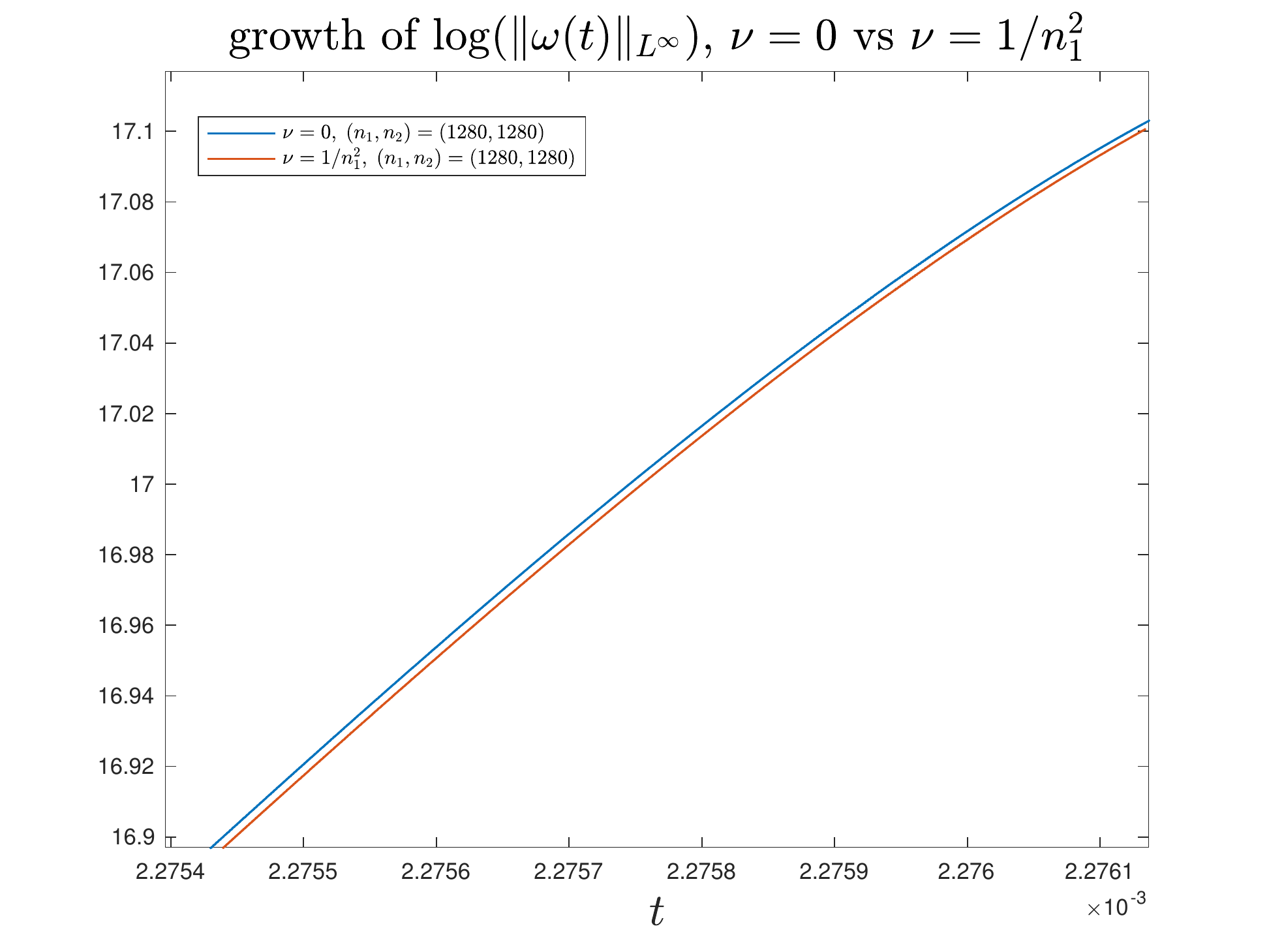}
    \caption{Comparison of $\log(\|\omega\|_\infty)$, close-up view}
    \end{subfigure} 
    \caption[Contour compare w1]{(a) Contours of $\omega_1$ without using numerical viscosity after $90000$ time steps. $\|\omega(t)\|_{L^\infty}/\|\omega(0)\|_{L^\infty}$ has grown roughly $3500$. We observe that oscillations develop by this time. (b) Contours of $\omega_1$ using numerical viscosity $\nu =1/n_1^2$ after $90000$ time steps. 
    (c) Comparison of the $\log(\|\omega\|_{L^\infty})$ with and without using numerical viscosity. (d) same as (c), a close-up view.}   
    \label{fig:numerical_viscosity}
       \vspace{-0.05in}
\end{figure}

\begin{figure}[!ht]
\centering
	\begin{subfigure}[b]{0.38\textwidth}
    \includegraphics[width=1\textwidth]{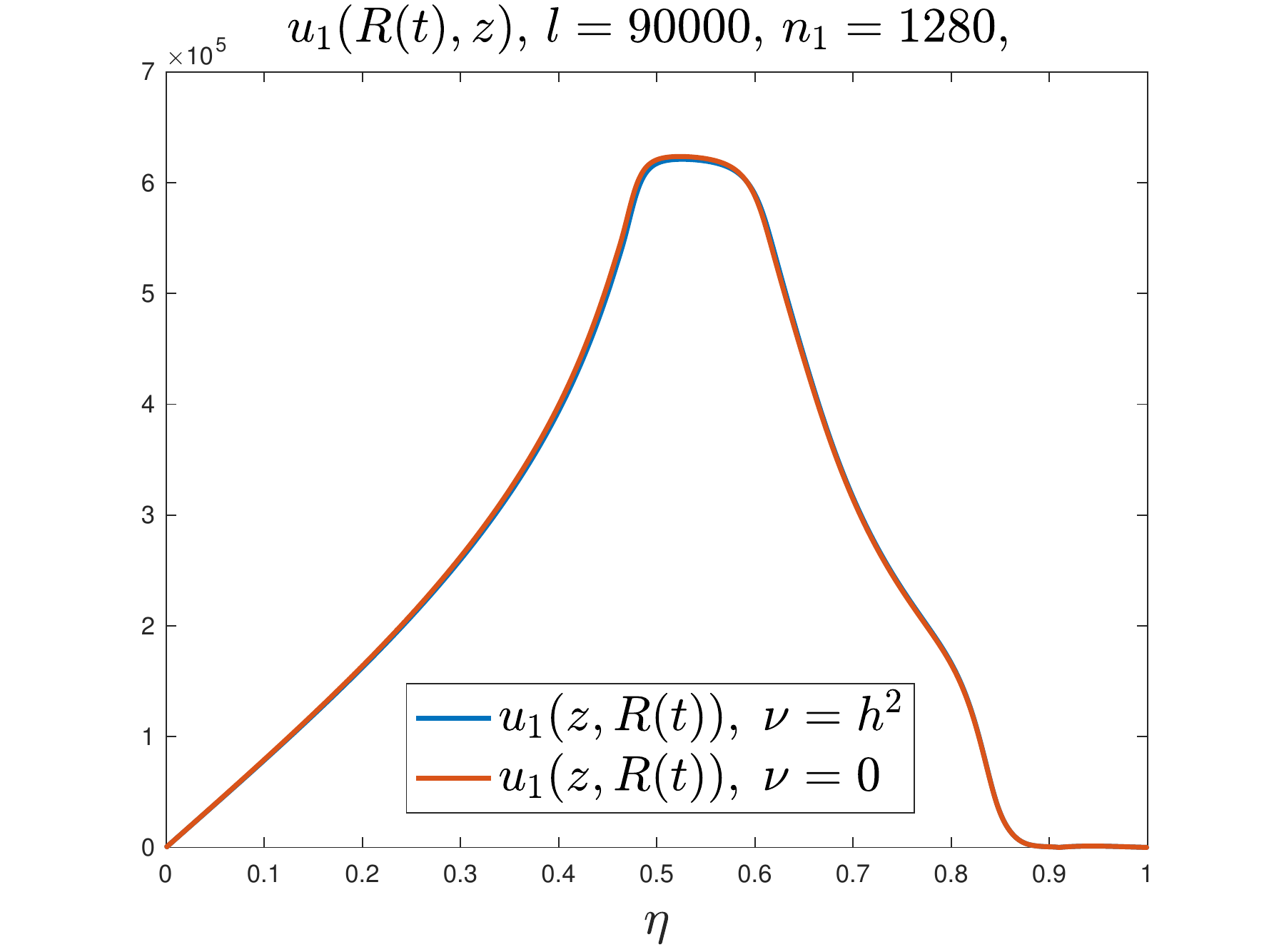}
    \caption{$z$-cross section $u_1$ at $r=R(t)$}
    \end{subfigure}
  	\begin{subfigure}[b]{0.38\textwidth} 
    \includegraphics[width=1\textwidth]{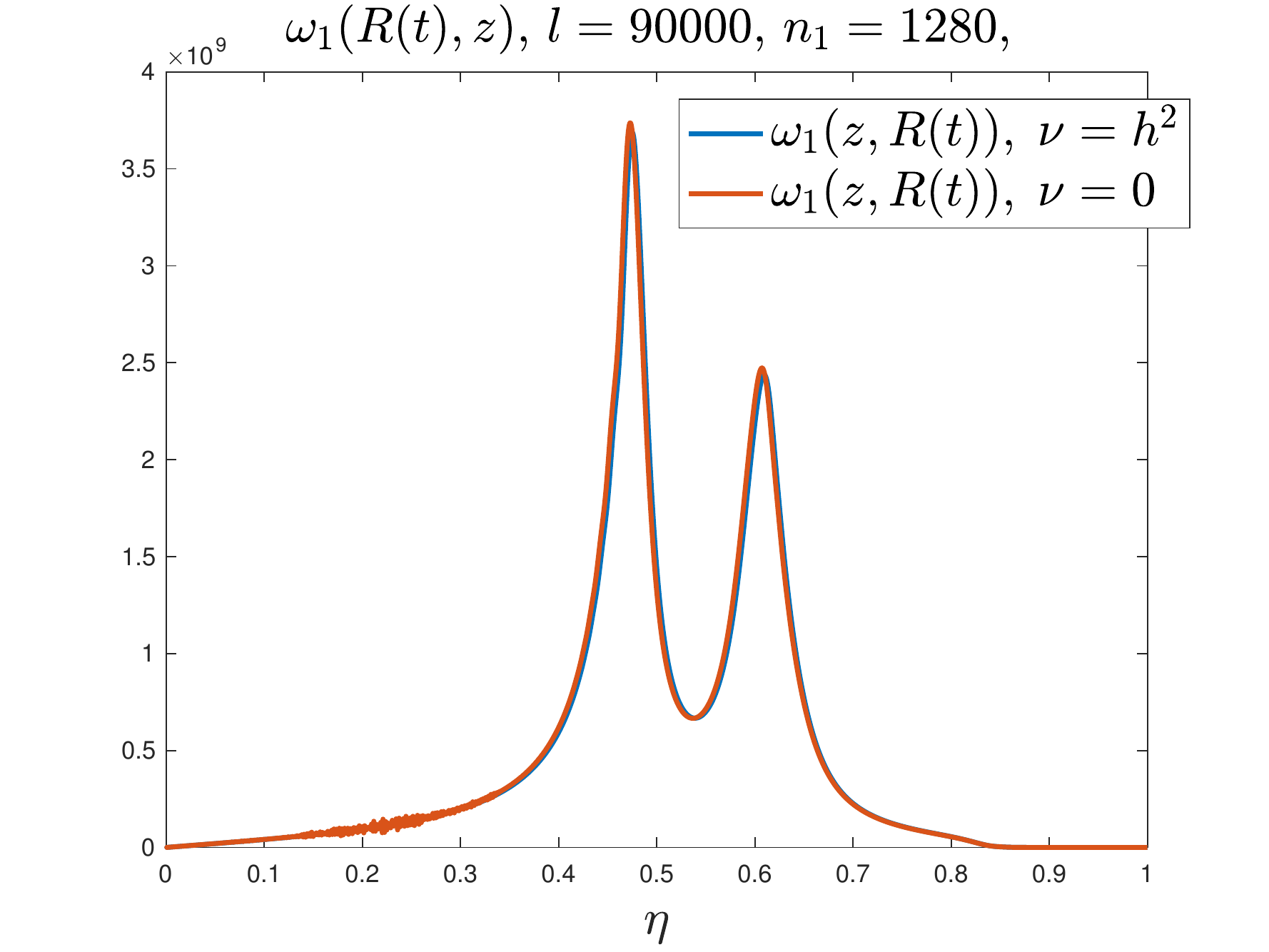}
    \caption{$z$-cross section $\omega_1$ at $r=R(t)$}
    \end{subfigure} 
    \vspace{0.05in}
    \begin{subfigure}[b]{0.38\textwidth} 
    \includegraphics[width=1\textwidth]{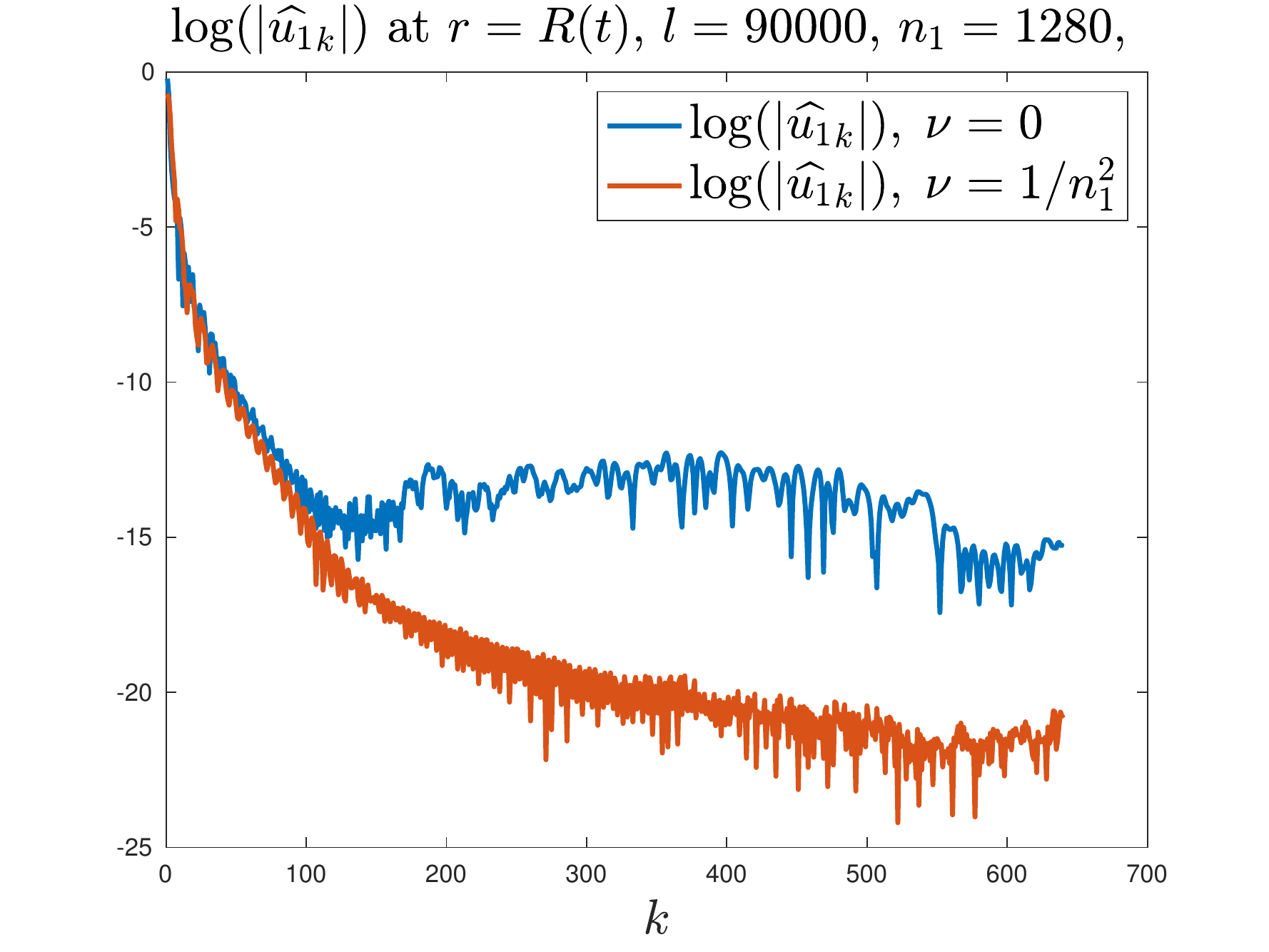}
    \caption{Fourier spectra of $u_1(R(t),z)$}
    \end{subfigure}
    \vspace{0.05in}
  	\begin{subfigure}[b]{0.38\textwidth} 
    \includegraphics[width=1\textwidth]{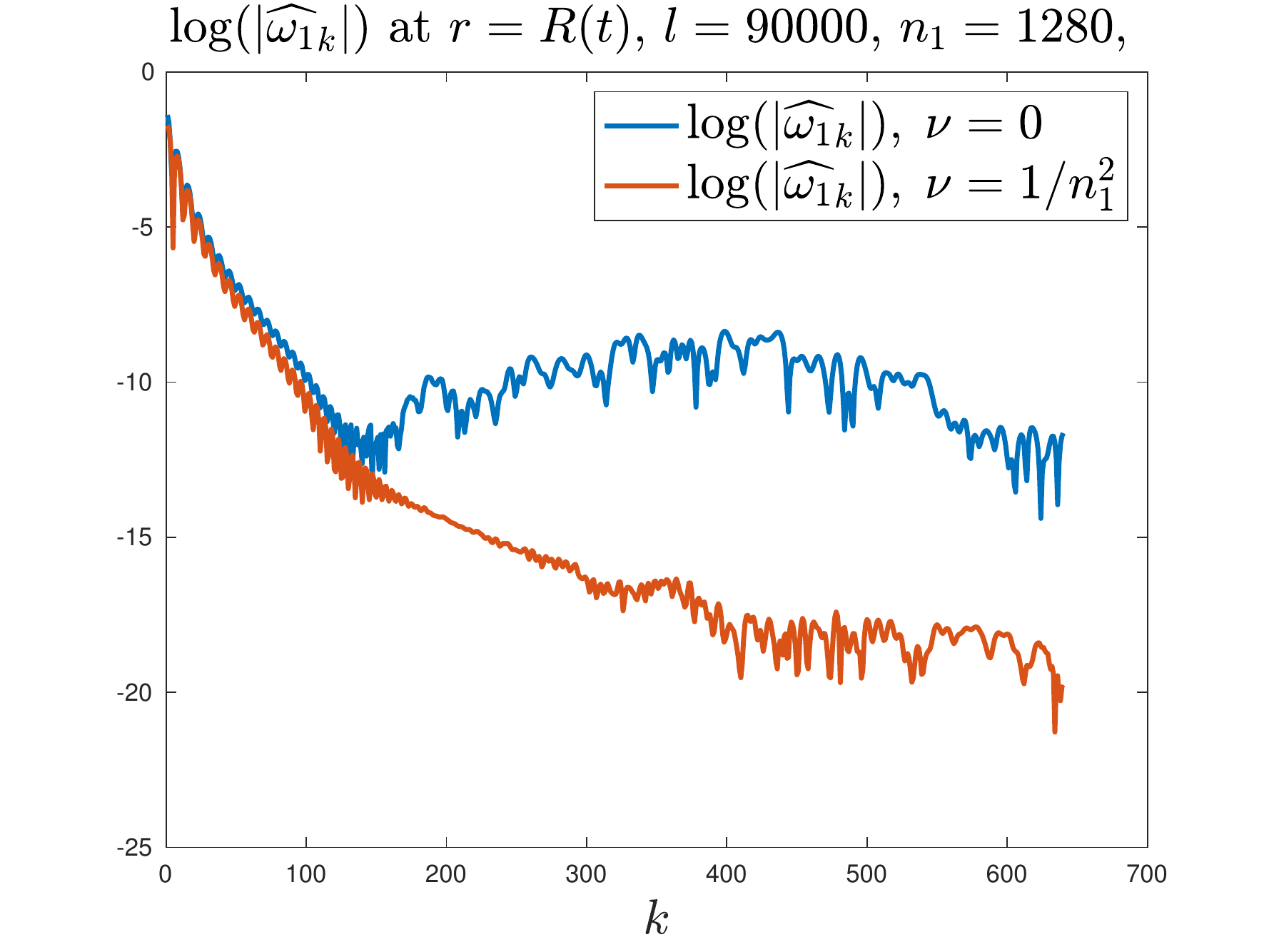}
    \caption{Fourier spectra of $\omega_1(R(t),z)$}
    \end{subfigure} 
    \caption[Cross section compare w1]{(a) Comparison of $z$-cross section of $u_1$ at $r=R(t)$ as a function of $\eta$ with and without using numerical viscosity after $90000$ time steps. (b) Comparison of $z$-cross section of $\omega_1$ at $r=R(t)$ with and without using numerical viscosity. (c) Comparison of the Fourier spectrum of the $z$-cross section $u_1/\max(u_1(R(t),z))$ at $r=R(t)$ as a function of $\eta$ with and without using numerical viscosity. (d) Comparison of the Fourier spectrum of the $z$-cross section $\omega_1/\max(\omega_1(R(t),z))$ at $r=R(t)$.  }
    \label{fig:numerical_viscosity2}
       \vspace{-0.05in}
\end{figure}

\subsection{The need for applying a second order numerical viscosity}

In this subsection, we provide some justification for the need to apply a second order numerical viscosity to compute the Euler equations. We will compare the numerical solution of the $3$D Euler equations using a second order numerical viscosity with that without using any numerical viscosity. 

In Figure \ref{fig:numerical_viscosity}(a), we plot the contours of $\omega_1$ using grid $1280\times1280$ to solve the Euler equations without using any numerical viscosity. The solution is obtained after $90,000$ time steps at $t=0.0022761539$. We observe some high frequency oscillations developed in the contour where $\omega_1$ is smooth and relatively small in amplitude. 
Judging from Figure \ref{fig:numerical_viscosity}(a), these high frequency oscillations occur near the interface between the two adaptive mesh phases. It is most likely that this high frequency instability is triggered by the frequent changes of adaptive mesh in the late stage. 

By the design of our adaptive mesh strategy, we change the adaptive mesh when the most singular part of solution is about to depart from phase $1$ of the adaptive mesh, which has the finest grid, and enter phase $0$ of the adaptive mesh, which has a relatively coarse grid. Thus the interpolation from the old mesh to the new mesh introduces some non-smooth high frequency errors that are largest in amplitude near the interface between the phase $1$ mesh and the phase $0$ mesh. 

Another contributing factor is that
the Poisson solver for the stream function becomes very ill-conditioned in the late stage due to the extremely large ratio between the coarsest mesh and the finest mesh. As a result, the high frequency errors introduced by the frequent changes of mesh in the late stage get amplified through the very ill-conditioned Poisson solver. 
We remark that these high frequency oscillations are different from the oscillations generated by fluid dynamic instability, which typically occurs near the sharp front where the shearing instability is strongest. A second order vanishing numerical viscosity would not be strong enough to suppress such instability.

In Figure \ref{fig:numerical_viscosity}(b), we plot the contours of $\omega_1$ using a second order numerical viscosity with $\nu = 1/n_1^2$ after $90,000$ time steps at $t=0.0022761342$. 
We can see that the contours are very smooth and the high frequency oscillations that we observed without using any numerical viscosity have been regularized by the second order numerical viscosity. 

In Figure \ref{fig:numerical_viscosity}(c), we compare the growth of maximum vorticity with and without using numerical viscosity. To better visualize the difference between the two solutions, we plot $\log(\|\vom\|_\infty)$. These two solutions are almost indistinguishable. We show a close-up view of the same picture in the late stage of the computation in Figure \ref{fig:numerical_viscosity}(d). We can see that the maximum vorticity computed with a second order numerical viscosity grows slightly slower than that without using any numerical viscosity. The difference is really small.

In Figure \ref{fig:numerical_viscosity2}(a)-(b), we plot the $z$-cross section of $u_1$ and $\omega_1$ at $r=R(t)$ as a function of $\eta$ with and without using numerical viscosity. Recall that $z=z(\eta)$ is the mesh map in the $z$ direction. We can see that two solutions are almost indistinguishable except for some mild oscillation in $\omega_1$ without using the second order numerical viscosity. To further illustrate the source of numerical instability and the effect of using a second order numerical viscosity, we plot the Fourier spectrum of $u_1/\max_z(u_1(t,R(t),z)$ and $\omega_1/\max_z(\omega_1(t,R(t),z)$ as a function of $\eta$ in Figure \ref{fig:numerical_viscosity2}(c)-(d). In order to reduce the boundary effect at $\eta=0,1$, we apply a soft cut-off $f_c(\eta)$ that is approximately equal to $1$ for $\eta \in (0.1,0.9)$ and goes to zero smoothly at the boundary $\eta=0,1$. We can see that without using numerical viscosity, the solution develops some high frequency instability. On the other hand, the solution obtained by applying a second order numerical viscosity does not suffer from this high frequency instability.


\vspace{-0.1in}
\subsection{Remark on ``hydrodynamic instability of blow-up solutions to $3$D Euler''}
\label{blowup-stability}

In \cite{vasseur2020,vasseur2021}, the authors showed that blow-up solutions to $3$D Euler are hydrodynamically unstable. 
Assume that ${\bf u}(t)$ is a smooth solution of the $3$D Euler equations with smooth initial value ${\bf u}^0 \in H^s(\Omega)$ ( $s > 9/2$) for $t < T^*$ but becomes singular at $t=T^*$ (possibly infinite). The authors considered the following linearized equation for the perturbation ${\bf v}$ around the blow-up solution ${\bf u}$:
\begin{equation}
\label{linear-euler}
{\bf v}_t + {\bf u} \cdot \nabla {\bf v} + ({\bf v} \cdot \nabla){\bf u} = - \nabla P' , \quad \nabla \cdot {\bf v} = 0, \quad {\bf v}\cdot {\bf n}|_{\partial \Omega} = 0.
\end{equation}
Denote $\gamma_p (T)$ as the growth in the $L^p$ norm of the perturbation ${\bf v}$ as follows:
\[
\gamma_p (T) = \sup_{{\bf v}^0 \in H^1(\Omega),\|{\bf v}^0\|_{L^p} \leq 1}\|{\bf v}(T)\|_{L^p} \; .
\]
One of the main results obtained in \cite{vasseur2020} is that 
\[
\lim \sup_{T \rightarrow T^*} \gamma_p (T) = \infty.
\]

In a recent preprint \cite{chenhou2022}, we use a Riccati equation, the inviscid Burgers equation and the $3$D Euler equations to demonstrate that using the linearized equation \eqref{linear-euler} around a blow-up solution may not be suitable to study stability of a blow-up solution. Here we summarize the main findings in \cite{chenhou2022}. We begin with a Riccati equation
\[
u_t = u^2, \quad u(0) = 1.
\]
We know that $u(t) = 1/(1-t)$ for $0 \leq t < 1$ is the exact solution that blows up at $T^*=1$. Let $\epsilon$ be a small perturbation to the initial value, i.e. $v(0) = \epsilon$. The linearized equation for $v$ is given by
\[
v_t = 2 u v, \quad v(0)=\epsilon ,\]
which can be solved analytically to give 
\[
v(t)/v(0) = 1/(T^*-t)^2 \rightarrow \infty , \quad \mbox{as}\; t \rightarrow T^*. 
\]
By the notion of stability described in \cite{vasseur2020}, the blow-up of this Riccati equation would be considered unstable.
On the other hand, it is clear that the blow-up of the Riccati equation with a perturbed initial value $u^\epsilon(0)=1+\epsilon$ is stable in the sense that the blow-up solution has the same form $u^\epsilon(t) = 1/(T^\epsilon - t)$ with $T^\epsilon = 1/(1+\epsilon )$ for any small $\epsilon$ and we have $|T^\epsilon -T^*| = \epsilon/(1+\epsilon)$. 

We can make a similar argument for the blow-up solution of the inviscid Burgers equation
\begin{equation}
\label{Burger-Eq}
u_t + (u^2/2)_x = 0, \quad
u(0,x) = u_0(x) = - \frac{x}{1+x^2}\;, x \in \Omega \equiv (-\infty, \infty),
\end{equation}
where $u_0 \in C^\infty$ satisfies the following properties: (i) $\partial_x u_0$ is minimal at $x=0$, (ii) $u_0(0) = 0$, (iii) $\partial_x u_0(0) = -1 < 0$, (iv) $\partial_x^2 u_0(0) = 0$ and $\partial_x^3 u_0(0) = 6 > 0$. Moreover, $u_0(x)$ is strictly monotonically decreasing for $x \in [-0.5, 0.5]$.
It is well known that the inviscid Burgers equation has an implicit solution formula of the form 
\begin{equation}
\label{Burger-solution}
u(t,x) = u_0(x-tu(t,x)).
\end{equation}
It is easy to see that $u(t,x)$ is odd and monotonically decreasing over $[-0.5,0.5]$ and develops a finite time singularity in $u_x$ at $x^*=0$ and $T^* = -1/u_{0,x}(0) = 1$. 
The linearized equation is given by
\[
v_t + u v_x + u_x v = 0\;, \quad v(0,x) = v_0^\epsilon (x) \in V\;.
\]
Here $V$ is a subspace of $H^1(\Omega)\cap L^p(\Omega)$ that has the properties that any $v_0^\epsilon(x) \in V$ is an odd smooth function with compact support in $[-\epsilon,\epsilon ]$ ($\epsilon \leq 0.5$) and
satisfies $\|\partial_x^j  v_0^\epsilon \|_{L^\infty} \leq c_j \epsilon^{4-j}$ for $ 0 \leq j \leq 3$. Further by choosing $\epsilon$ sufficiently small, we have $\partial_x^3(u_0+v_0^\epsilon)(0) >0$. Since $v$ is transported by velocity $u$, $v(t,x)$ has compact support in $[-X^{\epsilon}(t), X^\epsilon(t)]$, here $X^\epsilon(t) = X(t,\epsilon)$ and $X(t,\alpha)$ is the characteristics defined by $\frac{dX(t,\alpha)}{dt} = u(t,X(t,\alpha))$ with $X(0,\alpha)=\alpha$.  A simple calculation shows that 
\begin{eqnarray*}
\frac{d}{dt} \| v(t)\|_{L^p}^p &= & (p-1) \int (-u_x(t,x)) v(t,x)^p dx = (p-1) \int_{-X^{\epsilon}(t)}^{X^\epsilon(t)} (-u_x(t,x)) v(t,x)^p dx \\
& \geq & (p-1)\min_{X^{-\epsilon}(t) \leq x \leq X^\epsilon(t)} (-u_x(t,x)) \int_{-X^{\epsilon}(t)}^{X^\epsilon(t)} v(t,x)^p dx\;.
\end{eqnarray*}
By using the implicit solution formula \eqref{Burger-solution} for $u(t,x)$, we can show that 
\[
\min_{X^{-\epsilon}(t) \leq x \leq X^\epsilon(t)} (-u_x(t,x)) \geq \frac{1-\epsilon^2}{(1+\epsilon^2)^2} \left (\frac{1}{T^*-t+ 4 \epsilon^2}\right).
\]
For any fixed $0<T < T^*$, by taking  $\epsilon = \min\{1/3, \sqrt{(T^* - t)/20}\}$, we obtain 
\[
\frac{d}{dt} \| v(t)\|_{L^p}^p \geq \frac{(p-1)}{2(T^*-t)} \| v(t)\|_{L^p}^p, \quad \mbox{for any} \; t \in [0, T], \quad 1< p <\infty \;.
\]
Thus, by choosing a different initial condition $v_0^\epsilon \neq 0$ with a smaller $\epsilon$ as $T$ approaches to $T^*$ we obtain 
\[
\sup_{v \in V,\; v \neq 0} \frac{\| v(t)\|_{L^p}}{\| v(0)\|_{L^p}} \geq \frac{(T^*)^\mu}{(T^* - t)^\mu} , \quad \mu = \frac{1}{2} - \frac{1}{2p}, \quad \mbox{for any} \;t \in [0, T]\;, \; p > 1 .
\]
Since $V$ is a subspace of $H^1(\Omega)\cap L^p(\Omega )$, we conclude that 
\[
\gamma_p (T) \geq \sup_{v \in V,\; v \neq 0} \frac{\| v(T)\|_{L^p}}{\| v(0)\|_{L^p}} \;,
\]
which implies 
$\lim \sup_{T \rightarrow T^*} \gamma_p (T) = \infty \;$.
On the other hand, using the implicit solution formula, we can solve for the solution of the inviscid Burgers equation with a perturbed initial condition, $u^\epsilon_0(x) = u_0(x)+v_0^\epsilon(x)$. The solution $u^\epsilon(t,x)$ is defined implicitly by $u^\epsilon_0(x)$ via \eqref{Burger-solution} with the blow-up time given by $T^\epsilon = 1/\max_x (-(u_0^\epsilon)_x)$. By construction, we have $(v^\epsilon_0)_x = O(\epsilon^3)$, therefore we can show that the perturbed solution $u^\epsilon(t,x)$ develops a blow-up at time $T^\epsilon = T^* +O(\epsilon^3)$ with singularity located at $x^\epsilon = x^* + O(\epsilon^3)$. By construction, our perturbed initial condition $u^\epsilon_0$ satisfies the condition stated in Proposition $9$ in \cite{masmoudi2020} with $i=1$. By applying Proposition $9$ in \cite{masmoudi2020}, we conclude that the inviscid Burgers equation develops stable asymptotically self-similar blow-up solutions. However, by the notion of stability from \cite{vasseur2020}, the blow-up of the inviscid Burgers equation would be considered hydrodynamically unstable even within this class of perturbation.

In \cite{chenhou2022}, we will show that
$\lim \sup_{T \rightarrow T^*} \gamma_p (T) = \infty $ ($1<p < \infty $)
for the $3$D axisymmetric Euler equations with $C^{1,\alpha}$ initial velocity and boundary in the same setting as the Hou-Luo blow-up scenario \cite{luo2014toward} except for the regularity of the initial data. On the other hand, by using the dynamic rescaling formulation \cite{mclaughlin1986focusing,chen2019finite2,chen2020finite,chen2021finite3}, we have been able to prove in \cite{chen2019finite2} that the $3$D axisymmetric Euler equations with boundary develop a stable blow-up solution for a class of $C^{1,\alpha}$  initial velocity with finite energy. The blow-up solution is stable in the sense that it can be rescaled dynamically, and the rescaled solution is close to an approximate self-similar profile in some weighted Sobolev space.

The key difference between our approach to stability and the approach adopted in \cite{vasseur2020,vasseur2021} is that we first reformulate the Euler equations using the dynamic rescaling formulation (or using the modulation technique as in \cite{elgindi2019finite,Elg19}), and then perform linearization around an approximate self-similar profile. The linearized equation that we obtained is very different from \eqref{linear-euler}. By choosing appropriate normalization conditions and using the odd symmetry of the solution along the $z$ direction, we can eliminate some potentially unstable modes and prove stability of the linearized equation in an appropriately chosen functional space. 
The dependence of the blow-up time on the perturbation is naturally accounted in our analysis. The odd symmetry of the solution and the perturbation plays a crucial role in establishing linear stability of the approximate self-similar profile. The results in \cite{vasseur2020,vasseur2021} still apply with minor modification when this odd symmetry is imposed. 

\vspace{-0.12in}

\subsection{Stability of the self-similar profile to small perturbation of initial data} 

In this subsection, we study whether the observed nearly self-similar profile is stable with respect to a small perturbation of the initial condition. Our study shows that the approximate self-similar profile seems to be very stable to small perturbation of the initial data. On the other hand, the solution structure could have a very different behavior if we make an $O(1)$ perturbation to the initial data.

In our study, we solve the $3$D Euler equations using four different initial data given below. 

\vspace{0.05in}
\noindent
{\bf Case $1$.} We choose the same initial condition given in \eqref{eq:initial-data}, 

\begin{equation}
\label{IC:case1}
u_1 (0,r,z) =\frac{12000(1-r^2)^{18}
\sin(2 \pi z)}{1+12.5(\sin(\pi z))^2}, \quad \om_1(0,r,z)=0.
\end{equation}
 
 \vspace{0.05in}
 \noindent
{\bf Case $2$.} We choose the initial condition as a small perturbation to the initial condition defined in \eqref{IC:case1}, 

 \begin{equation}
\label{IC:case2}
u_1 (0,r,z) =\frac{12000(1-r^2)^{18}
\sin(2 \pi z)}{1+12.5(\sin(\pi z))^2} + 
\frac{(1-r^2)^{10}
\sin(6 \pi z)}{1+12.5(\sin(\pi z))^2} 
, \quad \om_1(0,r,z)=0.
\end{equation}
 
 \vspace{0.05in}
 \noindent
{\bf Case $3$.} We choose the initial condition as a larger perturbation to the initial condition defined in \eqref{IC:case1}, 

 \begin{equation}
\label{IC:case3}
u_1 (0,r,z) =\frac{12000(1-r^2)^{18}
\sin(2 \pi z)}{1+12.5(\sin(\pi z))^2} + 
\frac{42(1-r^2)^{6}
\sin(10 \pi z)}{1+12.5(\sin(\pi z))^2} 
, \quad \om_1(0,r,z)=0.
\end{equation}
 
\vspace{0.05in}
 \noindent
{\bf Case $4$.} We choose the initial condition that is $O(1)$ perturbation to the initial condition \eqref{IC:case1}, 

 \begin{equation}
\label{IC:case4}
u_1 (0,r,z) =\frac{12000(1-r^2)^{18}
\sin(4 \pi z)}{1+12.5(\sin(\pi z))^2}
, \quad \om_1(0,r,z)=0.
\end{equation}

  
 \begin{figure}[!ht]
\centering
	\begin{subfigure}[b]{0.38\textwidth}
    \includegraphics[width=1\textwidth]{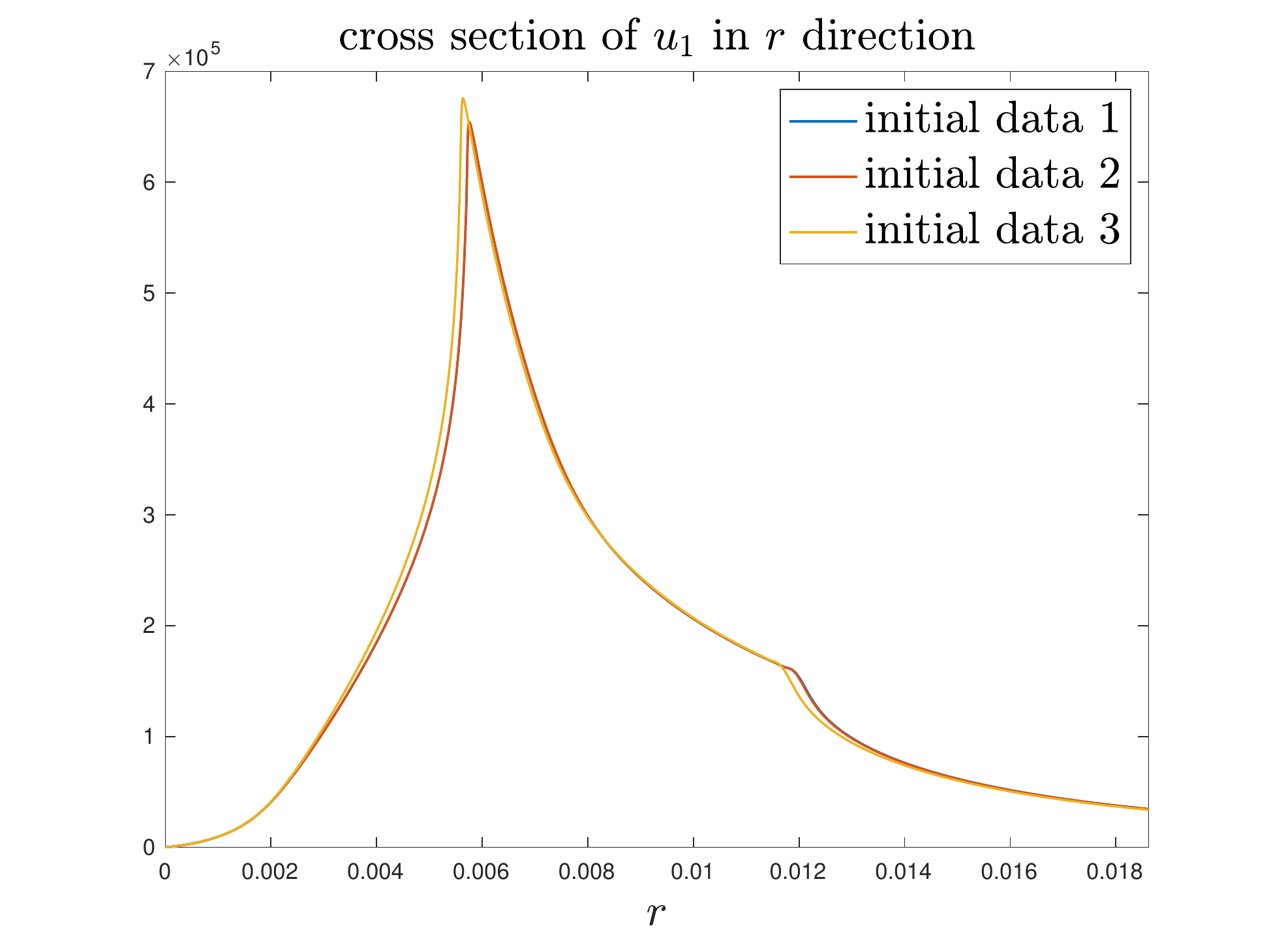}
    \caption{un-normalized $r$-cross section of $u_1$}
    \end{subfigure}
  	\begin{subfigure}[b]{0.38\textwidth} 
    \includegraphics[width=1\textwidth]{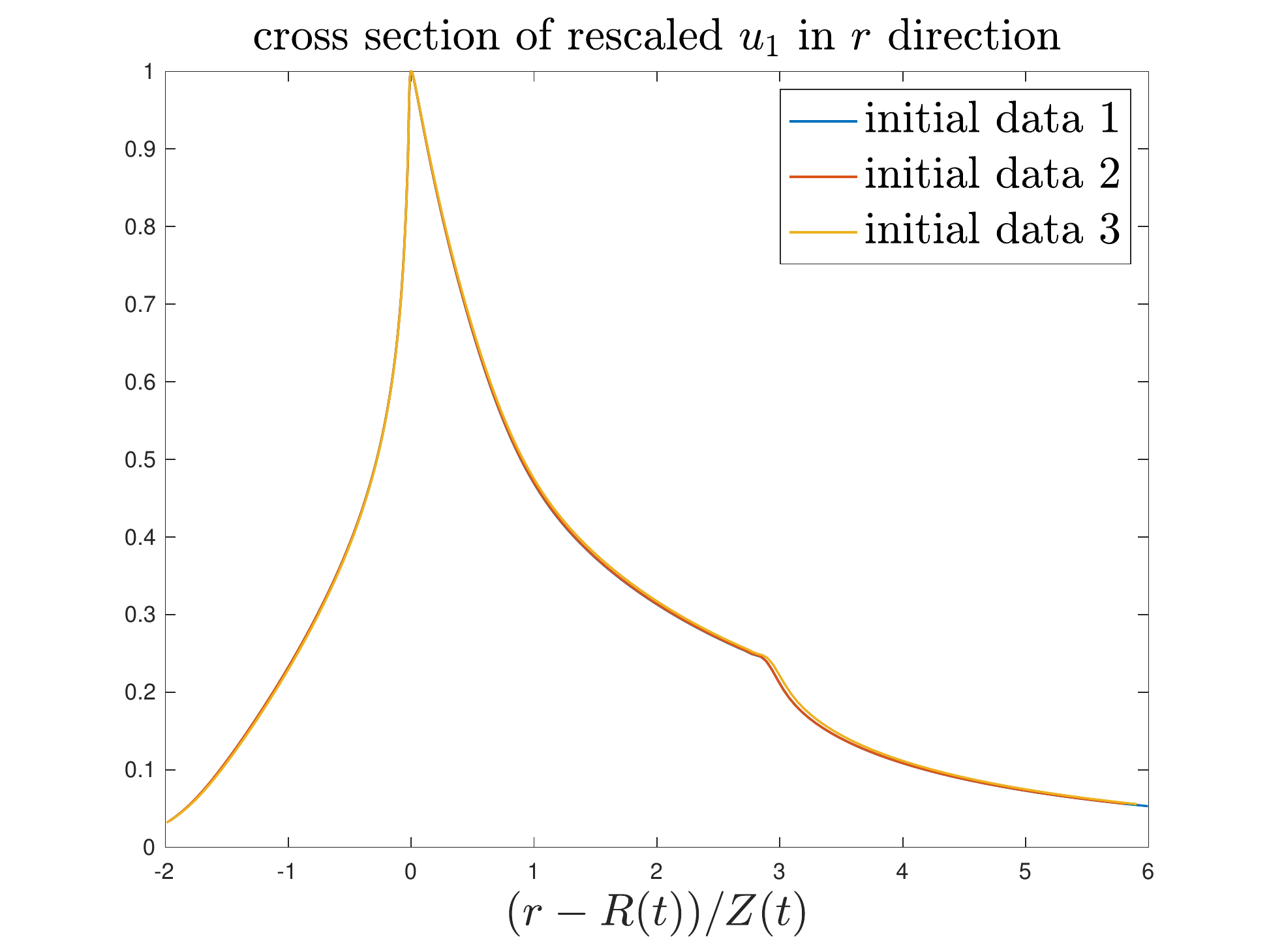}
    \caption{rescaled $r$-cross section of $u_1$}
    \end{subfigure} 
    \begin{subfigure}[b]{0.38\textwidth} 
    \includegraphics[width=1\textwidth]{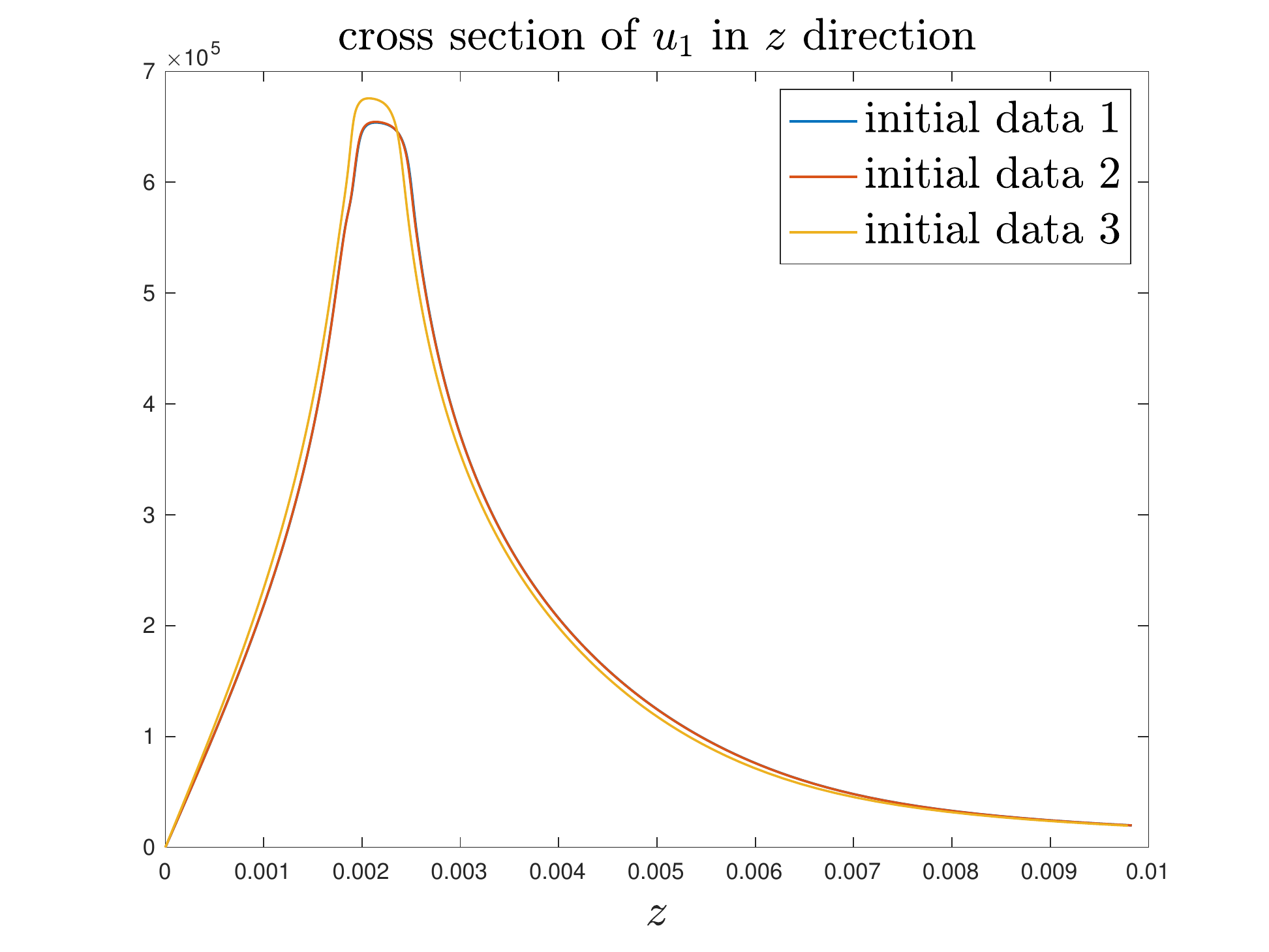}
    \caption{un-normalized $z$-cross section of $u_1$}
    \end{subfigure}
  	\begin{subfigure}[b]{0.38\textwidth} 
    \includegraphics[width=1\textwidth]{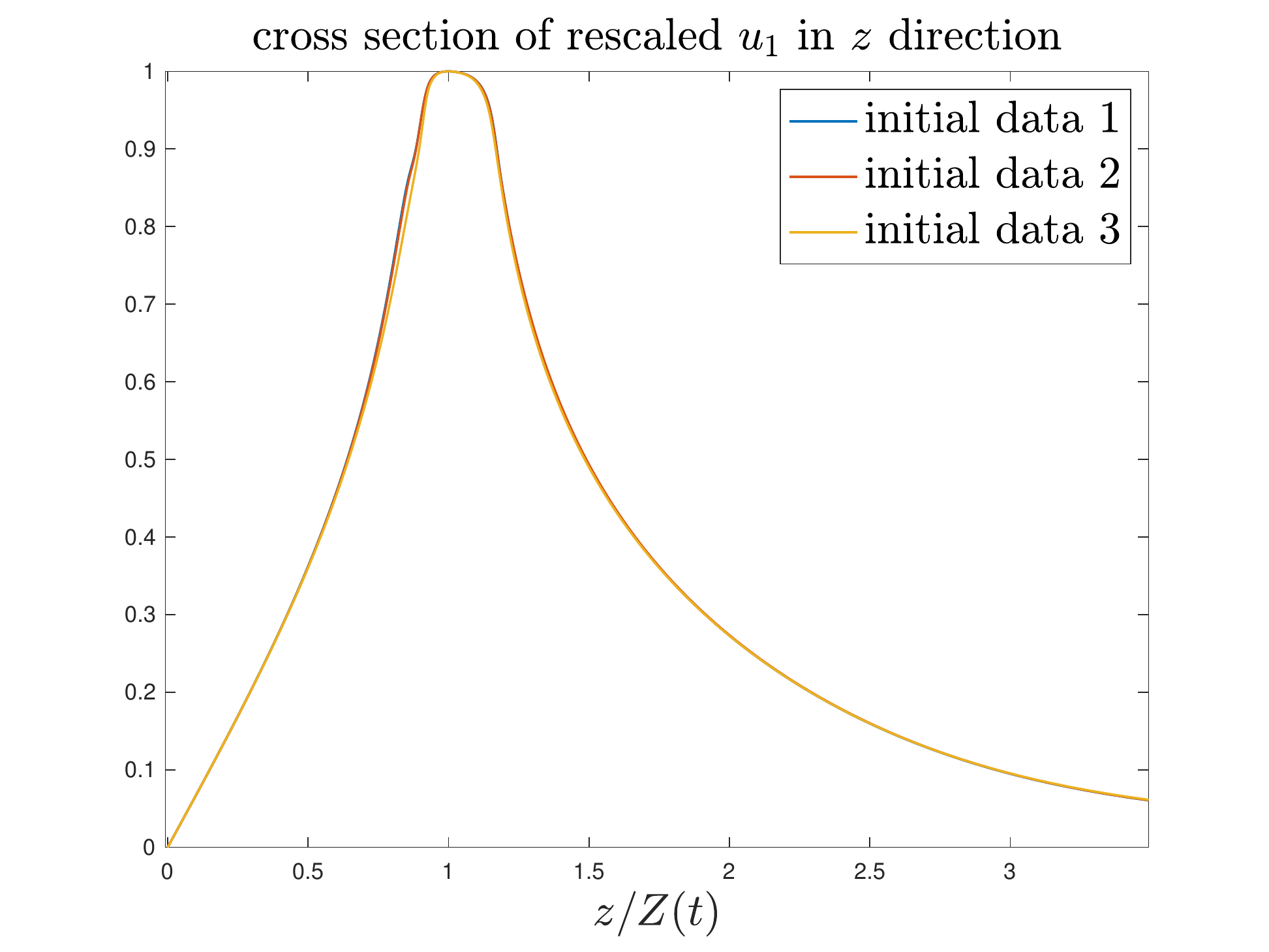}
    \caption{rescaled $z$-cross section of $u_1$}
    \end{subfigure} 
    \caption[compare u1-profile]{(a) The $r$-cross section of $u_1$ at $z=Z(t)$ at $75000$ time steps for the first three initial data in the original physical space. (b) Rescaled $\xi$-cross section of $u_1$ at $75000$ time steps for the first three initial data. (c) The $z$-cross section of $u_1$ at $r=R(t)$ at $75000$ time steps for the first three initial data in the original physical space. (d) Rescaled $\zeta$-cross section of $u_1$ at $75000$ time steps for the first three initial data.}   
    \label{fig:stability_data}
       \vspace{-0.1in}
\end{figure}

 \vspace{0.05in}
The relative size of the perturbation in Case $2$ is approximately $1.928\cdot 10^{-4}$ while the relative size of the perturbation in Case $3$ is about $10^{-2}$. Notice that the decay of the perturbation along the $r$-direction is slower than the original unperturbed initial condition. Moreover, the perturbation along the $z$-direction is more oscillatory. One would expect that the blow-up process may be sensitive to the perturbation of the initial data due to the strong shearing induced by the Euler equations. To our surprise, the blow-up profile is very stable with respect to the small perturbation of the initial data. 

In Case $4$, we just change $\sin(2 \pi z)$ in the numerator in the original initial condition to $\sin(4 \pi z)$. Everything else is the same. With this change, $u_1$ is no longer non-negative for $z \in [0, 0.5]$. In fact, $u_1$ becomes negative for $z$ near $0.5$. This introduces an $O(1)$ structural change to $u_1$ and the solution develops a two-scale structural similar to that observed in \cite{Hou-Huang-2021,Hou-Huang-2022}. In Cases $2$ and $3$, the perturbation is also oscillatory and becomes negative for $z \in [0,0.5]$. But in these two cases, the size of the perturbation is relatively small. The dominating part still comes from the original initial condition. Thus, the solution behavior is dominated by the unperturbed part of the initial condition and the solution behaves qualitatively similar to that of the original initial condition given by \eqref{IC:case1}.

 \begin{figure}[!ht]
\centering
  	\begin{subfigure}[b]{0.38\textwidth} 
    \includegraphics[width=1\textwidth]{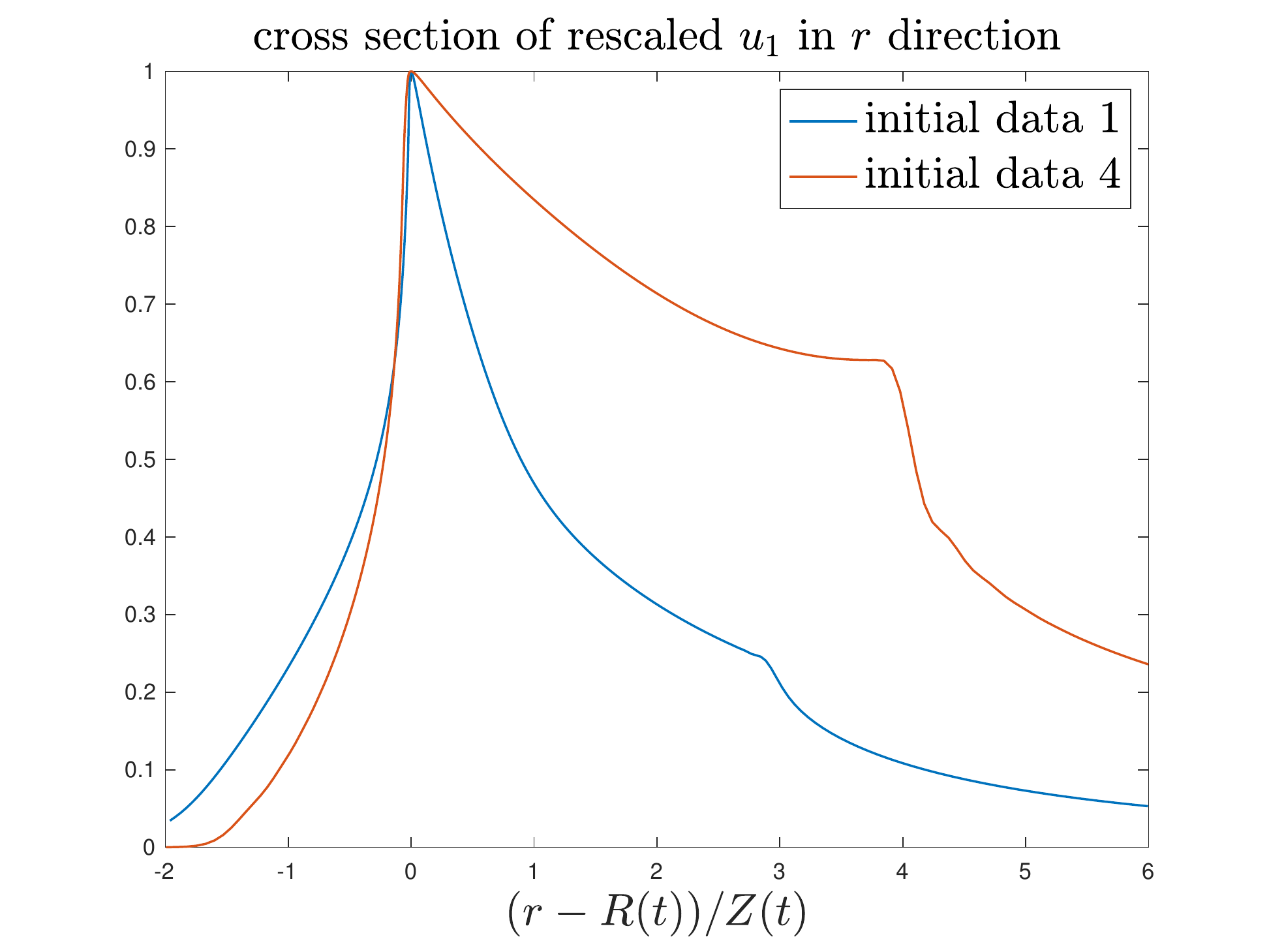}
    \caption{rescaled $r$-cross section of $u_1$}
    \end{subfigure} 
  	\begin{subfigure}[b]{0.38\textwidth} 
    \includegraphics[width=1\textwidth]{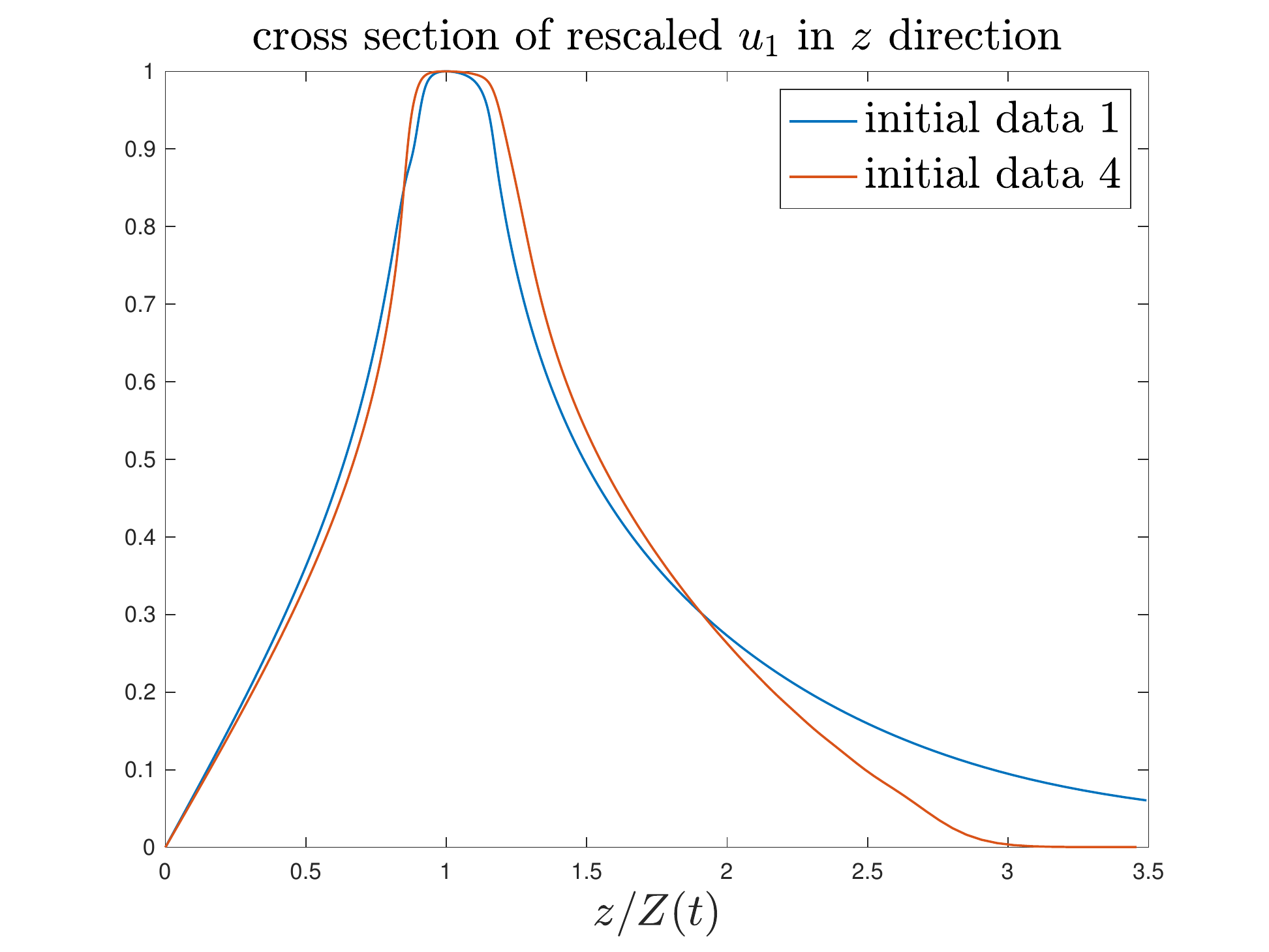}
    \caption{rescaled $z$-cross section of $u_1$}
    \end{subfigure} 
    \caption[compare u1-profile]{
    (a) The $r$-cross section of $u_1$ at $z=Z(t)$ at $75000$ time steps for the Case $1$ and the Case $4$ initial data in the original physical space. 
    (b) Rescaled $\rho$-cross section of $u_1$ at $75000$ time steps for the Case $1$ and Case $4$ initial data. }
    \label{fig:stability_data_4}
       \vspace{-0.15in}
\end{figure}

In Figure \ref{fig:stability_data}, we plot the solution $u_1$ after computing $75000$ time steps using the three initial conditions defined above. The time for the first initial condition after $75000$ time steps is $t_1=0.002276527323$, the time for the second initial condition is $t_2 = 0.002276019761$, and the time for the third initial condition is $t_3=0.002229117298$. We can see that $t_1 > t_2 > t_3$ and the gap between $t_1$ and $t_3$ is the largest. However, the growth of the third initial condition is the fastest in time.

From the $r$ and $z$ cross sections of $u_1$ in the original physical space, we can see that the solution obtained from the first initial condition agrees very well with that obtained from the second initial condition. This is also expected since the perturbation of the second initial condition is quite small. On the other hand, there is a noticeable difference between the solution obtained from the third initial condition and that from the first initial condition. However, when we plot the rescaled solution $U$ as a function of $\xi$ and $\zeta$ defined below: 
\[
u_1 = \max(u_1) U (t,\xi, \zeta ), \quad 
\xi = (r-R(t))/Z(t), \quad \zeta = z/Z(t)\;,
\]
we observe that the rescaled profiles $U$ as a function of the dynamically rescaled variables $(\xi,\zeta)$ for the first three initial conditions almost collapse on each other. This seems to suggest that the nearly self-similar profile of the solution is very stable to the small perturbation of the initial condition.

\begin{figure}[!ht]
\centering
	\begin{subfigure}[b]{0.38\textwidth}
    \includegraphics[width=1\textwidth]{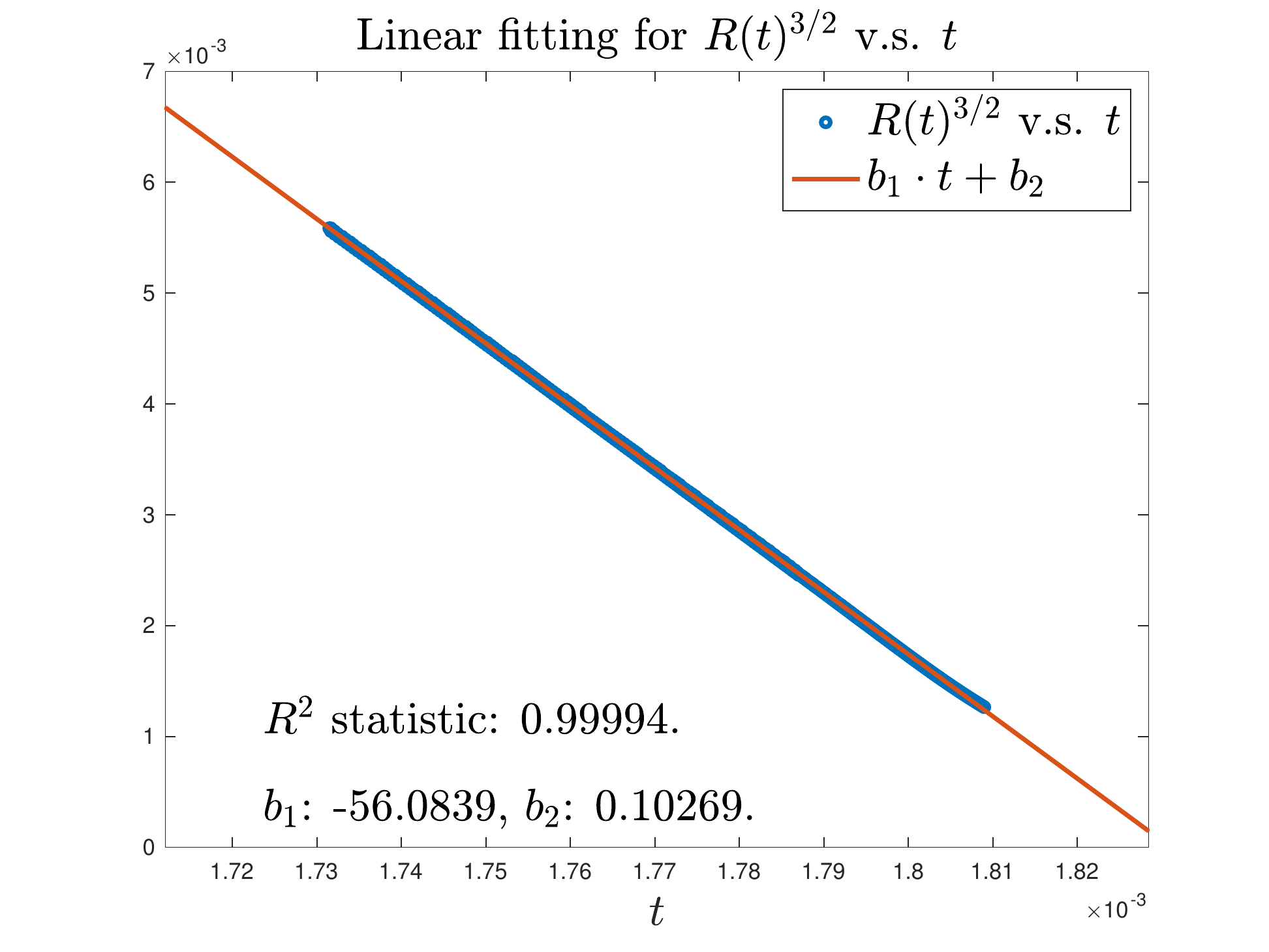}
    \caption{Linear fitting for $R(t)^{3/2}$}
    \end{subfigure}
  	\begin{subfigure}[b]{0.38\textwidth} 
    \includegraphics[width=1\textwidth]{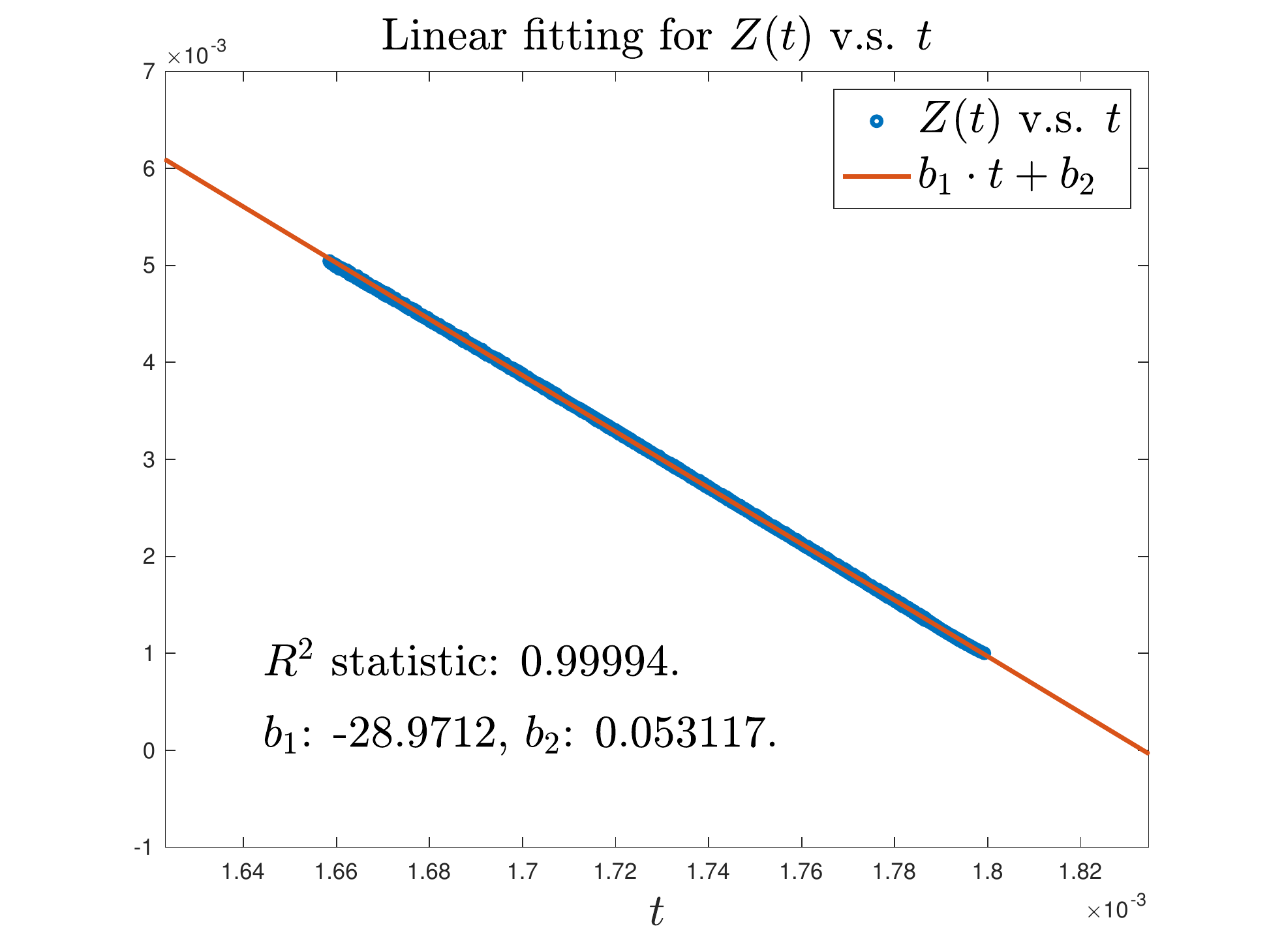}
    \caption{Linear fitting for $Z(t)$}
    \end{subfigure} 
    \begin{subfigure}[b]{0.38\textwidth} 
    \includegraphics[width=1\textwidth]{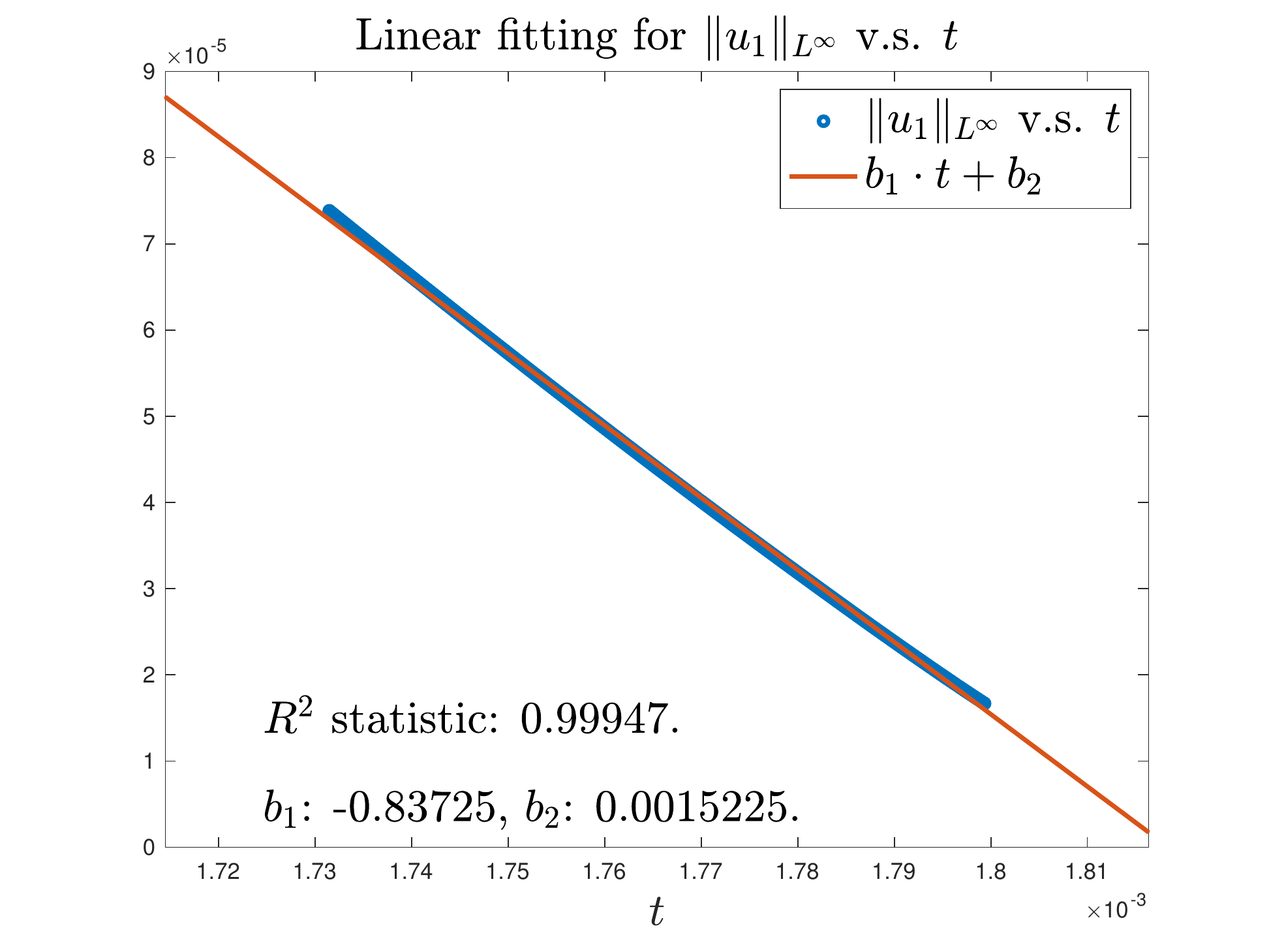}
    \caption{Linear fitting for $\|u_1\|_{L^\infty}^{-1}$}
    \end{subfigure}
  	\begin{subfigure}[b]{0.38\textwidth} 
    \includegraphics[width=1\textwidth]{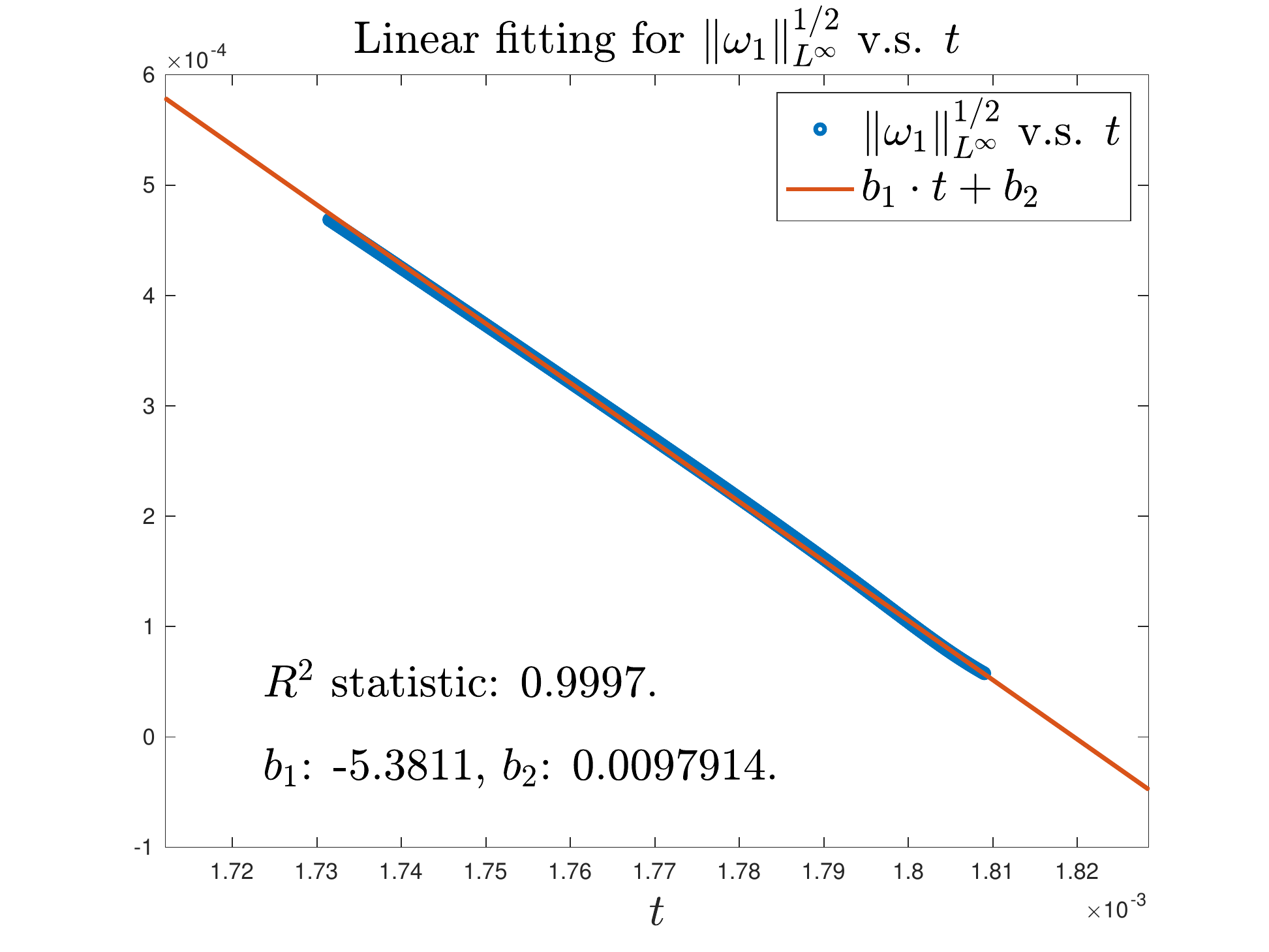}
    \caption{Linear fitting for $\|\omega_1\|_{L^\infty}^{-1/2}$}
    \end{subfigure} 
    \caption[compare u1-profile]{Linear fitting for the solution using the Case $4$ initial condition. (a) Linear fitting for $R(t)^{3/2}$. (b) Linear fitting for $Z(t)$. (c) Linear fitting for $\|u_1\|_{L^\infty}^{-1}$. (d) Linear fitting for $\|\omega_1\|_{L^\infty}^{-1/2}$.}   
    \label{fig:Fitting_data_4}
       \vspace{-0.15in}
\end{figure}

Next, we study the comparison between the Case $1$ and Case $4$ initial data. As we mentioned before, the perturbation in the Case $4$ initial data is an $O(1)$ perturbation to the Case $1$ initial data and the structures of the two initial data for $u_1$ are also different. In Figure \ref{fig:stability_data_4}, we plot the cross sections of the solution $u_1$ after computing $75000$ time steps using these two initial conditions. The time for the Case $1$ initial condition after $75000$ time steps is $t_1=0.002276527323$, the time for the Case $4$ initial condition is $t_4 = 0.001818245904$. We observe that the solution for the Case $4$ initial condition develops a two-scale solution structure with a sharp front. 

In Figure \ref{fig:stability_data_4}(a), we plot the $r$ cross section of $u_1$ in the original physical space. We can see that the solution obtained from the Case $1$ initial condition is significantly different from that obtained from the Case $4$ initial condition. Interestingly, we observe a no-spinning region developed between the sharp front and the symmetry axis for the Case $4$ initial data as we observed in the two-scale traveling wave singularity in \cite{Hou-Huang-2021,Hou-Huang-2022}. In Figure \ref{fig:stability_data_4}(b), we plot the cross section using the dynamically rescaled variables defined above.
We observe that the rescaled profiles $U$ as a function of the dynamically rescaled variables $(\xi,\zeta)$ for the two initial conditions do not agree with each other. They seem to converge to two different self-similar profiles.

To illustrate the two-scale solution structure, we perform linear fitting for $R(t)$, $Z(t)$, $\|u_1\|_{L^\infty}$ and $\|\omega_1\|_{L^\infty}$ in Figure \ref{fig:Fitting_data_4}. The linear fitting results obtained in Figure \ref{fig:Fitting_data_4} seem to imply the following scaling properties for the solution using the Case $4$ initial condition:
\[
R(t) \sim (T-t)^{2/3}, \quad Z(t) \sim (T-t),\quad
 \|u_1\|_{L^\infty} \sim \frac{1}{T-t}\;, \quad
\|\omega_1\|_{L^\infty} \sim \frac{1}{(T-t)^2} \;.
\]
These scaling properties suggest a two-scale blow-up structure with  $Z(t)/R(t) \sim (T-t)^{1/3}$. The fact that $\|\omega_1\|_{L^\infty} \sim \frac{1}{(T-t)^2} $ is also consistent with the scaling property $Z(t) \sim (T-t)$. In comparison, we have $Z(t) \sim (T-t)^{1/2}$ and $\|\omega_1\|_{L^\infty} \sim \frac{1}{(T-t)^{3/2}} $ for the solution using the Case $1$ initial condition.

\vspace{-0.1in}
\section{Concluding Remarks}
\label{sec:conclude}

In this paper, we presented strong numerical evidences that the $3$D incompressible axisymmetric Euler equations develop a finite time singularity at the origin. An important feature of this potential singularity is that the solution develops an essentially one-scale traveling wave solution with scaling properties compatible with those of the Navier--Stokes equations. This property is critical for the potentially singular behavior of the $3$D Navier--Stokes equations.  We presented some qualitative numerical evidence that the Euler solution seems to develop nearly self-similar scaling properties.
Moreover, the nearly self-similar profile seems to be very stable to small perturbation of the initial data. 
On the other hand, the specific form of our initial data is important in generating the observed scaling properties. 
If we make an $O(1)$ change to the initial data, the solution could behave very differently. We provided one such example in which a simple change of the initial data leads to a two-scale traveling wave solution, which is similar in spirit to what we observed in \cite{Hou-Huang-2021,Hou-Huang-2022}.

We plan to use the dynamic rescaling formulation \cite{mclaughlin1986focusing,chen2019finite2,chen2020finite,chen2021finite3} to investigate  the potential blow-up of the $3$D Euler and Navier--Stokes equations in our future work. One important advantage of using the dynamic rescaling approach is that we can use a fixed mesh to solve the dynamically rescaled Euler or Navier--Stokes equations. This avoids the numerical dissipation introduced by the frequent changes of adaptive mesh in the late stage when we solve the $3$D Euler or Navier--Stokes equations in the physical domain. Another benefit is that the use of the dynamic rescaling formulation may eliminate the mild high frequency instability that we observed when we solve the Euler equation in the physical domain without using a second order numerical diffusion.

\vspace{0.1in}
{\bf Acknowledgments.} The research was in part supported by NSF Grants DMS-$1907977$ and DMS-$1912654$, and the Choi Family Gift Fund.
I would like to thank Dr. De Huang for very helpful discussions regarding the design of the adaptive mesh strategy. I would also like to thank Professor Vladimir Sverak, Jiajie Chen, Dr. De Huang, and the referees for their very constructive comments and suggestions, which significantly improves the quality of this paper. 
Finally, I have benefited a lot from the AIM SQarRE ``Towards a $3$D Euler singularity''.

\appendix

\section{Construction of the adaptive mesh}\label{apdx:adaptive_mesh}
In this appendix, we describe our adaptive mesh strategy to study the singularity formation near the origin $(r,z) = (0,0)$. 
We will use the method described in Appendix B of \cite{Hou-Huang-2021} to construct our adaptive mesh maps $r=r(\rho)$ and $z=z(\eta)$. We will discretize the equations in the transformed variables $(\rho,\eta)$ with $n_1$ grid points along the $z$ direction and $n_2$ grid points along the $r$-direction. 

The adaptive mesh strategy described in \cite{Hou-Huang-2021} was inspired by the adaptive mesh strategy introduced in \cite{luo2014toward}. On the other hand, the adaptive mesh strategy developed in \cite{luo2014toward} is simpler due to the fact that the singularity is located at a fixed stagnation point on the boundary $(r,z)=(1,0)$, $\omega_1 >0$ for $z>0$ and has a bell-shaped structure near the singularity. In our case, we have a traveling wave solution that approaches the origin with $\omega_1$ changing sign in the most singular region. The singular region has a more complicated shape since the potentially singular solution produces a strong shearing flow traveling downstream. Thus we need to introduce a moving frame and design our adaptive mesh map $r(\rho)$ and $z(\eta)$ to resolve the solution in different regions. Our adaptive mesh strategy does not require that the solution has a bell-shaped structure in the most singular region. 
To construct $r(\rho)$, we use the distance  $dr$ between the location at which $u_1$  achieves its maximum and the location at which $u_{1r}$ achieves its maximum to define the boundary $r_i$ for different singular regions (phases). Similarly, we construct $z(\eta)$ by using the distance $dz$ between the location at which $\omega_1$  achieves its maximum and the location at which $\omega_{1z}$ achieves its maximum to define the boundary $z_j$ for different phases.

\vspace{-0.05in}
\subsection{The adaptive (moving) mesh algorithm}
To effectively and accurately compute the potential blow-up, we have designed a special meshing strategy that is dynamically adaptive to the singular structure of the solution. The adaptive mesh covering the half-period computational domain $\mathcal{D}_1 = \{(r,z):0\leq r\leq 1,0\leq z\leq 1/2\}$ is characterized by a pair of analytic mesh mapping functions
\[r = r(\rho),\quad \rho\in [0,1];\quad z = z(\eta),\quad \eta\in[0,1].\]
These mesh mapping functions are both monotonically increasing and infinitely differentiable on $[0,1]$, and satisfy 
$r(0) = 0,\; r(1) = 1,\; z(0) = 0,\; z(1) = 1/2.$
In particular, we construct these mapping functions by carefully designing their Jacobians/densities
$r_\rho = r'(\rho),\quad z_\eta = z'(\eta)$
using analytic functions that are even functions at $0$. The even symmetries ensure that the resulting mesh can be smoothly extended to the full-period cylinder $ \{(r,z):0\leq r\leq 1,-1/2\leq z\leq 1/2\}$. The density functions contain a small number of parameters, which are dynamically adjusted to the solution. Once the mesh mapping functions are constructed, the computational domain is covered with a tensor-product mesh:
\begin{equation}\label{eq:mesh}
\mathcal{G} = \{(r_i,z_j): 0\leq i\leq n_2,\ 0\leq j\leq n_1\},
\end{equation}
where $r_i^h = r(ih_\rho),\; h_\rho = 1/n_2;\; z_j^h = z(jh_\eta),\; h_\eta = 1/n_1.$
The precise definition and construction of the mesh mapping functions are described in Appendix B of \cite{Hou-Huang-2021}. 

Figure \ref{fig:map_density} gives an example of the densities $r_\rho,z_\eta$ (in log scale) we use in the computation. 
We design the densities $r_\rho,z_\eta$ to have three phases: 
\begin{itemize}
\item Phase $1$ covers the inner profile of the smaller scale near the sharp front;
\item Phase $2$ covers the outer profile of the larger scale of the solution;
\item Phase $3$ covers the (far-field) solution away from the symmetry axis $r=0$.
\end{itemize} 
We add a phase $0$ in the density $r_\rho$ to cover the region near $r=0$ and also add a phase $0$ in the density $z_\eta$ to cover the region near $z=0$ in the late stage. In our computation, the number (percentage) of mesh points in each phase are fixed, but the physical location of each phase will change in time, dynamically adaptive to the structure of the solution. Between every two neighboring phases, there is also a smooth transition region that occupies a fixed percentage of mesh points. 

\begin{figure}[!ht]
\centering
    \includegraphics[width=0.38\textwidth]{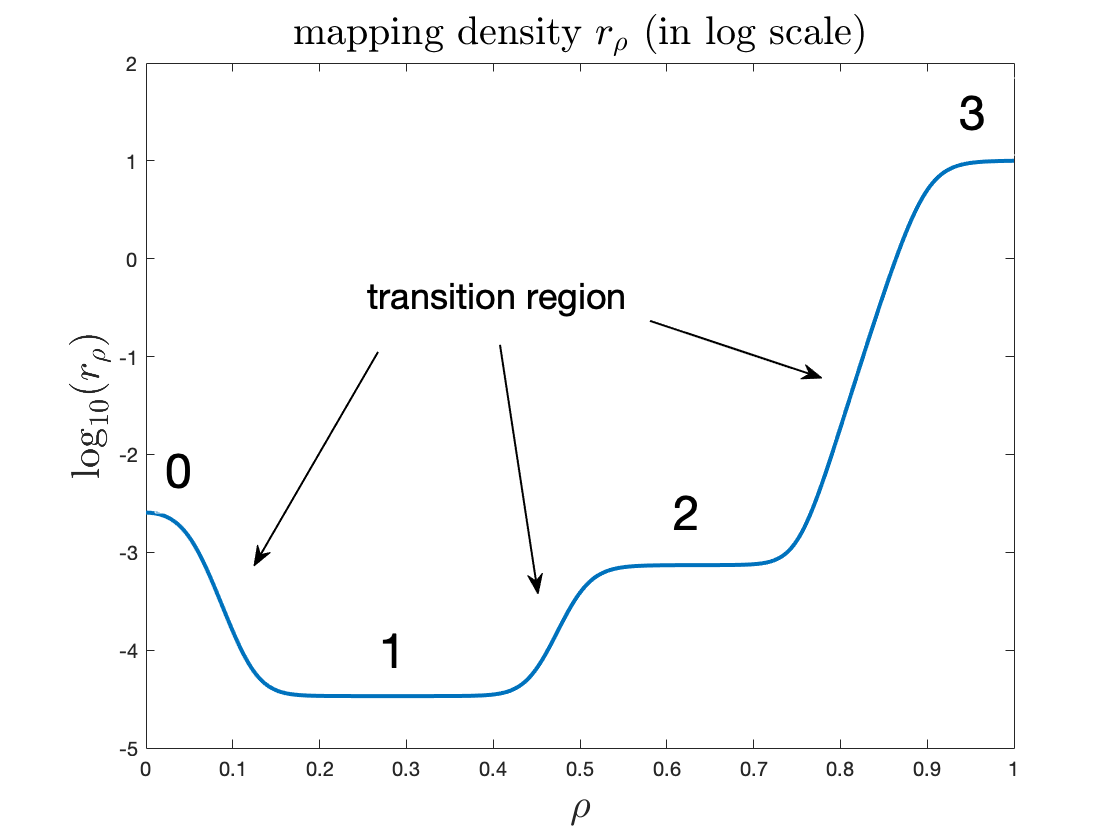}
    \includegraphics[width=0.38\textwidth]{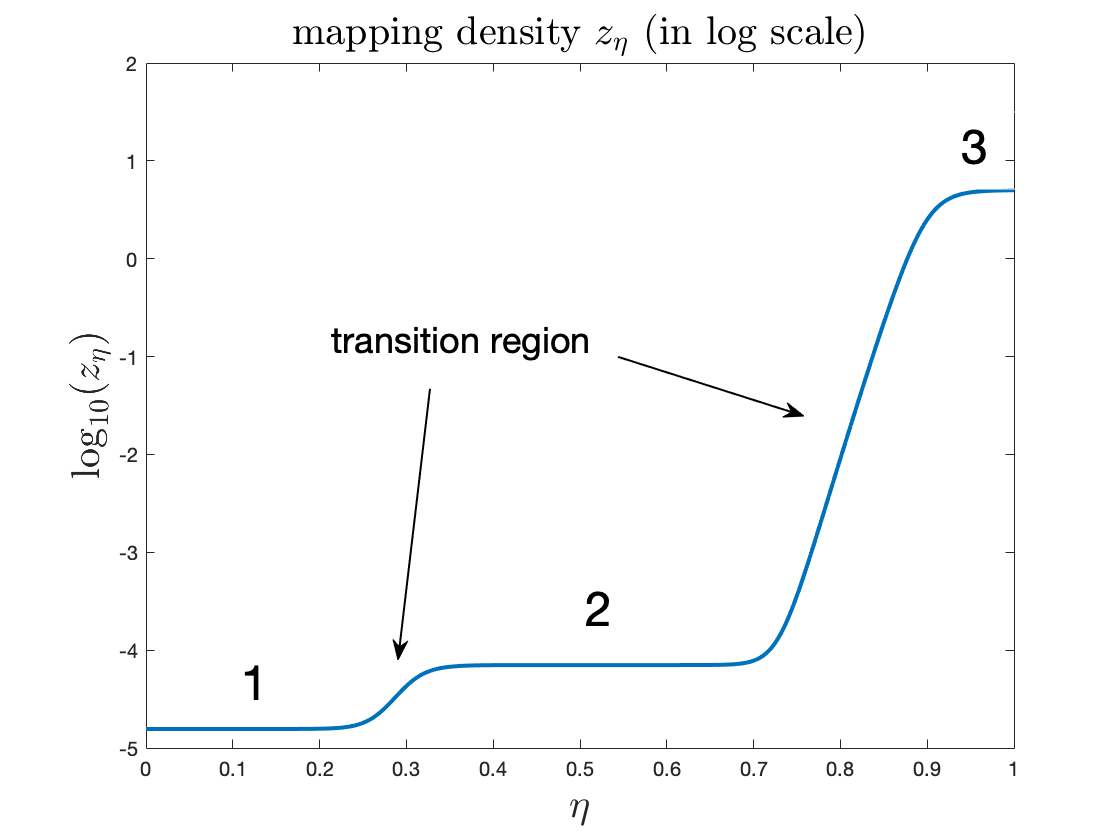}
    \caption[Map densities]{The mapping densities $r_\rho$ (left) and $z_\eta$ (right) with phase numbers labeled. This figure is for illustration only. The parameters do not reflect the adaptive mesh used in our computation.}  
    \label{fig:map_density}
       \vspace{-0.05in}
\end{figure}

\vspace{-0.1in}
\subsection{Adaptive mesh for the $3$D Euler equations}
We use three different adaptive mesh strategies for three different time periods. The first time period corresponds to the time interval between $t=0$ and $t_1=0.002231338$ with $\|\omega(t_1)\|_{L^\infty}/\|\omega(0)\|_{L^\infty} \approx 46.54325$ for the $1536\times 1536$ grid and the number of time steps equal to $45000$. The second time period corresponds to the time interval between $t_1=0.002231338$ and $t_2 = 0.002264353$ with $\|\omega(t_2)\|_{L^\infty}/\|\omega(0)\|_{L^\infty} \approx 295.39986$ for the $1536\times 1536$ grid and the number of time steps equal to $60000$. The third time period is for $ t \geq t_2$. 

In the first time period, since we use a very smooth initial condition whose support covers the whole domain, we use the following parameters $r_1=0.001,\; r_2=0.05, \; r_3=0.2$ and $s_{\rho_1}=0.001$, $s_{\rho_2}=0.5$, $s_{\rho_3}=0.85$ to construct the mapping $r=r(\rho)$ using a four-phase map. Similarly, we use the following parameters $z_1=0.1,\; z_2=0.25$ and 
$s_{\eta_1}=0.5$, $s_{\eta_2}=0.85$ to construct the mapping $z=z(\eta)$ using a three-phase map. We then update the mesh $z=z(\eta)$ dynamically using $z_1=2 z(I)$ and $z_2=10 z(I)$ with $s_{\eta_1}=0.6$, $s_{\eta_2}=0.9$ when $I < 0.25n_1$, but keep $r=r(\rho)$ unchanged during this early stage. Here $I$ is the grid point index along the $z$-direction at which $\omega_1$ achieves its maximum. 

In the second time period, we use the following parameters 
$s_{\rho_1}=0.05$, $s_{\rho_2}=0.6$, $s_{\rho_3}=0.9$, 
$r_2 = r(J) + 2dr$,
$r_1=\max((s_{\rho_1}/s_{\rho_2})r_2,r(J_r) - 5dr)$, and 
$r_3 = \max(3r(J),(r_2 - r_1)(s_{\rho_3}-s_{\rho_2})/(s_{\rho_2}-s_{\rho_1}) + r_2)$, where $J$ is the grid index at which $u_1$ achieves its maximum along the $r$-direction, $J_r$ is the grid index at which $u_{1,r}$ achieves its maximum along the $r$-direction, and $dr = r(J) - r(J_r)$. We update the mapping $r(\rho)$ dynamically when $J_r < 0.2n_2$. The adaptive mesh map for $z(\eta)$ in the second time period remains the same as in the first time period.

In the third time period, we need to allocate more grid points to resolve the sharp front. We use the following parameters 
$s_{\rho_1}=0.05$, $s_{\rho_2}=0.65$, $s_{\rho_3}=0.9$, 
$r_2 = r(J) + 10dr$,
$r_1=\max((s_{\rho_1}/s_{\rho_2})r_2,r(J_r) - 3dr)$, and 
$r_3 = \max(2.3r(J),(r_2 - r_1)(s_{\rho_3}-s_{\rho_2})/(s_{\rho_2}-s_{\rho_1}) + r_2)$.
To construct the mesh map $z(\eta$, we use the following parameters 
$s_{\eta_1}=0.05$, $s_{\eta_2}=0.65$, $s_{\eta_3}=0.9$, 
$z_2 = z(I_w) + 2dz$,
$z_1=\max((s_{\eta_1}/s_{\eta_2})z_2,z(I_{wz}) - 16dz)$, and 
$r_3 = \max(2.3z(I_w),(z_2 - z_1)(s_{\eta_3}-s_{\eta_2})/(s_{\eta_2}-s_{\eta_1}) + z_2)$, where $I_w$ is the grid index at which $\omega_1$ achieves its maximum along the $z$-direction, $I_{wz}$ is the grid index at which $\omega_{1,z}$ achieves its maximum along the $z$-direction, and $dz = z(I_w) - r(I_{wz})$. We will update $r(\rho)$ dynamically when $J_r < 0.2 n_2$ and update $z(\eta)$ when $I_z < 0.23 n_1$.

The three time periods for different grids will be defined similarly through the number of time steps by using a linear scaling relationship. For example, for the $1024\times1024$ grid, the first time period will be between $t=0$ and the time step $30000=45000 \cdot (1024/1536)$. The second time period will be between time steps $30000$ and $40000 = 60000\cdot (1024/1536)$. The third time period will be beyond time step $40000$. The three time periods for other grids are defined similarly.

As we mentioned in the Introduction, we use a second order Runge-Kutta method to discretize the Euler and Navier--Stokes equations in time with an adaptive time stepping strategy. We discretize the Euler and Navier--Stokes equations in the transformed domain $(\eta,\rho)$ using a uniform mesh with $h_1=1/n_1$ and $h_2=1/n_2$. We choose our adaptive time step as follows:
\begin{eqnarray*}
    k &=& \min(k_1,k_2)\;, \;\;k_1 = \min(\;0.2\min(h1,h_2)/\mbox{umax},\;10^{-3}\|u_1\|_{L^\infty}^{-1},\;2.5\cdot 10^{-7})\;,\\       
k_2 &=& 0.1\min(\min(h_1z_\eta)^2/\nu,\min(h_2r_\rho)^2/\nu )\;,
\end{eqnarray*}
where $\mbox{umax} = \max(\|u^r/r_\rho\|_{L^\infty},\|u^z/z_\eta\|_{L^\infty})$ is the maxixmum velocity in the transformed domain.
\bibliographystyle{abbrv}
\bibliography{reference}

\end{document}